\documentclass[11pt, openany, notitlepage]{book}

\usepackage[a4paper, hcentering, hmargin=2.5cm, vmargin=2.5cm]{geometry}
\usepackage[T1]{fontenc}
\usepackage[english]{babel}
\usepackage{lmodern}
\usepackage{microtype}

\usepackage{setspace}
\usepackage{indentfirst}
\usepackage{tocloft}
\usepackage{titlesec}
\usepackage{fancyhdr}
\usepackage{enumitem}
\usepackage{longtable}
\usepackage[bottom]{footmisc}

\usepackage[intlimits]{amsmath}
\usepackage{amsthm, amssymb, amsfonts}
\usepackage{mathtools}
\usepackage{thmtools}
\usepackage{mathrsfs}
\usepackage{bbm}
\usepackage[new]{old-arrows}
\usepackage{moreenum}

\newenvironment{abstract}{}{}
\usepackage{abstract}
\usepackage{tikz-cd}

\usepackage{csquotes}
\usepackage[backend=biber, style=numeric-comp, giveninits=true, sorting=nyt, maxbibnames=4]{biblatex}
\usepackage{hyperref}
\usepackage{xurl}
\usepackage[nameinlink, capitalize]{cleveref}
\hypersetup{pdfstartview={FitH}, colorlinks=true, allcolors=black, linktoc=all, bookmarksnumbered=true, breaklinks=true}

\setdisplayskipstretch{1.5}
\allowdisplaybreaks
\numberwithin{equation}{chapter}

\renewcommand{\thechapter}{\Roman{chapter}}
\renewcommand{\thesubsection}{\thesection.\alph{subsection}}
\titleformat{\subsection}{\normalfont\large\bfseries}{\thesubsection}{1em}{\itshape}
\titlespacing*{\section}{0pt}{0.75cm}{0.5cm}
\titlespacing*{\subsection}{0pt}{0.75cm}{0.5cm}

\renewcommand{\headrulewidth}{0pt}
\setlist{parsep=0pt}
\setlist[1]{labelindent=\parindent, leftmargin=*, align=left, topsep=8pt}
\setlist[2]{leftmargin=*, align=left}
\setlist[enumerate]{label=\normalfont(\alph*)}
\setlist[itemize]{topsep=12pt}

\SetEnumitemKey{num}{label=\normalfont(\arabic*)}
\SetEnumitemKey{cond}{label=\normalfont(\greek*)}
\SetEnumitemKey{env}{wide=\parindent, topsep=8pt, label=\normalfont\bfseries(\arabic*), ref=\normalfont(\arabic*)}
\SetEnumitemKey{equiv}{topsep=8pt, label=\normalfont(\roman*), align=right}

\DeclareNameAlias{default}{family-given}

\DeclareFieldFormat[article]{volume}{\mkbibbold{#1}}
\DeclareFieldFormat[article]{pages}{#1}
\DeclareFieldFormat[article]{year}{\mkbibparens{#1}}
\DeclareFieldFormat[book]{number}{\mkbibbold{#1}}
\DeclareFieldFormat[book]{date}{\mkbibparens{#1}}
\DeclareFieldFormat[incollection]{number}{\mkbibbold{#1}}
\DeclareFieldFormat[incollection]{date}{\mkbibparens{#1}}

\renewbibmacro{in:}{\ifentrytype{article}{}{\printtext{\bibstring{in}\intitlepunct}}}

\renewbibmacro*{issue+date}{}

\renewbibmacro*{volume+number+eid}{%
  \printfield{volume}%
  \printfield[parens]{number}%
  \setunit{\addcomma\space}%
}

\renewbibmacro*{note+pages}{%
  \setunit{\bibpagespunct}%
  \printfield{pages}%
  \setunit{\addspace}%
  \printfield{year}%
  \setunit{\addperiod\addspace}%
  \printfield{note}%
  \newunit
}

\renewbibmacro*{publisher+location+date}{%
  \printlist{location}%
  \iflistundef{publisher}
    {\setunit*{\addcomma\space}}
    {\setunit*{\addcolon\space}}%
  \printlist{publisher}%
  \setunit{\addspace}%
  \usebibmacro{date}%
  \newunit
}

\DeclareFieldFormat{doilink}{%
  \iffieldundef{doi}{%
    \iffieldundef{url}{%
      #1%
    }{%
      \href{\thefield{url}}{#1}%
    }%
  }{%
    \href{https://doi.org/\thefield{doi}}{#1}%
  }%
}

\DeclareBibliographyDriver{article}{%
  \usebibmacro{bibindex}%
  \usebibmacro{begentry}%
  \usebibmacro{author/translator+others}%
  \setunit{\printdelim{nametitledelim}}\newblock
  \usebibmacro{title}%
  \newunit
  \printlist{language}%
  \newunit\newblock
  \usebibmacro{byauthor}%
  \newunit\newblock
  \usebibmacro{bytranslator+others}%
  \newunit\newblock
  \printfield{version}%
  \newunit\newblock
  \usebibmacro{in:}%
  \printtext[doilink]{%
    \usebibmacro{journal+issuetitle}%
    \newunit
    \usebibmacro{byeditor+others}%
    \newunit
    \usebibmacro{note+pages}%
  }%
  \newunit\newblock
  \iftoggle{bbx:isbn}
    {\printfield{issn}}
    {}%
  \newunit\newblock
  \usebibmacro{doi+eprint+url}%
  \newunit\newblock
  \usebibmacro{addendum+pubstate}%
  \setunit{\bibpagerefpunct}\newblock
  \usebibmacro{pageref}%
  \newunit\newblock
  \iftoggle{bbx:related}
    {\usebibmacro{related:init}%
     \usebibmacro{related}}
    {}%
  \usebibmacro{finentry}%
}

\DeclareBibliographyDriver{book}{%
  \usebibmacro{bibindex}%
  \usebibmacro{begentry}%
  \usebibmacro{author/editor+others/translator+others}%
  \setunit{\printdelim{nametitledelim}}\newblock
  \usebibmacro{maintitle+title}%
  \newunit
  \printlist{language}%
  \newunit\newblock
  \usebibmacro{byauthor}%
  \newunit\newblock
  \usebibmacro{byeditor+others}%
  \newunit\newblock
  \printfield{edition}%
  \newunit
  \iffieldundef{maintitle}
    {\printfield{volume}%
     \printfield{part}}
    {}%
  \newunit
  \printfield{volumes}%
  \newunit\newblock
  \usebibmacro{series+number}%
  \newunit\newblock
  \printfield{note}%
  \newunit\newblock
  \printtext[doilink]{%
    \usebibmacro{publisher+location+date}%
  }%
  \newunit\newblock
  \usebibmacro{chapter+pages}%
  \newunit
  \printfield{pagetotal}%
  \newunit\newblock
  \iftoggle{bbx:isbn}
    {\printfield{isbn}}
    {}%
  \newunit\newblock
  \usebibmacro{doi+eprint+url}%
  \newunit\newblock
  \usebibmacro{addendum+pubstate}%
  \setunit{\bibpagerefpunct}\newblock
  \usebibmacro{pageref}%
  \newunit\newblock
  \iftoggle{bbx:related}
    {\usebibmacro{related:init}%
     \usebibmacro{related}}
    {}%
  \usebibmacro{finentry}
}

\DeclareBibliographyDriver{incollection}{%
  \usebibmacro{bibindex}%
  \usebibmacro{begentry}%
  \usebibmacro{author/translator+others}%
  \setunit{\printdelim{nametitledelim}}\newblock
  \usebibmacro{title}%
  \newunit
  \printlist{language}%
  \newunit\newblock
  \usebibmacro{byauthor}%
  \newunit\newblock
  \usebibmacro{in:}%
  \usebibmacro{maintitle+booktitle}%
  \newunit\newblock
  \usebibmacro{byeditor+others}%
  \newunit\newblock
  \printfield{edition}%
  \newunit
  \iffieldundef{maintitle}
    {\printfield{volume}%
     \printfield{part}}
    {}%
  \newunit
  \printfield{volumes}%
  \newunit\newblock
  \usebibmacro{series+number}%
  \newunit\newblock
  \printfield{note}%
  \newunit\newblock
  \printtext[doilink]{%
    \usebibmacro{publisher+location+date}%
    \newunit\newblock
    \usebibmacro{chapter+pages}%
  }%
  \newunit\newblock
  \iftoggle{bbx:isbn}
    {\printfield{isbn}}
    {}%
  \newunit\newblock
  \usebibmacro{doi+eprint+url}%
  \newunit\newblock
  \usebibmacro{addendum+pubstate}%
  \setunit{\bibpagerefpunct}\newblock
  \usebibmacro{pageref}%
  \newunit\newblock
  \iftoggle{bbx:related}
    {\usebibmacro{related:init}%
     \usebibmacro{related}}
    {}%
  \usebibmacro{finentry}%
}

\Crefname{section}{Sect.}{Sects.}
\Crefname{subsection}{Sect.}{Sects.}
\Crefname{para}{Para.}{Paras.}
\Crefname{definitions}{Definition}{Definitions}
\Crefname{examples}{Example}{Examples}
\Crefname{remarks}{Remark}{Remarks}
\Crefname{lemma}{Lemma}{Lemmata}
\Crefname{ntheorem}{Theorem}{Theorems}
\Crefname{nproposition}{Proposition}{Propositions}
\Crefname{page}{p.}{pp.}

\declaretheoremstyle[
    spaceabove=12pt, spacebelow=15pt,
    headfont=\normalfont\bfseries,
    notefont=\bfseries,
	bodyfont=\normalfont,
    notebraces={}{}
]{para}

\declaretheoremstyle[
    spaceabove=12pt, spacebelow=15pt,
    headfont=\normalfont\bfseries,
    notefont=\bfseries,
    bodyfont=\normalfont,
    headformat=swapnumber
]{def}

\declaretheoremstyle[
    spaceabove=12pt, spacebelow=15pt,
    headfont=\normalfont\bfseries,
    notefont=\bfseries,
    bodyfont=\normalfont,
    headindent=\parindent
]{ref}

\declaretheoremstyle[
    spaceabove=12pt, spacebelow=15pt,
    headfont=\normalfont\bfseries,
    notefont=\bfseries,
    bodyfont=\itshape,
    headformat=swapnumber
]{thm}

\declaretheoremstyle[
    spaceabove=12pt, spacebelow=15pt,
    headfont=\normalfont\bfseries,
    notefont=\bfseries,
    bodyfont=\itshape,
    headpunct={.\vspace{0.25\topsep}\nopagebreak\newline},
    postheadspace=0pt,
    headformat=swapnumber
]{thm_name}

\declaretheoremstyle[
    spaceabove=12pt, spacebelow=12pt,
    headfont=\normalfont\bfseries,
    notefont=\bfseries,
    bodyfont=\itshape,
    headpunct=.,
    headformat={\NAME~\NUMBER\NOTE}
]{stmnt}

\declaretheoremstyle[
    spaceabove=12pt, spacebelow=12pt,
    headfont=\normalfont\bfseries,
    notefont=\bfseries,
    bodyfont=\itshape,
    headpunct=.,
    headformat={\NAME\NOTE}
]{nstmnt}

\declaretheoremstyle[
    spaceabove=12pt, spacebelow=15pt,
    headfont=\normalfont\bfseries,
    notefont=\bfseries,
    bodyfont=\normalfont,
    notebraces={}{},
    qed=$\blacksquare$
]{proof}

\declaretheoremstyle[
    spacebelow=8pt,
    headfont=\normalfont\itshape,
    notefont=\itshape,
    bodyfont=\normalfont,
    notebraces={}{},
    qed=$\square$
]{subproof}

\declaretheorem[name=, style=para, numberwithin=section]{para}
\declaretheorem[name=Definition, style=def, numberlike=para]{definition}
\declaretheorem[name=Notation, style=ref, numbered=no]{notation}
\declaretheorem[name=References, style=ref, numbered=no]{references}
\declaretheorem[name=Example, style=def, numberlike=para, qed=$\blacklozenge$]{example}
\declaretheorem[name=Examples, style=def, numberlike=para, qed=$\blacklozenge$]{examples}
\declaretheorem[name=Remark, style=def, numberlike=para]{remark}
\declaretheorem[name=Remarks, style=def, numberlike=para]{remarks}

\declaretheorem[name=Proposition, style=thm, numberlike=para]{proposition}
\declaretheorem[name=Lemma, style=thm, numberlike=para]{lemma}
\declaretheorem[name=Corollary, style=thm, numberlike=para]{corollary}
\declaretheorem[name=Statement, style=stmnt, numberwithin=para]{statement}
\declaretheorem[name=Statement, style=nstmnt, numbered=no]{nstatement}
\declaretheorem[name=Proof, style=proof, numbered=no]{Proof}
\declaretheorem[name=Proof, style=subproof, numbered=no, preheadhook=\begingroup, postfoothook=\endgroup]{subproof}
\declaretheorem[name=Theorem, style=thm_name, numberlike=para]{ntheorem}

\makeatletter
\preto\notation{
  \patchcmd\cref@thmnoarg
    {\trivlist}
    {\list{}{\leftmargin\parindent}}
    {}{}
  \patchcmd\cref@thmoptarg
    {\trivlist}
    {\list{}{\leftmargin\parindent}}
    {}{}
  \patchcmd\thmt@original@endnotation{\endtrivlist}{\endlist}{}{}
}
\makeatother

\makeatletter
\preto\statement{
  \patchcmd\cref@thmnoarg
    {\trivlist}
    {\list{}{\leftmargin\parindent}}
    {}{}
  \patchcmd\cref@thmoptarg
    {\trivlist}
    {\list{}{\leftmargin\parindent}}
    {}{}
  \patchcmd\thmt@original@endstatement{\endtrivlist}{\endlist}{}{}
}
\makeatother

\makeatletter
\preto\nstatement{
  \patchcmd\cref@thmnoarg
    {\trivlist}
    {\list{}{\leftmargin\parindent}}
    {}{}
  \patchcmd\cref@thmoptarg
    {\trivlist}
    {\list{}{\leftmargin\parindent}}
    {}{}
  \patchcmd\thmt@original@endnstatement{\endtrivlist}{\endlist}{}{}
}
\makeatother

\renewcommand{\implies}{\Rightarrow}

\renewcommand{\iff}{\Leftrightarrow}
\renewcommand{\emptyset}{\varnothing}
\renewcommand{\subset}{\subseteq}
\renewcommand{\supset}{\supseteq}

\renewcommand{\epsilon}{\varepsilon}
\renewcommand{\Re}{\operatorname{Re}}
\renewcommand{\Im}{\operatorname{Im}}
\let\originalleft\left
\let\originalright\right
\renewcommand{\left}{\mathopen{}\mathclose\bgroup\originalleft}
\renewcommand{\right}{\aftergroup\egroup\originalright}
\renewcommand{\Gamma}{\varGamma}
\renewcommand{\Delta}{\varDelta}
\renewcommand{\Theta}{\varTheta}
\renewcommand{\Lambda}{\varLambda}
\renewcommand{\Xi}{\varXi}
\renewcommand{\Sigma}{\varSigma}
\renewcommand{\Phi}{\varPhi}
\renewcommand{\Psi}{\varPsi}
\renewcommand{\Omega}{\varOmega}
\DeclareTextFontCommand{\bemph}{\bfseries\em}

\newcommand{\pp}[1]{\textsl{#1}}
\newcommand{\ce}{\coloneqq}
\newcommand{\ec}{\eqqcolon}
\newcommand{\longto}{\varlongrightarrow}
\newcommand{\longmto}{\varlongmapsto}
\newcommand{\ol}[1]{\overline{#1}}
\newcommand{\wt}[1]{\widetilde{#1}}
\newcommand{\wh}[1]{\widehat{#1}}
\newcommand{\N}{\mathbb{N}}

\newcommand{\R}{\mathbb{R}}
\newcommand{\C}{\mathbb{C}}
\newcommand{\K}{\mathbb{K}}

\newcommand{\I}{\mathbb{I}}

\newcommand{\MB}{\mathcal{B}}

\newcommand{\MD}{\mathcal{D}}
\newcommand{\MF}{\mathcal{F}}
\newcommand{\MG}{\mathcal{G}}
\newcommand{\MH}{\mathcal{H}}
\newcommand{\MK}{\mathcal{K}}

\newcommand{\MO}{\mathcal{O}}
\newcommand{\MP}{\mathcal{P}}

\newcommand{\MR}{\mathcal{R}}
\newcommand{\MS}{\mathcal{S}}

\newcommand{\MFA}{\mathfrak{A}}
\newcommand{\MFB}{\mathfrak{B}}
\newcommand{\MFC}{\mathfrak{C}}

\newcommand{\MFJ}{\mathfrak{J}}

\newcommand{\MFM}{\mathfrak{M}}
\newcommand{\MFN}{\mathfrak{N}}
\newcommand{\MFP}{\mathfrak{P}}

\newcommand{\MFS}{\mathfrak{S}}
\newcommand{\MFT}{\mathfrak{T}}

\newcommand{\MFq}{\mathfrak{q}}

\newcommand{\MFt}{\mathfrak{t}}

\newcommand{\MSP}{\mathscr{P}}
\newcommand{\MST}{\mathscr{T}}

\newcommand{\ii}{\mathrm{i}}
\newcommand{\ee}{\mathrm{e}}

\newcommand{\Mat}{\mathrm{Mat}}
\newcommand{\GL}{\mathrm{GL}}

\newcommand{\BO}{\mathscr{B}}
\newcommand{\KO}{\mathscr{K}}
\newcommand{\NO}{\BO_1}
\newcommand{\HS}{\BO_2}
\newcommand{\PO}{\mathscr{P}}
\newcommand{\UO}{\mathscr{U}}
\newcommand{\DM}{\mathscr{S}}

\newcommand{\cdual}[1]{#1^{\ast}}
\newcommand{\pdual}[1]{#1_{\ast}}

\newcommand{\comm}[1]{#1^{\prime}}
\newcommand{\bicomm}[1]{#1^{\prime\prime}}
\newcommand{\SE}[1]{#1^{\mathrm{sa}}}
\newcommand{\PE}[1]{#1^{+}}
\renewcommand{\AE}[1]{#1^{(\eta)}}

\newcommand{\NF}[1]{\PE{\pdual{#1}}}
\newcommand{\NPC}{\mathcal{P}}
\renewcommand{\SS}{\Sigma}
\newcommand{\NS}{\pdual{\SS}}

\newcommand{\id}{\mathrm{Id}}
\newcommand{\idm}{\mathbbm{1}}
\newcommand{\1}{\mathbf{1}}
\newcommand{\diff}{\mathop{}\!\mathrm{d}}

\newcommand{\dlim}[2]{\operatorname*{\text{$#1$}-lim}_{#2}\,}
\newcommand{\expan}{\mathrm{E}}
\newcommand{\sbullet}[1][.75]{\mathbin{\vcenter{\hbox{\scalebox{#1}{$\bullet$}}}}}

\newcommand{\nsupp}{\mathbf{n}}
\newcommand{\lsupp}{\mathbf{l}}
\newcommand{\rsupp}{\mathbf{r}}
\newcommand{\ssupp}{\mathbf{s}}

\newcommand{\EV}{\mathbb{E}}

\newcommand{\inv}[1]{#1^{-1}}
\newcommand{\tran}[1]{#1^{\mathrm{T}}}

\newcommand{\dadj}[1]{#1^{\ast\ast}}

\newcommand{\set}[2][]{#1\{{#2}#1\}}

\newcommand{\abs}[1]{\vert #1 \vert}

\newcommand{\vabs}[2][]{#1\vert{#2}#1\vert}
\newcommand{\norm}[1]{\Vert #1 \Vert}

\newcommand{\vnorm}[2][]{#1\Vert{#2}#1\Vert}
\newcommand{\braket}[1]{\langle{#1}\rangle}

\newcommand{\vbraket}[2][]{#1\langle{#2}#1\rangle}
\newcommand{\sumbra}[2]{#1(\,{#2}#1)}
\newcommand{\od}[3][]{\frac{\diff^{#1} #2}{\diff {#3}^{#1}}}

\newcommand{\odbound}[4][]{\frac{\diff^{#1} #2}{\diff {#3}^{#1}}_{\upharpoonright \, #4}}
\newcommand{\pmat}[1]{\ensuremath{\begin{pmatrix} #1 \end{pmatrix}}}
\newcommand{\bmat}[1]{\ensuremath{\begin{bmatrix} #1 \end{bmatrix}}}

\newcommand{\tand}{\quad\text{and}\quad}
\newcommand{\tas}{\quad\text{as}\quad}
\newcommand{\tfor}{\quad\text{for}\quad}

\newcommand{\tor}{\quad\text{or}\quad}
\newcommand{\twhere}{\quad\text{where}\quad}
\newcommand{\twith}{\quad\text{with}\quad}
\newcommand{\com}{{\ , \quad}}
\newcommand{\ie}{\textit{i.e.}}
\newcommand{\cf}{\textit{cf.}}
\newcommand{\eg}{\textit{e.g.}}
\newcommand{\tAd}{\textit{Ad}}
\newcommand{\ndot}{\norm{\cdot}}
\newcommand{\bdot}{\braket{\cdot, \cdot}}
\newcommand{\mop}{\mathrm{op}}
\DeclareMathOperator{\dom}{dom}
\DeclareMathOperator{\ran}{ran}
\DeclareMathOperator{\gr}{gr}
\DeclareMathOperator{\clos}{clos}
\DeclareMathOperator{\supp}{supp}
\DeclareMathOperator{\lin}{lin}
\DeclareMathOperator{\Hom}{Hom}
\DeclareMathOperator{\Aut}{Aut}
\DeclareMathOperator{\Eig}{Eig}
\DeclareMathOperator{\tr}{tr}

\begin{document}

\frontmatter

\thispagestyle{empty}

\vspace*{\fill}

\begin{center}
    {\LARGE\bfseries Monotonicity of the Relative Entropy and the Two-sided Bogoliubov Inequality in von Neumann Algebras\par}
    \vspace{0.5cm}
    {\Large\scshape Benedikt M. Reible\footnote{\emph{Current address:} Institut für Mathematik, Freie Universität Berlin}\textsuperscript{,}\footnote{\emph{Electronic mail:} \href{mailto:benedikt.reible@fu-berlin.de}{benedikt.reible@fu-berlin.de}}\par}
    \vspace{0.5cm}
    {\large\itshape Institut für Theoretische Physik, Universität Leipzig\par}
    \vspace{0.5cm}
    {\normalsize June, 2024\par}
\end{center}

\vspace{2cm}

\begin{abstract}
    This text studies, on the one hand, certain monotonicity properties of the Araki-Uhlmann relative entropy and, on the other hand, unbounded perturbation theory of KMS-states which facilitates a proof of the two-sided Bogoliubov inequality in general von Neumann algebras. After introducing the necessary background from the theory of operator algebras and Tomita-Takesaki modular theory, the relative entropy functional is defined and its basic properties are studied. In particular, a full and detailed proof of Uhlmann's important monotonicity theorem for the relative entropy is provided. This theorem will then be used to derive a number of monotonicity inequalities for the relative entropy of normal functionals induced by vectors of the form $V \Omega, V \Phi \in \MH$, where $V \in \BO(\MH)$ is a suitable transformation. After that, an introduction to perturbation theory in von Neumann algebras is given, with an emphasis on unbounded perturbations of KMS-states following the framework of Derezi\'{n}ski-Jak\v{s}i\'{c}-Pillet. This mathematical apparatus will then be used to extend the two-sided Bogoliubov inequality for the relative free energy, which was very recently proved for quantum-mechanical systems, to arbitrary von Neumann algebras.
\end{abstract}

\vspace*{\fill}

\newpage
\tableofcontents
\addtocontents{toc}{\protect{\pdfbookmark{\contentsname}{toc}}}

\chapter{Notation}
\markboth{Notation}{Notation}

The end of a proof is indicated by the symbol $\blacksquare$, the end of a \enquote{sub-proof}, that is, of a proof of a statement contained in another proof, is indicated by the symbol $\square$, and the end of an example is marked by the symbol $\blacklozenge$. The phrase \enquote{iff} is an abbreviation of the logical connective \enquote{if and only if}.

In the following, a compilation of mathematical symbols which are frequently used in the text is provided for convenience of the reader since these symbols are not always used consistently in the existing literature.

\vspace{0.5cm}

{\def\arraystretch{1.1}
\begin{longtable}{@{}p{2.5cm}p{12cm}}
    $\K$ & Field of real numbers $\R$ or complex numbers $\C$ \\
    $\MH, \MK, \dotsc$ & Complex Hilbert spaces with inner product $\bdot$ and norm $\ndot$ \\
    $\MFA, \MFB, \dotsc$ & $C^\ast$-algebras \\
    $\MFM, \MFN, \dotsc$ & von Neumann algebras \\
    & \\
    $\BO(\MH)$ & Space of bounded linear operators \\
    $\NO(\MH)$ & Space of trace-class operators \\
    $\HS(\MH)$ & Space of Hilbert-Schmidt operators \\
    $\KO(\MH)$ & Space of compact operators \\
    $\UO(\MH)$ & Group of unitary operators \\
    $\PO(\MH)$ & Set of projection operators \\
    $\DM(\MH)$ & Set of density matrices \\
    & \\
    $\SE{\MFA}$ & Set of self-adjoint elements in $\MFA$ \\
    $\PE{\MFA}$ & Set of positive elements in $\MFA$ \\
    $\comm{\MFA}$ & Commutant of the set $\MFA \subset \BO(\MH)$ \\
    $\GL(\MFA)$ & Group of invertible elements of $\MFA$ \\
    $\Aut(\MFA)$ & Group of $\ast$-automorphisms of $\MFA$ \\
    $\Mat(n; \MFA)$ & $C^\ast$-algebra of $(n \times n)$-matrices with entries in $\MFA$ \\
    $\PO(\MFM)$ & Set of projections in $\MFM$ \\
    $\AE{\MFM}$ & Set of closed operators affiliated with $\MFM$ \\
    $\MFM^\tau$ & Set of $\tau$-entire elements of a $W^\ast$-dynamical system $(\MFM, \tau)$ \\
    $\pdual{\MFM}$ & Predual space of $\MFM$ \\
    $\NF{\MFM}$ & Space of normal functionals on $\MFM$ \\
    $\SS(\MFA)$ & Space of states on $\MFA$ \\
    $\NS(\MFM)$ & Space of normal states on $\MFM$ \\
    & \\
    $\clos_\MFT(U)$ & Closure of the set $U \subset X$ in the topological space $(X, \MFT)$ \\
    $\lin_\K(A)$ & Linear hull of the set $A \subset E$ in the $\K$-vector space $E$ \\
    $\MFA \MK$ & $\lin_\K \set{A \xi \, : \, A \in \MFA, \, \xi \in \MK}$ for $\MFA \subset \BO(\MH)$, $\MK \subset \MH$ \\
    $[\MFA \MK]$ & $\clos_{\ndot}(\MFA \MK)$ \\
    $\MB(X)$ & Borel $\sigma$-algebra of the topological space $X$ \\
    $\cdual{E}$ & Continuous dual space of the topological vector space $(E, \MFT)$ \\
    & \\
    $\dom(T)$ & Domain of the operator $T$ \\
    $\ran(T)$ & Range of the operator $T$ \\
    $\ker(T)$ & Kernel of the operator $T$ \\
    $\Eig(T, \lambda)$ & Eigenspace of $T$ corresponding to the eigenvalue $\lambda \in \K$ \\
    $\ssupp(T)$ & Support projection of the operator $T$ \\
    $E_T$ & Spectral measure $\MB(\sigma(T)) \owns A \longmto E_T(A) \in \PO(\MH)$ of the operator $T$ \\
    $\mu_{\xi, \eta}^T$ & Complex measure $\MB(\sigma(T)) \owns A \longmto \braket{\xi, E_T(A) \eta} \in \C$ for $\xi, \eta \in \MH$ \\
    $\mu_{\xi}^T$ & Positive measure $\mu_{\xi, \xi}^T$ for $\xi \in \MH$ \\
    & \\
    $\id_X$ & Identity mapping $\id_X : X \longto X$, $x \longmto x$, on the set $X$ \\
    $\1_X$ & Constant function $\1_X : X \longto \K$, $x \longmto 1$, on the set $X$ \\
    $\idm_\MFA$ & Unit element of the algebra $\MFA$ \\
    $[U]$ & Orthogonal projection onto the closure of $U \subset \MH$ \\
    $M_\xi$ & Multiplication operator with the element $\xi \in \MH$ on the Hilbert space $\MH$ \\
    & \\
    $\ndot_\infty$ & Supremum norm on $C^0(X; \K)$ \\
    $\ndot_\mop$ & Operator norm on $\BO(\MH)$ \\
    $\bdot_\mathrm{HS}$ & Hilbert-Schmidt inner product on $\HS(\MH)$ \\
    $\tr$ & Trace functional on $\NO(\MH)$ \\
    & \\
    $\dlim{\MH}{n \to \infty}$ & Limit in the norm topology of $\MH$ \\
    $\dlim{w}{n \to \infty}$ & Limit in the weak topology of $\MH$ \\
    $\dlim{so}{n \to \infty}$ & Limit in the strong operator topology of a von Neumann algebra $\MFM$ \\
    & \\
    $\omega_\mu$ & Functional $f \longmto \int_X f \diff \mu$ induced by the measure $\mu$ \\
    $\omega_\rho$ & Functional $A \longmto \tr(\rho A)$ induced by $\rho \in \NO(\MH)$ \\
    $\omega_\xi$ & Vector functional $A \longmto \braket{\xi, A \xi}$ induced by $\xi \in \MH$ \\
    $\comm{\omega_\xi}$ & Vector functional $A \longmto \braket{\xi, A \xi}$ on the commutant $\comm{\MFA} \subset \BO(\MH)$ \\
    $\xi_\omega$ & Vector representative of $\omega \in \NF{\MFM}$ in the natural positive cone $\NPC$ \\
    $\ssupp(\varphi), \ssupp_\MFM(\varphi), \ssupp_\varphi$ & Support projection of the normal functional $\varphi \in \NF{\MFM}$ \\
    $\ssupp_{\comm{\MFM}}(\varphi), \comm{\ssupp_\varphi}$ & Support projection of the normal functional $\comm{\varphi} \in \NF{(\comm{\MFM})}$ \\
    $\supp(\varphi)$ & Range $\ran\bigl(\ssupp(\varphi)\bigr)$ of the support projection of $\varphi \in \NF{\MFM}$ \\
    & \\
    $S_\omega$, $S_\Omega$ & Tomita operator with respect to the functional $\omega = \omega_\Omega$ \\
    $\Delta_\omega$, $\Delta_\Omega$ & Modular operator with respect to $\omega = \omega_\Omega$ \\
    $\sigma^\omega$, $\sigma^\Omega$ & Modular automorphism group on $\MFM$ with respect to $\omega = \omega_\Omega$ \\
    $(\MFM, \MH, J, \NPC)$ & Standard form representation of $\MFM$ \\
    $j(A)$ & Conjugation of $A \in \MFM$ with $J$: $j(A) = J A J$ \\
    $S_{\psi, \varphi}$, $S_{\Psi, \Phi}$ & Relative Tomita operator with respect to $\psi = \omega_\Psi$ and $\varphi = \omega_\Phi$ \\
    $\Delta_{\psi, \varphi}$, $\Delta_{\Psi, \Phi}$ & Relative modular operator with respect to $\psi = \omega_\Psi$ and $\varphi = \omega_\Phi$ \\
    $D(\MH, \psi)$ & Lineal of the semi-finite normal weight $\psi$ on $\MFM \subset \BO(\MH)$ \\
    $\Delta(\varphi / \comm{\psi})$ & Spatial derivative of $\varphi \in \NF{\MFM}$ and $\comm{\psi} \in \NF{(\comm{\MFM})}$ \\
    $\MFq_\varphi$ & Quadratic form $\xi \longmto \varphi\bigl(\Theta^{\comm{\psi}}(\xi)\bigr)$ associated with $\Delta(\varphi / \comm{\psi})$ on $\MH$ \\
    & \\
    $D(P, Q)$ & Kullback-Leibler divergence of the probability measures $P, Q$ \\
    $S(\rho, \sigma)$ & Umegaki relative entropy of the density matrices $\rho, \sigma \in \DM(\MH)$ \\
    $\MS_\MFM^\mathrm{std}(\omega, \varphi)$ & Araki-Uhlmann relative entropy of $\omega, \varphi \in \NF{\MFM}$ for $\MFM$ in standard form \\
    $\MS_\MFM^\mathrm{spa}(\omega, \varphi)$ & Araki-Uhlmann relative entropy for arbitrary representation of $\MFM$ \\
    $\MR_\MFM(\Omega, \Phi)$ & Araki-Uhlmann relative entropy $\MS_\MFM^\mathrm{spa}(\omega_\Omega, \omega_\Phi)$ for $\Omega, \Phi \in \MH$ \\
    & \\
    $\beta$ & Inverse temperature \\
    $\tau^V$ & Perturbation of the $W^\ast$-dynamics $\tau$ with the element $V$ \\
    $\expan_\tau^V$ & Araki-Dyson expansional \\
    $L_V$ & Perturbation of the Liouvillian $L$ with $V$ \\
    $\Omega_V$ & Perturbation of the vector $\Omega \in \NPC$ with $V$ \\
    $\omega_V$ & Perturbed KMS-state corresponding to $\Omega_V$ \\
    $\MF(\omega_V, \omega_0)$ & Relative free energy of perturbed and unperturbed KMS-states
\end{longtable}
}

\mainmatter

\chapter{Introduction}

In recent years, a substantial interest in information-theoretic aspects of quantum field theory has developed. An indispensible tool in corresponding studies is the \emph{relative entropy functional} which was already defined in 1977 in its most general form for arbitrary von Neumann algebras. Lately, it has experienced a renaissance in mathematical physics, especially in algebraic quantum field theory, with many new applications being discovered and investigated. To motivate the specific problems studied in the present text, the use of the relative entropy in mathematical physics is reviewed briefly in \cref{sec:introduction_history}. Following this, \cref{sec:introduction_outline} provides a concise outline of the text.

\section{The Relative Entropy in Mathematical Physics}\label{sec:introduction_history}

\begin{para}[The Kullback-Leibler divergence]\label{para:inroduction_KL}
    The notion of relative entropy was first introduced in mathematical statistics by \textsc{S. Kullback} and \textsc{R. A. Leibler} in 1951 \cite{KullbackLeibler51}. In modern terminology, it can be defined as follows: let $(\Omega, \Sigma)$ be a measurable space and $P, Q$ be two probability measures on $\Sigma$ such that $P$ is absolutely continuous with respect to $Q$. In this case, the Radon-Nikodým derivative $\od{P}{Q}$ exists \cite[Thm. 4.2.2]{Cohn13}, and the \bemph{classical relative entropy} or \bemph{Kullback-Leibler divergence} between $P$ and $Q$ is defined to be \cite[p. 184]{Jaksic19}
    \begin{equation}\label{eq:introduction_KL}
        D(P, Q) \ce \int_\Omega \log\left(\od{P}{Q}\right) \diff P = \int_\Omega \od{P}{Q} \, \log\left(\od{P}{Q}\right) \diff Q \ .
    \end{equation}
    If $P$ is not absolutely continuous with respect to $Q$, one can extend the above definition by setting $D(P, Q) \ce + \infty$. \eqref{eq:introduction_KL} has a number of interesting properties, for example \cite[Prop. 4.1, 4.4 \& 4.8]{Jaksic19}:
    \begin{enumerate}[num]
        \item \label{enu:introduction_KLposdef} $D(P, Q) \ge 0$ and $D(P, Q) = 0$ if and only if $P = Q$;
        \item \label{enu:introduction_KLmonotonicity} $D\bigl(\alpha(P), \alpha(Q)\bigr) \le D(P, Q)$ for every stochastic mapping $\alpha$;
        \item \label{enu:introduction_KLconvex} $(P, Q) \longmto D(P, Q)$ is jointly convex.
    \end{enumerate}

    Based on \ref{enu:introduction_KLposdef}, the relative entropy $D(P, Q)$ can be interpreted as a measure for the \emph{distinguishability} of the probability distributions $P$ and $Q$ \cite[pp. 8 \& 84]{OP04}. In fact, it can be shown that if one draws $N$ samples according to $Q$, then the probability that they will look as if they were drawn from $P$ is proportional to $\exp(- N D(P, Q))$ for large $N$ \cite[p. 200]{Vedral02}. In light of this interpretation, \ref{enu:introduction_KLmonotonicity} implies that the probability distributions $P$ and $Q$ can only become less distinguishable under stochastic evolution \cite[p. 198]{Vedral02}, and property \ref{enu:introduction_KLconvex} asserts that mixing probability distributions decreases distinguishability \cite[p. 203]{Vedral02}.

    The notion of distinguishability is of fundamental importance for information theory because the degree to which one can distinguish between different physical states of a system determines the amount of information that can be encoded and manipulated \cite[p. 197]{Vedral02}. This is the main reason why generalizations of the relative entropy \eqref{eq:introduction_KL} to non-commutative probability spaces have become so important in quantum information theory and quantum field theory, as will be outlined in the following.
\end{para}

\begin{para}[The Umegaki relative entropy]\label{para:introduction_Umegaki}
    Let $\MH$ be a Hilbert space. The Kullback-Leibler divergence \eqref{eq:introduction_KL} was extended to the von Neumann algebra $\MFM = \BO(\MH)$ by \textsc{H. Umegaki} in 1962 \cite{Umegaki62}: if $\omega$ and $\varphi$ are two normal states on $\MFM$ represented by density matrices $\rho_\omega, \rho_\varphi \in \DM(\MH)$ such that $\ker(\rho_\omega)^\perp \subset \ker(\rho_\varphi)^\perp$, then the \bemph{quantum} or \bemph{Umegaki relative entropy} of $\omega$ with respect to $\varphi$ is defined to be \cite[p. 16]{OP04}
    \begin{equation}\label{eq:introduction_Umegaki}
        S(\rho_\omega, \rho_\varphi) \ce \tr\bigl(\rho_\omega (\log \rho_\omega - \log \rho_\varphi)\bigr) \ .
    \end{equation}
    One can show that if $\rho_\omega$ and $\rho_\varphi$ commute with each other, then \eqref{eq:introduction_Umegaki} reduces to \eqref{eq:introduction_KL} \cite[p. 16]{OP04}. Furthermore, analoga of the properties \ref{enu:introduction_KLposdef} -- \ref{enu:introduction_KLconvex} from \cref{para:inroduction_KL} can be proved for the Umegaki entropy as well \cite[Prop. 1.1 \& Thm 1.4, 1.5]{OP04}. Therefore, one can interpret the quantity $S(\omega, \varphi)$ as a measure for the distinguishability of $\omega$ and $\varphi$ \cite[p. 207]{Vedral02}.

    The relative entropy \eqref{eq:introduction_Umegaki} was introduced to mathematical physics by \textsc{G. Lindblad} in 1973 -- 1975 \cite{Lindblad73, Lindblad74, Lindblad75}. He used this quantity to analyze the quantum-mechanical measurement process, and he showed that $S$ is non-increasing under completely positive maps; see \cite{MHReeb17, Petz03} for recent discussions and refinements of this monotonicity property. Since then, the quantum relative entropy has become an indispensible tool in two areas of mathematical physics: \emph{statistical mechanics} \cite{Wehrl78} and \emph{quantum information theory} \cite{Schumacher00, Vedral02}. In the former, Lindblad's monotonicity result can be seen as a \enquote{generalized $H$-theorem} for open quantum systems \cite[p. 272]{GustafsonSigal20}, and the analogue of \cref{para:inroduction_KL} \ref{enu:introduction_KLposdef} can be used to prove that the canonical Gibbs state is the unique minimizer of the free energy \cite[Prop. 1.10]{OP04}. Recently, there have also been studies which suggest reformulations of the second and third law of thermodynamics using the relative entropy \cite{DowlingFloerchingerHaas20, FloerchingerHaas20, Sagawa23}. In quantum information theory, the Umegaki relative entropy plays a role analogous to that of the Kullback-Leibler divergence in classical statistics \cite[p. 208]{Vedral02}. Among other reasons, this is due to the fact that the relative entropy can be used to quantify the amount of \emph{entanglement} present in a single quantum state \cite[p. 2]{Schumacher00}, \cite[p. 219]{Vedral02}.

    Lately, also an interest in entanglement measures for relativistic quantum field theory has developed, see the overviews \cite{HollandsSanders18, Witten18} and other references given below in \cref{para:introduction_ArakiUhlmann}. Since the Umegaki relative entropy is not available for general quantum field theories \cite[p. 114]{LongoXu18}, a more general definition of a relative entropy functional is needed.
\end{para}

\begin{para}[The Araki-Uhlmann relative entropy]\label{para:introduction_ArakiUhlmann}
    A generalization of the Umegaki relative entropy to arbitrary von Neumann algebras was discovered by \textsc{H. Araki} in 1977, who had previously already made substantial contributions to the mathematical theory of von Neumann algebras (\eg{}, \cite{Araki74a, ArakiWoods68}) and algebraic quantum field theory (\eg{}, \cite{Araki63, Araki64b, Araki64a}). He developed an extension of \emph{Tomita-Takesaki modular theory} and used it to define a notion of relative entropy $\MS_\MFM(\omega, \varphi)$ between normal functionals $\omega, \varphi$ on an arbitrary von Neumann algebra $\MFM$ \cite{Araki76, Araki77}. To some extend, he was motivated by questions of statistical mechanics regarding certain properties of \emph{KMS-states}, where he had already applied methods from modular theory and perturbation theory in operator algebras \cite{Araki74b, Araki75}.
    
    Simultaneously and independently, \textsc{A. Uhlmann} also introduced a notion of relative entropy between positive linear functionals on $\ast$-algebras using an abstract quadratic interpolation theory \cite{Uhlmann77}. Later, it was shown by \textsc{D. Petz} that \textsc{Uhlmann}'s and \textsc{Araki}'s definition actually coincide on von Neumann algebras \cite[Lem. 2]{Petz86a}, \cite[pp. 79 f.]{OP04}, whence this functional, which is the protagonist of this text, is termed the \emph{Araki-Uhlmann relative entropy}; it will be formally introduced in \cref{def:relativeEntropy_stdRelativeEntropy,def:relativeEntropy_spaRelativeEntropy}. In \cite[Prop. 18]{Uhlmann77}, \textsc{Uhlmann} also proved a very general \emph{monotonicity theorem} for the relative entropy which generalizes Lindblad's result mentioned in \cref{para:introduction_Umegaki}. Namely, he showed that for a unital normal Schwarz mapping $\alpha : \MFM_1 \longto \MFM_2$, the following inequality holds true \cite[Cor. 5.12 (iii)]{OP04}
    \begin{equation}\label{eq:introduction_UhlmannMonotonicity}
        \MS_{\MFM_1}(\omega \circ \alpha, \varphi \circ \alpha) \le \MS_{\MFM_2}(\omega, \varphi) \ .
    \end{equation}
    This property of the relative entropy is highly important for many reasons, \eg{}, the fact that analoga of \ref{enu:introduction_KLposdef} and \ref{enu:introduction_KLconvex} from \cref{para:inroduction_KL} can be proved with the help of \eqref{eq:introduction_UhlmannMonotonicity}; \cref{sec:relativeEntropy_monotonicity} of this text is devoted to this topic. Further properties of the Araki-Uhlmann relative entropy were studied in the 80s and 90s by \textsc{M. J. Donald} \cite{Donald86, Donald87b}, \textsc{H. Kosaki} \cite{Kosaki86}, and \textsc{D. Petz} \cite{Petz86a, Petz86b, Petz88, Petz91, Petz92, Petz94}, to name just a selected few.

    After the introduction of the relative entropy on general von Neumann algebras, it was primarily applied in statistical mechanics to study KMS-states \cite{ArakiSewell77, BR2}. The interest of quantum field theorists in the relative entropy started only a couple of years ago. To the best knowledge of the author, this can be traced back to a paper of \textsc{H. Casini} from 2008 \cite{Casini08} who used the Umegaki relative entropy \eqref{eq:introduction_Umegaki} to give a proof of the so-called \emph{Bekenstein bound} which provides an upper limit for the entropy of a system in terms of its energy and size. Ten years later, \textsc{R. Longo} and \textsc{F. Xu} \cite{LongoXu18} pointed out that \textsc{Casini}'s argument breaks down for general quantum field theories, and they gave a more rigorous proof of the Bekenstein bound using the Araki-Uhlmann relative entropy. Since then, this quantity has been applied to various problems in (algebraic) quantum field theory (in curved spacetimes). To begin with, \textsc{Casini} \emph{et al.} \cite{CGP19} and \textsc{Longo} \cite{Longo19} independently computed, in the setting of a free scalar field, the relative entropy  between the vacuum state and a coherent excitation thereof; generalizations of these computations were considered afterwards by various authors \cite{Hollands20, BCV22, BCM23, Galanda23}. Furthermore, a connection between the relative entropy and \emph{quantum energy inequalities} was discovered by \textsc{F. Ceyhan} and \textsc{T. Faulkner} \cite{CeyhanFaulkner2020} (see also \cite{CLRR22, MTW22, Wall17}), and more general investigations of \emph{relativistic quantum information theory} were conducted by \textsc{F. Hiai}, \textsc{S. Hollands}, \textsc{R. Longo}, and others \cite{CLR20, FaulknerHollands2022, FHSW2022, FMP23, Hiai18a, Hiai19, Hiai21a, Hollands21, Hollands23, HollandsIshibashi19, HollandsSanders18, Longo20}.
\end{para}

\section{Outline of the Text}\label{sec:introduction_outline}

As illustrated above, the study of the relative entropy in von Neumann algebras is a highly active field of research. It is the goal of this text to contribute to this field by providing a detailed introduction to its mathematical foundations, and by discussing certain properties of the relative entropy related to monotonicity and perturbation theory of KMS-states. To this end, \cref{ch:operatorAlgebras} gives an overview of the theory of $C^\ast$-algebras and von Neumann algebras, and \cref{ch:vonNeumann} contains a fairly detailed discussion of Tomita-Takesaki modular theory.

The Araki-Uhlmann relative entropy is introduced in \cref{ch:relativeEntropy}, first for von Neumann algebras in standard form, and second for arbitrary representations. The main part of this chapter is devoted to the monotonicity theorem of \textsc{A. Uhlmann}. Fully detailed proofs of this result and of all the auxiliary lemmata required to establish it are presented; in this form, this is not readily available in the literature. \cref{sec:relativeEntropy_vectorMonotonicity} contains some original results regarding the following question: given von Neumann algebras $\MFM_1, \MFM_2 \subset \BO(\MH)$ and vectors $\Omega, \Phi \in \MH$, for which Hilbert-space transformations $V \in \BO(\MH)$ can one prove monotonicity inequalities of the form $\MS_{\MFM_1}(V \Omega, V \Phi) \le \MS_{\MFM_2}(\Omega, V)$? Results along these lines will be obtained by applying Uhlmann's general theorem to specific situations (see \cref{pro:relativeEntropy_vectorMonotonicityIsometry,eq:relativeEntropy_vectorMonotonicityIsometry2,pro:relativeEntropy_vectorMonotonicityPartialIsometry,pro:relativeEntropy_vectorMonotonicityUnitary,pro:vectorMonotonicity_contraction}).

In \cref{ch:perturbationTheory}, perturbation theory in operator algebras is discussed. Primarily, the chapter deals with unbounded perturbations of KMS-states by reviewing relatively recent results of \textsc{J. Derezi\'{n}ski}, \textsc{V. Jak\v{s}i\'{c}}, and \textsc{C.-A. Pillet} \cite{DJP03}, and extending them slightly in \cref{subsec:perturbationTheory_unboundedPerturbationKMS} (see \cref{lem:perturbationTheory_unboundedConvergenceLiouvillian} \ref{enu:perturbationTheory_unbdConvergenceSum}, \cref{thm:perturbationTheory_unboundedPerturbation} \ref{enu:perturbationTheory_unboundedRelModOp2}, and \cref{pro:perturbationTheory_relEntUnboundedPerturbedState}). The final \cref{sec:perturbationTheory_twoSidedBogoliubov} is concerned with the so-called \emph{two-sided Bogoliubov inequality}. This inequality was first discussed in the setting of classical statistical mechanics as a tool to estimate finite-size effects in molecular simulations \cite{DS17}. It was extended to quantum statistical mechanics by the author of this text in his bachelor's thesis; the findings were published in \cite{Reible22} (they will be reviewed very briefly in \cref{subsec:perturbationTheory_twoSidedBogoliubovQM}), and concrete applications of this inequality to different problems were subsequently presented in \cite{Reible23, DS24}. In this work, using the extended framework of unbounded perturbation theory of KMS-states, the results of \cite{Reible22} are generalized to the setting of arbitrary von Neumann algebras. First, a version of the two-sided Bogoliubov inequality is obtained (\cref{pro:perturbationTheory_BogoliubovVN}) which reduces to the known inequality of \cite[Thm. 4.1]{Reible22} on the von Neumann algebra $\BO(\MH)$. Second, the variational bounds of \cite[Sect. 4.1]{Reible22} are extended to certain classes of unbounded perturbations in \cref{subsec:perturbationTheory_variationalBoundsVN} (see \cref{pro:perturbationTheory_unboundedGibbsVP,pro:perturabtionTheory_unboundedDVVP} as well as \cref{cor:perturbationTheory_unboundedVariationalBounds}).

The main part of text will be concluded with \cref{ch:conclusion} which summarizes the results and provides a brief outlook. Finally, \cref{app:topologicalVectorSpaces,app:operators,app:forms} contain some selected definitions and results from the fields of topological vector spaces, (un)bounded linear operators, and quadratic forms on Hilbert spaces which are needed throughout the text.

\chapter{Aspects of the Theory of Operator Algebras}\label{ch:operatorAlgebras}

To commence the study of the Araki-Uhlmann relative entropy, this chapter presents some basic definitions and results from the mathematically rich theory of operator algebras. First, \cref{sec:operatorAlgebras_CStar} introduces \emph{$C^\ast$-algebras} and their \emph{positive elements}. Following that, \cref{sec:operatorAlgebras_positiveMaps} studies properties of \emph{positive mappings} between $C^\ast$-algebras. In \cref{sec:operatorAlgebras_vonNeumann}, the important subclass of \emph{von Neumann algebras} and some of their basic properties are discussed. Finally, \cref{sec:operatorAlgebras_representationsStates} presents the important \emph{GNS-construction}. Due to the vastness of the theory, the limited amount of space available, and the wide range of literature on the topic, only a few selected proofs will be given.

\section{\texorpdfstring{$C^\ast$-Algebras}{C*-Algebras}}\label{sec:operatorAlgebras_CStar}

\subsection{Definitions and Examples}

The properties characterizing $C^\ast$-algebras are, on the one hand, an abstraction of properties displayed by the space $(\BO(\MH), \ndot_\mop)$ of bounded linear operators on a Hilbert space $\MH$. On the other hand, the class of $C^\ast$-algebras forms a subclass of the more general Banach algebras, characterized by additional requirements on the algebra norm. Here, the definition of a $C^\ast$-algebra shall be presented directly without a deductive tour through the theory of Banach algebras; for the latter, see \cite{BonsallDuncan73, KadisonRingrose97, Takesaki79}. First, however, some purely algebraic concepts are recalled.

\begin{para}[$\ast$-algebras]\label{para:operatorAlgebras_astAlgebras}
    (\cite[Ch. I]{BonsallDuncan73}, \cite[Sect. 2.1]{Schmüdgen20})
    A \bemph{$\ast$-algebra} is a $\K$-algebra together with an \emph{algebra involution} $\MFA \longto \MFA$, $A \longmto A^\ast$, that is, a mapping satisfying for all $A, B \in \MFA$ and $\lambda, \mu \in \K$:
    \begin{equation*}
        (A^\ast)^\ast = A \com (\lambda \cdot A + \mu \cdot B)^\ast = \ol{\lambda} \cdot A^\ast + \ol{\mu} \cdot B^\ast \ , \tand (AB)^\ast = B^\ast A^\ast \ .
    \end{equation*}
    A subset $\MFJ \subset \MFA$ is called \emph{self-adjoint} iff $A^\ast \in \MFJ$ for every $A \in \MFJ$, and a self-adjoint subalgebra (ideal) of $\MFA$ is called a \emph{$\ast$-subalgebra} (\emph{$\ast$-ideal}). An element $A \in \MFA$ is called \bemph{self-adjoint} iff $A^\ast = A$, and the set of all self-adjoint elements is denoted by $\SE{\MFA}$. This set is a real linear subspace of $\MFA$ such that $\MFA = \SE{\MFA} \oplus_\R \ii \, \SE{\MFA}$ \cite[Lem. I.12.3]{BonsallDuncan73}. A linear mapping $f : \MFA \longto \MFB$ between two $\ast$-algebras is called a \bemph{$\ast$-homomorphism} iff $f(AB) = f(A) f(B)$ and $f(A^\ast) = f(A)^\ast$ for all $A, B \in \MFA$. Finally, an element $A \in \MFA$ of a unital $\ast$-algebra $(\MFA, \idm)$ is called \bemph{invertible} iff there exists an element $\inv{A} \in \MFA$ such that $A \inv{A} = \inv{A} A = \idm$. The set of all invertible elements of $\MFA$ is denoted by $\GL(\MFA)$.
\end{para}

\begin{definition}[$C^\ast$-algebra]\label{def:operatorAlgebras_CStar}
    Let $\MFA$ be a $\ast$-algebra. If there exists a norm $\ndot : \MFA \longto [0, + \infty)$ such that (1) $\MFA$ is $\ndot$-complete, (2) $\ndot$ is \emph{sub-multiplicative}, that is,
    \begin{equation}\label{eq:operatorAlgebras_algebraNorm}
        \forall A, B \in \MFA \ : \ \norm{AB} \le \norm{A} \, \norm{B} \ ,
    \end{equation}
    and (3) $\ndot$ satisfies the the so-called \emph{$C^\ast$-property}
    \begin{equation}\label{eq:operatorAlgebras_CStarProperty}
        \forall A \in \MFA \ : \ \norm{A^\ast A} = \norm{A}^2 \ ,
    \end{equation}
    then $(\MFA, \ndot)$ is called an \bemph{abstract} \bemph{$C^\ast$-algebra},\footnotemark
    \footnotetext{The terminology \enquote{abstract} will be apparent after \emph{concrete $C^\ast$-algebras} are introduced in \cref{exa:operatorAlgebras_CStar} \ref{enu:operatorAlgebras_exaConcreteCstar}. The name \enquote{$C^\ast$-algebra} was coined by \textsc{I. E. Segal} in 1947 \cite{Segal47}.}
    and $\ndot$ is called a \emph{$C^\ast$-algebra norm}. If there is no risk of confusion, such an algebra will be called $C^\ast$-algebra and denoted by $\MFA$ for short. The algebra norm $\ndot$ defines a metric topology on $\MFA$ which is referred to as the \bemph{uniform} or \bemph{norm topology} $\MFT_\mathrm{norm}$. Finally, a \bemph{$C^\ast$-subalgebra} of $\MFA$ is defined to be a uniformly closed $\ast$-subalgebra $\MFB \subset \MFA$.
\end{definition}

\begin{examples}\label{exa:operatorAlgebras_CStar}
    \leavevmode
    \begin{enumerate}[env]
        \item The complex numbers $\C$ form a commutative $C^\ast$-algebra with involution given by complex conjugation $z = a + \ii b \longmto \ol{z} \ce a - \ii b$ and norm given by the absolute value.
        
        \item \label{enu:operatorAlgebras_exaCCStar} Let $X$ be a compact Hausdorff space and $C^0(X) \ce C^0(X; \C)$ be the space of continuous functions $f : X \longto \C$. With multiplication defined pointwise, $(fg)(x) \ce f(x) g(x)$ for $f, g \in C^0(X)$ and $x \in X$, involution given by $f^\ast(x) \ce \ol{f(x)}$, and norm defined to be
        \begin{equation*}
            \norm{f}_\infty \ce \sup_{x \in X} \abs{f(x)} \ ,
        \end{equation*}
        this space becomes a commutative $C^\ast$-algebra \cite[Exa. 2.1.4]{BR1}. It is unital with unit given by the constant function $\1_X : X \longto \K$, $x \longmto 1$.

        \item \label{enu:operatorAlgebras_exaLInfCStar} Similarly to the previous example, for a measure space $(\Omega, \Sigma, \mu)$ it holds that the essentially bounded functions $L^\infty(\Omega, \mu)$ with pointwise multiplication, the involution from \ref{enu:operatorAlgebras_exaCCStar}, and the essential supremum norm $\ndot_{L^\infty}$ form a unital commutative $C^\ast$-algebra.

        \item \label{enu:operatorAlgebras_exaMatCStar} Let $n \in \N$. The space $\Mat(n; \K)$ consisting of $(n \times n)$-matrices with entries in $\K$ is a unital $\K$-algebra with multiplication given by the matrix product. If $\K = \R$, one can define an involution in terms of the transposition $A \longmto \tran{A}$, and if $\K = \C$, an involution is given by $A \longmto \tran{(\ol{A})}$. There exist various sub-multiplicative, complete norms on $\Mat(n; \K)$, \eg{}, the Frobenius norm \cite[p. 211]{Werner22}, but only the operator norm satisfies the $C^\ast$-property \cite[Cor. 2.2.6]{BR1}.

        \item \label{enu:operatorAlgebras_exaBOCStar} Let $(\MH, \bdot)$ be a Hilbert space. The space $\BO(\MH)$ of bounded linear operators $T : \MH \longto \MH$ is a unital $\ast$-algebra with multiplication $(T, S) \longmto T \circ S \equiv TS$ and involution given by the operation $T \longmto T^\ast$ of taking adjoints. Furthermore, the \bemph{operator norm} \cite[p. 50]{Werner18}
        \begin{equation*}
            \norm{T}_\mop \ce \norm{T}_{\BO(\MH)} \ce \sup_{\xi \in \MH \setminus \{0\}} \frac{\norm{T \xi}}{\norm{\xi}} = \sup_{\norm{\xi} \le 1} \norm{T \xi} = \sup_{\norm{\xi} = 1} \norm{T \xi}
        \end{equation*}
        is sub-multiplicative because for $T, S \in \BO(\MH)$ and $\xi \in \MH$, the inequality
        \begin{equation*}
            \norm{T (S \xi)} \le \norm{T}_\mop \norm{S \xi} \le \norm{T}_\mop \norm{S}_\mop \norm{\xi} \ ,
        \end{equation*}
        which follows from boundedness of $T$ and $S$, implies $\norm{T S}_\mop \le \norm{T}_\mop \norm{S}_\mop$ \cite[Lem. II.1.6]{Werner18}. It holds that $(\BO(\MH), \ndot_\mop)$ is complete \cite[Thm. II.1.4]{Werner18}, and the involution satisfies $\norm{T^\ast}_\mop = \norm{T}_\mop$ \cite[Thm. V.5.2]{Werner18}. Finally, for all $T \in \BO(\MH)$ and $\xi \in \MH$, the Cauchy-Schwarz inequality implies
        \begin{align*}
            \norm{T \xi}^2 = \braket{T \xi, T \xi} = \braket{\xi, T^\ast T \xi} \le \norm{\xi} \, \norm{T^\ast T \xi} \ .
        \end{align*}
        With this inequality and \eqref{eq:operatorAlgebras_algebraNorm}, it follows that
        \begin{align*}
            \norm{T}_\mop^2 &= \sup_{\norm{\xi} \le 1} \norm{T \xi}^2 \le \sup_{\norm{\xi} \le 1} \norm{\xi} \, \norm{T^\ast T \xi} \le \norm{T^\ast T}_\mop \le \norm{T^\ast}_\mop \norm{T}_\mop = \norm{T}_\mop^2 \ .
        \end{align*}
        This chain of inequalities shows that there is actually equality, \ie{}, $\norm{T^\ast T}_\mop = \norm{T}_\mop^2$ for all $T \in \BO(\MH)$, hence $(\BO(\MH), \ndot_\mop)$ is a unital $C^\ast$-algebra \cite[Thm. V.5.2 (f)]{Werner18}.

        \item \label{enu:operatorAlgebras_exaKOCStar} Let $\KO(\MH) \subset \BO(\MH)$ be the space of compact operators on $\MH$. It holds that $\KO(\MH)$ is a $\ndot_\mop$-closed subalgebra of $\BO(\MH)$ \cite[Thm. II.3.2]{Werner18} which is not unital if $\dim(\MH) = + \infty$ \cite[p. 506]{Werner18}. Furthermore, $T \in \KO(\MH)$ if and only if $T^\ast \in \KO(\MH)$ \cite[Thm. III.4.4]{Werner18}, hence $\KO(\MH)$ is a normed $\ast$-ideal of $\BO(\MH)$, and thus, in particular, a $C^\ast$-algebra as well.
        
        \item \label{enu:operatorAlgebras_exaNOHSCStar} The space $\NO(\MH)$ of trace-class operators and the space $\HS(\MH)$ of Hilbert-Schmidt operators (\cref{def:operators_NOHS}) are $\ast$-ideals in $\BO(\MH)$, and with respect to the trace norm $\ndot_\mathrm{tr}$, respectively, the Hilbert-Schmidt norm $\ndot_\mathrm{HS}$, they are complete (\cref{pro:operators_NOHS}). However, they are not $C^\ast$-algebras because $\norm{A}_\mathrm{HS}^2 = \norm{A^\ast A}_\mathrm{tr}$ for all $A \in \NO(\MH) \subset \HS(\MH)$.
        
        \item \label{enu:operatorAlgebras_exaConcreteCstar} More generally, it holds that any uniformly closed $\ast$-subalgebra $\MFA \subset \BO(\MH)$ is a $C^\ast$-algebra, referred to as a \bemph{concrete $C^\ast$-algebra}. To see this, let $A \in \MFA$ be arbitrary. Then there exists a sequence $(A_n)_{n \in \N} \subset \MFA$ such that $\norm{A_n - A}_\mop \to 0$ as $n \to + \infty$. It follows from continuity of the multiplication and involution, and the $C^\ast$-property \eqref{eq:operatorAlgebras_CStarProperty} that
        \begin{equation*}
            \norm{A^\ast A}_\mop = \lim_{n \to \infty} \norm{A_n^\ast A_n}_\mop = \lim_{n \to \infty} \norm{A_n}_\mop^2 = \norm{A}_\mop^2 \ . \tag*{\qedhere}
        \end{equation*}
    \end{enumerate}
\end{examples}

\subsection{Positive Elements}

In the following, the class of positive elements in a $C^\ast$-algebra will be introduced. It shall be assumed that all appearing algebras are considered over the field $\K = \C$ of complex numbers. Recall that the \bemph{spectrum} of an element $A \in \MFA$ of a unital $\C$-algebra is given by
\begin{equation*}
    \sigma(A) \ce \set{\lambda \in \C \, : \, A - \lambda \idm \notin \GL(\MFA)} \ .
\end{equation*}

\begin{definition}[Positive element]\label{def:operatorAlgebras_positiveElement}
    Let $\MFA$ be a unital $\ast$-algebra. An element $A \in \MFA$ is called \bemph{positive}, denoted by $A \ge 0$, iff $A$ is self-adjoint and satisfies $\sigma(A) \subset [0, + \infty)$. The set of all positive elements will be denoted by the symbol $\PE{\MFA}$.
\end{definition}

The following theorem characterizes the positive elements in a $C^\ast$-algebra. The proof can be found in \cite[Thm. 2.2.10 \& 2.2.12]{BR1} or \cite[Thm. VIII.3.6]{Conway85}.

\begin{ntheorem}[Characterization of positive elements]\label{thm:operatorAlgebras_characterizationPE}
    Let $\MFA$ be a unital $C^\ast$-algebra. For an element $A \in \MFA$, the following assertions are equivalent:
    \begin{enumerate}[equiv]
        \item $A \in \PE{\MFA}$;
        \item there exists a unique $B \in \PE{\MFA}$ such that $A = B^2$;
        \item there exists $B \in \MFA$ such that $A = B^\ast B$.
    \end{enumerate}
\end{ntheorem}

\begin{para}[Self-dual convex cones]\label{para:operatorAlgebras_cones}
    Recall the following algebraic definitions \cite[p. 38]{SchaeferWolff99}: a \bemph{cone} in an $\R$-vector space $V$ is a subset $C \subset V$ which is closed under multiplication by non-negative scalars, \ie{}, $x \in C, \, \lambda \ge 0 \, \implies \, \lambda x \in C$. The cone is called \bemph{convex} iff $\lambda, \mu \ge 0, \, x, y \in C \, \implies \, \lambda x + \mu y \in C$. For later use, introduce also the following notion \cite[p. 218]{SchaeferWolff99}: a cone $C \subset \MH$ in a Hilbert space $(\MH, \bdot)$ is called \bemph{self-dual} iff $C = \set{y \in \MH \, : \, \text{$\braket{y, x} \ge 0$ for all $x \in C$}}$.
\end{para}

\begin{proposition}[\protect{\cite[Prop. 2.2.11]{BR1}}]\label{pro:operatorAlgebras_propertiesPE}
    Let $\MFA$ be a unital $C^\ast$-algebra. The set $\PE{\MFA}$ is a uniformly closed convex cone satisfying $\PE{\MFA} \cap (- \PE{\MFA}) = \{0\}$. Moreover, if $A \in \SE{\MFA}$, then the elements $A_\pm \ce \frac{1}{2} (\abs{A} \pm A)$ are positive, and they satisfy $A_+ A_- = 0$ and $A = A_+ - A_-$.
\end{proposition}

According to this result, the sum of two positive elements is again positive. Regarding the product of two such elements, the following property is found which will be needed later in \cref{ch:relativeEntropy}.

\begin{lemma}\label{lem:operatorAlgebras_productPositiveElements}
    Let $\MFA$ be a unital $C^\ast$-algebra and $A, B \in \PE{\MFA}$. Then $AB \in \PE{\MFA}$ if and only if $AB = BA$.
\end{lemma}

\begin{Proof}
    \pp{1.} Assume first that $A$ and $B$ commute. Using \cref{thm:operatorAlgebras_characterizationPE}, one can define the \emph{square root} of a positive element $C \in \PE{\MFA}$: it is the unique element $\sqrt{C} \ce C^{1/2} \ce D \in \PE{\MFA}$ such that $D^2 = C$ \cite[p. 34]{BR1}. With this, the following decomposition holds true:
    \begin{equation*}
        AB = \sqrt{A} \sqrt{A} \sqrt{B} \sqrt{B} = \sqrt{B} \sqrt{A} \sqrt{A} \sqrt{B} = \bigl(\sqrt{B}\bigr)^\ast \bigl(\sqrt{A}\bigr)^\ast \sqrt{A} \sqrt{B} = \bigl(\sqrt{A} \sqrt{B}\bigr)^\ast \sqrt{A} \sqrt{B} \ .
    \end{equation*}
    Note that $\sqrt{A}$ and $\sqrt{B}$ commute according to the assumption and the continuous functional calculus in $C^\ast$-algebras \cite[p. 22]{StratilaZsido19}. From \cref{thm:operatorAlgebras_characterizationPE}, it now follows that $AB$ is a positive element.
    
    \pp{2.} Conversely, if $AB \ge 0$, then $\sigma(AB) \subset [0, + \infty)$ and $(AB)^\ast = AB$ by definition. But $(AB)^\ast = B^\ast A^\ast = BA$ since $A$ and $B$ are self-adjoint, hence $AB = BA$.
\end{Proof}

\begin{para}[Order relation for self-adjoint elements]\label{para:operatorAlgebras_orderSE}
    (\cite[p. 36]{BR1})
    Let $\MFA$ be a unital $C^\ast$-algebra. The closed convex cone $\PE{\MFA}$ induces an order structure on the real linear subspace $\SE{\MFA}$ of self-adjoint elements: for $A, B \in \SE{\MFA}$ define
    \begin{equation*}
        A \le B \quad \text{(or $B \ge A$)} \quad :\iff \quad B - A \in \PE{\MFA} \ .
    \end{equation*}
    This relation defines a partial order on $\SE{\MFA}$: (1) it is clear that $A \le A$ for all $A \in \SE{\MFA}$ because $0 \in \PE{\MFA}$; (2) if $A, B \in \SE{\MFA}$ such that $A \le B$ and $B \le A$, it follows from \cref{pro:operatorAlgebras_propertiesPE} that $B - A \in \PE{\MFA} \cap (- \PE{\MFA}) = \{0\}$, hence $B = A$; (3) finally, if there is a third element $C \in \SE{\MFA}$ such that $A \le B$ and $B \le C$, then $C - A = C - B + B - A$ is the sum of two positive elements, hence $C - A \ge 0$ as $\PE{\MFA}$ is a convex cone.
\end{para}

\begin{lemma}\label{lem:operatorAlgebras_invertibleLowerBound}
    Let $(\MFA, \idm)$ be a unital $C^\ast$-algebra and $A \in \PE{\MFA}$ be a positive element. Then $A$ is invertible if and only if there exists $\epsilon > 0$ such that $A \ge \epsilon \, \idm$.
\end{lemma}

\begin{Proof}
    \pp{1.} Assume that $A$ is invertible. From positivity of $A$, it follows that $\sigma(A) \subset [0, + \infty)$, and from invertibility of $A$, one can conclude that $0 \notin \sigma(A)$. Therefore, there exists $\epsilon > 0$ such that $\sigma(A) \subset [\epsilon, + \infty)$. Since $\sigma(A - \epsilon \, \idm) = \set{\lambda - \epsilon : \lambda \in \sigma(A)}$ \cite[Prop. 2.2.3]{BR1} and $\lambda \ge \epsilon$ for all $\lambda \in \sigma(A)$, it follows that $A - \epsilon \, \idm \ge 0$.

    \pp{2.} Conversely, assume that $A \ge \epsilon \, \idm$ for some $\epsilon > 0$. Employing an analogous argument as before, one obtains $\lambda \ge \epsilon$ for all $\lambda \in \sigma(A)$. Hence, $\sigma(A) \subset [\epsilon, + \infty)$, that is, $0 \notin \sigma(A)$, so $A$ must be invertible.
\end{Proof}

\begin{example}[\protect{Positive elements in $\BO(\MH)$}]\label{exa:operatorAlgebras_positiveBO}
    Let $(\MH, \bdot)$ be a Hilbert space and $\MFA = \BO(\MH)$ be the $C^\ast$-algebra of bounded linear operators on $\MH$ from \cref{exa:operatorAlgebras_CStar} \ref{enu:operatorAlgebras_exaBOCStar}. Then for all $A \in \MFA$, the following characterization holds true \cite[Thm. VIII.3.8]{Conway85}:
    \begin{equation*}
        A \in \PE{\MFA} \quad \iff \quad \forall \xi \in \MH \, : \, \braket{\xi, A \xi} \ge 0 \ . \tag*{\qedhere}
    \end{equation*}
\end{example}

\section{Positive Linear Maps}\label{sec:operatorAlgebras_positiveMaps}

Having introduced the closed convex cone $\PE{\MFA}$ of positive elements, now, a class of linear mappings between $C^\ast$-algebras shall be defined which preserve these cones.

\subsection{Basic Properties}

\begin{definition}[Positive mapping]\label{def:operatorAlgebras_positiveMapping}
    A linear mapping $f : \MFA \longto \MFB$ between two unital $C^\ast$-algebras $\MFA$ and $\MFB$ is called \bemph{positive} iff for every $A \in \MFA$, $A \ge 0$ implies $f(A) \ge 0$, that is, iff $f$ maps the cone $\PE{\MFA}$ into the cone $\PE{\MFB}$.
\end{definition}

\begin{remark}\label{rem:operatorAlgebras_positiveMappingOrderPreserving}
    It follows immediately from this definition that positive linear mappings $f : \MFA \longto \MFB$ are \emph{order-preserving} with respect to the order relation introduced in \cref{para:operatorAlgebras_orderSE}: if $A, B \in \SE{\MFA}$ satisfy $A \le B$, then $B - A \ge 0$ and hence $f(B) - f(A) \ge 0$ in $\MFB$, \ie{}, $A \le B \implies f(A) \le f(B)$.
\end{remark}

Next, it will be shown that positive mappings automatically satisfy an important algebraic property related to the algebra involution \cite[p. 2]{Stormer13}.

\begin{proposition}\label{pro:operatorAlgebras_positiveImpliesStarPreserving}
    Let $f : \MFA \longto \MFB$ be a positive linear mapping between unital $C^\ast$-algebras $\MFA$ and $\MFB$. Then $f$ is $\ast$-preserving, \ie{}, $f(A^\ast) = f(A)^\ast$ for all $A \in \MFA$.
\end{proposition}

\begin{Proof}
    Any element $A \in \MFA$ can be written as $A = A_1 - A_2 + \ii (A_3 - A_4)$ for $\{A_i\}_{i=1}^{4} \subset \PE{\MFA}$ (\cf{} \cref{para:operatorAlgebras_astAlgebras} and \cref{pro:operatorAlgebras_propertiesPE}). It holds that $A^\ast = A_1 - A_2 - \ii (A_3 - A_4)$ since positive elements are, by definition, self-adjoint. From positivity of $f$, it follows that $f(A_i) \in \PE{\MFB} \subset \SE{\MFB}$ for all $i \in \{1, \dotsc, 4\}$, hence linearity of $f$ implies $f(A^\ast) = f(A)^\ast$.
\end{Proof}

Positive mappings are also well-behaved in regards to the topological structure of the algebras: they are automatically continuous. A simple proof can be found in \cite[Prop. 2.1]{Paulsen03} or \cite[p. 2]{Stormer13}.

\begin{proposition}\label{pro:operatorAlgebras_positiveMappingContinuous}
    Let $f : \MFA \longto \MFB$ be a positive linear mapping between unital $C^\ast$-algebras. Then $f$ is bounded with $\norm{f}_\mop \le 2 \norm{f(\idm)}_\MFB$.
\end{proposition}

\subsection{Completely Positive Maps}\label{subsec:operatorAlgebras_CP}

For applications, one often requires mappings which satisfy a stronger positivity-preserving property than the one introduced in \cref{def:operatorAlgebras_positiveMapping}. (See, for example, \cref{sec:relativeEntropy_monotonicity}.) These maps will be introduced in this subsection.

\begin{para}[Matrices over a $C^\ast$-algebra]\label{para:operatorAlgebras_matricesCStar}
    (\cite[Sect. 30.4]{BlanchardBrüning15}, \cite[pp. 2 f.]{Paulsen03})
    Let $\MFA$ be a $C^\ast$-algebra and $n \in \N$ be arbitrary. Denote by $\Mat(n; \MFA)$ the space of all $(n \times n)$-matrices $A = [A_{ij}]$ with entries $A_{ij}$, $i, j \in \set{1, \dotsc, n}$, in the algebra $\MFA$; one also writes $A$ as
    \begin{equation*}
        A =
        \begin{bmatrix}
            A_{11} & \dots & A_{1n} \\
            \vdots & \ddots & \vdots \\
            A_{n1} & \dots & A_{nn}
        \end{bmatrix} \ .
    \end{equation*}
    It is clear that the space $\Mat(n; \MFA)$ is a $\ast$-algebra with respect to the usual matrix multiplication and involution, both being defined in terms of the respective algebraic operations on $\MFA$:
    \begin{equation*}
        [A_{ij}] [B_{ij}] \ce \biggl[\,\sum_{k=1}^{n} A_{ik} B_{kj}\biggr] \tand [A_{ij}]^\ast \ce [A_{ji}^\ast] \ .
    \end{equation*}
    Furthermore, one can show that there exists a canonical algebra norm on $\Mat(n; \MFA)$ with respect to which this space becomes a $C^\ast$-algebra \cite[p. 143]{BR1}, \cite[p. 192]{Takesaki79}. Finally, one can show that \cite[p. 192]{Conway00}
    \begin{equation*}
        \Mat(n; \MFA) = \MFA \otimes \Mat(n; \C) \ . 
    \end{equation*}

    Let $\MFB$ be another $C^\ast$-algebra and $f : \MFA \longto \MFB$ be a linear map. For all $n \in \N$, define linear mappings $f_{(n)} : \Mat(n; \MFA) \longto \Mat(n; \MFB)$ by setting for all $A = [A_{ij}] \in \Mat(n; \MFA)$:
    \begin{equation}\label{eq:operatorAlgebras_nMatrixMap}
        f_{(n)}\bigl([A_{ij}]\bigr) \ce \bigl[f(A_{ij})\bigr] \ .
    \end{equation}
\end{para}

\begin{definition}[Completely positive mapping]\label{def:operatorAlgebras_CP}
    Let $f : \MFA \longto \MFB$ be a linear map between unital $C^\ast$-algebras $\MFA$ and $\MFB$. For some $n \in \N$, $f$ is called \bemph{$n$-positive} iff $f_{(n)} : \Mat(n; \MFA) \longto \Mat(n; \MFB)$ is a positive mapping in the sense of \cref{def:operatorAlgebras_positiveMapping}. Furthermore, $f$ is called \bemph{completely positive} iff for every $n \in \N$, $f_{(n)}$ is positive.
\end{definition}

\begin{remark}[\protect{Physical interpretation}]\label{rem:operatorAlgebras_physicalInterpretationCP}
    (\cite[p. 148]{Lindblad75})
    Consider a physical system $S_1$ described by a concrete $C^\ast$-algebra $\MFA \subset \BO(\MH_1)$ over a Hilbert space $\MH_1$, and let $S_2$ be a second system, independent of the system $S_1$, which is described by the algebra $\Mat(n; \C)$, $n \in \N$. According to the postulates of quantum mechanics \cite[Sect. 2.1]{Petz08}, the Hilbert space for the combined system $S_1 + S_2$ is given by the Hilbert-space tensor product $\MH = \MH_1 \otimes \MH_2$, where $\dim(\MH_2) = n$.
    
    Let $T : \MFA \longto \MFA$ be a linear mapping. By the duality between the bounded operators $\BO(\MH_1)$ and the trace-class operators $\NO(\MH_1)$, \cf{} \cref{pro:operatorAlgebras_sigmaWeakTopology} below, $T$ induces a linear mapping $T^\ast$ on the space $\NO(\MH_1)$ \cite[Lem. 31.3]{BlanchardBrüning15}. For $T^\ast$ to send quantum states of $S_1$ (that is, density matrices on $\MH_1$) again to quantum states, it is necessary that the mapping $T^\ast$ is positive.\footnotemark
    \footnotetext{Although $\NO(\MH_1)$ is not a $C^\ast$-algebra itself (\cref{exa:operatorAlgebras_CStar} \ref{enu:operatorAlgebras_exaBOCStar}), it is a $\ast$-ideal in $\BO(\MH)$ (\cref{pro:operators_NOHS}). Therefore, one can define notions such as positivity of elements and complete positivity of mappings similarly as before \cite[p. 486]{BlanchardBrüning15}.}
    Furthermore, under the assumption that $T$ does not influence the system $S_2$ directly, the mapping $T_{(n)} = T \otimes \id : \MFA \otimes \Mat(n; \C) \longto \MFA \otimes \Mat(n; \C)$, defined as in \cref{eq:operatorAlgebras_nMatrixMap}, describes the resulting transformation on the combined system $S_1 + S_2$: $T_{(n)}(A_1 \otimes A_2) = T(A_1) \otimes A_2$. Therefore, in order for $(T_{(n)})^\ast$ to map quantum states of the combined system $S_1 + S_2$ to quantum states, $(T_{(n)})^\ast$ has to be positive, \ie{}, $T^\ast$ has to be completely positive which is the case if and only if $T$ is completely positive \cite[Lem. 31.4]{BlanchardBrüning15}.
\end{remark}

\begin{examples}\label{exa:operatorAlgebras_CP}
    \leavevmode
    \begin{enumerate}[env]
        \item \label{enu:operatorAlgebras_exaConjugationCP} Let $\MFA$ and $\MFB$ be unital $C^\ast$-algebras, let $f : \MFA \longto \C$ be a positive linear functional, let $\pi : \MFA \longto \MFB$ be a $\ast$-homomorphism, and let $\alpha : \MFA \longto \MFB$ be defined by $\alpha(A) \ce V^\ast \pi(A) V$ for some $V \in \MFB$. All of these three mappings are completely positive \cite[pp. 474 f.]{BlanchardBrüning15}. To show this, it suffices to check, according to \cite[Cor. 30.1]{BlanchardBrüning15}, that for all $n \in \N$, $\{A_i\}_{i=1}^{n} \subset \MFA$, and $\{B_j\}_{j=1}^{n} \subset \MFB$, there holds $\sum_{i,j=1}^{n} B_i^\ast \alpha(A_i^\ast A_j) B_j \ge 0$ in the algebra $\MFB$. Indeed, for the mapping $\alpha$ one finds that
        \begin{align*}
            \sum_{i,j=1}^{n} B_i^\ast \alpha(A_i^\ast A_j) B_j &= \sum_{i,j=1}^{n} B_i^\ast V^\ast \pi(A_i)^\ast \pi(A_j) V B_j = \sumbra{\bigg}{\sum_{i=1}^{n} \pi(A_i) V B_i }^\ast \sumbra{\bigg}{\sum_{j=1}^{n} \pi(A_j) V B_j}
        \end{align*}
        is a positive element in $\MFB$ by \cref{thm:operatorAlgebras_characterizationPE}, hence $\alpha$ is completely positive.

        \item The notion of complete positivity can be used to distinguish between classical and quantum situations \cite[p. 260]{GustafsonSigal20}, \cite[p. 449]{Keyl2002}: a positive mapping $f : \MFA \longto \MFB$ between $C^\ast$-algebras is automatically completely positive if either $\MFA$ or $\MFB$ is of the form $C^0(X; \K)$ for a compact Hausdorff space $X$ \cite[Thm. 3.9 \& 3.11]{Paulsen03}, \cite[Thm. 1.2.4 \& 1.2.5]{Stormer13}.
        
        \item \label{enu:operatorAlgebras_exaTrCP} Let $\MFA$ be a unital $C^\ast$-algebra and $\tr : \Mat(n; \MFA) \longto \MFA$ be defined by $\tr([A_{ij}]) \ce \sum_{i=1}^{n} A_{ii}$. Then $\tr$ is a positive linear functional because for all $A = [A_{ij}] \in \Mat(n; \MFA)$, there holds
        \begin{align*}
            \tr(A^\ast A) = \tr\Biggl(\biggl[\,\sum_{k=1}^{n} A_{ki}^\ast A_{kj}\biggr]\Biggr) = \sum_{i=1}^{n} \sum_{k=1}^{n} A_{ki}^\ast A_{ki} \ge 0
        \end{align*}
        since $A_{ki}^\ast A_{ki} \ge 0$ in $\MFA$ for all $i, k \in \set{1, \dotsc, n}$ (\cref{thm:operatorAlgebras_characterizationPE}) and $\PE{\MFA}$ is a convex cone (\cref{pro:operatorAlgebras_propertiesPE}). One can show that $\tr$ is actually completely positive \cite[p. 40]{Paulsen03}.

        \item Let $\MH_1$ and $\MH_2$ be two Hilbert spaces and $\MH \ce \MH_1 \otimes \MH_2$. Define a mapping $\alpha : \BO(\MH_1) \longto \BO(\MH)$, $A \longmto A \otimes \id_{\MH_2}$, which is completely positive \cite[Exa. 9.1]{Petz08}. Its dual map $\alpha^\ast : \NO(\MH) \longto \NO(\MH_1)$, in the sense discussed in \cref{rem:operatorAlgebras_physicalInterpretationCP}, is the so-called \emph{partial trace} $\alpha^\ast(A \otimes B) = A \tr(B)$ which is physically very relevant \cite[pp. 11 f.]{Petz08}. By \cite[Lem. 31.4]{BlanchardBrüning15}, $\alpha^\ast$ is completely positive.
        
        \item Consider the $(2 \times 2)$-matrices $\MFA = \Mat(2; \C)$ and the transposition $f : \MFA \longto \MFA$, $A \longmto A^T$. This is a positive mapping: for all $A \in \MFA$, it holds that
        \begin{equation*}
            f(A^\ast A) = (A^\ast A)^T = A^T (A^\ast)^T = A^T (A^T)^\ast \ge 0 \ .
        \end{equation*}
        However, $f$ fails to be 2-positive \cite[Exa. 2.48]{HiaiPetz14}. To see this, consider the following element of the algebra $\Mat(2; \MFA) = \Mat(2; \C) \otimes \Mat(2; \C) = \Mat(4; \C)$:
        \begin{equation*}
            A \ce \bmat{\pmat{2 & 0 \\ 0 & 1} & \pmat{0 & 2 \\ 1 & 0} \\ \pmat{0 & 1 \\ 2 & 0} & \pmat{1 & 0 \\ 0 & 2}} = \bmat{2 & 0 & 0 & 2 \\ 0 & 1 & 1 & 0 \\ 0 & 1 & 1 & 0 \\ 2 & 0 & 0 & 2} \ .
        \end{equation*}
        Its eigenvalues are $\sigma(A) = \set{4, 2, 0}$ (the last eigenvalue appears twice), hence it is positive. Yet,
        \begin{equation*}
            f_{(2)}(A) = \bmat{\pmat{2 & 0 \\ 0 & 1} & \pmat{0 & 1 \\ 2 & 0} \\ \pmat{0 & 2 \\ 1 & 0} & \pmat{1 & 0 \\ 0 & 2}} = \bmat{2 & 0 & 0 & 1 \\ 0 & 1 & 2 & 0 \\ 0 & 2 & 1 & 0 \\ 1 & 0 & 0 & 2}
        \end{equation*}
        has eigenvalues $\sigma(f_2(A)) = \set{3, 1, -1}$ (the first eigenvalue appears twice), \ie{}, it is not a positive element in $\Mat(2; \MFA)$. One can even show that for all $n > 1$, $f_{(n)}$ is not positive \cite[p. 192]{Conway00}.\qedhere
    \end{enumerate}
\end{examples}

The following result of \textsc{M.-D. Choi} \cite[Prop. 4.1]{Choi80} characterizes 2-positive mappings.

\begin{proposition}\label{pro:operatorAlgebras_2positivity}
    Let $f : \MFA \longto \MFB$ be a positive linear mapping between two unital $C^\ast$-algebras, and suppose that $f(\idm_\MFA) \in \PE{\MFB} \cap \GL(\MFB)$ is positive and invertible. Furthermore, let $R, S, T \in \MFA$ be arbitrary with $T \in \PE{\MFA} \cap \GL(\MFA)$. Then the following properties are equivalent:
    \begin{enumerate}[equiv]
        \item $f$ is 2-positive;
        \item If $R \ge S \inv{T} S^\ast$ in $\MFA$, then $f(R) \ge f(S) f(T)^{-1} f(S)^\ast$ in $\MFB$;
        \item $f(S \inv{T} S^\ast) \ge f(S) f(T)^{-1} f(S)^\ast$ in $\MFB$.
    \end{enumerate}
\end{proposition}

\subsection{Schwarz Mappings}

To conclude this section, another class of positive mappings between $C^\ast$-algebras is introduced which will be of particular relevance for the study of monotonicity of the relative entropy in \cref{sec:relativeEntropy_monotonicity}. The following definition subscribes to the convention of \cite[p. 81]{OP04}.

\begin{definition}[Schwarz map]\label{def:operatorAlgebras_SchwarzMap}
    Let $(\MFA, \idm_\MFA)$ and $(\MFB, \idm_\MFB)$ be two unital $C^\ast$-algebras. A linear mapping $\alpha : \MFA \longto \MFB$ is called a \bemph{Schwarz mapping} iff $\alpha$ is unital, that is, $\alpha(\idm_\MFA) = \idm_\MFB$, and satisfies for all $A \in \MFA$ the following inequality, referred to as \bemph{Schwarz inequality}:
    \begin{equation}\label{eq:operatorAlgebras_SchwarzInequality}
        \alpha(A)^\ast \alpha(A) \le \alpha(A^\ast A) \ .
    \end{equation}
\end{definition}

\begin{lemma}\label{lem:operatorAlgebras_basicPropertiesSchwarzMap}
    Every Schwarz mapping $\alpha : \MFA \longto \MFB$ is positive and satisfies the inequality $\alpha(A) \alpha(A)^\ast \le \alpha(A A^\ast)$ for all $A \in \MFA$.
\end{lemma}

\begin{Proof}
    \pp{1.} Let $A = B^\ast B$, $B \in \MFA$, be an arbitrary positive element in the $C^\ast$-algebra $\MFA$ (\cf{} \cref{thm:operatorAlgebras_characterizationPE}). Then $\alpha(A) = \alpha(B^\ast B) \ge \alpha(B)^\ast \alpha(B) \ge 0$ in $\MFB$ by \eqref{eq:operatorAlgebras_SchwarzInequality}.

    \pp{2.} Since $\MFA$ is closed under the algebra involution, one can choose $A^\ast$ in \eqref{eq:operatorAlgebras_SchwarzInequality} to obtain $\alpha(A^\ast)^\ast \alpha(A^\ast) \le \alpha(A A^\ast)$. As was just shown, $\alpha$ is a positive mapping, hence it is especially $\ast$-preserving by \cref{pro:operatorAlgebras_positiveImpliesStarPreserving}. This proves the asserted inequality.
\end{Proof}

According to \cref{lem:operatorAlgebras_basicPropertiesSchwarzMap}, Schwarz mappings are stronger than positive mappings. Similarly, the next proposition shows that 2-positive mappings are stronger than Schwarz mappings; the proof can be found in \cite[Prop. 3.3]{Paulsen03}.

\begin{proposition}\label{pro:operatorAlgebras_2PmapsAreSchwarzMaps}
    Let $\MFA$ and $\MFB$ be two unital $C^\ast$-algebras, and let $\alpha : \MFA \longto \MFB$ be a unital 2-positive mapping. Then $\alpha$ satisfies the Schwarz inequality \eqref{eq:operatorAlgebras_SchwarzInequality}.
\end{proposition}

\section{von Neumann Algebras}\label{sec:operatorAlgebras_vonNeumann}

In \cref{exa:operatorAlgebras_CStar} \ref{enu:operatorAlgebras_exaConcreteCstar}, the notion of a concrete $C^\ast$-algebra was defined to be a uniformly closed $\ast$-subalgebra of $\BO(\MH)$. Considering further classes of $\ast$-subalgebras of $\BO(\MH)$ which satisfy stronger closure conditions leads to the concept of \emph{von Neumann algebras}.

\subsection{Operator Topologies}

Let $(\MH, \bdot)$ be a Hilbert space and $\ndot$ be the norm on $\MH$ induced by $\bdot$. In the following, three important locally convex topologies on $\BO(\MH)$ shall be defined by specifying a family of semi-norms (\cf{} \cref{pro:topologicalVectorSpaces_LCS}). The definitions can be found in \cite[Sect. 26.3]{BlanchardBrüning15}, \cite[Sect. 2.4.1]{BR1}, or \cite[Ch. 1]{StratilaZsido19}.

\begin{para}[Strong topology]\label{para:operatorAlgebras_strongTop}
    For $\xi \in \MH$, the mapping $A \longmto \norm{A \xi}$ on $\BO(\MH)$ is clearly a semi-norm. The topology $\MFT_{so} \ce \MFT_{P_{so}}$ induced by the system $P_{so} \ce \set{p_\xi : \xi \in \MH}$ of semi-norms
    \begin{equation*}
        p_\xi(A) \ce \norm{A \xi} \com A \in \BO(\MH) \ ,
    \end{equation*}
    is called the \bemph{strong operator topology} (also: \bemph{$so$-topology}) on $\BO(\MH)$. Since by \cref{pro:topologicalVectorSpaces_LCS} a net $(A_i)_{i \in I} \subset \BO(\MH)$ converges to an operator $A \in \BO(\MH)$ with respect to $\MFT_{so}$ if and only if $\lim_{i \in I} p_x(A_i - A) = \lim_{i \in I} \norm{(A_i - A) \xi} \to 0$ for all $\xi \in \MH$, it follows that the strong operator topology is the topology of pointwise convergence on $\BO(\MH)$.
\end{para}

\begin{para}[Weak topology]\label{para:operatorAlgebras_weakTop}
    For $\xi, \eta \in \MH$, the mapping $A \longmto \abs{\braket{\xi, A \eta}}$ on $\BO(\MH)$ is a semi-norm. The topology $\MFT_{wo} \ce \MFT_{P_{wo}}$ induced by the family $P_{wo} \ce \set{p_{\xi,\eta} : \xi, \eta \in \MH}$, where
    \begin{equation*}
        p_{\xi, \eta}(A) \ce \abs{\braket{\xi, A \eta}} \com A \in \BO(\MH) \ ,
    \end{equation*}
    is called the \bemph{weak operator topology} (also: \bemph{$wo$-topology}) on $\BO(\MH)$. As for the strong operator topology, it follows from \cref{pro:topologicalVectorSpaces_LCS} that the weak operator topology is the topology of pointwise weak convergence on $\BO(\MH)$.
\end{para}

\begin{para}[$\sigma$-weak topology]\label{para:operatorAlgebras_sigmaWeakTop}
    Introduce the notation
    \begin{equation*}
        \ell^2(\N; \MH) \ce \set[\Bigg]{(\xi_n)_{n \in \N} \ : \ \xi_n \in \MH, \, \sum_{n=1}^{\infty} \norm{\xi_n}^2 < + \infty}
    \end{equation*}
    for the space of all sequences $(\xi_n)_{n \in \N} \subset \MH$ which are square-summable, \ie{}, $\sum_{n=1}^{\infty} \norm{\xi_n}^2 < + \infty$.
    
    Let $(\xi_n)_{n \in \N}, (\eta_n)_{n \in \N} \in \ell^2(\N; \MH)$ be arbitrary. For every $A \in \BO(\MH)$, it follows by using first the Cauchy-Schwarz and then Hölder's inequality that \cite[p. 67]{BR1}
    \begin{equation*}
        \sum_{n=1}^{\infty} \abs{\braket{\xi_n, A \eta_n}} \le \sum_{n=1}^{\infty} \norm{A}_\mop \norm{\xi_n} \norm{\eta_n} \le \norm{A}_\mop \sumbra{\Bigg}{\sum_{n=1}^{\infty} \norm{\xi_n}^2}^{1/2} \sumbra{\Bigg}{\sum_{n=1}^{\infty} \norm{\eta_n}^2}^{1/2} < + \infty \ .
    \end{equation*}
    This shows that $A \longmto \sum_{n=1}^{\infty} \abs{\braket{\xi_n, A \eta_n}}$ is a semi-norm on $\BO(\MH)$. The topology $\MFT_{uw} \ce \MFT_{P_{uw}}$ induced by the system $P_{uw} \ce \set[\big]{q_{(\xi_n), (\eta_n)} : (\xi_n)_{n \in \N}, (\eta_n)_{n \in \N} \in \ell^2(\N; \MH)}$ of semi-norms
    \begin{equation*}
        q_{(\xi_n), (\eta_n)}(A) \ce \sum_{n=1}^{\infty} \abs{\braket{\xi_n, A \eta_n}} \com A \in \BO(\MH) \ ,
    \end{equation*}
    is called the \bemph{$\sigma$-weak operator topology} (or: \bemph{ultraweak operator topology}; \bemph{$uw$-topology}). It holds that the mappings $A \longmto A B$ and $A \longmto B A$ on $\BO(\MH)$ for fixed $B \in \BO(\MH)$, as well as the mapping $A \longmto A^\ast$, are continuous with respect to both $\MFT_{uw}$ and $\MFT_{wo}$ \cite[Prop. 2.4.2]{BR1}.
\end{para}

The next proposition characterizes the $\sigma$-weak topology in terms of trace-class operators. For proofs of the assertions, see \cite[Thm. 26.5 and pp. 378 f.]{BlanchardBrüning15} and \cite[Prop. 2.4.3]{BR1}.

\begin{proposition}[$\sigma$-weak topology and $\NO(\MH)$]\label{pro:operatorAlgebras_sigmaWeakTopology}
    It holds that $\cdual{\bigl(\NO(\MH), \ndot_\mathrm{tr}\bigr)} \cong \bigl(\BO(\MH), \ndot_\mop\bigr)$ with respect to the duality mapping
    \begin{equation*}
        \Phi :
        \begin{cases}
            \begin{aligned}
                \BO(\MH) &\longto \cdual{\NO(\MH)} \ , \\
                A &\longmto \bigl(\Phi_A : \rho \longmto \tr(A \rho)\bigr) \ .
            \end{aligned}
        \end{cases}
    \end{equation*}
    The weak-$\ast$ topology $\sigma\bigl(\BO(\MH), \NO(\MH)\bigr) = \sigma\bigl(\cdual{\NO(\MH)}, \NO(\MH)\bigr)$ on $\BO(\MH)$ arising from the above duality, \ie{}, the topology generated by the semi-norms $p_\rho(A) \ce \abs{\tr(A \rho)}$, $\rho \in \NO(\MH)$, on $\BO(\MH)$ (\cf{} \cref{para:topologicalVectorSpaces_weakStar}) coincides with the $\sigma$-weak operator topology:
    \begin{equation*}
        \MFT_{uw} = \sigma\bigl(\BO(\MH), \NO(\MH)\bigr) \ .
    \end{equation*}
    
    In particular, for every $\omega \in \cdual{\bigl(\BO(\MH), \MFT_{uw}\bigr)}$, there exist sequences $(\xi_n)_{n \in \N}, (\eta_n)_{n \in \N} \in \ell^2(\N; \MH)$, equivalently, a trace-class operator $\rho \in \NO(\MH)$, such that
    \begin{equation*}
        \omega = \sum_{n=1}^{\infty} \braket{\xi_n, \, \sbullet \ \eta_n} \com \text{equivalently} \com \omega = \tr(\rho \ \sbullet\,) \ .
    \end{equation*}
\end{proposition}

The different topologies introduced in \ref{para:operatorAlgebras_strongTop}, \ref{para:operatorAlgebras_weakTop} and \ref{para:operatorAlgebras_sigmaWeakTop} are in the following relation (\cf{} \cref{para:topologicalVectorSpaces_comparisonTopologies}) to one another \cite[p. 379]{BlanchardBrüning15}, \cite[p. 70]{BR1}.

\begin{proposition}\label{pro:operatorAlgebras_comparisonTopologies}
    The following relations hold true:
    \begin{equation*}
        \begin{tikzcd}[row sep=large, column sep=large]
            \MFT_{wo} \arrow[draw=none]{r}[sloped,auto=false]{\le} \arrow[draw=none]{d}[sloped,auto=false]{\le} & \MFT_{so} \arrow[draw=none]{d}[sloped,auto=false]{\le} \\
            \MFT_{uw} \arrow[draw=none]{r}[sloped,auto=false]{\le} & \MFT_{\mathrm{norm}}
        \end{tikzcd}
    \end{equation*}
\end{proposition}

\subsection{Definition, Bicommutant Theorem, Examples}

\begin{definition}[von Neumann algebra]\label{def:operatorAlgebras_vonNeumann}
    Let $\MH$ be a Hilbert space. A $\ast$-subalgebra of operators $\MFM \subset \BO(\MH)$ which contains the identity $\id_\MH$ and which is closed in the strong operator topology is called a \bemph{von Neumann algebra}\footnote{This terminology was introduced by \textsc{J. Dixmier} in 1957 \cite{Dixmier81}. In the works of \textsc{J. von Neumann} \cite{vonNeumann30} (and \textsc{von Neumann} \& \textsc{F. J. Murray} \cite{MurrayVonNeumann36}), they were called \enquote{rings of operators}.} (or: \emph{concrete $W^\ast$-algebra}).
\end{definition}

\begin{remark}\label{rem:operatorAlgebras_vonNeumann}
    From Mazur's theorem \cite[Cor. 2.11]{Voigt20}, it follows that a von Neumann algebra is also closed with respect to the weak operator topology. Furthermore, since $\MFT_{so} \subset \MFT_\mathrm{norm}$ by \cref{pro:operatorAlgebras_comparisonTopologies}, it holds that $\MFM$ is uniformly closed as well, \cf{} \cref{para:topologicalVectorSpaces_comparisonTopologies}, and hence every von Neumann algebra is a concrete $C^\ast$-algebra according to \cref{exa:operatorAlgebras_CStar} \ref{enu:operatorAlgebras_exaConcreteCstar}.
\end{remark}

\begin{definition}[Commutant]\label{def:operatorAlgebras_commutant}
    For a subset $\MFA \subset \BO(\MH)$, the \bemph{commutant} is defined to be the set
    \begin{equation*}
        \comm{\MFA} \ce \set{A \in \BO(\MH) \ : \ \text{$A B = B A$ for all $B \in \MFA$}}
    \end{equation*}
    of all bounded linear operators which commute with every element from the set $\MFA$. Similarly, the \bemph{bicommutant} (or: \emph{double commutant}) of $\MFA$ is defined as $\bicomm{\MFA} \ce \comm{(\comm{\MFA})}$.
\end{definition}

\begin{lemma}[\protect{\cite[Prop. 6.3]{Moretti19}}]\label{lem:operatorAlgebras_propComm}
    Let $\MFA \subset \BO(\MH)$ be a subset. Then $\comm{\MFA}$ is a strongly closed, unital subalgebra of $\BO(\MH)$. If $\MFA$ is self-adjoint, then $\comm{\MFA}$ is even a von Neumann algebra.
\end{lemma}

The next theorem is one of the most important and deep results in the theory of von Neumann algebras. It was discovered by \textsc{J. von Neumann} in 1930 \cite[Thm. 8]{vonNeumann30}, and a modern proof can be found in \cite[Thm. 2.4.11]{BR1}, \cite[Thm. 6.4]{Moretti19}, or \cite[Thm. 3.2 \& Cor. 3.11]{StratilaZsido19}.

\begin{ntheorem}[von Neumann's bicommutant theorem]\label{thm:operatorAlgebras_bicommutant}
    Let $\MFM \subset \BO(\MH)$ be a $\ast$-subalgebra which contains the identity $\id_\MH$. Then the following assertions are equivalent:
    \begin{enumerate}[equiv]
        \item $\MFM$ is a von Neumann algebra;
        \item $\MFM$ is weakly closed;
        \item $\MFM$ is $\sigma$-weakly closed;
        \item $\MFM = \bicomm{\MFM}$.
    \end{enumerate}
\end{ntheorem}

\begin{examples}\label{exa:operatorAlgebras_vonNeumann}
    \leavevmode
    \begin{enumerate}[env]
        \item \label{enu:operatorAlgebras_exaLInf} Let $(X, \Sigma, \mu)$ be a $\sigma$-finite measure space. Consider the Hilbert space $\MH = L^2(X, \mu)$ and the $C^\ast$-algebra $\MFM = L^\infty(X, \mu)$ from \cref{exa:operatorAlgebras_CStar} \ref{enu:operatorAlgebras_exaLInfCStar}. One can view $\MFM$ as a $\ast$-subalgebra of $\BO(\MH)$ by virtue of the identification $\MFM \owns f \longmto M_f \in \BO(\MH)$, where $M_f g \ce f g$, $g \in \MH$, is the multiplication operator with the function $f$; that is, one identities
        \begin{equation*}
            \MFM \cong \set[\big]{M_f \in \BO(\MH) \ : \ f \in L^\infty(X, \mu)} \ .
        \end{equation*}
        In this identification, $\MFM$ becomes a von Neumann algebra acting on $\MH$. This can be proved by showing that $\comm{\MFM} \subset \MFM$, see \cite[Thm. IX.6.6]{Conway85}, \cite[Prop. 12.4]{Conway00}, or \cite[Sect. 6.10]{StratilaZsido19} for detailed arguments.
        
        The $C^\ast$-algebra of continuous functions $\MFA = C^0(X; \K)$ from \cref{exa:operatorAlgebras_CStar} \ref{enu:operatorAlgebras_exaCCStar}, which, for $X$ a compact topological space, can also be viewed as $\MFA \subset \BO(\MH)$ in terms of the above identification, does not form a von Neumann algebra because $\MFA$ is only uniformly, but neither strongly nor weakly closed in $\BO(\MH)$ \cite[Exa. 5.1.6]{KadisonRingrose97}. ($L^\infty$-functions are pointwise almost-everywhere limits of continuous functions.)

        \item \label{enu:operatorAlgebras_exaBOvN} Let $\MH$ be a Hilbert space. Then the $C^\ast$-algebra $\BO(\MH)$ from \cref{exa:operatorAlgebras_CStar} \ref{enu:operatorAlgebras_exaBOCStar} is a von Neumann algebra since $\comm{\BO(\MH)} = \C \cdot \id_\MH$ and $\comm{\{\C \cdot \id_\MH\}} = \BO(\MH)$, hence $\bicomm{\BO(\MH)} = \BO(\MH)$. To see this, note that if $A \in \comm{\BO(\MH)}$, then $A$ commutes, in particular, with every one-dimensional projection $P_\xi = \braket{\xi, \,\sbullet\ } \, \xi$, $\xi \in \MH$. This implies that there exists a constant $\alpha \in \C$ such that $A \xi = \alpha \xi$. Since $A$ is linear, it follows for all $\xi, \eta \in \MH$ that $A(\xi + \eta) = \alpha \xi + \beta \eta$ for some $\alpha, \beta \in \C$, but one also has $A(\xi + \eta) = \gamma (\xi + \eta)$ for $\gamma \in \C$, hence one may conclude $\alpha = \beta = \gamma$, that is, $A = \alpha \, \id_\MH$.
        
        The $C^\ast$-algebra of compact operators $\KO(\MH)$ from \cref{exa:operatorAlgebras_CStar} \ref{enu:operatorAlgebras_exaKOCStar} is not a von Neumann algebra because $\comm{\KO(\MH)} = \C \cdot \id_\MH$ as before, and hence $\bicomm{\KO(\MH)} = \BO(\MH)$. This follows since finite-rank operators are compact, and they are dense in $\BO(\MH)$.\qedhere
    \end{enumerate}
\end{examples}

Many applications in mathematical physics require working with \emph{unbounded operators} (see, for example, \cref{ch:perturbationTheory}). Due to the flexibility of von Neumann algebras, it turns out that one can \emph{affiliate} unbounded operators with this class of $C^\ast$-algebras.

\begin{definition}[Affiliated operator]\label{def:operatorAlgebras_affiliatedElement}
    Let $\MFM \subset \BO(\MH)$ be a von Neumann algebra. A closed linear operator $T : \MH \supset \dom(T) \longto \MH$ on $\MH$ (\cref{def:operators_closed}) is said to be \bemph{affiliated with $\MFM$}, and the set of these operators is denoted by $\AE{\MFM}$, iff for all $\comm{A} \in \comm{\MFM}$ it holds that
    \begin{equation*}
        \comm{A} \bigl(\dom(T)\bigr) \subset \dom(T) \tand \comm{A} T \subset T \comm{A} \ .
    \end{equation*}
\end{definition}

Essentially, a closed operator $T$ is affiliated with a von Neumann algebra $\MFM$ if and only if $T$ strongly commutes \cite[Prop. 3.70]{Moretti19} with the elements of the commutant $\comm{\MFM}$. The next important result shows that the above definition is somewhat natural \cite[Thm. 2.1]{DJP03}.

\begin{proposition}\label{pro:operatorAlgebras_characterizationAffiliation}
    Let $\MFM \subset \BO(\MH)$ be a von Neumann algebra and $T : \MH \supset \dom(T) \longto \MH$ be self-adjoint. It holds that $T \in \AE{\MFM}$ if and only if all bounded Borel functions of $T$ belong to the algebra $\MFM$.
\end{proposition}

\subsection{Cyclic and Separating Vectors}

The concept of a cyclic and separating vector $\xi \in \MH$ for a von Neumann algebra $\MFM \subset \BO(\MH)$ is of great importance for the theory and will be used extensively in the following three chapters.

\begin{para}[Closed subspaces and projections]\label{para:operatorAlgebras_closedSubspaceProjection}
    For arbitrary subsets $\MFA \subset \BO(\MH)$ and $\MK \subset \MH$, introduce the following notation:
    \begin{equation*}
        \MFA \MK \ce \lin_\K \set{A \xi \, : \, A \in \MFA,\, \xi \in \MK} \tand [\MFA \MK] \ce \clos_{\ndot} (\MFA \MK) \ .
    \end{equation*}
    By the correspondence between orthogonal projections and closed subspaces (\cref{rem:operators_orthogonalProjectionClosedSubspace}), there exists a unique orthogonal projection $P \in \PO(\MH)$ onto $[\MFA \MK]$. Often, the projection $P$ will also be denoted by the symbol $[\MFA \MK]$. In case that $\MK = \{\xi\}$ for some $\xi \in \MH$, one writes $\MFA \MK \equiv \MFA \xi$ and $[\MFA \MK] \equiv [\MFA \xi]$.
\end{para}

\begin{definition}[Cyclic and separating vectors]\label{def:operatorAlgebras_cyclicSeparating}
    Let $\MFM \subset \BO(\MH)$ be a von Neumann algebra. A set of vectors $\MS \subset \MH$ is said to be \bemph{cyclic} for the algebra $\MFM$ iff $[\MFM \MS] = \MH$, that is, if the subspace generated by $\MFM \MS$ lies dense in $\MH$ with respect to the norm topology. The set $\MS$ is called \bemph{separating} for $\MFM$ iff for every $A \in \MFM$, $A \xi = 0$ for all $\xi \in \MS$ implies that $A = 0$.
\end{definition}

\begin{proposition}\label{pro:operatorAlgebras_cyclicSeparating}
    Let $\MFM \subset \BO(\MH)$ be a von Neumann algebra and $\MS \subset \MH$ be a subset.
    \begin{enumerate}
        \item \label{enu:operatorAlgebras_cyclicMseparatingCommM} $\MS$ is cyclic for $\MFM$ if and only if $\MS$ is separating for $\comm{\MFM}$.
        
        \item \label{enu:operatorAlgebras_separatingMcyclicCommM} $\MS$ is separating for $\MFM$ if and only if $\MS$ is cyclic for $\comm{\MFM}$.
    \end{enumerate}
\end{proposition}

\begin{Proof}
    \tAd{} \ref{enu:operatorAlgebras_cyclicMseparatingCommM}. \cite[Prop. 2.5.3]{BR1} Assume first that $\MS$ is cyclic for $\MFM$. Let $\comm{A} \in \comm{\MFM}$ be arbitrary such that for every $\xi \in \MS$, $\comm{A} \xi = 0$. For all $A \in \MFM$, it follows that $ \comm{A} A \xi = A \comm{A} \xi = 0$, and hence $\comm{A} [\MFM \MS] = 0$ by linearity and continuity of $\comm{A}$. Since $[\MFM \MS] = \MH$ by assumption, one can conclude $\comm{A} = 0$ which shows that $\MS$ is separating for $\comm{\MFM}$.

    Conversely, assume that for every $\comm{A} \in \comm{\MFM}$ which satisfies $\comm{A} \xi = 0$ for all $\xi \in \MS$, it follows that $\comm{A} = 0$. Let $P \ce [\MFM \MS] \in \PO(\MH)$, then the claim $[\MFM \MS] = \MH$ is equivalent to the operator identity $P = \id_\MH$, see \cref{cor:operators_characterizationPerp} \ref{enu:operators_densePerp}. On the set $\MS$, it clearly holds that $\id_\MH - P = 0$ since $\id_\MH \in \MFM$ and hence $\MS \subset [\MFM \MS]$. Moreover, $P \in \comm{\MFM}$ \cite[p. 73]{BR1}: for all $A \in \MFM$ and $\eta \in \MH$, one has $A P \eta \in [\MFM \MS]$ by definition of this space, hence $A P = P A P$ (\cref{lem:operators_rangeProjection}). This gives $P A  = (A^\ast P)^\ast = (P A^\ast P)^\ast = P A P = A P$. Therefore, it follows that $\id_\MH - P \in \comm{\MFM}$ by \cref{lem:operatorAlgebras_propComm} which, together with the previous observation and the assumption, implies that $P = \id_\MH$.
    
    \tAd{} \ref{enu:operatorAlgebras_separatingMcyclicCommM}. If $\MS$ is separating for the algebra $\MFM$, it is for $\bicomm{\MFM}$ as well by \cref{thm:operatorAlgebras_bicommutant}. Using assertion \ref{enu:operatorAlgebras_cyclicMseparatingCommM} together with the fact that $\comm{\MFM}$ is a von Neumann algebra (\cref{lem:operatorAlgebras_propComm}), it follows that $\MS$ is cyclic for $\comm{\MFM}$. Conversely, if $\MS$ is cyclic for $\comm{\MFM}$, it is separating for $\bicomm{\MFM}$ by \ref{enu:operatorAlgebras_cyclicMseparatingCommM}, and hence also for $\MFM$.
\end{Proof}

\begin{para}[Cyclic projections]\label{para:operatorAlgebras_cyclicProjections}
    Let $\MFM \subset \BO(\MH)$ be a von Neumann algebra and $\xi \in \MH$ be arbitrary. Projections of the form $P_\xi \ce [\comm{\MFM} \xi]$ and $\comm{P}_\xi \ce [\MFM \xi]$ are called \bemph{cyclic projections}. It holds that $P_\xi \in \MFM$ and $\comm{P}_\xi \in \comm{\MFM}$. Indeed, in the proof of \cref{pro:operatorAlgebras_cyclicSeparating}, it was already shown that $\comm{P}_\xi \in \comm{\MFM}$. Therefore, it also follows that $P_\xi = [\comm{\MFM} \xi] \in \bicomm{\MFM} = \MFM$.

    A vector $\xi \in \MH$ is cyclic for the algebra $\MFM$ if and only if $\comm{P}_\xi = \id_\MH$, and it is separating for $\MFM$ if and only if $P_\xi = \id_\MH$. This can be seen as follows: $\xi$ is cyclic for $\MFM$ $\iff$ $[\MFM \xi] = \MH$ $\iff$ $\comm{P}_\xi = \id_\MH$ according to \cref{cor:operators_characterizationPerp} \ref{enu:operators_densePerp}. Similarly, $\xi$ is separating for $\MFM$ $\iff$ $\xi$ is separating for $\bicomm{\MFM}$ $\iff$ $\xi$ is cyclic for $\comm{\MFM}$ $\iff$ $[\comm{\MFM} \xi] = \MH$ $\iff$ $P_\xi = \id_\MH$ by \cref{pro:operatorAlgebras_cyclicSeparating}.
\end{para}

\section{States and Representations}\label{sec:operatorAlgebras_representationsStates}

The goal of this section is to argue that every abstract $C^\ast$-algebra can be realized as a subalgebra of $\BO(\MH)$ for a suitable Hilbert space $\MH$. To this end, the important notions of a \emph{state} on a $C^\ast$-algebra and a \emph{normal state} on a von Neumann algebra will be introduced and studied in some detail. The former will be employed in the so-called \emph{GNS-representation} of a $C^\ast$-algebra.

\subsection{\texorpdfstring{General States on a $C^\ast$-Algebra}{General States on a C*-Algebra}}

\begin{definition}[State]\label{def:operatorAlgebras_state}
    Let $(\MFA, \idm)$ be a unital $C^\ast$-algebra. A linear functional $\omega : \MFA \longto \C$ is called a \bemph{state} iff it is positive and $\omega(\idm) = 1$. The set of all states on the algebra $\MFA$ will be denoted by the symbol $\SS(\MFA)$. If, for every $A \in \MFA$, $\omega(A^\ast A) = 0$ implies $A = 0$, then the state $\omega$ is called \bemph{faithful}.
\end{definition}

The following proposition contains fundamental properties of positive linear functionals on $C^\ast$-algebras. Its proof can be found, for example, in \cite[Prop. 5]{BärFredenhagen09} or \cite[Lem. 2.3.10]{BR1}.

\begin{proposition}\label{pro:operatorAlgebras_formFromState}
    Let $\MFA$ be a $C^\ast$-algebra and $\omega : \MFA \longto \C$ be a positive linear functional. Then the mapping $\MFq_\omega : \MFA \times \MFA \longto \C$, $\MFq_\omega(A, B) \ce \omega(A^\ast B)$, is a positive semi-definite Hermitian sesquilinear form on $\MFA$. In particular, it satisfies the \textbf{Cauchy-Schwarz inequality}:
    \begin{equation}\label{eq:operatorAlgebras_CS}
        \forall A, B \in \MFA \ : \ \abs{\omega(A^\ast B)}^2 \le \omega(A^\ast A) \, \omega(B^\ast B) \ .
    \end{equation}
\end{proposition}

\begin{examples}\label{exa:operatorAlgebras_states}
    \leavevmode
    \begin{enumerate}[env]
        \item Let $n \in \N$ and $\MFA = \Mat(n; \C)$ from \cref{exa:operatorAlgebras_CStar} \ref{enu:operatorAlgebras_exaMatCStar}. Consider the trace functional $\tr : \MFA \longto \C$, $A \longmto \tr(A)$, which was shown to be a positive mapping in \cref{exa:operatorAlgebras_CP} \ref{enu:operatorAlgebras_exaTrCP}. A corresponding state $\tau : \MFA \longto \C$ can be obtained by setting $\tau(A) \ce \frac{1}{n} \tr(A)$ for all $A \in \MFA$.
        
        \item \label{enu:operatorAlgebras_exaMeasureState} Let $(X, \MFT)$ be a compact Hausdorff space, let $\MFA = C^0(X; \C)$ be the $C^\ast$-algebra of continuous functions on $X$ (\cref{exa:operatorAlgebras_CStar} \ref{enu:operatorAlgebras_exaCCStar}), and let $\mu$ be a measure on the Borel $\sigma$-algebra $\MB(X)$ of $X$ (that is, the smallest $\sigma$-algebra containing $\MFT$). Define a functional $\omega_\mu : \MFA \longto \C$ by
        \begin{equation*}
            \omega_\mu(f) \ce \int\nolimits_X f \diff \mu \quad (f \in \MFA) \ .
        \end{equation*}
        From the properties of the $\mu$-integral, it follows that $\omega_\mu$ is positive ($f \ge 0 \, \implies \, \int_X f \diff \mu \ge 0$). Furthermore, if $\mu$ is a probability measure, \ie{}, $\mu(X) = 1$, then $\omega_\mu$ is a state on $\MFA$.

        \item \label{enu:operatorAlgebras_exaVectorFunctional} Let $\MFA \subset \BO(\MH)$ be a unital concrete $C^\ast$-algebra (\cref{exa:operatorAlgebras_CStar} \ref{enu:operatorAlgebras_exaConcreteCstar}) and $\Omega \in \MH$ be a vector. Define a functional $\omega_\Omega : \MFA \longto \C$ by
        \begin{equation*}
            \omega_\Omega(A) \ce \braket{\Omega, A \Omega} \quad (A \in \MFA) \ .
        \end{equation*}
        This functional is positive since for all $A \in \MFA$, it holds that $\omega_\Omega(A^\ast A) = \braket{A \Omega, A \Omega} = \norm{A \Omega}^2 \ge 0$. One refers to $\omega_\Omega$ as the \bemph{vector functional} induced by the vector $\Omega$. If $\norm{\Omega} = 1$, then it follows that $\omega_\Omega(\id_\MH) = 1$, hence $\omega_\Omega$ defines a state on $\MFA$.

        \item \label{enu:operatorAlgebras_exaNormalState} Let $\MH$ be a separable Hilbert space, let $\MFA = \BO(\MH)$ (\cref{exa:operatorAlgebras_CStar} \ref{enu:operatorAlgebras_exaBOCStar}), and let $\rho \in \PE{\NO(\MH)}$ be a positive trace-class operator. Defining $\omega_\rho : \MFA \longto \C$ by
        \begin{equation*}
            \omega_\rho(A) \ce \tr(\rho A) \quad (A \in \MFA)
        \end{equation*}
        yields a positive linear functional: $\omega_\rho(A^\ast A) = \tr(\rho A^\ast A) = \tr(\rho^{1/2} A^\ast A \rho^{1/2}) \ge 0$ due to cyclicity and positivity of the trace. Furthermore, if $\rho$ is a density matrix (\cref{def:operators_densityMatrix}), \ie{}, 
        $\tr(\rho) = 1$, then $\omega_\rho$ is a state on $\MFA$. States of this form are called \bemph{normal states}. One can show that $\omega_\rho$ is faithful if and only if $\rho$ is invertible \cite[Exa. 2.5.5]{BR1}.
        
        If $\dim(\MH) < + \infty$, then every state on $\BO(\MH)$ is induced by a density matrix, \ie{}, $\SS(\BO(\MH)) = \DM(\MH)$ \cite[Thm. 2.7 \& 2.8]{Landsman17}. This is no longer true if $\dim(\MH) = + \infty$, see \cite[p. 109]{Landsman17}.\qedhere
    \end{enumerate}
\end{examples}

\subsection{Normal States on a von Neumann Algebra}\label{subsec:operatorAlgebras_normalStates}

In this subsection, normal states defined on an arbitrary von Neumann algebra will be investigated more closely.

\begin{definition}[Normal functional]\label{def:operatorAlgebras_normalState}
    Let $\MFM \subset \BO(\MH)$ be a von Neumann algebra and $\pdual{\MFM} \ce \cdual{(\MFM, \MFT_{uw})}$ be the set of all $\sigma$-weakly continuous linear functionals. A positive element of $\pdual{\MFM}$ is called a \bemph{normal functional}, and the set of all these is denoted by the symbol $\NF{\MFM}$. A normal functional $\omega$ satisfying $\omega(\id_\MH) = 1$ is called a \bemph{normal state}, and the set of these is denoted $\NS(\MFM)$.
\end{definition}

\begin{proposition}[\protect{\cite[Prop. 2.4.18]{BR1}}]
    Let $\MFM \subset \BO(\MH)$ be a von Neumann algebra. The space $\pdual{\MFM}$ is a Banach space with respect to the operator norm, and $\MFM$ is the dual of $\pdual{\MFM}$ with respect to the duality $\MFM \times \pdual{\MFM} \owns (A, \omega) \longto \omega(A) \in \C$, that is, $\MFM = \cdual{(\pdual{\MFM})}$.
\end{proposition}

\begin{remarks}
    \leavevmode
    \begin{enumerate}[env]
        \item The previous assertion establishes that every von Neumann algebra possesses a \bemph{predual space}. Generally, an abstract $C^\ast$-algebra which possesses a predual space is called a \bemph{$W^\ast$-algebra}.\footnotemark
        \footnotetext{The terminology \enquote{$W^\ast$-algebra} was also introduced by \textsc{I. E. Segal} in 1951 \cite[Def. 2.1]{Segal51b}.}
        It was shown by \textsc{S. Sakai} in 1956 that every $W^\ast$-algebra is a von Neumann algebra \cite{Sakai56}, \cite[Thm. 1.16.7]{Sakai98}, hence these two notions are actually equivalent.
        
        \item Let $\omega \in \pdual{\MFM}$ be $\sigma$-weakly continuous. From \cref{pro:operatorAlgebras_comparisonTopologies} and \cref{para:topologicalVectorSpaces_comparisonTopologies}, it follows that $\omega \in \cdual{\MFM}$. Furthermore, characterization of continuity in a locally convex space \cite[Cor. VIII.2.4]{Werner18} and the Hahn-Banach theorem \cite[Thm. III.1.4]{Werner18} imply that $\omega$ can be extended to a $\sigma$-weakly continuous linear functional on $\BO(\MH)$ so that \cref{pro:operatorAlgebras_sigmaWeakTopology} applies \cite[p. 75]{BR1}.
    \end{enumerate}
\end{remarks}

The next important result gives a characterization of normal states. Proofs of the assertions can be found in \cite[Thm. 2.4.21]{BR1}, \cite[Thm. 5.11]{StratilaZsido19}, and \cite[Prop. II.3.20]{Takesaki79}.

\begin{ntheorem}[Characterization of normal states]\label{thm:operatorAlgebras_characterizationNormalStates}
    Let $\omega \in \SS(\MFM)$. The following assertions are equivalent:
    \begin{enumerate}[equiv]
        \item $\omega$ is normal;
        
        \item there exist a sequence $(\xi_n)_{n \in \N} \in \ell^2(\N; \MH)$ such that $\omega = \sum_{n=1}^{\infty} \braket{\xi_n, \,\sbullet\ \xi_n}$;
        
        \item there exists a density matrix $\rho \in \DM(\MH)$ such that $\omega = \tr(\rho \ \sbullet\,)$;

        \item for every increasing net $(A_i)_{i \in I} \subset \PE{\MFM}$ with an upper bound (with respect to the order relation from \cref{para:operatorAlgebras_orderSE}), there holds $\omega\bigl(\sup_{i \in I} A_i\bigr) = \sup_{i \in I} \omega(A_i)$.
    \end{enumerate}
\end{ntheorem}

\begin{examples}\label{exa:operatorAlgebras_normalStates}
    \leavevmode
    \begin{enumerate}[env]
        \item \label{enu:operatorAlgebras_exaNormalStateAbelian} Let $\MFM \subset \BO(\MH)$ be a commutative von Neumann algebra on a separable Hilbert space $\MH$. According to a theorem of \textsc{I. E. Segal} \cite{Segal51b, Segal51}, there is a compact Hausdorff space $X$ and a measure $\mu$ on $\bigl(X, \MB(X)\bigr)$ such that $\MH \cong L^2(X, \mu)$ and $\MFM \cong L^\infty(X, \mu)$ \cite[Thm. IX.7.8]{Conway85}, \cite[Sect. I.7.3, Thm. 1]{Dixmier81}, \cite[Prop. 1.18.1]{Sakai98}. Since $\cdual{L^1(X, \mu)} = L^\infty(X, \mu)$ \cite[Thm. II.2.4]{Werner18}, it follows that the predual space of $\MFM$ is given by $\pdual{\MFM} = L^1(X, \mu)$. Explicitly, a function $h \in L^1(X, \mu)$ gives rise to the measure $\diff \nu \ce h \diff \mu$ which induces the functional $\omega_\nu(f) = \int_X f \diff \nu$, $f \in \MFM$ (\cref{exa:operatorAlgebras_states} \ref{enu:operatorAlgebras_exaMeasureState}). If, additionally, $h$ is assumed to be non-negative, then one can define $g \ce \sqrt{h} \in \MH$ to obtain
        \begin{equation*}
            \omega_\nu(f) = \int\nolimits_X f h \diff \mu = \int\nolimits_X \ol{g} f g \diff \mu = \braket{g, f g} = \omega_g(f) \ ,
        \end{equation*}
        that is, $\omega_\nu$ is equal to the vector functional $\omega_g$, $g \in \MH$ (\cref{exa:operatorAlgebras_states} \ref{enu:operatorAlgebras_exaVectorFunctional}). This functional is normal by \cref{thm:operatorAlgebras_characterizationNormalStates}, and, in fact, there is a one-to-one correspondence between $\omega \in \NF{\MFM}$ and $h_\omega \in L^1(X, \mu)_+$ \cite[p. 83]{BR1}. Since $g_\omega \ce \sqrt{h_\omega} \in L^2(X, \mu)$ is well-defined, one may also interpret this as a correspondence $\NF{\MFM} \cong L^2(X, \mu)_+$, where $\omega = \omega_{g_\omega} = \braket{g_\omega, \,\sbullet\ g_\omega}$ is a vector functional.

        \item For $\MFM = \BO(\MH)$, \cref{pro:operatorAlgebras_sigmaWeakTopology} showed that $\cdual{\NO(\MH)} \cong \BO(\MH)$ and $\pdual{\BO(\MH)} \cong \NO(\MH)$ (considering \cref{pro:topologicalVectorSpaces_dualSpaceWeakTopology}), hence $\cdual{(\pdual{\BO(\MH)})} \cong \BO(\MH)$. Thus, in this case the first equivalences of \cref{thm:operatorAlgebras_characterizationNormalStates} are clear. More generally, let $\MFM \subset \BO(\MH)$ be an arbitrary von Neumann algebra. \cref{thm:operatorAlgebras_characterizationNormalStates} implies that the functionals $\omega_\xi = \braket{\xi, \,\sbullet\ \xi}$, $\xi \in \MH$, and $\omega_\rho = \tr(\rho \ \sbullet\,)$, $\rho \in \DM(\MH)$, from \cref{exa:operatorAlgebras_states} \ref{enu:operatorAlgebras_exaVectorFunctional} and \ref{enu:operatorAlgebras_exaNormalState} are normal states on $\MFM$ in the sense of \cref{def:operatorAlgebras_normalState}.\qedhere
    \end{enumerate}
\end{examples}

\subsection{Support of a Normal Functional}

\begin{para}[Support projection]\label{para:operatorAlgebras_support}
    (\cite[Sect. 5.15]{StratilaZsido19})
    Let $\MFM \subset \BO(\MH)$ be a von Neumann algebra and $\varphi \in \NF{\MFM}$ be a normal functional on $\MFM$. Using spectral theory, one can show the following result \cite[p. 123]{StratilaZsido19}.

    \begin{nstatement}
        If $A \in \PE{\MFM}$ is a positive element in $\MFM$ such that $\varphi(A) = 0$, then it follows that $\varphi\bigl(\ssupp(A)\bigr) = 0$, where $\ssupp(A)$ denotes the support of $A$ (introduced in \cref{para:operators_supportPO}). Consequently, the family $\MF \ce \set{P \in \PO(\MFM) \, : \, \varphi(P) = 0}$ is increasingly directed.
    \end{nstatement}
   
    Denote the supremum of the family $\MF$ by $P_0 \ce \id_\MH - \ssupp(\varphi) \in \PO(\MFM)$, and let $(P_i)_{i \in I} \subset \MF$ be a net such that $P_i \uparrow P_0$. Using that $\varphi$ is normal (\cref{thm:operatorAlgebras_characterizationNormalStates}), it follows that
    \begin{equation*}
        \varphi(P_0) = \varphi\sumbra{\bigg}{\sup_{i \in I} P_i} = \sup_{i \in I} \varphi(P_i) = 0 \ .
    \end{equation*}
    This shows that $P_0$ is the largest projection in $\MFM$ on which $\varphi$ vanishes. Equivalently, the projection $\ssupp(\varphi) = \id_\MH - P_0$, which is called the \bemph{support} (\bemph{projection}) of $\varphi$, is the \emph{smallest projection in $\MFM$} satisfying the relation
    \begin{equation}\label{eq:operatorAlgebras_supportProjection}
        \varphi\bigl(\ssupp(\varphi)\bigr) = \varphi(\id_\MH) \ .
    \end{equation}
\end{para}

The following two identities are basic consequences of the definition of $\ssupp(\varphi)$ which are often useful. They are stated in \cite[Sect. 5.15]{StratilaZsido19} and \cite[Eq. (2.2)]{Stratila20}, respectively.

\begin{lemma}
    Let $\varphi \in \NF{\MFM}$. For all $A \in \MFM$, there holds
    \begin{equation}\label{eq:operatorAlgebras_supportInFunctional}
        \varphi(A) = \varphi\bigl(A \, \ssupp(\varphi)\bigr) = \varphi\bigl(\ssupp(\varphi) A\bigr) = \varphi\bigl(\ssupp(\varphi) A \, \ssupp(\varphi)\bigr) \ .
    \end{equation}
\end{lemma}

\begin{Proof}
    By linearity of $\varphi$, the definition of $\ssupp(\varphi)$, and the Cauchy-Schwarz inequality \eqref{eq:operatorAlgebras_CS} from \cref{pro:operatorAlgebras_formFromState}, it follows that
    \begin{align*}
        \vabs[\big]{\varphi(A) - \varphi\bigl(A \ssupp(\varphi)\bigr)}^2 = \vabs[\big]{\varphi\bigl(A (\id_\MH - \ssupp(\varphi))\bigr)} \le \varphi(A A^\ast) \, \varphi\bigl(\id_\MH - \ssupp(\varphi)\bigr) = 0
    \end{align*}
    which shows the first equality. The second identity is proved similarly, and third one follows from the first two identities.
\end{Proof}

\begin{lemma}\label{lem:operatorAlgebras_normalFunctionalZero}
    For all $A \in \MFM$, it holds that $\varphi(A^\ast A) = 0 \iff A \, \ssupp(\varphi) = 0$.
\end{lemma}

\begin{Proof}
    \textit{1.} Assume that $\varphi(A^\ast A) = 0$. From \cref{para:operatorAlgebras_support}, it follows that $\varphi\bigl(\ssupp(A^\ast A)\bigr) = 0$. Since $\ssupp(A^\ast A) = \rsupp(A)$ \cite[p. 26]{StratilaZsido19}, where $\rsupp(A) = \id_\MH - \bigl[\ker(A)\bigr]$ is the right support of $A$ (\cref{para:operators_supportPO}), and $\id_\MH - \ssupp(\varphi)$ is the largest projection in $\MFM$ contained in the kernel of $\varphi$, one obtains
    \begin{equation*}
        \id_\MH - \bigl[\ker(A)\bigr] \le \id_\MH - \ssupp(\varphi) \ .
    \end{equation*}
    This implies $\ssupp(\varphi) \le \bigl[\ker(A)\bigr]$, hence $\ran\bigl(\ssupp(\varphi)\bigr) \subset \ker(A)$ by definition of the order relation for projection operators (\cref{para:operators_orderPO}). It is now an immediate consequence that $A \ssupp(\varphi) = 0$.
    
    \textit{2.} Conversely, if one assumes that the last identity is satisfied for arbitrary $A \in \MFM$, one computes using \cref{eq:operatorAlgebras_supportInFunctional} that
    \begin{align*}
        \varphi(A^\ast A) = \varphi\bigl(A^\ast A \, \ssupp(\varphi)\bigr) = \varphi(0) = 0 \ . \tag*{\qedhere}
    \end{align*}
\end{Proof}

One of the most important applications of the support projection is to detect whether a normal functional is faithful. The following is stated in \cite[p. 10]{Stratila20} and \cite[p. 123]{StratilaZsido19}.

\begin{proposition}\label{pro:operatorAlgebras_normalFaithful}
    The normal functional $\varphi$ is faithful if and only if $\ssupp(\varphi) = \id_\MH$.
\end{proposition}

\begin{Proof}
    \pp{1.} Assume first that $\varphi$ is faithful, \ie{}, for all $A \in \MFM$ the implication $\varphi(A^\ast A) = 0 \, \implies \, A = 0$ holds true. By \cref{lem:operatorAlgebras_normalFunctionalZero}, this is equivalent to $A \, \ssupp(\varphi) = 0 \, \implies \, A = 0$ for all $A \in \MFM$. In particular, one can choose any projection $P \in \PO(\MFM)$ in place of $A$. But then the identity $P \, \ssupp(\varphi) = 0$ is equivalent to the projections $P$ and $\ssupp(\varphi)$ being orthogonal \cite[Prop. 3.17]{Moretti19}, hence the previous implication asserts that $\ssupp(\varphi)$ is only orthogonal to the zero projection, that is, $\ran(\ssupp(\varphi))^\perp = \{0\}$. \cref{cor:operators_characterizationPerp} \ref{enu:operators_densePerp} now implies that $\ssupp(\varphi) = \id_\MH$.
    
    \pp{2.} Conversely, assume that $\ssupp(\varphi) = \id_\MH$. Relying once again on \cref{lem:operatorAlgebras_normalFunctionalZero}, it follows that for all $A \in \MFM$, $\varphi(A^\ast A) = 0 \implies A = 0$, thus $\varphi$ is faithful. This completes the proof.
\end{Proof}

As another important application, the support of vector functionals, introduced in \cref{exa:operatorAlgebras_states} \ref{enu:operatorAlgebras_exaVectorFunctional}, shall be computed; this will yield cyclic projections (\cref{para:operatorAlgebras_cyclicProjections}) \cite[p. 131]{StratilaZsido19}.

\begin{proposition}\label{pro:operatorAlgebras_supportVectorFunctional}
    Let $\MFM \subset \BO(\MH)$ be a von Neumann algebra, let $\xi \in \MH$ be arbitrary, and let $\omega_\xi = \braket{\xi, \,\sbullet\ \xi}$ be the vector functional induced by $\xi$ on $\MFM$. Denote the same functional on the commutant $\comm{\MFM}$ of the algebra $\MFM$ by $\comm{\omega_\xi} = \braket{\xi, \,\sbullet\ \xi}$. Then
    \begin{equation*}
        \ssupp(\omega_\xi) = P_\xi = [\comm{\MFM} \xi] \tand \ssupp(\comm{\omega_\xi}) = \comm{P_\xi} = [\MFM \xi] \ .
    \end{equation*}
\end{proposition}

\begin{Proof}
    First, it is clear that $\omega_\xi(P_\xi) = \omega_\xi(\id_\MH)$ since $\xi \in [\comm{\MFM} \xi]$, hence $P_\xi$ satisfies the same property as the support $\ssupp(\omega_\xi)$. Therefore, it only remains to show that $P_\xi \le \ssupp(\omega_\xi)$ which, by the minimality of the support projection, then implies $P_\xi = \ssupp(\omega_\xi)$. Recall from \cref{para:operators_orderPO} that the relation $P_\xi \le \ssupp(\omega_\xi)$ is equivalent to $\ran(P_\xi) \subset \ran\bigl(\ssupp(\omega_\xi)\bigr)$. Let $\eta \in \ran(P_\xi)$ be arbitrary, that is, $\eta = \dlim{\MH}{n \to \infty} \comm{A_n} \xi$ for a sequence $(\comm{A_n})_{n \in \N} \subset \comm{\MFM}$. One computes
    \begin{align*}
        \ssupp(\omega_\xi) \eta = \dlim{\MH}{n \to \infty} \ssupp(\omega_\xi) \comm{A_n} \xi = \dlim{\MH}{n \to \infty} \comm{A_n} \, \ssupp(\omega_\xi) \xi = \dlim{\MH}{n \to \infty} \comm{A_n} \xi = \eta \ ,
    \end{align*}
    where first continuity of the support, then the fact that $\ssupp(\omega_\xi) \in \PO(\MFM) \subset \MFM$, and finally the identity $\ssupp(\omega_\xi) \xi = \xi$ were used. (Note that the last two properties hold by construction of the support projection.) From \cref{lem:operators_rangeProjection}, it now follows that $\eta \in \ran\bigl(\ssupp(\omega_\xi)\bigr)$, and thus the proof of the first identity is complete. For the second relation, one proceeds analogously.
\end{Proof}

\begin{corollary}\label{cor:operatorAlgebras_faithfulSeparating}
    Let $\MFM \subset \BO(\MH)$ be a von Neumann algebra, and let $\varphi = \omega_\xi$ be a vector functional on $\MFM$ induced by $\xi \in \MH$. Then $\varphi$ is faithful if and only if $\xi$ is separating for $\MFM$.
\end{corollary}

\begin{Proof}
    By \cref{pro:operatorAlgebras_normalFaithful}, the normal functional $\varphi$ is faithful if and only if $\ssupp(\varphi) = \id_\MH$, and by \cref{pro:operatorAlgebras_supportVectorFunctional}, the support takes the form $\ssupp(\varphi) = [\comm{\MFM} \xi]$. According to \cref{para:operatorAlgebras_cyclicProjections}, the condition $[\comm{\MFM} \xi] = \id_\MH$ is equivalent to $\xi$ being separating for $\MFM$.
\end{Proof}

\subsection{The GNS-Construction}\label{subsec:operatorAlgebras_GNS}

\begin{definition}\label{def:operatorAlgebras_representation}
    Let $\MFA$ be a $C^\ast$-algebra. A \bemph{$\ast$-representation} of $\MFA$ is pair $(\MH, \pi)$ consisting of a Hilbert space $\MH$ and a $\ast$-homomorphism $\pi : \MFA \longto \BO(\MH)$. It is called \bemph{faithful} iff $\pi$ is injective. Two $\ast$-representations $(\MH_1, \pi_1)$ and $(\MH_2, \pi_2)$ are said to be \bemph{unitarily equivalent} iff there exists a unitary operator $U : \MH_1 \longto \MH_2$ such that
    \begin{equation*}
        \forall A \in \MFA \ : \ \pi_2(A) = U \pi_1(A) U^\ast \ .
    \end{equation*}
    A \bemph{cyclic representation} of $\MFA$ is a triple $(\MH, \pi, \Omega)$ consisting of a $\ast$-representation $(\MH, \pi)$ of $\MFA$ and a cyclic vector $\Omega \in \MH$ for the set $\pi(\MFA) \subset \BO(\MH)$, that is, a vector such that $[\pi(\MFA) \Omega] = \MH$.
\end{definition}

Let $(\MFA, \idm)$ be a unital $C^\ast$-algebra and $\omega \in \SS(\MFA)$ be an arbitrary state. The goal of this subsection is to construct a specific cyclic representation $(\MH_\omega, \pi_\omega, \Omega_\omega)$ of $\MFA$ such that $\omega = \omega_{\Omega_\omega}$ for a suitable vector $\Omega_\omega \in \MH_\omega$. The following presentation closely follows \cite[Sect. 2.3.3 \& 2.3.4]{BR1}.

\begin{para}[Construction of the representation]\label{para:operatorAlgebras_constructionGNS}
    \emph{Step 1.} The idea to construct $\MH_\omega$ is the following: by definition, $\MFA$ is a Banach space, and according to \cref{pro:operatorAlgebras_formFromState}, the state $\omega$ may be used to define a semi-inner product $\bdot_\omega$ on $\MFA$, given for all $A, B \in \MFA$ by
    \begin{equation}\label{eq:operatorAlgebras_statePreIP}
        \braket{A, B}_\omega \ce \omega(A^\ast B) \ .
    \end{equation}
    This is not positive-definite because $\braket{A, A}_\omega = \omega(A^\ast A) = 0$ does not, in general, imply that $A = 0$. (Indeed, this is only the case for a \emph{faithful} state.) One has to factor out the set
    \begin{equation}\label{eq:operatorAlgebras_GNSideal}
        \MFJ_\omega \ce \set{A \in \MFA \ : \ \omega(A^\ast A) = 0}
    \end{equation}
    in order for $\bdot_\omega$ to become a proper inner product on $\MFA / \MFJ_\omega$. It turns out that $\MFJ_\omega$ is a left ideal of the algebra $\MFA$ \cite[p. 54]{BR1}, \cite[p. 19]{BärFredenhagen09}.

    For $A \in \MFA$, let $[A]_\omega \ce \set{A + B : B \in \MFJ_\omega}$ be the equivalence class of $A$ in $\MFA / \MFJ_\omega$. As usual, one defines a linear structure on this quotient by setting $[A]_\omega + [B]_\omega \ce [A + B]_\omega$ and $\lambda [A]_\omega \ce [\lambda A]_\omega$ which does not depend on the specific representatives of the classes. For all $A, B \in \MFA$, define
    \begin{equation*}
        \vbraket[\big]{[A]_\omega, [B]_\omega}_\omega \ce \braket{A, B}_\omega = \omega(A^\ast B) \ .
    \end{equation*}
    This is an inner product on $\MFA / \MFJ_\omega$: it is well-defined because $\MFJ_\omega$ is a left ideal of $\MFA$ \cite[p. 55]{BR1}, \cite[pp. 19 f.]{BärFredenhagen09}, and it is positive definite since $\vbraket[\big]{[A]_\omega, [A]_\omega}_\omega = 0$ implies $\omega(A^\ast A) = 0$, hence $A \in \MFJ_\omega$ and $[A]_\omega = [0]_\omega$. Let $\MH_\omega$ denote the completion of $\MFA / \MFJ_\omega$ with respect to the norm induced by $\bdot_\omega$. It is well-known that this space is a Hilbert space \cite[Thm. V.1.8 (c)]{Werner18}, and that $\MH_\omega$ contains the quotient $\MFA / \MFJ_\omega$ as a dense linear subspace \cite[p. 3 \& Cor. III.3.2]{Werner18}.

    \emph{Step 2.} The next step is to construct the representatives $\pi_\omega(A) \in \BO(\MH_\omega)$ for all $A \in \MFA$. First, define $\pi_\omega(A)$ on $\MFA / \MFJ_\omega$ by setting for all $B \in \MFA$:
    \begin{equation}\label{eq:operatorAlgebrasGNSpi}
        \pi_\omega(A) [B]_\omega \ce [AB]_\omega \ .
    \end{equation}
    It holds that $\pi_\omega : \MFA \longto \BO(\MFA / \MFJ_\omega)$ is a well-defined $\ast$-homomorphism \cite[p. 55]{BR1}, \cite[p. 20]{BärFredenhagen09}. From this and the fact that $\MFA / \MFJ_\omega$ lies dense in  $\MH_\omega$, it follows that $\pi_\omega(A) \in \BO(\MFA / \MFJ_\omega)$ can be uniquely extended to a bounded linear operator on $\MH_\omega$ (\cref{thm:operators_BLT}); this operator shall be denoted by the same symbol. Hence, one obtains a $\ast$-representation $\pi_\omega : \MFA \longto \BO(\MH_\omega)$.

    \emph{Step 3.} It remains to construct a vector $\Omega_\omega \in \MH_\omega$ such that $\omega = \omega_{\Omega_\omega}$ \cite[p. 55]{BR1}, \cite[p. 20]{BärFredenhagen09}. Define $\Omega_\omega \ce [\idm]_\omega \in \MH_\omega$. For all $A \in \MFA$, it follows that
    \begin{equation*}
        \braket{\Omega_\omega, \pi_\omega(A) \Omega_\omega}_\omega = \vbraket[\big]{[\idm]_\omega, \pi_\omega(A) [\idm]_\omega}_\omega = \vbraket[\big]{[\idm]_\omega, [A]_\omega}_\omega = \omega(\idm^\ast A) = \omega(A) \ .
    \end{equation*}
    Furthermore, since $\pi_\omega(A) \Omega_\omega = [A]_\omega$ for all $A \in \MFA$, one obtains that $\pi_\omega(\MFA) \Omega_\omega = \MFA / \MFJ_\omega$. By the very construction of the Hilbert space $\MH_\omega = \clos_{\ndot_\omega}(\MFA / \MFJ_\omega)$, this shows that the set $\pi_\omega(\MFA) \Omega_\omega$ lies dense in the latter. Therefore, $(\MH_\omega, \pi_\omega, \Omega_\omega)$ is a \emph{cyclic representation} of the $C^\ast$-algebra $\MFA$.
\end{para}

The following important theorem asserts that the representation constructed above is essentially unique \cite[Thm. 2.3.16]{BR1}. Motivated by an algebraic description of quantum theory, the theorem was originally obtained by \textsc{I. E. Segal} in 1947 \cite[Thm. 1]{Segal47}, building on the pioneering work of \textsc{I. M. Gelfand} and \textsc{M. A. Naimark} \cite{GelfandNaimark43}.

\begin{ntheorem}[Gelfand-Naimark-Segal construction]\label{thm:operatorAlgebras_GNS}
    Let $\MFA$ be a $C^\ast$-algebra and $\omega \in \SS(\MFA)$ be a state. Then there exists a cyclic representation $(\MH_\omega, \pi_\omega, \Omega_\omega)$ of $\MFA$ such that for every $A \in \MFA$, there holds $\omega(A) = \braket{\Omega_\omega, \pi_\omega(A) \Omega_\omega}_\omega$. This representation is, moreover, unique up to unitary equivalence.
\end{ntheorem}

\begin{Proof}
    The existence of the representation $(\MH_\omega, \pi_\omega, \Omega_\omega)$ was proved above in \cref{para:operatorAlgebras_constructionGNS}. Hence, only uniqueness must be shown. To this end, assume that $(\MH_\omega^\prime, \pi_\omega^\prime, \Omega_\omega^\prime)$ is another cyclic representation for $\MFA$ such that $\omega = \omega_{\Omega_\omega^\prime}$. Define an operator $U : \MFA / \MFJ_\omega \longto \MFA / \MFJ_\omega^\prime$ by
    \begin{equation}\label{eq:operatorAlgebras_constructionUnitaryEquiv}
        U \pi_\omega(A) \Omega_\omega \ce \pi_\omega^\prime(A) \Omega_\omega^\prime \quad (A \in \MFA) \ .
    \end{equation}
    Note that this mapping is clearly linear due to the linearity of $\pi_\omega$ and $\pi_\omega^\prime$, and it satisfies the following relation for all $A, B \in \MFA$:
    \begin{align*}
        \vbraket[\big]{U \pi_\omega(A) \Omega_\omega, U \pi_\omega(B) \Omega_\omega}_\omega &= \vbraket[\big]{\pi_\omega^\prime(A) \Omega_\omega^\prime, \pi_\omega^\prime(B) \Omega_\omega^\prime}_\omega = \omega(A^\ast B) = \vbraket[\big]{\pi_\omega(A) \Omega_\omega, \pi_\omega(B) \Omega_\omega}_\omega \ .
    \end{align*}
    Hence, $U$ is an isometry (so in particular injective and well-defined) with dense domain $\pi_\omega(\MFA) \Omega_\omega \subset \MH_\omega$ and with dense range $\pi_\omega^\prime(\MFA) \Omega_\omega^\prime \subset \MH_\omega^\prime$. Therefore, it can be extended uniquely to a unitary operator $U : \MH_\omega \longto \MH_\omega^\prime$ (\cref{thm:operators_BLT}). Finally, for $A, B \in \MFA$ one computes using $U^\ast U = \id_{\MH_\omega}$ that
    \begin{align*}
        \pi_\omega^\prime(A) \pi_\omega^\prime(B) \Omega_\omega^\prime &= U \pi_\omega(AB) \Omega_\omega = U \pi_\omega(A) U^\ast U \pi_\omega(B) \Omega_\omega = U \pi_\omega(A) U^\ast \pi_\omega^\prime(B) \Omega_\omega^\prime \ ,
    \end{align*}
    hence $\pi_\omega^\prime(A) = U \pi_\omega(A) U^\ast$ on the dense subspace $\pi_\omega^\prime(\MFA) \Omega_\omega^\prime$; by continuous extension, this identity holds true also on $\MH_\omega^\prime$ which shows that the two cyclic representations are unitarily equivalent.
\end{Proof}

\begin{definition}[GNS-representation]\label{def:operatorAlgebras_GNS}
    Let $\MFA$ be a $C^\ast$-algebra and $\omega \in \SS(\MFA)$ be a state. The unique cyclic representation $(\MH_\omega, \pi_\omega, \Omega_\omega)$ constructed from $\omega$ in \cref{thm:operatorAlgebras_GNS} is termed \bemph{GNS-representation} or \bemph{canonical cyclic representation} of $\MFA$ associated with $\omega$.
\end{definition}

The following two corollaries can be obtained immediately from the GNS-construction. They concern states with additional properties, namely \emph{faithful} and \emph{invariant states}, and they will be needed later on. The second corollary is from \cite[Cor. 2.3.17]{BR1}.

\begin{corollary}\label{cor:operatorAlgebras_GNSfaithful}
    Let $\MFA$ be a $C^\ast$-algebra and $\omega \in \SS(\MFA)$ be a state. If $\omega$ is faithful, then the GNS-representation $(\MH_\omega, \pi_\omega, \Omega_\omega)$ is faithful.
\end{corollary}

\begin{Proof}
    By \cref{def:operatorAlgebras_representation}, it has to be verified that $\pi_\omega$ is injective. To this end, let $A \in \ker(\pi) \subset \MFA$ be arbitrary. Then it follows that $\omega(A^\ast A) = \braket{\Omega_\omega, \pi(A^\ast A) \Omega_\omega}_\omega = \norm{\pi(A) \Omega_\omega}^2 = 0$, hence one may conclude that $A = 0$ because $\omega$ is faithful. This shows that $\ker(\pi) = \{0\}$.
\end{Proof}

\begin{corollary}\label{cor:operatorAlgebras_unitaryInvariantState}
    Let $\MFA$ be a $C^\ast$-algebra, let $\omega \in \SS(\MFA)$ be a state on $\MFA$, and let $\tau \in \Aut(\MFA)$ be a $\ast$-automorphism on $\MFA$ which leaves $\omega$ invariant, that is, $\omega\bigl(\tau(A)\bigr) = \omega(A)$ for all $A \in \MFA$. Then there exists a uniquely defined unitary operator $U_\omega : \MH_\omega \longto \MH_\omega$, where $\MH_\omega$ is the representation space of the cyclic representation $(\MH_\omega, \pi_\omega, \Omega_\omega)$ of $\MFA$, such that
    \begin{equation*}
        U_\omega \Omega_\omega = \Omega_\omega \tand U_\omega \pi_\omega(A) U_\omega^\ast = \pi_\omega\bigl(\tau(A)\bigr) \quad (A \in \MFA) \ .
    \end{equation*}
\end{corollary}

\begin{Proof}
    This is a consequence of the uniqueness of the cyclic representation, applied to the representation $(\MH_\omega, \pi_\omega^\prime, \Omega_\omega)$ with $\pi_\omega^\prime \ce \pi_\omega \circ \tau$. Note, in particular, that \cref{eq:operatorAlgebras_constructionUnitaryEquiv} determines the unitary operator $U_\omega : \MH_\omega \longto \MH_\omega$ uniquely due to the fact that $\pi_\omega(\MFA) \Omega_\omega$ lies dense in $\MH_\omega$.
\end{Proof}

The following theorem, which was first obtained in \cite[Thm. 1]{GelfandNaimark43} and which is a cornerstone of $C^\ast$-algebra theory, shows that $\ast$-representations actually exist by combining the GNS-representation with the Hahn-Banach theorem. The proof can be found, \eg{}, in \cite[Thm. 2.1.10]{BR1} or \cite[Thm. IX.3.15]{Werner18}.

\begin{ntheorem}[Gelfand-Naimark; non-commutative version]\label{thm:operatorAlgebras_structureCStar}
    Let $\MFA$ be a $C^\ast$-algebra. Then there exists a Hilbert space $\MH$ such that $\MFA$ is isometrically $\ast$-isomorphic to a uniformly closed $\ast$-subalgebra of $\BO(\MH)$.
\end{ntheorem}

\begin{remark}\label{rem:operatorAlgebras_concreteCStar}
    In \cref{exa:operatorAlgebras_CStar} \ref{enu:operatorAlgebras_exaConcreteCstar}, a \emph{concrete $C^\ast$-algebra} was defined to be a uniformly closed $\ast$-subalgebra of $\BO(\MH)$, and it was shown that such an algebra is an abstract $C^\ast$-algebra in the sense of \cref{def:operatorAlgebras_CStar}. Conversely, the fundamental \cref{thm:operatorAlgebras_structureCStar} asserts that {\itshape every abstract $C^\ast$-algebra is $\ast$-isomorphic to a concrete $C^\ast$-algebra}. Therefore, combining these results, it follows that abstract and concrete $C^\ast$-algebras are equivalent notions.
\end{remark}

\begin{examples}\label{exa:operatorAlgebras_GNS}
    \leavevmode
    \begin{enumerate}[env]
        \item \label{enu:operatorAlgebras_exaGNSC} Let $\MFA$ be a commutative unital $C^\ast$-algebra. According to the commutative version of the Gelfand-Naimark theorem (see \cite[Lem. 1]{GelfandNaimark43} as well as \cite[Thm. 2.1.11]{BR1} and \cite[Thm. IX.3.4]{Werner18}), there exists a compact Hausdorff space $X$ such that $\MFA \cong C^0(X)$. Let $\omega \in \SS(\MFA)$ be an arbitrary state on $\MFA$. Then from the Riesz-Markov-Kakutani representation theorem \cite[Thm. C.17]{Conway85}, \cite[Thm. II.2.5]{Werner18}, one obtains a unique measure $\mu$ on the space $\bigl(X, \MB(X)\bigr)$ such that $\mu(X) = 1$ and $\omega = \omega_\mu$ (\cf{} \cref{exa:operatorAlgebras_states} \ref{enu:operatorAlgebras_exaMeasureState}), that is
        \begin{equation*}
            \omega(f) = \int\nolimits_X f \diff \mu \quad (f \in \MFA) \ .
        \end{equation*}
        
        To obtain the GNS-representation $(\MH_\omega, \pi_\omega, \Omega_\omega)$ of $\MFA$, note that the ideal $\MFJ_\omega$ from \cref{eq:operatorAlgebras_GNSideal} is given by $\MFJ_\omega = \set{f \in \MFA \, : \, \int_X \abs{f}^2 \diff \mu = 0}$, and that the inner product $\bdot_\omega$ on $\MFA / \MFJ_\omega$ takes the form $\braket{f, g}_\omega = \int_X \ol{f} g \diff \mu$. Since this is the usual $L^2$-inner product, it follows that $\MH_\omega = L^2(X, \mu)$. According to \cref{eq:operatorAlgebrasGNSpi}, the representation $\pi_\omega : \MFA \longto \BO(\MH_\omega)$ acts like $\pi_\omega(f) g = fg$ for $f \in \MFA$ and $g \in \MH_\omega$, that is, $\pi_\omega(f) = M_f$ is the multiplication operator with the function $f$. Finally, the vector $\Omega_\omega \in \MH_\omega$ is identified with the constant function $\1_X \in \MFA$. Thus, the GNS-representation of $(C^0(X), \omega_\mu)$ is given by
        \begin{equation*}
            \bigl(L^2(X, \mu),\, M_{\sbullet},\, \1_X\bigr) \ .
        \end{equation*}

        \item \label{enu:operatorAlgebras_exaGNSBO} Let $\MH$ be a separable Hilbert space and $\MFA = \BO(\MH)$. Consider a state $\omega \in \SS(\MFA)$ given by a density matrix $\rho \in \DM(\MH)$ as in \cref{exa:operatorAlgebras_states} \ref{enu:operatorAlgebras_exaNormalState}, that is,
        \begin{equation*}
            \omega(A) = \tr(\rho A) \quad (A \in \MFA) \ .
        \end{equation*}
        Assume, for simplicity, that $\rho$ is invertible, \ie{}, that $\omega$ is faithful. The more general case of singular $\rho$ is discussed, \eg{}, in \cite[Exa. 4.32]{AJP06}.
        
        The pre-inner product $\bdot_\omega$ from \cref{eq:operatorAlgebras_statePreIP} takes the form
        \begin{equation}\label{eq:operatorAlgebras_GNSBO}
            \braket{A, B}_\omega = \omega(A^\ast B) = \tr(\rho A^\ast B) = \tr\bigl((A \rho^{1/2})^\ast \, B \rho^{1/2}\bigr) = \vbraket[\big]{A \rho^{1/2}, B \rho^{1/2}}_\mathrm{HS}
        \end{equation}
        for all $A, B \in \MFA$, where in the last step \eqref{eq:operators_HSinnerProduct} was used. Note that since $\rho$ is trace-class, $\rho^{1/2} \in \HS(\MH)$ is a Hilbert-Schmidt operator, hence the above expression is well-defined because $\HS(\MH) \subset \BO(\MH)$ is an ideal (\cref{pro:operators_NOHS}). The left ideal $\MFJ_\omega$ from \cref{eq:operatorAlgebras_GNSideal} used in the GNS-construction is trivial since $\omega$ is faithful. Therefore, the GNS-representation space $\MH_\omega$ is given by the completion of $\MFA$ with respect to the norm induced by $\bdot_\omega$, the $\ast$-homomorphism $\pi_\omega : \MFA \longto \BO(\MH_\omega)$ can be defined by $\pi_\omega(A) B = AB$ for all $A, B \in \MFA$, and finally the vector representative of $\omega$ is given by $\Omega_\omega = \id_\MH \in \MH_\omega$.

        Define the following linear and unitary mapping:
        \begin{equation*}
            U :
            \begin{cases}
                \begin{aligned}
                    \bigl(\MFA, \bdot_\omega\bigr) &\longto \bigl(\HS(\MH), \bdot_\mathrm{HS}\bigr) \ , \\
                    A &\longmto A \rho^{1/2} \ .
                \end{aligned}
            \end{cases}
        \end{equation*}
        The inverse $\inv{U} : \HS(\MH) \longto \MFA$ is given by $\inv{U} B = B \rho^{-1/2}$, and $U$ is an isometry by \eqref{eq:operatorAlgebras_GNSBO}. Hence, by uniqueness of the GNS-representation (\cref{thm:operatorAlgebras_GNS}), one may identify $(\MH_\omega, \bdot_\omega)$ with $(\HS(\MH), \bdot_\mathrm{HS})$, the transformed $\ast$-automorphism $\wt{\pi}_\omega : \MFA \longto \BO\bigl(\HS(\MH)\bigr)$, $\wt{\pi}_\omega(A) = U \pi_\omega(A) \inv{U}$, still takes the form $\wt{\pi}_\omega(A) B = A B$ for all $A \in \MFA$, $B \in \HS(\MH)$, and the vector representative $\Omega_\omega$ of $\omega$ is transformed to $U \Omega_\omega = \rho^{1/2}$. Thus, the GNS-representation of $(\BO(\MH), \omega_\rho)$ is given by
        \begin{equation*}
            \bigl(\HS(\MH),\, M_{\sbullet},\, \rho^{1/2}\bigr) \ . \tag*{\qedhere}
        \end{equation*}
    \end{enumerate}
\end{examples}

\chapter{Modular Theory in von Neumann Algebras}\label{ch:vonNeumann}

This chapter is devoted to \emph{Tomita-Takesaki modular theory} and related notions which have revolutionized both the purely mathematical study of von Neumann algebras \cite[p. 276]{StratilaZsido19} as well their applications to algebraic quantum field theory \cite{Borchers00}, \cite[p. 5]{BuchholzFredenhagen23} and quantum statistical mechanics \cite{BR2}. First, \cref{sec:vonNeumann_TomitaTakesaki} introduces the \emph{modular data} for $\sigma$-finite von Neumann algebras. This framework is then used to define the \emph{standard form} representation in \cref{sec:vonNeumann_standardForm}. Building on that, \cref{sec:vonNeumann_relativeModularOperator} discusses \emph{relative modular operators} and some of their properties. Finally, \cref{sec:vonNeumann_spatialDerivative} treats the \emph{spatial derivative} which is a generalization of the relative modular operator.

\section{Tomita-Takesaki Modular Theory}\label{sec:vonNeumann_TomitaTakesaki}

In this section, the modular theory of von Neumann algebras, which was developed first by \textsc{M. Tomita} in two unpublished notes in 1967 and then presented in a didactically enriched and mathematically extended version by \textsc{M. Takesaki} in 1970 \cite{Takesaki70}, will be outlined.

\subsection{\texorpdfstring{$\sigma$-finite von Neumann Algebras}{σ-finite von Neumann Algebras}}

For applications to mathematical physics, the following type of von Neumann algebras turns out to be particularly relevant \cite[p. 84]{BR1}.

\begin{definition}[$\sigma$-finite von Neumann algebra]\label{def:vonNeumann_sigmaFinite}
    A von Neumann algebra $\MFM \subset \BO(\MH)$ is said to be \bemph{$\sigma$-finite} iff any family $(P_i)_{i \in I} \subset \PO(\MFM) \setminus \{0\}$ of mutually orthogonal non-zero projections in $\MFM$ (that is, $P_i P_j = 0$ for all $i \neq j$) is at most countable.
\end{definition}

$\sigma$-finite von Neumann algebras are characterized by the following theorem \cite[Prop. 2.5.6]{BR1}.

\begin{ntheorem}[Characterization of $\sigma$-finite von Neumann algebras]\label{thm:vonNeumann_characterizationSigmaFinite}
    Let $\MFM \subset \BO(\MH)$ be a von Neumann algebra. The following assertions are equivalent:
    \begin{enumerate}[equiv]
        \item $\MFM$ is $\sigma$-finite;
        \item there exists a countable subset which is separating for $\MFM$;
        \item there exists a faithful normal state $\omega \in \NS(\MFM)$;
        \item $\MFM$ is isomorphic to a von Neumann algebra which contains a cyclic and separating vector.
    \end{enumerate}
\end{ntheorem}

\begin{example}\label{exa:vonNeumann_BOnormalState}
    Consider the von Neumann algebra $\MFM = \BO(\MH)$ over a separable Hilbert space $\MH$. Every normal state $\omega \in \NS(\MFM)$ is of the form $\omega = \tr(\rho \ \sbullet\,)$ for some $\rho \in \DM(\MH)$ (\cref{thm:operatorAlgebras_characterizationNormalStates}). Furthermore, in \cref{exa:operatorAlgebras_states} \ref{enu:operatorAlgebras_exaNormalState} it was mentioned that $\omega$ is faithful if and only if $\rho$ is invertible. One can proceed to show that such a state exists, \ie{}, that $\MFM$ is $\sigma$-finite, if and only if $\MH$ is separable \cite[Exa. 2.5.5]{BR1}.
\end{example}

\subsection{The Tomita Operator}\label{subsec:vonNeumann_TomitaOperator}

Let $\MFM \subset \BO(\MH)$ be a $\sigma$-finite von Neumann algebra. According to \cref{thm:vonNeumann_characterizationSigmaFinite}, there exists a cyclic and separating vector $\Omega \in \MH$, and thus the mapping $\MFM \owns A \longmto A \Omega \in \MH$ is a linear bijection. One may use it to transfer the algebra involution $\MFM \owns A \longmto A^\ast \in \MFM$, which is an isometry, to an operation $A \Omega \longmto A^\ast \Omega$ on $\MFM \Omega$ which is, in general, not an isometry. The analysis of this mapping is the starting point of Tomita-Takesaki modular theory.

\begin{para}[Construction of the Tomita operator]\label{para:vonNeumann_preclosedTomitaOperator}
    (\cite[pp. 13 f.]{Hiai21})
    Consider two anti-linear operators (\cf{} \cref{para:operators_antiLinear}) $S_0 : \MH \supset \dom(S_0) \longto \MH$ and $F_0 : \MH \supset \dom(F_0) \longto \MH$ with domains $\dom(S_0) \ce \MFM \Omega$ and $\dom(F_0) \ce \comm{\MFM} \Omega$; note that by assumption on $\Omega$ (see also \cref{pro:operatorAlgebras_cyclicSeparating}), both $S_0$ and $F_0$ are densely defined operators. For every $A \in \MFM$ and every $\comm{A} \in \comm{\MFM}$, define their action by
    \begin{equation}\label{eq:vonNeumann_preclosedTomitaOperator}
        S_0 A \Omega \ce A^\ast \Omega \tand F_0 \comm{A} \Omega \ce (\comm{A})^\ast \Omega \ .
    \end{equation}

    \begin{nstatement}
        $S_0$ and $F_0$ are closable anti-linear operators on $\MH$.
    \end{nstatement}

    \begin{Proof}
        From the definition \eqref{eq:vonNeumann_preclosedTomitaOperator} and the properties of the Hilbert-space adjoint, anti-linearity immediately follows. Let $A \in \MFM$ and $\comm{B} \in \comm{\MFM}$ be arbitrary. One computes
        \begin{align*}
            \braket{A \Omega, F_0 \comm{B} \Omega} &= \braket{A \Omega, (\comm{B})^\ast \Omega} = \braket{\comm{B} A \Omega, \Omega} \\
            &= \braket{A \comm{B} \Omega, \Omega} = \braket{\comm{B} \Omega, A^\ast \Omega} = \braket{\comm{B} \Omega, S_0 A \Omega} \ .
        \end{align*}
        By the definition of the adjoint of an anti-linear operator (\cref{para:operators_antiLinear}), this shows that $A \Omega \in \dom(F_0^\ast)$ and $F_0^\ast A \Omega = S_0 A \Omega$; since $A \in \MFM$ was arbitrary, this is to say that $S_0 \subset F_0^\ast$. As the adjoint is always a closed operator (\cref{pro:operators_propertiesAdjoint} \ref{enu:operators_adjointClosed}), it follows that $S_0$ is closable. By reading the above chain of equations backwards, one sees that $F_0 \subset S_0^\ast$, hence $F_0$ is closable as well.
    \end{Proof}

    Consider the anti-linear operators
    \begin{equation}\label{eq:vonNeumann_closedTomitaOperator}
        S \ce \ol{S_0} = \dadj{S_0} \tand F \ce S^\ast = S_0^\ast \ ,
    \end{equation}
    where the bar denotes the closure of the operator $S_0$, see \cref{para:operators_closure}. (Recall also \cref{pro:operators_propertiesAdjoint} \ref{enu:operators_adjointClosure} for the second identities in both definitions.) The operator $S$ is called the \bemph{Tomita operator} of the von Neumann algebra $\MFM$ with respect to the cyclic and separating vector $\Omega$, and it is also denoted by $S_\Omega$ or $S_\omega$ to emphasize the dependence on $\Omega$ or $\omega = \omega_\Omega$.
\end{para}

\begin{lemma}\label{lem:vonNeumann_propertiesTomitaOperator}
    The operator $S$ and its adjoint $S^\ast$ enjoy the following properties:
    \begin{enumerate}
        \item \label{enu:vonNeumann_Sinvertible} $S$ is invertible with $\inv{S} = S$.
        \item \label{enu:vonNeumann_ranSdense} The ranges of $S$ and $S^\ast$ lie dense in $\MH$
    \end{enumerate}
\end{lemma}

\begin{Proof}
    \tAd{} \ref{enu:vonNeumann_Sinvertible}. By definition, it holds that $\ran(S_0) = \dom(S_0)$ and $S_0^2 = \id_\MH$. Therefore, $S_0$ is invertible with $\inv{S_0} = S_0$. From \cref{pro:operators_propertiesAdjoint} \ref{enu:operators_inverseClosure}, it follows that $\inv{S_0}$ is closable with
    \begin{equation*}
        \inv{S} = \inv{(\ol{S_0})} = \ol{(\inv{S_0})} = S \ .
    \end{equation*}
    
    \tAd{} \ref{enu:vonNeumann_ranSdense}. First, it will be shown that $\dom(S_0) \subset \ran(S)$. To this end, let $\xi \in \dom(S_0)$ be arbitrary, \ie{}, $\xi = A \Omega$ for some $A \in \MFM$. Since $S$ is an extension of $S_0$, it follows that $\dom(S_0) \subset \dom(S)$ and $S B \Omega = S_0 B \Omega = B^\ast \Omega$ for all $B \in \MFM$. Thus, since $\eta = A^\ast \Omega \in \dom(S_0)$, one obtains $S \eta = S_0 \eta = A \Omega = \xi$, hence $\xi \in \ran(S)$.

    To see that the range of the adjoint $S^\ast$ is dense in $\MH$ as well, first note that since $S_0$ is closable, it follows from \cref{pro:operators_propertiesAdjoint} \ref{enu:operators_closability} that $\dom(S_0^\ast) = \dom(S^\ast)$ is dense in $\MH$. Hence, one can use \ref{enu:operators_ranTperp} of the same proposition and the definition of $S$ to conclude that $\ran(S^\ast)^\perp = \ker(S^{\ast \ast}) = \ker{S} = \set{0}$ by assertion \ref{enu:vonNeumann_Sinvertible} proved above. From \cref{cor:operators_characterizationPerp} \ref{enu:operators_densePerp}, the claim now follows.
\end{Proof}

\begin{para}[Polar decomposition of the Tomita operator]\label{para:vonNeumann_polarDecompositionTomita}
    One can compute the polar decomposition of the closed operator $S$ (\cf{} \cref{thm:operators_polarDecompositionCO}) to obtain a unique partial isometry $J : \MH \longto \MH$ and a unique positive self-adjoint linear operator $\Delta : \MH \supset \dom(\Delta) \longto \MH$ such that
    \begin{equation}\label{eq:vonNeumann_modularOperator}
        S = J \Delta^{1/2} \twith \Delta \ce S^\ast S = F S \ .
    \end{equation}

    \begin{nstatement}
        $J$ is an anti-linear unitary operator, and $\Delta$ is invertible.
    \end{nstatement}

    \begin{Proof}
        The operator $J$ is anti-linear and maps its initial space $\ker(S)^\perp = [\ran(S^\ast)]$ isometrically onto its final space $\ker(S^\ast)^\perp = [\ran(S)]$, and it satisfies $J \equiv 0$ on $\ker(S)^{\perp\perp} = [\ker(S)]$. From \cref{lem:vonNeumann_propertiesTomitaOperator}, it follows that both initial and final space of $J$ are given by $\MH$, and that $\ker(S) = \set{0}$. Therefore, $J$ maps $\MH$ isometrically to $\MH$ and is thus a unitary operator \cite[p. 20]{Conway85}. To see that $\Delta = S^\ast S$ is invertible, note that from \cref{pro:operators_propertiesAdjoint} \ref{enu:operators_inverseAdjoint}, it follows that $S^\ast$ is invertible with $(S^\ast)^{-1} = (S^{-1})^\ast$. Hence, $\Delta$ is invertible with $\inv{\Delta} = \inv{S} \inv{(S^\ast)} = \inv{S} (\inv{S})^\ast = S S^\ast$.
    \end{Proof}
\end{para}

\subsection{Modular Data and Tomita's Theorem}\label{subsec:vonNeumann_modularData}

The observations made in the previous \cref{para:vonNeumann_polarDecompositionTomita} motivate the following definition.

\begin{definition}[Modular data]\label{def:vonNeumann_modularData}
    The operator $J$ is called the \bemph{modular conjugation} with respect to the von Neumann algebra $\MFM \subset \BO(\MH)$ and the cyclic and separating vector $\Omega \in \MH$, and the operator $\Delta$ is called the \bemph{modular operator} with respect to $(\MFM, \Omega)$. To emphasize the dependence on the vector $\Omega$ or the functional $\omega = \omega_\Omega$, one also writes $J_\Omega, J_\omega$ and $\Delta_\Omega, \Delta_\omega$.
\end{definition}

The next proposition contains some of the most important properties of the modular data; they will be used very often in the following discussions (sometimes implicitly). The proof is an elaborated version of \cite[Lem. 2.1]{Hiai21} with some insights from \cite[p. 278]{StratilaZsido19}.

\begin{proposition}\label{pro:vonNeumann_propertiesModularData}
    Let $\MFM \subset \BO(\MH)$ be a von Neumann algebra with cyclic and separating vector $\Omega \in \MH$, and let $(J, \Delta)$ denote the modular data with respect to $(\MFM, \Omega)$.
    \begin{enumerate}
        \item \label{enu:vonNeumann_propertiesJ} $J^\ast = J$ and $J^2 = \id_\MH$.
        
        \item \label{enu:vonNeumann_deltaFS} $\Delta = FS$ and $\inv{\Delta} = SF$.
        
        \item \label{enu:vonNeumann_modDataTomitaOp} $S = J \Delta^{1/2} = \Delta^{-1/2} J$ and $F = J \Delta^{-1/2} = \Delta^{1/2} J$.
        
        \item \label{enu:vonNeumann_inverseModOp} $\inv{\Delta} = J \Delta J$ and $J \Delta^{\ii t} = \Delta^{\ii t} J$ for all $t \in \R$.
        
        \item \label{enu:vonNeumann_actionModDatOmega} $J \Omega = \Omega$ and $\Delta \Omega = \Omega$.
        
        \item \label{enu:vonNeumann_funcCalcModOp} If $f : [0, + \infty) \longto \C$ is Borel-measurable, then $f(\Delta) \Omega = f(1) \Omega$.
    \end{enumerate}
\end{proposition}

\begin{Proof}
    Since $S = \inv{S}$ by \cref{lem:vonNeumann_propertiesTomitaOperator} \ref{enu:vonNeumann_Sinvertible} and $\inv{J} = J^\ast$ by \cref{para:vonNeumann_polarDecompositionTomita}, it follows that
    \begin{equation*}
        S = J \Delta^{1/2} = \bigl(J \Delta^{1/2}\bigr)^{-1} = \Delta^{-1/2} \inv{J} = \Delta^{-1/2} J^\ast = J^\ast J \Delta^{- 1/2} J^\ast \ .
    \end{equation*}
    From this relation, uniqueness of the polar decomposition, and the fact that $J^\ast$ is a partial isometry as well \cite[p. 162]{Schmüdgen12}, one can conclude that $J = J^\ast$ and $\Delta^{1/2} = J \Delta^{- 1/2} J^\ast$. The first identity (together with $J^\ast = J^{-1}$) establishes \ref{enu:vonNeumann_propertiesJ}, and from the second, one can conclude the first part of \ref{enu:vonNeumann_inverseModOp}: using that $\Delta^{1/2} = J \Delta^{-1/2} J$ and $J^2 = \id_\MH$, one computes
    \begin{equation*}
        J \Delta J = J (J \Delta^{-1/2} J) (J \Delta^{-1/2} J) J = \inv{\Delta} \ .
    \end{equation*}
    
    Let $f : [0, + \infty) \longto \C$ be any Borel-measurable function. From the previous identity and \cref{lem:operators_funcCalcUnitaryConj} (together with the observation  $J (\lambda \Delta) J = \ol{\lambda} J \Delta J = \ol{\lambda} \inv{\Delta}$, $\lambda \in \C$), it follows that
    \begin{equation}\label{eq:vonNeumann_funcCalcModOp}
        \ol{f}(\inv{\Delta}) = J f(\Delta) J \ .
    \end{equation}
    In particular, using this relation with $f(t) \ce \ee^{\ii t}$, one obtains $J \Delta^{\ii t} = \Delta^{\ii t} J$ for all $t \in \R$ which completes the argument to establish \ref{enu:vonNeumann_inverseModOp}.

    Next, assertion \ref{enu:vonNeumann_deltaFS} was already established in \cref{para:vonNeumann_polarDecompositionTomita} by noting that $F = S^\ast$, \cf{} \cref{eq:vonNeumann_closedTomitaOperator}. Similarly, \ref{enu:vonNeumann_modDataTomitaOp} is an easy consequence of \cref{eq:vonNeumann_modularOperator,eq:vonNeumann_funcCalcModOp}: $S = J \Delta^{1/2} = J \Delta^{1/2} J^2 = \Delta^{-1/2} J$ and $F = S^\ast = \Delta^{1/2} J = J^2 \Delta^{1/2} J = J \Delta^{-1/2}$. For assertion \ref{enu:vonNeumann_actionModDatOmega}, first note that $S \Omega = F \Omega = \Omega$ by \eqref{eq:vonNeumann_preclosedTomitaOperator}, hence $\Delta \Omega = FS \Omega = \Omega$ by \ref{enu:vonNeumann_deltaFS}. With this observation, \ref{enu:vonNeumann_funcCalcModOp} follows from \cref{lem:operators_funcCalcEigenvalue}. Finally, using this and \ref{enu:vonNeumann_modDataTomitaOp}, one obtains also the first relation in \ref{enu:vonNeumann_actionModDatOmega}: $J \Omega = \Delta^{1/2} S \Omega = \Delta^{1/2} \Omega = \Omega$.
\end{Proof}

The following theorem is the fundamental result of the modular theory in von Neumann algebras. The proof is very technical and involved, and thus cannot be reproduced here. For detailed arguments, see \cite[91 -- 95]{BR1}, \cite[15 -- 20]{Hiai21}, \cite[278 -- 283]{StratilaZsido19}, \cite[54 -- 56]{Takesaki70}, or \cite[13 -- 17]{Takesaki03}.

\begin{ntheorem}[Tomita]\label{thm:vonNeumann_Tomita}
    Let $J$ and $\Delta$ be the modular data of a $\sigma$-finite von Neumann algebra $\MFM \subset \BO(\MH)$ with cyclic and separating vector $\Omega \in \MH$. Then the following identities hold true.
    \begin{enumerate}
        \item \label{enu:vonNeumann_Tomita1} $J \MFM J = \comm{\MFM}$.
        
        \item \label{enu:vonNeumann_Tomita2} $\Delta^{\ii t} \MFM \Delta^{- \ii t} = \MFM$ for all $t \in \R$.
    \end{enumerate}
\end{ntheorem}

\begin{definition}[Modular group]\label{def:vonNeumann_modularGroup}
    Let $\MFM \subset \BO(\MH)$ be a von Neumann algebra with a faithful normal state $\omega \in \NS(\MFM)$, and let $\Delta$ be the modular operator associated with $(\MFM, \omega)$. The one-parameter group $\R \owns t \longmto \sigma_t^\omega \in \Aut(\MFM)$ of $\ast$-automorphisms of $\MFM$, defined by
    \begin{equation*}
        \sigma_t^\omega(A) \ce \Delta^{\ii t} A \Delta^{- \ii t}
    \end{equation*}
    for all $A \in \MFM$ and $t \in \R$, is called the \bemph{modular (automorphism) group} of $(\MFM, \omega)$.
\end{definition}

\begin{examples}\label{exa:vonNeumann_modularData}
    \leavevmode
    \begin{enumerate}[env]
        \item \label{enu:vonNeumann_exaModOpLInf} (\cite[Exa. 2.4]{Hiai21}) Let $\MFM \subset \BO(\MH)$ be a von Neumann algebra possessing a faithful normal finite \bemph{trace} $\omega \in \NF{\MFM}$, that is, a functional satisfying $\omega(AB) = \omega(BA)$ for all $A, B \in \MFM$. Assume that $\omega = \omega_\Omega$ for a cyclic and separating vector $\Omega \in \MH$. For all $A \in \MFM$, it follows that
        \begin{equation*}
            \norm{A \Omega}^2 = \braket{\Omega, A^\ast A \Omega} = \omega(A^\ast A) = \omega(A A^\ast) = \braket{\Omega, A A^\ast \Omega} = \norm{A^\ast \Omega}^2 = \norm{S A \Omega}^2 \ .
        \end{equation*}
        This shows that the Tomita operator $S$ with respect to $(\MFM, \omega)$ is an isometry, hence also anti-unitary by \cref{lem:vonNeumann_propertiesTomitaOperator} \ref{enu:vonNeumann_Sinvertible}. Therefore, the modular data is given by
        \begin{equation*}
            J = S \tand \Delta = \id_\MH \ .
        \end{equation*}
        Note that if $\MFM$ were a commutative von Neumann algebra and $\omega$ an arbitrary faithful normal state, the same result would follow because $\omega$ is automatically tracial. This discussion suggests that in a sense, the modular operator $\Delta$ measures the \emph{non-tracial character} of the state $\omega$ on the algebra $\MFM$ \cite[p. 90]{BR1}. In particular, Tomita-Takesaki modular theory is trivial in the commutative case. (There has, however, been a suggestion of a \enquote{commutative version} of Tomita's theorem in the context of classical statistical mechanics \cite{Gallavotti76}.)

        \item \label{enu:vonNeumann_exaModOpBO} Let $\MH$ be a separable Hilbert space and $\MFM = \BO(\MH)$. By \cref{exa:vonNeumann_BOnormalState}, there exists a faithful normal state $\omega$ on $\MFA$ which is represented by an invertible density matrix $\rho \in \DM(\MH)$ in the form $\omega(A) = \tr(\rho A)$ for all $A \in \MFM$. Furthermore, recall from \cref{exa:operatorAlgebras_GNS} \ref{enu:operatorAlgebras_exaGNSBO} that the GNS-representation space of $\MFM$ is given by $\MH_\omega = \HS(\MH)$, with Hilbert-Schmidt inner product $\bdot_\mathrm{HS}$, on which $\MFM$ acts by multiplication, the GNS-vector representative of $\omega$ is $\Omega = \rho^{1/2} \in \MH_\omega$, and $\omega = \omega_\Omega = \braket{\Omega, \,\sbullet\ \Omega}_\mathrm{HS}$. In this case, the modular operator is given by \cite[Eq. (2.6)]{HollandsSanders18}:
        \begin{equation}\label{eq:vonNeumann_exaModOpBO}
            \Delta X = \rho X \rho^{-1} \quad (X \in \MH_\omega) \ .
        \end{equation}
        
        To see this, note that $\MFM \Omega$ lies dense in $\MH_\omega$ by the GNS-construction, hence it suffices to check the identity on this subspace. Let $A, B \in \MFM$ be arbitrary. Then
        \begin{align*}
            \braket{B \Omega, \Delta A \Omega}_\mathrm{HS} &= \braket{B \rho^{1/2}, \rho A \rho^{-1/2}}_\mathrm{HS} = \tr\bigl((B \rho^{1/2})^\ast \rho A \rho^{-1/2}\bigr) = \tr(\rho A B^\ast) \\
            &= \omega(A B^\ast) = \braket{A^\ast \Omega, B^\ast \Omega}_\mathrm{HS} = \braket{S A \Omega, S B \Omega}_\mathrm{HS} = \braket{B \Omega, S^\ast S A \Omega}_\mathrm{HS} \ .
        \end{align*}
        In the second to last step, the Tomita operator \eqref{eq:vonNeumann_preclosedTomitaOperator} was introduced, and in the last step, the definition of the adjoint of an anti-linear operator was used (\cf{} \cref{para:operators_antiLinear}). It follows that $\Delta$, as given in \eqref{eq:vonNeumann_exaModOpBO}, agrees with the operator $S^\ast S$, hence, by uniqueness of the polar decomposition of the Tomita operator, $\Delta$ must be the modular operator.\qedhere
    \end{enumerate}
\end{examples}

\section{Standard Form Representation}\label{sec:vonNeumann_standardForm}

For the discussion of the relative modular operator in the next \cref{sec:vonNeumann_relativeModularOperator} as well as the development of perturbation theory in \cref{ch:perturbationTheory}, the theory of the standard form representation of von Neumann algebras is an indispensible tool. It was developed, independently, by \textsc{H. Araki} \cite{Araki74a} and \textsc{A. Connes} \cite{Connes74} for $\sigma$-finite von Neumann algebras in 1974, and generalized by \textsc{U. Haagerup} \cite{Haagerup75} in 1975.

\subsection{Definition and Examples}\label{subsec:vonNeumann_standardFormDef}

\begin{para}[Natural positive cone]\label{para:vonNeumann_NPC}
    (\cite[p. 25]{Hiai21})
    Let $\MFM \subset \BO(\MH)$ be a $\sigma$-finite von Neumann algebra with cyclic and separating vector $\Omega \in \MH$. Denote by $\Delta$ and $J$ the modular data with respect to $(\MFM, \Omega)$, and let $j : \MFM \longto \comm{\MFM}$ be the anti-linear $\ast$-isomorphism defined by
    \begin{equation*}
        j(A) \ce J A J \quad (A \in \MFM) \ .
    \end{equation*}
    (\cref{thm:vonNeumann_Tomita} \ref{enu:vonNeumann_Tomita1} shows that the range of $j$ is indeed the commutant $\comm{\MFM}$.) The \bemph{natural positive cone} $\NPC$ in $\MH$ associated with $(\MFM, \Omega)$ is defined to be
    \begin{equation}\label{eq:vonNeumann_NPC}
        \NPC \ce \clos_{\ndot} \! \set[\big]{A j(A) \Omega \ : \ A \in \MFM} = \clos_{\ndot} \! \set[\big]{A J A \Omega \ : \ A \in \MFM} \ .
    \end{equation}
\end{para}

The set $\NPC$ displays a number of very useful intrinsic properties as well as relations to the modular data. The proof of the next proposition can be found, \eg{}, in \cite[Prop. 2.5.26]{BR1} or \cite[Thm. 3.2]{Hiai21}.

\begin{proposition}\label{pro:vonNeumann_propertiesNPC}
    The natural positive cone $\NPC \subset \MH$ enjoys the following properties:
    \begin{enumerate}
        \item \label{enu:vonNeumann_NPCclosedConvex} $\NPC = \clos_{\ndot}(\Delta^{1/4} \PE{\MFM} \Omega)$. In particular, $\NPC$ is a closed convex cone inside $\MH$.
        
        \item \label{enu:vonNeumann_JonNPC} $J \xi = \xi$ for all $\xi \in \NPC$.
        
        \item \label{enu:vonNeumann_DeltaNPC} $\Delta^{\ii t} \NPC = \NPC$ for all $t \in \R$.
        
        \item \label{enu:vonNeumann_jNPC} $A j(A) \NPC \subset \NPC$ for all $A \in \MFM$.
        
        \item If $f : \R \longto (0, + \infty)$ is a positive-definite function on $\R$, then $f(\log \Delta) \NPC \subset \NPC$.
        
        \item $\NPC$ is self-dual in the sense that $\NPC = \set{\eta \in \MH \, : \, \text{$\braket{\xi, \eta} \ge 0$ for all $\xi \in \NPC$}}$.
    \end{enumerate}
\end{proposition}

\begin{definition}[Standard form]\label{def:vonNeumann_standardForm}
    A quadruple $(\MFM, \MH, J, \NPC)$ consisting of a von Neumann algebra $\MFM$ represented faithfully on a Hilbert space $\MH$, an anti-unitary involution $J : \MH \longto \MH$, and a self-dual closed convex cone $\NPC \subset \MH$ is called a \bemph{standard form} of the von Neumann algebra $\MFM$ iff the following properties are satisfied:
    \begin{enumerate}[num]
        \item \label{enu:vonNeumann_standardFormTomita} $J \MFM J = \comm{\MFM}$;
        
        \item \label{enu:vonNeumann_standardFormJCenter} $J A J = A^\ast$ for all $A \in \MFM \cap \comm{\MFM}$;
        
        \item \label{enu:vonNeumann_standardFormJIdentityNPC} $J \xi = \xi$ for all $\xi \in \NPC$;
        
        \item \label{enu:vonNeumann_standardFormJConjugationNPC} $A j(A) \NPC \subset \NPC$ for all $A \in \MFM$.
    \end{enumerate}
\end{definition}

\cref{para:vonNeumann_NPC} and \cref{pro:vonNeumann_propertiesNPC} show that every $\sigma$-finite von Neumann algebra can be represented in standard form. Regarding non-$\sigma$-finite von Neumann algebras, \textsc{U. Haagerup} proved the following theorem in 1975 \cite[Thm. 1.6]{Haagerup75}. (See also \cite[Sect. 10.14 \& 10.23]{StratilaZsido19}.)

\begin{ntheorem}[Haagerup]\label{thm:vonNeumann_existenceStandardForm}
    Every von Neumann algebra possesses a standard form representation.
\end{ntheorem}

\begin{examples}\label{exa:vonNeumann_standardForm}
    \leavevmode
    \begin{enumerate}[env]
        \item \label{enu:vonNeumann_exaStandardFormAbelian} Let $(X, \Sigma, \mu)$ be a $\sigma$-finite measure space. The von Neumann algebra $\MFM = L^\infty(X, \mu)$ acts on $\MH = L^2(X, \mu)$ by multiplication (\cref{exa:operatorAlgebras_vonNeumann} \ref{enu:operatorAlgebras_exaLInf}), and this action is faithful. Moreover, by assumption there exists a measurable function $g : X \longto \C$ such that $g > 0$ and $\int_X g \diff \mu < + \infty$ \cite[p. 21]{Kallenberg17}. This gives rise to the faithful normal functional $\omega(f) \ce \int_X f g \diff \mu$ on $\MFM$ (\cf{} \cref{exa:operatorAlgebras_normalStates} \ref{enu:operatorAlgebras_exaNormalStateAbelian}) with cyclic and separating vector representative $h \ce \sqrt{g} \in \MH$.
        
        In \cref{exa:vonNeumann_modularData} \ref{enu:vonNeumann_exaModOpLInf}, it was shown that $J = S$ is given by the Tomita operator. Furthermore, it holds that $j(f) = J f J$ implements the $\ast$-operation on $\MFM$: for all $f, g \in \MFM$, one has
        \begin{equation*}
            J f J (g h) = J (f g^\ast h) = g f^\ast h = \ol{f} (g h) \ .
        \end{equation*}
        Therefore, the natural positive cone \eqref{eq:vonNeumann_NPC} is given by $\NPC = \clos_{\ndot} \set{\abs{f}^2 h \, : \, f \in \MFM} = L^2(X, \mu)_+$, and hence the standard form representation of $\MFM$ is \cite[Exa. 3.6 (1)]{Hiai21}
        \begin{equation*}
            \bigl(L^\infty(X, \mu),\, L^2(X, \mu),\, J \xi = \ol{\xi},\, L^2(X, \mu)_+\bigr) \ .
        \end{equation*}

        \item \label{enu:vonNeumann_exaStandardFormTypeI} Let $\MH$ be a separable Hilbert space and $\MFM = \BO(\MH)$. This algebra acts faithfully on $\HS(\MH)$ by multiplication, and one can show that the standard form representation of $\MFM$ is given by \cite[Exa. 3.6 (2)]{Hiai21}
        \begin{equation*}
            \bigl(\BO(\MH), \HS(\MH), J A = A^\ast, \HS(\MH)_+\bigr) \ ,
        \end{equation*}
        where $J : \HS(\MH) \longto \HS(\MH)$, $A \longmto A^\ast$, is given by the operation of taking adjoints, and $\HS(\MH)_+$ is the cone of positive Hilbert-Schmidt operators.\qedhere
    \end{enumerate}
\end{examples}

\subsection{Properties of the Standard Form Representation}\label{subsec:vonNeumann_standardFormProperties}

The following proposition, proved in \cite[Prop. 2.5.28]{BR1} or \cite[Prop. 3.7]{Hiai21}, concerns further geometric properties of the cone $\NPC$. Note the close analogy between the following result and \cref{pro:operatorAlgebras_propertiesPE} which characterized the cone $\PE{\MFA}$ of positive elements in a $C^\ast$-algebra $\MFA$.

\begin{proposition}\label{pro:vonNeumann_geometricPropertiesNPC}
    Let $(\MFM, \MH, J, \NPC)$ be a standard form. The closed convex cone $\NPC$ satisfies $\NPC \cap (- \NPC) = \{0\}$. Furthermore, if $\xi \in \MH$ such that $J \xi = \xi$, then $\xi$ can be uniquely decomposed as $\xi = \xi_1 - \xi_2$, where $\xi_1, \xi_2 \in \NPC$ with $\xi_1 \perp \xi_2$. Finally, $\MH$ is linearly spanned by $\NPC$.
\end{proposition}

One can show that the data $(J, \NPC)$ does not depend on the choice of a cyclic and separating vector $\Omega \in \MH$, hence the objects $J$ and $\NPC$ are \emph{universal} \cite[Prop. 2.5.30]{BR1}, \cite[Prop. 3.10]{Hiai21}. Furthermore, the elements of $\NPC$ enjoy some special properties.

\begin{lemma}\label{lem:vonNeumann_vectorsNPC}
    Let $(\MFM, \MH, J, \NPC)$ be a standard form. For any $\xi \in \NPC$, it holds that $\xi$ is cyclic for $\MFM$ if and only if $\xi$ is separating for $\MFM$.
\end{lemma}

\begin{Proof}
    (\cite[p. 30]{Hiai21})
    Let $\xi \in \NPC$ be arbitrary. If $\xi$ is assumed to be cyclic for the algebra $\MFM$, then $\xi = J \xi$ is cyclic for $J \MFM J = \comm{\MFM}$, hence it is also separating for $\MFM$ by \cref{pro:operatorAlgebras_cyclicSeparating} \ref{enu:operatorAlgebras_separatingMcyclicCommM}. The other direction is proved analogously.
\end{Proof}

The following theorem, which is the main result of this section, establishes a relationship between the natural positive cone of a von Neumann algebra and the space of normal functionals. The proof can be found in \cite[Thm. 2.5.31]{BR1} or \cite[Thm. 3.12]{Hiai21}.

\begin{ntheorem}[Correspondence between $\NF{\MFM}$ and $\NPC$]\label{thm:vonNeumann_vectorRepresentativesNormalStates}
    Let $(\MFM, \MH, J, \NPC)$ be a standard form. For every $\varphi \in \NF{\MFM}$, there exists a unique vector $\xi \in \NPC$ such that $\varphi = \omega_\xi$. In fact, the mapping $\NPC \longto \NF{\MFM}$, $\xi \longmto \omega_\xi$, is a homeomorphism with respect to the norm topology on both spaces.
\end{ntheorem}

\begin{remark}\label{rem:vonNeumann_cyclicSeparatingFaithful}
    The inverse of the mapping $\NPC \longto \NF{\MFM}$ of the previous theorem will be denoted by $\NF{\MFM} \longto \NPC$, $\varphi \longmto \xi_\varphi$, and $\xi_\varphi$ will be called the \bemph{standard vector representative} of $\varphi$. Combining \cref{cor:operatorAlgebras_faithfulSeparating} and \cref{lem:vonNeumann_vectorsNPC}, it follows that
    \begin{equation*}
        \forall \xi \in \NPC \quad : \quad \xi \ \text{cyclic} \ \iff \ \xi \ \text{separating} \ \iff \ \omega_\xi \ \text{faithful} \ .
    \end{equation*}
\end{remark}

An important consequence of \cref{thm:vonNeumann_vectorRepresentativesNormalStates} is that all $\ast$-automorphisms of the algebra $\MFM$ are implemented by unitary elements from $\MH$. This will be a key technical tool for the investigations in \cref{ch:perturbationTheory}. The proof can be found in \cite[Cor. 2.5.32]{BR1} and \cite[Thm. 4.43]{AJP06}.

\begin{corollary}\label{cor:vonNeumann_unitaryReprAut}
    Let $(\MFM, \MH, J, \NPC)$ be a von Neumann algebra in standard form. Then there exists a unique unitary representation $\Aut(\MFM) \owns \alpha \longmto U(\alpha) \in \UO(\MH)$ of the group $\Aut(\MFM)$ of $\ast$-automorphisms of the von Neumann algebra $\MFM$ on the Hilbert space $\MH$ such that
    \begin{enumerate}[num]
        \item $U(\alpha) \MFM U(\alpha)^\ast = \MFM$ and $U(\alpha) \comm{\MFM} U(\alpha)^\ast = \comm{\MFM}$;
        \item $U(\alpha) \NPC \subset \NPC$;
        \item $\forall A \in \MFM \, : \, \alpha(A) = U(\alpha) A U(\alpha)^\ast$;
        \item $\forall \varphi \in \NF{\MFM} \, : \, U(\alpha)^\ast \xi_\varphi = \xi_{\varphi \circ \alpha}$;
        \item $[U(\alpha), J] = 0$.
    \end{enumerate}
    Furthermore, the mapping $\Aut(\MFM) \owns \alpha \longmto U(\alpha) \in \UO(\MH)$ is a homeomorphism when both spaces are equipped with the norm topology.
\end{corollary}

The final result of this section establishes uniqueness of the standard form representation up to unitary equivalence \cite[Thm. 3.13]{Hiai21}.

\begin{corollary}\label{cor:vonNeumann_uniquenessStandardForm}
    Let $(\MFM, \MH, J, \NPC)$ and $(\wt{\MFM}, \wt{\MH}, \wt{J}, \wt{\NPC})$ be two standard forms of von Neumann algebras $\MFM$ and $\wt{\MFM}$, respectively. If $\Phi : \MFM \longto \wt{\MFM}$ is a $\ast$-isomorphism, then there exists a uniquely defined unitary operator $U : \MH \longto \wt{\MH}$ such that
    \begin{enumerate}[label=\normalfont(\arabic*)]
        \item $\forall A \in \MFM \, : \, \Phi(A) = U A U^\ast$;
        \item $\wt{J} = U J U^\ast$;
        \item $\wt{\NPC} = U \NPC$.
    \end{enumerate}
\end{corollary}

\section{Relative Modular Operators}\label{sec:vonNeumann_relativeModularOperator}

The construction of the modular operator with respect to a von Neumann algebra $\MFM \subset \BO(\MH)$ and a cyclic and separating vector $\Omega \in \MH$ can be generalized to the situation in which two (not necessarily faithful) positive normal functionals $\varphi, \psi \in \NF{\MFM}$ are given. The resulting object can be seen as a \emph{non-commutative Radon-Nikodým derivative} of these functionals \cite[p. 163]{Hiai21}; it was invented by \textsc{H. Araki} \cite{Araki76, Araki77} who used it to extend the notion of relative entropy to general von Neumann algebras.

\subsection{Construction}

Let $(\MFM, \MH, J, \NPC)$ be a standard form and $\varphi, \psi \in \NF{\MFM}$. According to \cref{thm:vonNeumann_vectorRepresentativesNormalStates}, there exist vectors $\Phi, \Psi \in \NPC$ such that $\varphi = \omega_\Phi$ and $\psi = \omega_\Psi$. Recall from \cref{pro:operatorAlgebras_supportVectorFunctional} that the support projections of the functionals $\omega_\Phi$ on $\MFM$ and $\comm{\omega_\Phi}$ on $\comm{\MFM}$ are given by $[\comm{\MFM} \Phi]$ and $[\MFM \Phi]$, respectively. In this and the subsequent chapters, the following notations shall be used:
\begin{equation}\label{eq:vonNeumann_supportProjectionsMT}
    \ssupp_\MFM(\varphi) \equiv \ssupp_\varphi \equiv \ssupp_\Phi \ce [\comm{\MFM} \Phi] \tand \ssupp_{\comm{\MFM}}(\varphi) \equiv \comm{\ssupp_\varphi} \equiv \comm{\ssupp_\Phi} \ce [\MFM \Phi] \ .
\end{equation}

\begin{lemma}\label{lem:vonNeumann_commMsupport}
    $J \, \ssupp_\MFM(\varphi) J = \ssupp_{\comm{\MFM}}(\varphi)$.
\end{lemma}

\begin{Proof}
    (\cite[p. 163]{Hiai21})
    Since the Hilbert space $\MH$ is linearly spanned by the natural positive cone $\NPC$ (\cref{pro:vonNeumann_geometricPropertiesNPC}) and since $J \xi = \xi$ for all $\xi \in \NPC$ (\cref{def:vonNeumann_standardForm}), it follows that $J \ssupp_\MFM(\varphi) J \MH = J \ssupp_\MFM(\varphi) \MH$. Using that $J \comm{\MFM} J = \MFM$ one obtains
    \begin{align*}
        J \ssupp_\MFM(\varphi) J \MH &= J \ssupp_\MFM(\varphi) \MH = \clos_{\ndot}(J \comm{\MFM} \Phi) \\
        &= \clos_{\ndot}(J \comm{\MFM} J \Phi) = \clos_{\ndot}(\MFM \Phi) = \ssupp_{\comm{\MFM}}(\varphi) \MH \ . \tag*{\qedhere}
    \end{align*}
\end{Proof}

\begin{para}[The relative Tomita operator]\label{para:vonNeumann_relativeTomita}
    (\cite[p. 163]{Hiai21})
    For every pair of positive normal functionals $\varphi, \psi \in \NF{\MFM}$, define two operators $S_{\psi, \varphi}^0 : \MH \supset \dom(S_{\psi, \varphi}^0) \longto \MH$ and $F_{\psi, \varphi}^0 : \MH \supset \dom(F_{\psi, \varphi}^0) \longto \MH$ by
    \begin{equation}\label{eq:vonNeumann_relativeTomitaOperator}
        S_{\psi, \varphi}^0(A \Phi + \eta) \ce \ssupp_\MFM(\varphi) A^\ast \Psi \tand F_{\psi, \varphi}^0(\comm{A} \Phi + \zeta) \ce \ssupp_{\comm{\MFM}}(\varphi) (\comm{A})^\ast \Psi
    \end{equation}
    for all $A \in \MFM$, $\comm{A} \in \comm{\MFM}$, $\eta \in (\id_\MH - \comm{\ssupp_\varphi}) \MH$ and $\zeta \in (\id_\MH - \ssupp_\varphi) \MH$. Since $\id_\MH - \comm{\ssupp_\varphi}$ is the orthogonal projection onto $[\MFM \Phi]^\perp$ and $\id_\MH - \ssupp_\varphi$ onto $[\comm{\MFM} \Phi]^\perp$ (\cref{thm:operators_orthogonalProjection}), it follows that
    \begin{equation*}
        \dom(S_{\psi, \varphi}^0) = \MFM \Phi + [\MFM \Phi]^\perp \tand \dom(F_{\psi, \varphi}^0) = \comm{\MFM} \Phi + [\comm{\MFM} \Phi]^\perp \ .
    \end{equation*}

    \begin{nstatement}
        The operators $S_{\psi, \varphi}^0$ and $F_{\psi, \varphi}^0$ are well-defined anti-linear operators. Moreover, their domains lie dense in $\MH$, and they are closable.
    \end{nstatement}

    \begin{Proof}
        (\cite[Lem. 10.2]{Hiai21})
        \pp{1. The operators are well-defined}. Let $A_i \in \MFM$ and $\eta_i \in [\MFM \Phi]^\perp$, $i \in \set{1, 2}$, be arbitrary such that $A_1 \Phi + \eta_1 = A_2 \Phi + \eta_2$. Observe that $(A_1 - A_2) \Phi = \eta_2 - \eta_1 \in [\MFM \Phi] \cap [\MFM \Phi]^\perp = \set{0}$, hence $A_1 \Phi = A_2 \Phi$. By definition of $\ssupp_\MFM(\varphi)$, this implies $A_1 \, \ssupp_\MFM(\varphi) = A_2 \, \ssupp_\MFM(\varphi)$ on $\MH$. In particular, one may conclude that $\ssupp_\MFM(\varphi) A_1^\ast \Psi = \ssupp_\MFM(\varphi) A_2^\ast \Psi$. Thus,
        \begin{equation*}
            S_{\psi, \varphi}^0 (A_1 \Phi + \eta_1) = \ssupp_\MFM(\varphi) A_1^\ast \Psi = \ssupp_\MFM(\varphi) A_2^\ast \Psi = S_{\psi, \varphi}^0 (A_2 \Phi + \eta_2)
        \end{equation*}
        which shows that the operator $S_{\psi, \varphi}^0$ is well-defined; an analogous argument establishes the same for the operator $F_{\psi, \varphi}^0$.

        \pp{2. The domains lie dense}. It holds that $\MH \supset \clos_{\ndot} \dom(S_{\psi, \varphi}^0) = \clos_{\ndot}(\MFM \Phi + [\MFM \Phi]^\perp) \supset \clos_{\ndot}(\MFM \Phi) + \clos_{\ndot}([\MFM \Phi]^\perp) = [\MFM \Phi] + [\MFM \Phi]^\perp = \MH$ (\cf{} \cref{thm:operators_orthogonalProjection} and note that the orthogonal complement is always closed). This shows that $S_{\psi, \varphi}^0$ is densely defined; for the operator $F_{\psi, \varphi}^0$ defined on $\comm{\MFM} \Phi + [\comm{\MFM} \Phi]^\perp$, the argument is analogous. The anti-linearity of these operators is clear from their definition.

        \pp{3. Closability}. Let $A \in \MFM$, $\eta \in [\MFM \Phi]^\perp$ and $\comm{A} \in \comm{\MFM}$, $\zeta \in [\comm{\MFM} \Phi]^\perp$ be arbitrary. One computes
        \begin{align*}
            \vbraket[\big]{A \Phi + \eta, F_{\psi, \varphi}^0 (\comm{A} \Phi + \zeta)} &= \vbraket[\big]{A \Phi + \eta, \ssupp_{\comm{\MFM}}(\varphi) (\comm{A})^\ast \Psi} = \vbraket[\big]{A \Phi, (\comm{A})^\ast \Psi} = \vbraket[\big]{\comm{A} \Phi, A^\ast \Psi} \\
            &= \vbraket[\big]{\ssupp_\MFM(\varphi) (\comm{A} \Phi + \zeta), A^\ast \Psi} = \vbraket[\big]{\comm{A} \Phi + \zeta, \ssupp_\MFM(\varphi) A^\ast \Psi} \\
            &= \vbraket[\big]{\comm{A} \Phi + \zeta, S_{\psi, \varphi}^0(A \Phi + \eta)} \ .
        \end{align*}
        From the definition of the adjoint of an anti-linear operator (\cref{para:operators_antiLinear}), it follows that $S_{\psi, \varphi}^0 \subset (F_{\psi, \varphi}^0)^\ast$ and $F_{\psi, \varphi}^0 \subset (S_{\psi, \varphi}^0)^\ast$, hence \cref{pro:operators_propertiesAdjoint} shows that both $S_{\psi, \varphi}^0$ and $F_{\psi, \varphi}^0$ are closable.
    \end{Proof}

    Let $S_{\psi, \varphi}$ and $F_{\psi, \varphi}$ denote the closures of the operators $S_{\psi, \varphi}^0$ and $F_{\psi, \varphi}^0$, respectively. One refers to $S_{\psi, \varphi}$ as the \bemph{relative Tomita operator} associated with the functionals $\varphi, \psi \in \NF{\MFM}$. If $\psi = \varphi$, one simply writes $S_\varphi$ and $F_\varphi$ for these operators. Note that if $\psi = \varphi$ is faithful, then $S_\varphi$ agrees with the Tomita operator defined in \cref{para:vonNeumann_preclosedTomitaOperator}: indeed, the vector representative $\Phi$ of $\varphi$ is cyclic and separating (\cref{rem:vonNeumann_cyclicSeparatingFaithful}), hence $[\MFM \Phi] = \MH$ and $[\MFM \Phi]^\perp = \{0\}$ so that $S_\varphi^0 A \Phi = A^\ast \Phi = S A \Phi$ for all $A \in \MFM$.
\end{para}

\begin{definition}[Relative modular operator]\label{def:vonNeumann_relativeModularOperator}
    The \bemph{relative modular operator} $\Delta_{\psi, \varphi}$ with respect to the positive normal functionals $\psi$ and $\varphi$ on $\MFM$ is defined to be
    \begin{equation}\label{eq:vonNeumann_relativeModularOperator}
        \Delta_{\psi, \varphi} \ce S_{\psi, \varphi}^\ast S_{\psi, \varphi} \ .
    \end{equation}
    In the case that $\psi = \varphi$, one simply writes $\Delta_\varphi$ for $\Delta_{\varphi, \varphi}$ and calls this the \bemph{modular operator} of $\varphi$. (Note that by the above remark, it holds that if $\varphi$ is faithful, then $\Delta_\varphi$ agrees with the operator $\Delta$ from \cref{def:vonNeumann_modularData}.) One also writes
    \begin{equation*}
        \Delta_{\Psi, \Phi} \equiv \Delta_{\psi, \varphi} \tand S_{\Psi, \Phi} \equiv S_{\psi, \phi} \ .
    \end{equation*}
\end{definition}

The following examples \cite[Exa. 10.4]{Hiai21} shed some light on the abstract concept of the relative modular operator. Note, in particular, that \ref{enu:vonNeumann_relativeModularOperatorAbelian} shows that this concept is non-trivial for commutative von Neumann algebras which is in contrast to the situation for the modular operator found in \cref{exa:vonNeumann_modularData} \ref{enu:vonNeumann_exaModOpLInf}.

\begin{examples}\label{exa:vonNeumann_relativeModularOperator}
    \leavevmode
    \begin{enumerate}[env]
        \item \label{enu:vonNeumann_relativeModularOperatorAbelian} Let $(X, \Sigma, \mu)$ be a $\sigma$-finite measure space and $\MFM = L^\infty(X, \mu)$, acting on $\MH = L^2(X, \mu)$. Recall from \cref{exa:operatorAlgebras_normalStates} \ref{enu:operatorAlgebras_exaNormalStateAbelian} that there is a one-to-one correspondence $\NF{\MFM} \cong L^1(X, \mu)_+$; in particular, every $\psi \in \NF{\MFM}$ gives rise to a measure $\diff \nu_\psi = h_\psi \diff \mu$. For every pair of normal functionals $\varphi, \psi \in \NF{\MFM}$, it holds that the relative modular operator $\Delta_{\psi, \varphi}$ is given by the multiplication operator with the function
        \begin{equation*}
            \1_{\set{h_\varphi > 0}} \cdot \frac{h_\psi}{h_\varphi} \ .
        \end{equation*}
        
        To see this, let $f, g \in \MFM$ be arbitrary and note that the vector representative of $\varphi$ in the natural positive cone is $\Phi \ce \sqrt{h_\varphi} \in L^2(X, \mu)_+$. The support projection of $\varphi$ is given by $\ssupp_\MFM(\varphi) = \1_{\set{\varphi > 0}} = \1_{\set{h_\varphi > 0}}$. Indeed, if there were $P \in \PO(\MFM)$ such that $\varphi(P) = \varphi(\1_X)$ and $P < \ssupp_\MFM(\varphi)$, then
        \begin{equation*}
            \int_X h_\varphi \diff \mu = \varphi(\1_X) = \varphi(P) = \int_X P h_\varphi \diff \mu < \int_X \1_{\set{h_\varphi > 0}} h_\varphi \diff \mu = \int_X h_\varphi \diff \mu \ ,
        \end{equation*}
        a contradiction. Since $\MFM$ is commutative, it holds that $\ssupp_{\comm{\MFM}}(\varphi) = \1_{\set{h_\varphi > 0}}$ as well. Let $\eta, \zeta \in (\1_X - \comm{\ssupp_\varphi}) \MH$ be arbitrary, and observe that $\1_{\set{h_\varphi > 0}} \eta = \1_{\set{h_\varphi > 0}} \zeta = 0$. One then computes
        \begin{align*}
            \vbraket[\big]{f \, \Phi + \eta, \Delta_{\psi, \varphi} (g \, \Phi + \zeta)} &= \int\nolimits_X \ol{(f \, \Phi + \eta)} \, \1_{\set{h_\varphi > 0}} \, \frac{h_\psi}{h_\varphi} \, (g \, \Phi + \zeta) \diff \mu = \int\nolimits_X \1_{\set{h_\varphi > 0}} \ol{f} g h_\psi \diff \mu \\
            &= \vbraket[\big]{\1_{\set{h_\varphi > 0}} \ol{g} \, \Psi, \1_{\set{h_\varphi > 0}} \ol{f} \, \Psi} = \vbraket[\big]{S_{\psi, \varphi} (g \, \Phi + \zeta), S_{\psi, \varphi} (f \, \Phi + \eta)} \\[5pt]
            &= \vbraket[\big]{f \, \Phi + \eta, S_{\psi, \varphi}^\ast S_{\psi, \varphi} (g \, \Phi + \zeta)} \ ,
        \end{align*}
        where $\Psi \ce \sqrt{h_\psi} \in L^2(X, \mu)_+$ is the vector representative of $\psi$. From the uniqueness of the polar decomposition of the relative Tomita operator, the claim now follows.
        
        It holds that $\Delta_{\psi, \varphi} = \1_{\set{h_\varphi > 0}} \, \frac{h_\psi}{h_\varphi}$ is equal to the Radon-Nikodým derivative $\od{\nu_\psi}{\nu_\varphi}$ of the measure $\diff \nu_\psi = h_\psi \diff \mu$ with respect to the measure $\diff \nu_\varphi = h_\varphi \diff \mu$. (This is the reason for calling the relative modular operator a non-commutative Radon-Nikodým derivative.) Indeed, for all $A \in \Sigma$, one easily computes that
        \begin{equation*}
            \nu_\psi(A) = \int_A h_\psi \diff \mu = \int_A \1_{\set{h_\varphi > 0}} \, \frac{h_\psi}{h_\varphi} \, h_\varphi \diff \mu = \int_A \1_{\set{h_\varphi > 0}} \, \frac{h_\psi}{h_\varphi} \diff \nu_\varphi \ .
        \end{equation*}       

        \item \label{enu:vonNeumann_relativeModularOperatorTypeI} Let $\MH$ be a separable Hilbert space. Consider the von Neumann algebra $\MFM = \BO(\MH)$ which is represented in standard form on the Hilbert space $\HS(\MH)$ (\cref{exa:vonNeumann_standardForm} \ref{enu:vonNeumann_exaStandardFormTypeI}). For every $\varphi, \psi \in \NF{\MFM}$, there exist positive $\rho_\varphi, \rho_\psi \in \NO(\MH)$ such that $\varphi = \tr(\rho_\varphi \ \sbullet\,)$ and $\psi = \tr(\rho_\psi \ \sbullet\,)$. Let $\rho_\psi = \sum_{i \in \N} \lambda_i P_i$ and $\rho_\varphi = \sum_{j \in \N} \mu_j Q_j$ be the spectral decompositions of $\rho_\varphi$ and $\rho_\psi$, where $\lambda_i, \mu_j > 0$ and $P_i, Q_j \in \PO(\MH)$ are finite-rank projections for all $i, j \in \N$. Then one shows, similarly to the previous example and \cref{exa:vonNeumann_modularData} \ref{enu:vonNeumann_exaModOpBO} that the relative modular operator $\Delta_{\psi, \varphi}$ on $\HS(\MH)$ is given by\footnote{The operator $\inv{\rho_\varphi}$ is defined with restriction to the support $\supp(\rho_\varphi)$ of $\rho_\varphi$, \ie{}, it is the generalized inverse of $\rho_\varphi$; see \cite[p. 167]{Hiai21}, \cite[p. 164]{Hiai21a}, \cite[Sect. 9.3]{BenIsraelGreville03}.}
        \begin{equation*}
            \Delta_{\psi, \varphi} X = \rho_\psi X \rho_\varphi^{-1} = \sum_{i,j=1}^{\infty} \lambda_i \inv{\mu_j} P_i X Q_j \ . \tag*{\qedhere}
        \end{equation*}
    \end{enumerate}
\end{examples}

\subsection{Properties of Relative Modular Operators}\label{subsec:vonNeumann_relativeModularOperatorProperties}

The most important properties of the relative modular operator are contained in the next proposition whose proof shall be skipped as it is quite involved; see \cite[Prop. 10.3]{Hiai21}.

\begin{proposition}\label{pro:vonNeumann_propertiesRelativeModularOp}
    Let $\psi, \varphi \in \NF{\MFM}$ be two positive normal functionals on a von Neumann algebra $(\MFM, \MH, J, \NPC)$ in standard form.
    \begin{enumerate}
        \item \label{enu:vonNeumann_supportRelModOp} The support projection (\cf{} \cref{para:operators_supportPO}) of the relative modular operator $\Delta_{\psi, \varphi}$ is
        \begin{equation*}
            \ssupp(\Delta_{\psi, \varphi}) = \ssupp_\MFM(\psi) \, \ssupp_{\comm{\MFM}}(\varphi) = \ssupp_\MFM(\psi) J \ssupp_\MFM(\varphi) J \ .
        \end{equation*}
        
        \item \label{enu:vonNeumann_polarDecompRelTomita} The polar decomposition of $S_{\psi, \varphi}$ takes the form $S_{\psi, \varphi} = J \Delta_{\psi, \varphi}^{1/2}$.
        
        \item \label{enu:vonNeumann_inverseRelModOp} The inverse operator $\Delta_{\varphi, \psi}^{-1}$, defined on the support of $\Delta_{\varphi, \psi}$, is given by $\Delta_{\varphi, \psi}^{-1} = J \Delta_{\psi, \varphi} J$.
    \end{enumerate}
\end{proposition}

The next two lemmata will be useful for arguments in the following chapters. These properties are mentioned, for example, in \cite[Thm. 4.1]{DJP03}. Let again $(\MFM, \MH, J, \NPC)$ be a von Neumann algebra in standard form and $\varphi, \psi \in \NF{\MFM}$.

\begin{lemma}\label{lem:vonNeumann_relModOpCommutesCenter}
    If $B \in \MFM \cap \comm{\MFM}$ belongs to the center of $\MFM$, then $B$ commutes with $\Delta_{\varphi, \psi}$.
\end{lemma}

\begin{Proof}
    Let $A, C \in \MFM$ and $\eta, \zeta \in [\MFM \Psi]^\perp$ be arbitrary. Then
    \begin{align*}
        \vbraket[\big]{C \Psi + \zeta, \Delta_{\varphi, \psi} B (A \Psi + \eta)} &= \vbraket[\big]{S_{\varphi, \psi} B (A \Psi + \eta), S_{\varphi, \psi} (C \Psi + \zeta)} \\
        &= \vbraket[\big]{\ssupp_\MFM(\psi) (BA)^\ast \Phi, \ssupp_\MFM(\psi) C^\ast \Phi} \\
        &= \vbraket[\big]{\ssupp_\MFM(\psi) A^\ast \Phi, \ssupp_\MFM(\psi) (B^\ast C)^\ast \Phi} \\
        &= \vbraket[\big]{S_{\varphi, \psi} (A \Psi + \eta), S_{\varphi, \psi} (B^\ast C \Psi + B^\ast \zeta)} \\
        &= \vbraket[\big]{C \Psi + \zeta, B \Delta_{\varphi, \psi} (A \Psi + \eta)} \ .
    \end{align*}
    First, the definition of the adjoint of an anti-linear operator was employed (\cref{para:operators_antiLinear}). In the second line, it was used that $\eta \in [\MFM \Psi]^\perp = \ker(\comm{\ssupp_\psi})$ implies $\ssupp_{\comm{\MFM}}(\psi) B \eta = B \, \ssupp_{\comm{\MFM}}(\psi) \eta = 0$, which shows that $B \eta \in \ker(\comm{\ssupp_\psi})$. Since the space $\MFM \Psi + [\MFM \Psi]^\perp$ lies dense in $\MH$, the claim follows from the above identity.
\end{Proof}

\begin{lemma}\label{lem:vonNeumann_relativeModularOpScaling}
    Let $\lambda, \mu \in (0, + \infty)$ be positive real numbers. Then
    \begin{equation*}
        \Delta_{\mu \varphi, \lambda \psi} = \frac{\mu}{\lambda} \, \Delta_{\varphi, \psi} \ .
    \end{equation*}
\end{lemma}

\begin{Proof}
    The vector representative in $\NPC$ of the functional $\lambda \psi$, which is defined by $(\lambda \psi)(A) \ce \lambda \psi(A)$ for all $A \in \MFM$, is given by $\sqrt{\lambda} \, \Psi$, and similarly the vector representative of $\mu \varphi$ is $\sqrt{\mu} \, \Phi$. 
    
    Observe that the support projection of the functional $\lambda \psi$ satisfies the identity $\ssupp_\MFM(\lambda \psi) = \ssupp_\MFM(\psi)$. To see this, note that on the one hand,
    \begin{equation*}
        \psi\bigl(\ssupp_\MFM(\lambda \psi)\bigr) = \frac{1}{\lambda} \, (\lambda \psi)\bigl(\ssupp_\MFM(\lambda \psi)\bigr) \overset{\eqref{eq:operatorAlgebras_supportProjection}}{=} \frac{1}{\lambda} (\lambda \psi)(\id_\MH) = \psi(\id_\MH) \ ,
    \end{equation*}
    hence $\ssupp_\MFM(\lambda \psi) \ge \ssupp_\MFM(\psi)$ by definition of the support projection (\cref{para:operatorAlgebras_support}). On the other hand, $(\lambda \psi)\bigl(\ssupp_\MFM(\psi)\bigr) = (\lambda \psi)(\id_\MH)$, implying $\ssupp_\MFM(\psi) \ge \ssupp_\MFM(\lambda \psi)$.
    
    With these two observations, one can now compute the relative Tomita operator for the functionals $\mu \varphi$ and $\lambda \psi$ as follows: for $A \in \MFM$ and $\eta \in [\MFM \Psi]^\perp$, there holds
    \begin{align*}
        S_{\mu \varphi, \lambda \psi}^0(A \Psi + \eta) &= \frac{1}{\sqrt{\lambda}} \, S_{\mu \varphi, \lambda \psi} \bigl(A \sqrt{\lambda} \, \Psi + \sqrt{\lambda} \, \eta\bigr) \\
        &= \frac{1}{\sqrt{\lambda}} \, \ssupp_\MFM(\lambda \psi) A^\ast \sqrt{\mu} \, \Phi = \sqrt{\frac{\mu}{\lambda}} \, \ssupp_\MFM(\psi) A^\ast \Phi \\
        &= \sqrt{\frac{\mu}{\lambda}} \, S_{\varphi, \psi}^0(A \Psi + \eta) \ .
    \end{align*}
    Observe that $\sqrt{\lambda} \, \eta \in (\id_\MH - \comm{\ssupp_{\lambda \psi}}) \MH$ was tacitly used which holds true since $\ssupp_{\comm{\MFM}}(\lambda \psi) = J \ssupp_\MFM(\lambda \psi) J = J \ssupp_\MFM(\psi) J = \ssupp_{\comm{\MFM}}(\psi)$ according to \cref{lem:vonNeumann_commMsupport}. The above relation implies $S_{\mu \varphi, \lambda \psi}^0 = \sqrt{\mu / \lambda} \, S_{\varphi, \psi}^0$; since this identity is stable under taking closures, it follows that
    \begin{equation*}
        \Delta_{\mu \varphi, \lambda \psi} = S_{\mu \varphi, \lambda \psi}^\ast S_{\mu \varphi, \lambda \psi} = \frac{\mu}{\lambda} \, S_{\varphi, \psi}^\ast S_{\varphi, \psi} = \frac{\mu}{\lambda} \, \Delta_{\varphi, \psi} \ . \tag*{\qedhere}
    \end{equation*}
\end{Proof}

The following result is also standard in modular theory; it can be found, for example, in \cite[Rem. 3.3]{Araki77}, and the proof given below is a more detailed version of the one given in this reference.

\begin{lemma}\label{lem:vonNeumann_relativeModularOpUnitary}
    Let $(\MFM, \MH, J, \NPC)$ be a von Neumann algebra in standard form, let $\varphi, \psi \in \NF{\MFM}$ be positive normal functionals with vector representatives $\Phi, \Psi \in \NPC$, and let $U \in \UO(\MH)$ be a unitary operator such that $U A U^\ast = A$ for all $A \in \MFM$ and $U \NPC \subset \NPC$. Then
    \begin{equation*}
        \Delta_{U \Phi, U \Psi} = U \Delta_{\Phi, \Psi} U^\ast \ .
    \end{equation*}
\end{lemma}

\begin{Proof}
    Define $\Phi^\prime \ce U \Phi$ and $\Psi^\prime \ce U \Psi$. By assumption, these are still vectors in the natural positive cone. Let $\psi^\prime \ce \omega_{\Psi^\prime}$ be the vector functional induced by $\Psi^\prime$ on $\MFM$, and let $\eta^\prime \in (\id_\MH - \comm{\ssupp}_{\psi^\prime}) \MH$ be arbitrary. The first assumption on $U$ implies that $U \in \comm{\MFM}$; since $U$ is also continuous, it follows that $\ran(\comm{\ssupp}_{\psi^\prime}) = [\MFM U \Psi] = U [\MFM \Psi]$, and hence that\footnotemark
    \footnotetext{Let $(\MH, \bdot)$ be a Hilbert space, $\MK \subset \MH$ be a subset, $U \in \UO(\MH)$ be a unitary operator, and $\xi \in \MK$ be arbitrary. Then $\eta \in (U \MK)^\perp$ $\iff$ $\braket{U \xi, \eta} = 0$ $\iff$ $\braket{\xi, U^\ast \eta} = 0$ $\iff$ $U^\ast \eta \in \MK^\perp$ $\iff$ $\eta \in U \MK^\perp$. Therefore, $(U \MK)^\perp = U \MK^\perp$.}
    $\eta^\prime \in U [\MFM \Psi]^\perp$. Consequently, there is $\eta \in [\MFM \Psi]^\perp$ such that $\eta^\prime = U \eta$. Furthermore, observe that for all $A \in \MFM$, there holds
    \begin{equation*}
        \psi^\prime(A) = \braket{U \Psi, A U \Psi} = \braket{\Psi, U^\ast A U \Psi} = \psi(A) \ ,
    \end{equation*}
    hence $\psi^\prime = \psi$ on $\MFM$ and, in particular, $\ssupp(\psi^\prime) = \ssupp(\psi)$. Thus, one obtains
    \begin{equation*}
        S_{\Phi^\prime, \Psi^\prime} (A \Psi^\prime + \eta^\prime) = S_{U \Phi, U \Psi} (A U \Psi + U \eta) = S_{U \Phi, U \Psi} U (A \Psi + \eta)
    \end{equation*}
    on the one hand, and
    \begin{equation*}
        S_{\Phi^\prime, \Psi^\prime} (A \Psi^\prime + \eta^\prime) = \ssupp(\psi^\prime) A^\ast \Phi^\prime = \ssupp(\psi) A^\ast U \Phi = U \ssupp(\psi) A^\ast \Phi = U S_{\Phi, \Psi} (A \Psi + \eta)
    \end{equation*}
    on the other hand. Combining these two relations, it follows that $S_{U \Phi, U \Psi} U = U S_{\Phi, \Psi}$ and, consequently, $S_{U \Phi, U \Psi} = U S_{\Phi, \Psi} U^\ast$ on $\dom(S_{\Phi, \Psi})$. Therefore, the modular operator $\Delta_{\Phi^\prime, \Psi^\prime}$ takes the following form:
    \begin{equation*}
        \Delta_{U \Phi, U \Psi} = S_{U \Phi, U \Psi}^\ast S_{U \Phi, U \Psi} = (U S_{\Phi, \Psi}^\ast U^\ast) (U S_{\Phi, \Psi} U^\ast) = U S_{\Phi, \Psi}^\ast S_{\Phi, \Psi} U^\ast = U \Delta_{\Phi, \Psi} U^\ast \ . \tag*{\qedhere}
    \end{equation*}
\end{Proof}

This section shall be concluded with a result concerning convergence properties of the relative modular operator which will be used in the perturbation theory of KMS-states in \cref{ch:perturbationTheory}. The proof presented below is a more detailed version of \cite[Thm. 4.2]{DJP03}.

\begin{lemma}\label{lem:vonNeumann_convergenceRelModOp}
    Let $\MFM$ be a von Neumann algebra represented in standard form $(\MFM, \MH, J, \NPC)$, and let $(\Phi_n)_{n \in \N} \subset \NPC$ and $(\Psi_n)_{n \in \N} \subset \NPC$ be two sequences. Assume that the following properties are satisfied:
    \begin{enumerate}[cond]
        \item \label{enu:vonNeumann_convergenceRelModOpA1} $\Delta_{\Phi_n, \Psi_n} \to M$ in the strong resolvent sense (\cref{def:operators_SRConvergence});
        \item \label{enu:vonNeumann_convergenceRelModOpA2} $\ssupp_{\Psi_n} \to \ssupp_\Psi$ in the strong operator topology;
        \item \label{enu:vonNeumann_convergenceRelModOpA4} $\Phi_n \to \Phi$ weakly in $\MH$;
        \item \label{enu:vonNeumann_convergenceRelModOpA3} $\Psi_n \to \Psi$ weakly in $\MH$.
    \end{enumerate}
    Then it follows that $M = \Delta_{\Phi, \Psi}$.
\end{lemma}

\begin{Proof}
    Let $A \in \MFM$ be arbitrary. From the representation $S_{\Phi, \Psi} = J \Delta_{\Phi, \Psi}^{1/2}$ of the relative Tomita operator (\cref{pro:vonNeumann_propertiesRelativeModularOp} \ref{enu:vonNeumann_polarDecompRelTomita}) and the definition \eqref{eq:vonNeumann_relativeTomitaOperator} of the latter, it follows that
    \begin{equation*}
        \Delta_{\Phi_n, \Psi_n}^{1/2} A \Psi_n = J S_{\Phi_n, \Psi_n} A \Psi_n = J \, \ssupp_{\Psi_n} A^\ast \Phi_n \ .
    \end{equation*}
    Furthermore, $A \Psi_n \to A \Psi$ weakly in $\MH$ since $\braket{\xi, A \Psi_n} = \braket{A^\ast \xi, \Psi_n} \to \braket{A^\ast \xi, \Psi} = \braket{\xi, A \Psi}$ for all $\xi \in \MH$ according to assumption \ref{enu:vonNeumann_convergenceRelModOpA3}; similarly, the assumptions \ref{enu:vonNeumann_convergenceRelModOpA2} and \ref{enu:vonNeumann_convergenceRelModOpA4} together imply that $J \, \ssupp_{\Psi_n} A^\ast \Phi_n \to J \, \ssupp_\Psi A^\ast \Phi$ weakly in $\MH$. Finally, it holds that $\Delta_{\Phi_n, \Psi_n}^{1/2} \to M^{1/2}$ in the strong resolvent sense by \cref{lem:operators_SRimpliesFSR} and assumption \ref{enu:vonNeumann_convergenceRelModOpA1}.
    
    Therefore, applying \cref{pro:operators_SRimpliesW} (with $T_n \equiv \Delta_{\Phi_n, \Psi_n}^{1/2}$, $T \equiv M^{1/2}$, $\Omega_n \equiv A \Psi_n$, $\Omega \equiv A \Psi$ while noting that $\dlim{w}{n \to \infty} T_n \Omega_n$ exists so that \cref{rem:operators_SRimpliesW} applies), it follows that $A \Psi \in \dom(M^{1/2})$ and
    \begin{equation}\label{eq:vonNeumann_MbehavesRelModOp}
        M^{1/2} A \Psi = \dlim{w}{n \to \infty} \Delta_{\Phi_n, \Psi_n}^{1/2} A \Psi_n = \dlim{w}{n \to \infty} J \, \ssupp_{\Psi_n} A^\ast \Phi_n = J \ssupp_\Psi A^\ast \Phi \ .
    \end{equation}
    Let $\eta \in (\id_\MH - \comm{\ssupp_\Psi}) \, \MH = [\MFM \Psi]^\perp$ be arbitrary and define $\eta_n \ce (\id_\MH - \comm{\ssupp_{\Psi_n}}) \, \eta$ for all $n \in \N$. It holds that $\comm{\ssupp_{\Psi_n}} \to \comm{\ssupp_\Psi}$ strongly because $\comm{\ssupp_{\Psi_n}} = J \, \ssupp_{\Psi_n} J$ by \cref{lem:vonNeumann_commMsupport} and $\ssupp_{\Psi_n} \to \ssupp_{\Psi}$ strongly by assumption \ref{enu:vonNeumann_convergenceRelModOpA2}. Therefore, it follows that $\eta_n \to \eta$ strongly in $\MH$. Moreover, since
    \begin{equation*}
        \Delta_{\Phi_n, \Psi_n}^{1/2} \eta_n = 0 \quad (n \in \N)
    \end{equation*}
    according to the construction of the relative Tomita operator \eqref{eq:vonNeumann_relativeTomitaOperator}, one obtains $\eta \in \dom(M^{1/2})$ and $M^{1/2} \eta = 0$ by applying \cref{pro:operators_SRimpliesW} once more. Together with the above identity \eqref{eq:vonNeumann_MbehavesRelModOp} and uniqueness of the polar decomposition of $S_{\Phi, \Psi}$, this yields that $M^{1/2} = \Delta_{\Phi, \Psi}^{1/2}$.
\end{Proof}

\section{Spatial Derivatives}\label{sec:vonNeumann_spatialDerivative}

The relative modular operator $\Delta_{\psi, \varphi}$ generalized the modular operator $\Delta_\omega$ of a faithful normal functional $\omega \in \NF{\MFM}$, and \cref{sec:vonNeumann_relativeModularOperator} demonstrated that it always acts on the Hilbert space from the standard form representation. In this section, the spatial derivative operator, a generalization of the relative modular operator which was introduced by \textsc{A. Connes} in 1980 \cite{Connes80}, will be constructed. It is defined with respect to a \emph{semi-finite normal weight} $\varphi$ on $\MFM$ and a \emph{faithful semi-finite normal weight} $\psi$ on $\comm{\MFM}$. The advantage of the spatial derivative over the relative modular operator is that it is defined in \emph{any representation} space for $\MFM$, not just the one from the standard form representation.

\subsection{Weights}\label{subsec:vonNeumann_weights}

As mentioned, the spatial derivative is defined with respect to a so-called normal weight on the commutant; this is, in a sense, an \enquote{unbounded normal functional} \cite[p. 318]{StratilaZsido19}. The definitions and results presented here can be found, for example, in \cite[Sect. 6.1 \& 11.1]{Hiai21} or \cite[Sect. 1 \& 7]{Stratila20}.

\begin{definition}[Weight]\label{def:vonNeumann_weight}
    Let $\MFM \subset \BO(\MH)$ be a von Neumann algebra. A \bemph{weight} on $\MFM$ is a functional $\varphi : \PE{\MFM} \longto [0, + \infty]$ satisfying $\varphi(A + B) = \varphi(A) + \varphi(B)$ and $\varphi(\lambda A) = \lambda \varphi(A)$ for all $A, B \in \PE{\MFM}$ and $\lambda \ge 0$. Define the following set which is actually a left ideal of $\MFM$:
    \begin{equation*}
        \MFN_\varphi \ce \set[\big]{A \in \MFM \ : \ \varphi(A^\ast A) < + \infty} \ .
    \end{equation*}
    Then the weight $\varphi$ is called
    \begin{enumerate}
        \item \bemph{faithful} iff $\varphi(A^\ast A) = 0$ implies $A = 0$ for any $A \in \MFM$;
        
        \item \bemph{semi-finite} iff $\MFN_\varphi$ is $\sigma$-weakly dense in $\MFM$;

        \item \bemph{normal} iff $\varphi(\sup_{i \in I} A_i) = \sup_{i \in I} \varphi(A_i)$ for any increasing net $(A_i)_{i \in I} \subset \PE{\MFM}$.
    \end{enumerate}
    
    Analogously as for normal functionals, one defines the \bemph{support} of a normal weight $\varphi$ to be the smallest projection $\ssupp(\varphi) \in \PO(\MFM)$ such that $\varphi\bigl(\id_\MH - \ssupp(\varphi)\bigr) = 0$ \cite[p. 10]{Stratila20}. As in the previous case, it holds that $\varphi$ is faithful if and only if $\ssupp(\varphi) = \id_\MH$. Below, the following notation will be used:
    \begin{equation*}
        \supp(\varphi) \ce \ran\bigl(\ssupp(\varphi)\bigr) \ .
    \end{equation*}
\end{definition}

\begin{example}\label{exa:vonNeumann_weights}
    (\cite[Exa. 6.3]{Hiai21})
    Let $\MFM = L^\infty(X, \mu)$ for a $\sigma$-finite measure space $(X, \Sigma, \mu)$, and define $\varphi : \PE{\MFM} \longto [0, + \infty]$ by $\varphi(f) \ce \int_X f \diff \mu$. Then $\varphi$ is a faithful semi-finite normal weight on $\MFM$ with $\MFN_\varphi = \MFM \cap L^2(X, \mu)$. Next, consider $\MFM = \BO(\MH)$ for a Hilbert space $\MH$. Then the trace functional $\tr : \BO(\MH) \longto [0, + \infty]$, $A \longmto \tr(A)$, is a faithful semi-finite normal weight with $\MFN_{\tr} = \HS(\MH)$.
\end{example}

\begin{para}[The GNS-construction]\label{para:vonNeumann_weightGNS}
    (\cite[p. 88]{Hiai21}, \cite[p. 2]{Stratila20})
    Similarly as for a positive linear functional on a $C^\ast$-algebra (\cf{} \cref{subsec:operatorAlgebras_GNS}), one can also perform a GNS-construction for a weight $\varphi$ on a von Neumann algebra $\MFM \subset \BO(\MH)$. Define the inner product $\braket{X, Y}_\varphi \ce \varphi(X^\ast Y)$ on $\MFN_\varphi$, and let $\MH_\varphi$ denote the completion of $\MFN_\varphi$ with respect to $\bdot_\varphi$. Furthermore, let $\eta_\varphi : \MFN_\varphi \longto \MH_\varphi$, $A \longmto \eta_\varphi(A)$, be the canonical injection. Since $\MFN_\varphi$ is a left-ideal of $\MFM$, and since
    \begin{equation}\label{eq:vonNeumann_weightInequality}
        \varphi(X^\ast A^\ast A X) \le \norm{A}_\mop^2 \braket{X, X}_\varphi
    \end{equation}
    for all $A \in \MFM$, $X \in \MFN_\varphi$, it follows that there exists a $\ast$-representation $\pi_\varphi : \MFM \longto \BO(\MH_\varphi)$ which is uniquely determined by $\pi_\varphi(A) \eta_\varphi(X) \ce \eta_\varphi(AX)$. This representation will be faithful if the weight $\varphi$ is faithful, and in this case one may consider $\MFM$ as a von Neumann algebra on $\MH_\varphi$.
\end{para}

The next theorem, taken from \cite[Thm. 6.2]{Hiai21}, was originally discovered by \textsc{U. Haagerup} in 1975 \cite{Haagerup75b}. A proof may be found in \cite[Thm. 1.3 \& Cor. 5.9]{Stratila20}.

\begin{ntheorem}[Haagerup]\label{thm:vonNeumann_Haagerup}
    For a weight $\varphi$ on a von Neumann algebra $\MFM \subset \BO(\MH)$, the following properties are equivalent:
    \begin{enumerate}[equiv]
        \item $\varphi$ is normal;
        \item $\varphi$ is $\sigma$-weakly lower semi-continuous;
        \item $\varphi(A) = \sup \set{\omega(A) : \omega \in \NF{\MFM}, \omega \le \varphi}$ for all $A \in \PE{\MFM}$;
        \item $\varphi(A) = \sum_{i \in I} \braket{\xi_i, A \xi_i}$ for all $A \in \PE{\MFM}$, where $(\xi_i)_{i \in I} \subset \MH$.
    \end{enumerate}
\end{ntheorem}

The object to be introduced next is required for the construction of the spatial derivative operator; the terminology is borrowed from \cite[Ch. 4]{OP04}.

\begin{definition}[Lineal]\label{def:vonNeumann_linealWeight}
    Let $\MFM \subset \BO(\MH)$ be a von Neumann algebra, let $\psi$ be a semi-finite normal weight on $\MFM$, and let $(\MH_{\psi}, \pi_{\psi}, \eta_{\psi})$ be the GNS-representation of $\MFM$ with respect to $\psi$. Define the \bemph{lineal} of $\psi$ to be the set
    \begin{equation}\label{eq:vonNeumann_linealWeight}
        D(\MH, \psi) \ce \set[\big]{\xi \in \MH \ : \ \exists C_\xi > 0 \ \forall A \in \MFN_\psi : \norm{A \xi}_\MH \le C_\xi \, \norm{\eta_\psi(A)}_{\MH_\psi}}
    \end{equation}
    consisting of so-called \bemph{$\psi$-bounded vectors}; this terminology is due to the fact that $\norm{\eta_\psi(A)}_{\MH_\psi} = \psi(A^\ast A)^{1/2}$. Note that the lineal is a linear subspace of $\MH$. Furthermore, if $\psi \in \NF{\MFM}$ happens to be a normal functional on $\MFM$, then the lineal is given by
    \begin{equation}\label{eq:vonNeumann_linealNormalFunctional}
        D(\MH, \psi) = \set[\big]{\xi \in \MH \ : \ \exists C_\xi > 0 \ \forall A \in \MFM : \norm{A \xi}_\MH^2 \le C_\xi \, \psi(A^\ast A)} \ .
    \end{equation}
\end{definition}

The proof of the following result follows \cite[Lem. 11.2]{Hiai21} and \cite[Sect. 7.1]{Stratila20}.

\begin{proposition}\label{pro:vonNeumann_propertiesLineal}
    In the terminology of \cref{def:vonNeumann_linealWeight}, the following assertions hold true.
    \begin{enumerate}
        \item \label{enu:vonNeumann_linealInvariant} The lineal $D(\MH, \psi)$ is invariant under the commutant $\comm{\MFM}$.
        
        \item \label{enu:vonNeumann_linealDense} If $\psi$ is assumed to be faithful, then $D(\MH, \psi)$ lies dense in $\MH$.
    \end{enumerate}
\end{proposition}

\begin{Proof}
    \tAd{} \ref{enu:vonNeumann_linealInvariant}. Let $\xi \in D(\MH, \psi)$ and $\comm{A} \in \comm{\MFM}$ be arbitrary. For all $A \in \MFN_\psi$, there holds
    \begin{equation*}
        \norm{A (\comm{A} \xi)}_\MH = \norm{\comm{A} (A \xi)}_\MH \le \norm{\comm{A}}_\mop \norm{A \xi}_\MH \le \norm{\comm{A}}_\mop \, C_\xi \, \norm{\eta_\psi(A)}_{\MH_\psi} \ .
    \end{equation*}

    \tAd{} \ref{enu:vonNeumann_linealDense}. Let $P \in \PO(\MH)$ denote the projection onto $\clos_{\ndot} \bigl(D(\MH, \psi)\bigr)$. It will be shown that $P = \id_\MH$; the assertion then follows from \cref{cor:operators_characterizationPerp} \ref{enu:operators_densePerp}. First, observe that the invariance of the lineal under the commutant implies that $P \comm{A} P = \comm{A} P = \comm{A} P P$ for all $\comm{A} \in \comm{\MFM}$. This shows that actually $P \in \bicomm{\MFM} = \MFM$.  Next, because $\psi$ is a normal weight, \cref{thm:vonNeumann_Haagerup} implies that $\psi = \sum_{i \in I} \omega_{\xi_i}$ for some family $(\xi_i)_{i \in I} \subset \MH$. For every $i \in I$, there holds $\xi_i \in D(\MH, \psi)$; this follows from the representation of $\psi$ and noting that for all $A \in \MFN_\psi$,
    \begin{equation}\label{eq:vonNeumann_vectorReprInLineal}
        \norm{A \xi_i}_\MH^2 = \braket{\xi_i, A^\ast A \xi_i}_\MH \le \psi(A^\ast A) = \norm{\eta_\psi(A)}_{\MH_\psi}^2 \ .
    \end{equation}
    Therefore, $P \xi_i = \xi_i$ for all $i \in I$ by definition of $P$, and so one obtains that $\psi(\mathord{\id_\MH} - P) = 0$. As $\psi$ was assumed to be faithful, it follows that $P = \id_\MH$. (Note that the previous observation $P \in \MFM$ was used in order to apply the weight $\psi$ to the element $\id_\MH - P \in \PE{\MFM}$.)
\end{Proof}

The next result, mentioned but not proved in \cite[p. 69]{OP04}, will be required down the line to establish another result which will find crucial application in \cref{ch:relativeEntropy}.

\begin{lemma}\label{lem:vonNeumann_linealVectorState}
    If $\psi = \omega_\Psi$ is a vector functional on $\MFM$ induced by $\Psi \in \MH$, then $D(\MH, \psi) = \comm{\MFM} \Psi$.
\end{lemma}

\begin{Proof}
    \pp{1.} Note that \cref{eq:vonNeumann_vectorReprInLineal} in the previous proof shows that the vector $\Psi$ is contained in the lineal of the state $\psi$ it induces, and since $D(\MH, \psi)$ is invariant under the commutant of $\MFM$ by \cref{pro:vonNeumann_propertiesLineal} \ref{enu:vonNeumann_linealInvariant}, it follows that $\comm{\MFM} \Psi \subset D(\MH, \psi)$.

    \pp{2.} Let $\xi \in D(\MH, \psi)$ be arbitrary, that is, $\norm{A\xi}^2 \le C_\xi \, \norm{A \Psi}^2$ for all $A \in \MFM$ by \eqref{eq:vonNeumann_linealNormalFunctional}. Consider the subspaces $\MH_\Psi \ce [\MFM \Psi]$ and $\MH_\xi \ce [\MFM \xi]$ of $\MH$; they, as well as their orthogonal complements, are invariant under $\MFM$. Define a linear operator $\comm{T_\xi} : \MFM \Psi \longto \MFM \xi$ in $\MH$ by setting $\comm{T_\xi} (A \Psi) \ce A \xi$ for all $A \in \MFM$. The previous inequality shows that $\comm{T_\xi}$ is bounded, and that $A \Psi = 0$ implies $\comm{T_\xi} (A \Psi) = A \xi = 0$. Therefore, $\comm{T_\xi}$ is well-defined and extends uniquely to a linear operator $\comm{T_\xi} : \MH_\Psi \longto \MH_\xi$ (\cref{thm:operators_BLT}). Note that if $(A_n \Psi)_{n \in \N}$ and $(B_n \Psi)_{n \in \N}$ are two different sequences approximating the element $\eta \in \MH_\Psi$, \ie{}, $\dlim{\MH}{n \to \infty} A_n \Psi = \dlim{\MH}{n \to \infty} B_n \Psi = \eta$, then
    \begin{equation*}
        \norm{A_n \xi - B_n \xi} \le \sqrt{C_\xi} \, \norm{(A_n - B_n) \Psi} \longto 0 \tas n \longto + \infty \ ,
    \end{equation*}
    so the definition of $\comm{T_\xi}$ is indeed independent of the approximating sequence. This establishes that $\comm{T_\xi} : \MH_\Psi \longto \MH_\xi$ is a well-defined bounded linear operator. Finally, since $\MH = \MH_\Psi \oplus (\MH_\Psi)^\perp$ (\cref{thm:operators_orthogonalProjection}), one can extend the operator $\comm{T_\xi}$ to the whole space $\MH$ by setting $\comm{T_\xi} \zeta \ce 0$ for all $\zeta \in (\MH_\Psi)^\perp$.

    Let $A \in \MFM$ and $\zeta \in (\MH_\Psi)^\perp$ be arbitrary. As mentioned above, it holds that $A \zeta \in (\MH_\Psi)^\perp$ and, moreover, $A \comm{T_\xi} \zeta = 0 = \comm{T_\xi} (A \zeta)$ by definition of $\comm{T_\xi}$; this shows that $\comm{T_\xi}$ commutes with $A$ on the subspace $(\MH_\Psi)^\perp$. Let now $\eta = \dlim{\MH}{n \to \infty} B_n \Psi \in \MH_\Psi$ be arbitrary. Then
    \begin{align*}
        A \comm{T_\xi} \eta = \dlim{\MH}{n \to \infty} A B_n \xi = \dlim{\MH}{n \to \infty} \comm{T_\xi} (A B_n \Psi) = \comm{T_\xi} (A \eta) \ ,
    \end{align*}
    that is, $\comm{T_\xi}$ commutes with $A$ on $\MH_\Psi$ as well, \ie{}, $\comm{T_\xi} \in \comm{\MFM}$. By construction, there holds $\xi = \comm{T_\xi} \Psi$, whence one obtains that $\xi \in \comm{\MFM} \Psi$ which concludes the proof.
\end{Proof}

\subsection{Intertwining Representations}

As mentioned in the introduction of this section, the goal is to construct an analogue of the relative modular operator which can be employed in any representation of the von Neumann algebra. Therefore, it is fruitful to study how the GNS-representation intertwines with other representations.

\begin{para}[Intertwining operators]\label{para:vonNeumann_intertwiningOperator}
    Let $\psi$ be a faithful semi-finite normal weight on the von Neumann algebra $\MFM \subset \BO(\MH)$, and let $(\MH_\psi, \pi_\psi, \eta_\psi)$ be the GNS-representation with respect to $\psi$. For arbitrary $\xi \in D(\MH, \psi)$, define a linear operator $R_\xi^\psi : \MH_\psi \longto \MH$ by setting
    \begin{equation*}
        R_\xi^\psi \, \eta_\psi(A) \ce A \xi \ \quad (A \in \MFN_\psi) \ .
    \end{equation*}
    Since this defines a bounded operator by the definition \eqref{eq:vonNeumann_linealWeight} of the lineal of $\psi$, and since $\MFN_\psi$ is dense in $\MH_\psi$ by construction (\cref{para:vonNeumann_weightGNS}), it follows that $R_\xi^\psi$ extends uniquely to the whole space $\MH_\psi$ (\cref{thm:operators_BLT}). The adjoint of this operator is defined by the relation \cite[p. 77]{Stratila20}
    \begin{equation*}
        \vbraket[\big]{(R_\xi^\psi)^\ast \zeta, \eta_\psi(A)}_{\MH_\psi} = \braket{\zeta, A \xi}_\MH \quad (\zeta \in \MH, \, A \in \MFN_\psi) \ .
    \end{equation*}
    For $\xi \in D(\MH, \psi)$, one furthermore defines a bounded self-adjoint operator $\Theta^\psi(\xi) : \MH \longto \MH$ by
    \begin{equation*}
        \Theta^\psi(\xi) \ce R_\xi^\psi (R_\xi^\psi)^\ast \ .
    \end{equation*}
    
    Observe that if $\psi \in \NF{\MFM}$, then $R_\xi^\psi$ is again defined using the GNS-representation $(\MH_\psi, \pi_\psi, \Psi_\psi)$ associated with $\psi$: for all $\xi \in D(\MH, \psi)$ and $A \in \MFM$, one sets \cite[Eq. (4.3)]{OP04}
    \begin{equation*}
        R_\xi^\psi \bigl(\pi_\psi(A) \Psi_\psi\bigr) \ce A \xi \ .
    \end{equation*}
\end{para}

The following proposition collects basic properties of $R_\xi^\psi$ and $\Theta^\psi(\xi)$ \cite[Lem. 11.2]{Hiai21}.

\begin{proposition}\label{pro:vonNeumann_propertiesIntertwiningOp}
    \leavevmode
    \begin{enumerate}
        \item \label{enu:vonNeumann_intertwiningOpArgument} For every $\xi \in D(\MH, \psi)$ and $\comm{A} \in \comm{\MFM}$, it holds that $R_{\comm{A} \xi}^\psi = \comm{A} R_\xi^\psi$.
        
        \item \label{enu:vonNeumann_intertwining} The operator $R_\xi^\psi$ intertwines the GNS-representation of $\MFM$ on $\MH_\psi$ with the identity representation $\MFM \owns A \longmto A \in \BO(\MH)$ of $\MFM$ on $\MH$: $A R_\xi^\psi = R_\xi^\psi \pi_\psi(A)$ for all $A \in \MFM$.
        
        \item \label{enu:vonNeumann_thetaComm} If $\xi_1, \xi_2 \in D(\MH, \psi)$, then $R_{\xi_1}^\psi (R_{\xi_2}^\psi)^\ast \in \comm{\MFM}$. In particular, $\Theta^\psi(\xi) \in \comm{\MFM}$ for all $\xi \in D(\MH, \psi)$.
    \end{enumerate}
\end{proposition}

\begin{Proof}
    \tAd{} \ref{enu:vonNeumann_intertwiningOpArgument}. It was observed before in \cref{pro:vonNeumann_propertiesLineal} \ref{enu:vonNeumann_linealInvariant} that for $\xi \in D(\MH, \psi)$ and $\comm{A} \in \comm{\MFM}$, also $\comm{A} \xi \in D(\MH, \psi)$. For arbitrary $X \in \MFN_\psi$, one obtains $R_{\comm{A} \xi}^\psi \, \eta_\psi(X) = X \comm{A} \xi = \comm{A} X \xi = \comm{A} R_\xi^\psi \, \eta_\psi(X)$, and this identity extends, by continuity, to the whole Hilbert space $\MH_\psi$.

    \tAd{} \ref{enu:vonNeumann_intertwining}. For every $\xi \in D(\MH, \psi)$, $A \in \MFM$, and $X \in \MFN_\psi$, one finds using the definition of the representation $\pi_\psi$ that $R_\xi^\psi \pi_\psi(A) \eta_\psi(X) = R_\xi^\psi \eta_\psi(A X) = A X \xi = A R_\xi^\psi \eta_\psi(X)$, hence $R_\xi^\psi \pi_\psi(A) = A R_\xi^\psi$ on $\MFN_\psi$ which, again, extends by continuity to $\MH_\psi$.

    \tAd{} \ref{enu:vonNeumann_thetaComm}. Applying the previous observation multiple times, one computes for arbitrary $A \in \MFM$:
    \begin{align*}
        A R_{\xi_1}^\psi (R_{\xi_2}^\psi)^\ast &= R_{\xi_1}^\psi \pi_\psi(A) (R_{\xi_2}^\psi)^\ast = R_{\xi_1}^\psi \Bigl(R_{\xi_2}^\psi \pi_\psi(A)^\ast\Bigr)^\ast \\
        &= R_{\xi_1}^\psi \Bigl(R_{\xi_2}^\psi \pi_\psi(A^\ast)\Bigr)^\ast = R_{\xi_1}^\psi \Bigl(A^\ast R_{\xi_2}^\psi\Bigr)^\ast = R_{\xi_1}^\psi (R_{\xi_2}^\psi)^\ast A \ .\tag*{\qedhere}
    \end{align*}
\end{Proof}

The next result will be a key ingredient in the proof of Uhlmann's monotonicity theorem for the relative entropy in \cref{sec:relativeEntropy_monotonicity}. The assertion is mentioned but not proved in \cite[p. 70]{OP04}.

\begin{proposition}\label{pro:vonNeumann_ThetaVectorState}
    Let $\psi = \omega_\Psi$ be a vector functional on $\MFM \subset \BO(\MH)$ induced by a vector $\Psi \in \MH$. Then it follows that for every $\comm{A} \in \comm{\MFM}$, the operator $\Theta^\psi(\comm{A} \Psi)$ is given by
    \begin{equation*}
        \Theta^\psi(\comm{A} \Psi) = \comm{A} \, [\MFM \Psi] (\comm{A})^\ast \ .
    \end{equation*}
\end{proposition}

\begin{Proof}
    According to \cref{lem:vonNeumann_linealVectorState}, it holds that $D(\MH, \psi) = \comm{\MFM} \Psi$, so one can indeed consider the operator $\Theta^\psi(\xi)$ for vectors $\xi = \comm{A} \Psi$, $\comm{A} \in \comm{\MFM}$. (In particular, $\Psi$ is itself an element of $D(\MH, \psi)$.) Using \cref{pro:vonNeumann_propertiesIntertwiningOp} \ref{enu:vonNeumann_intertwiningOpArgument} and the definition of $\Theta^\psi(\xi)$, one finds
    \begin{align*}
        \Theta^\psi(\comm{A} \Psi) &= R_{\comm{A} \Psi}^\psi (R_{\comm{A} \Psi}^\psi)^\ast = \comm{A} R_\Psi^\psi \bigl(\comm{A} R_\Psi^\psi\bigr)^\ast = \comm{A} R_\Psi^\psi (R_\Psi^\psi)^\ast (\comm{A})^\ast \ .
    \end{align*}
    
    It remains to show that $P \ce R_\Psi^\psi (R_\Psi^\psi)^\ast$ is the orthogonal projection $[\MFM \Psi]$ onto $U \ce \clos_{\ndot}(\MFM \Psi)$. To this end, let $(\MH_\psi, \pi_\psi, \Psi_\psi)$ be the GNS-representation with respect to the functional $\psi$, and observe first that for all $A, B \in \MFM$, one obtains
    \begin{align*}
        \vbraket[\big]{(R_\Psi^\psi)^\ast A \Psi, \pi_\psi(B) \Psi_\psi}_{\MH_\psi} &= \vbraket[\big]{A \Psi,  R_\Psi^\psi \bigl(\pi_\psi(B) \Psi_\psi\bigr)}_\MH = \braket{A \Psi, B \Psi}_\MH \\
        &= \psi(A^\ast B) = \vbraket[\big]{\Psi_\psi, \pi_\psi(A^\ast B) \Psi_\psi}_{\MH_\psi} = \braket{\pi_\psi(A) \Psi_\psi, \pi_\psi(B) \Psi_\psi}_{\MH_\psi} \ .
    \end{align*}
    Since $\pi_\psi(\MFM) \Psi_\psi \subset \MH_\psi$ is dense, this shows that $(R_\Psi^\psi)^\ast A \Psi = \pi_\psi(A) \Psi_\psi$ for all $A \in \MFM$, and hence it follows that $(R_\Psi^\psi)^\ast R_\Psi^\psi = \id_\MH$. Therefore,
    \begin{equation*}
        P^2 = \bigl(R_\Psi^\psi (R_\Psi^\psi)^\ast\bigr)^2 = R_\Psi^\psi (R_\Psi^\psi)^\ast R_\Psi^\psi (R_\Psi^\psi)^\ast = R_\Psi^\psi (R_\Psi^\psi)^\ast = P
    \end{equation*}
    which shows that $P = R_\Psi^\psi (R_\Psi^\psi)^\ast$ is a projection. Since this operator is self-adjoint, it is automatically an orthogonal projection (\cref{def:operators_orthogonalProjections}). For every $\xi = \dlim{\MH}{n \to \infty} A_n \Psi \in [\MFM \Psi]$, it holds that
    \begin{equation*}
        R_\Psi^\psi (R_\Psi^\psi)^\ast \xi = \dlim{\MH}{n \to \infty} R_\Psi^\psi (R_\Psi^\psi)^\ast A_n \Psi = \dlim{\MH}{n \to \infty} A_n \Psi = \xi \ ,
    \end{equation*}
    so $P|_U = \id_\MH$ which implies $U = [\MFM \Psi] \subset \ran(P)$ by \cref{lem:operators_rangeProjection}. Moreover, for all $u^\perp \in U^\perp = [\MFM \Psi]^\perp$ and $B \in \MFM$, it follows that
    \begin{equation*}
        \vbraket[\big]{(R_\Psi^\psi)^\ast \, u^\perp, \pi_\psi(B) \Psi_\psi}_{\MH_\psi} = \vbraket[\big]{u^\perp, R_\Psi^\psi \bigl(\pi_\psi(B) \Psi_\psi\bigr)} = \braket{u^\perp, B \Psi} = 0 \ ,
    \end{equation*}
    thus $(R_\Psi^\psi)^\ast \, u^\perp = 0$, and this implies $[\MFM \Psi]^\perp \subset \ker(P)$. From this observation and \cref{lem:operators_projectionRestrictedID}, it follows that $\ran(P) = [\MFM \Psi]$, hence $P = R_\Psi^\psi (R_\Psi^\psi)^\ast$ is the orthogonal projection onto $U$.
\end{Proof}

\subsection{Construction of the Spatial Derivative}

The following lemma is the essential result needed for the definition of the spatial derivative as the unique self-adjoint operator associated with a suitable quadratic form \cite[Lem. 11.3]{Hiai21}, \cite[Sect. 7.3]{Stratila20}.

\begin{lemma}\label{lem:vonNeumann_quadraticFormSpatialDerivative}
    Let $\MFM \subset \BO(\MH)$ be a von Neumann algebra, let $\MFN \ce \comm{\MFM} \subset \BO(\MH)$ be its commutant, let $\varphi$ be a semi-finite normal weight on $\MFM$, and let $\psi$ be a faithful semi-finite normal weight on $\MFN$. Define a function $\MFq_\varphi : \MH \supset \dom(\MFq_\varphi) \longto \R$ by
    \begin{equation*}
        \MFq_\varphi(\xi) \ce \varphi \bigl(\Theta^\psi(\xi)\bigr) \com \xi \in \dom(\MFq_\varphi) \ce \set[\big]{\eta \in D(\MH, \psi) \ : \ \MFq_\varphi(\eta) < + \infty} \ .
    \end{equation*}
    Then $\MFq_\varphi$ is a densely defined, lower semi-continuous, positive quadratic form on $\MH$, \cf{} \cref{def:forms_quadraticForm}.
\end{lemma}

\begin{Proof}
    \pp{1. $\MFq_\varphi$ is a positive quadratic form}. Let $\lambda \in \C$ and $\xi, \eta \in \dom(\MFq_\varphi) = D(\MH, \psi)$ be arbitrary. Using \cref{pro:vonNeumann_propertiesIntertwiningOp} \ref{enu:vonNeumann_intertwiningOpArgument}, one finds that
    \begin{equation*}
        \MFq_\varphi(\lambda \xi) = \varphi \bigl(\Theta^\psi(\lambda \xi)\bigr) = \varphi \bigl(R_{\lambda \xi}^\psi (R_{\lambda \xi}^\psi)^\ast\bigr) = \varphi \bigl(\abs{\lambda}^2 R_\xi^\psi (R_\xi^\psi)^\ast\bigr) = \abs{\lambda}^2 \MFq_\varphi(\xi) \ .
    \end{equation*}
    Similarly, one computes
    \begin{align*}
        \MFq_\varphi(\xi + \eta) + \MFq_\varphi(\xi - \eta) &= \varphi \bigl(\Theta^\psi(\xi + \eta)\bigr) + \varphi \bigl(\Theta^\psi(\xi - \eta)\bigr) = \varphi \bigl(R_{\xi + \eta}^\psi (R_{\xi + \eta}^\psi)^\ast\bigr) + \varphi \bigl(R_{\xi - \eta}^\psi (R_{\xi - \eta}^\psi)^\ast\bigr) \\
        &= \varphi \bigl((R_\xi^\psi + R_\eta^\psi) (R_\xi^\psi + R_\eta^\psi)^\ast\bigr) + \varphi \bigl((R_\xi^\psi - R_\eta^\psi) (R_\xi^\psi - R_\eta^\psi)^\ast\bigr) \\
        &= 2 \varphi \bigl(R_\xi^\psi (R_\xi^\psi)^\ast\bigr) + 2 \varphi \bigl(R_\eta^\psi (R_\eta^\psi)^\ast\bigr) \\
        &= 2 \MFq_\varphi(\xi) + 2 \MFq_\varphi(\eta) \ .
    \end{align*}
    Furthermore, $\MFq_\varphi$ is positive because $\Theta^\psi(\xi)$ is a positive operator by construction, and $\varphi$ is, as a weight, positive as well. This shows that $\MFq_\varphi$ is a quadratic form in the sense of \cref{def:forms_quadraticForm}.

    \pp{2. $\MFq_\varphi$ is densely defined}. Since $\psi$ is faithful, it follows from \cref{pro:vonNeumann_propertiesLineal} \ref{enu:vonNeumann_linealDense} that $D(\MH, \psi) \subset \MH$ is dense. Moreover, because $\varphi$ is semi-finite, the left ideal $\MFN_\varphi$ is dense in $\MFM$ with respect to the $\sigma$-weak operator topology by \cref{def:vonNeumann_weight}. As the algebra involution is $\sigma$-weakly continuous (\cref{para:operatorAlgebras_sigmaWeakTop}), it follows that the set $\set[\big]{A \eta \, : \, A \in \MFN_\varphi^\ast, \, \eta \in D(\MH, \psi)}$ is dense in $\MH$ with respect to the norm topology \cite[p. 79]{Stratila20} and contained in $D(\MH, \psi)$ because the lineal of $\psi$ is invariant under $\comm{\MFN} = \bicomm{\MFM} = \MFM$ by \cref{pro:vonNeumann_propertiesLineal} \ref{enu:vonNeumann_linealInvariant}. Furthermore, for any $A \in \MFN_\varphi^\ast$ and $\xi \in D(\MH, \psi)$, one finds that
    \begin{align*}
        \MFq_\varphi(A \xi) = \varphi\bigl(R_{A \xi}^\psi (R_{A \xi}^\psi)^\ast\bigr) = \varphi\bigl(A R_{\xi}^\psi (R_{\xi}^\psi)^\ast A^\ast\bigr) \le \norm{R_\xi^\psi (R_\xi^\psi)^\ast}_\mop \, \varphi(A A^\ast) < + \infty
    \end{align*}
    by using \cref{pro:vonNeumann_propertiesIntertwiningOp} \ref{enu:vonNeumann_intertwiningOpArgument} and \cref{eq:vonNeumann_weightInequality}, hence $A \xi \in \dom(\MFq_\varphi)$. This shows that $\dom(\MFq_\varphi)$ contains a norm-dense subset, and so it is itself dense in $\MH$.

    \pp{3. $\MFq_\varphi$ is lower semi-continuous on $D(\MH, \psi)$}. By \cref{thm:vonNeumann_Haagerup}, there exists a family of vectors $(\xi_i)_{i \in I} \subset \MH$ such that $\varphi = \sum_{i \in I} \omega_{\xi_i}$. Furthermore, note that by the Hahn-Banach theorem \cite[Cor. III.1.7]{Werner18} and the representation theorem of Fréchet-Riesz \cite[Thm. V.3.6]{Werner18},
    \begin{equation*}
        \norm{\zeta}_{\MH_\psi} = \sup_{f \in \cdual{(\MH_\psi)},\, \norm{f} \le 1} \abs{f(\zeta)} = \sup_{\xi \in \MH_\psi,\, \norm{\xi} \le 1} \abs{\braket{\xi, \zeta}} = \sup_{A \in\MFN_\psi,\, \norm{A} \le 1} \abs{\braket{\eta_\psi(A), \zeta}}
    \end{equation*}
    for all $\zeta \in \MH_\psi$, where in the last step, it was used that $\MH_\psi = \clos_{\ndot}(\MFN_\psi)$ by the GNS-construction, hence it suffices to take the supremum over elements from $\MFN_\psi$. With this, one obtains for all $\eta \in D(\MH, \psi)$ and fixed $i \in I$ that
    \begin{equation*}
        \omega_{\xi_i}\bigl(\Theta^\psi(\eta)\bigr) = \vbraket[\big]{\xi_i, R_\eta^\psi (R_\eta^\psi)^\ast \xi_i}_\MH = \vnorm[\big]{(R_\eta^\psi)^\ast \xi_i}_{\MH_\psi}^2 = \sup_{A \in\MFN_\psi,\, \norm{A} \le 1} \vabs[\big]{\vbraket[\big]{\eta_\psi(A), (R_\eta^\psi)^\ast \xi_i}_{\MH_\psi}}^2 \ .
    \end{equation*}
    Since $\braket{\eta_\psi(A), (R_\eta^\psi)^\ast \xi_i}_{\MH_\psi} = \braket{A \eta, \xi_i}_\MH = \braket{\eta, A^\ast \xi_i}_\MH$, it follows that $\MFq_\varphi(\eta)$ takes the form
    \begin{align*}
        \MFq_\varphi(\eta) = \sum_{i \in I} \braket{\xi_i, R_\eta^\psi (R_\eta^\psi)^\ast \xi_i}_\MH = \sum_{i \in I} \sup_{A \in\MFN_\psi,\, \norm{A} \le 1} \vabs[\big]{\braket{\eta, A^\ast \xi_i}_\MH}^2 \ .
    \end{align*}
    The function $\MH \supset D(\MH, \psi) \owns \eta \longmto \abs{\braket{\eta, A^\ast \xi_i}_\MH}^2 \in \R$ is continuous, so in particular lower semi-continuous. As the supremum and the sum of lower semi-continuous functions are again lower semi-continuous \cite[Cor. 3.2.8]{NiculescuPersson18}, it follows that $\MFq_\varphi$ is lower semi-continuous.
\end{Proof}

\begin{para}[The spatial derivative operator]\label{para:vonNeumann_spatialDerivative}
    (\cite[Def. 11.4]{Hiai21})
    Let $\MFM \subset \BO(\MH)$ be a von Neumann algebra, $\MFN \ce \comm{\MFM} \subset \BO(\MH)$ be its commutant, and $\psi$ be a faithful semi-finite normal weight on $\MFN$. For any semi-finite normal weight $\varphi$ on $\MFM$, let $\MFq_\varphi$ be the densely defined positive quadratic form from above.
    
    By \cref{lem:vonNeumann_quadraticFormSpatialDerivative} and \cref{pro:forms_closableForms}, it follows that $\MFq_\varphi$ is a closable form, hence the closure $\ol{\MFq_\varphi}$ of $\MFq_\varphi$ exists which is a closed positive quadratic form and the smallest closed extension of $\MFq_\varphi$ (\cf{} \cref{para:forms_closure}). Applying the form representation theorem, \cref{thm:forms_formReprThm} (see also \cref{para:forms_fromFormToOp}), one obtains a uniquely defined positive self-adjoint operator $Q_{\MFq_\varphi}$ associated with $\MFq_\varphi$ which satisfies
    \begin{equation*}
        \dom\bigl(Q_{\MFq_\varphi}^{1/2}\bigr) = \dom(\ol{\MFq_\varphi}) \tand \ol{\MFq_\varphi}(\xi) = \vnorm[\big]{Q_{\MFq_\varphi}^{1/2} \xi}^2 \com \xi \in \dom(\ol{\MFq_\varphi}) \ .
    \end{equation*}
    This $Q_{\MFq_\varphi}$ is called the \bemph{spatial derivative operator} of the semi-finite normal weight $\varphi$ on $\MFM$ with respect to the faithful semi-finite normal weight $\psi$ on $\MFN = \comm{\MFM}$, and it will be denoted by
    \begin{equation*}
        \Delta(\varphi / \psi) \tor \od{\varphi}{\psi} \ .
    \end{equation*}
    
    According to the discussion in \cref{para:forms_fromFormToOp}, it holds that $\Delta(\varphi / \psi)$ is the largest positive self-adjoint operator, with respect to the ordering defined in \cref{para:forms_orderRelation}, such that
    \begin{equation}\label{eq:vonNeumann_spatialDerivativeForm}
        \dom\bigl(\Delta(\varphi / \psi)^{1/2}\bigr) \supset \dom(\MFq_\varphi) \tand \MFq_\varphi(\xi) = \vnorm[\big]{\Delta(\varphi / \psi)^{1/2} \xi}^2 \com \xi \in \dom(\MFq_\varphi) \ .
    \end{equation}
    Moreover, the form domain $\dom(\MFq_\varphi)$ is a core for $\Delta(\varphi / \psi)^{1/2}$, and if $\varphi \in \NF{\MFM}$ is a normal functional, then
    \begin{equation*}
        D(\MH, \psi) = \dom(\MFq_\varphi) \subset \dom\bigl(\Delta(\varphi / \psi)^{1/2}\bigr) \ .
    \end{equation*}
    Finally, one can show that the support of $\Delta(\varphi / \psi)$ is given by $\ssupp(\varphi)$ \cite[Thm. 11.7 (3)]{Hiai21}.
\end{para}

The following proposition shows that the spatial derivative is a generalization of the relative modular operator to the setting in which the von Neumann algebra $\MFM \subset \BO(\MH)$ is represented on an arbitrary Hilbert space $\MH$, not necessarily the standard representation \cite[p. 183]{Hiai21}. In particular, it follows that the examples from \cref{exa:vonNeumann_relativeModularOperator} also apply to the spatial derivative. The proof of the assertion may be found in \cite[Prop. 11.6]{Hiai21}.

\begin{proposition}\label{pro:vonNeumann_spatialDerivativeRelativeModOp}
    Let $(\MFM, \MH, J, \MP)$ be a von Neumann algebra in standard form, and let $\varphi, \psi \in \NF{\MFM}$ with $\psi$ faithful. Define a faithful $\comm{\psi} \in \NF{(\comm{\MFM})}$ by $\comm{\psi}(\comm{A}) \ce \psi(J (\comm{A})^\ast J)$ for $\comm{A} \in \comm{\MFM}$, that is, $\comm{\psi}(\comm{A}) = \braket{\Psi, \comm{A} \Psi}$, where $\Psi \in \MP$ is the vector representative of $\psi$. Then
    \begin{equation*}
        \Delta(\varphi / \comm{\psi}) = \Delta_{\varphi, \psi} \ .
    \end{equation*}
\end{proposition}

\chapter{The Araki-Uhlmann Relative Entropy}\label{ch:relativeEntropy}

In the previous two chapters, the theory of operator algebras and modular theory in von Neumann algebras were outlined very roughly. Now, the focus of the presentation shifts towards applications of this theory to the study of the relative entropy functional. In \cref{sec:relativeEntropy_standardForm}, the \emph{Araki-Uhlmann relative entropy} is defined for von Neumann algebras in standard form, and its most important properties are proved. \cref{sec:relativeEntropy_spatialForm} discusses the same functional for arbitrary representations of the von Neumann algebra; both definitions are useful in different situations, and both will be employed in later investigations in this text. The main part of this chapter is \cref{sec:relativeEntropy_monotonicity}, where \emph{Uhlmann's monotonicity theorem} for the relative entropy is proved in great detail, and some of its corollaries are discussed. Finally, \cref{sec:relativeEntropy_vectorMonotonicity} contains some original results: Uhlmann's theorem is applied to the specific situation of two vector functionals in order to find monotonicity inequalities with respect to certain Hilbert-space transformations.

\section{von Neumann Algebras in Standard Form}\label{sec:relativeEntropy_standardForm}

In the following, a generalization of the Umegaki relative entropy, \cf{} \cref{eq:introduction_Umegaki}, to normal functionals on general von Neumann algebras will be given. It will be assumed that the von Neumann algebra is represented in standard form (\cf{} \cref{sec:vonNeumann_standardForm}). In this case, the relative entropy was defined for faithful normal states by \textsc{H. Araki} in 1976 \cite{Araki76} using the relative modular operator (\cf{} \cref{sec:vonNeumann_relativeModularOperator}). Subsequently, in 1977 he generalized his formula to non-faithful states in the seminal work \cite{Araki77}.

\subsection{Definition and Examples}

\begin{definition}[Relative entropy -- standard form]\label{def:relativeEntropy_stdRelativeEntropy}
    Let $(\MFM, \MH, J, \NPC)$ be a von Neumann algebra in standard form, let $\psi, \varphi \in \NF{\MFM}$ be two positive normal functionals, and let $\Psi, \Phi \in \NPC$ be the corresponding vector representatives of these functionals in the natural positive cone (\cref{thm:vonNeumann_vectorRepresentativesNormalStates}). The \bemph{Araki-Uhlmann relative entropy} of $\psi$ and $\varphi$ is defined to be
    \begin{equation}\label{eq:relativeEntropy_stdRelativeEntropy}
        \MS_\MFM^\mathrm{std}(\psi, \varphi) \ce
        \begin{cases}
            - \vbraket[\big]{\Psi, \log (\Delta_{\Phi, \Psi}) \Psi} & \text{if $\ssupp(\psi) \le \ssupp(\varphi)$} \ , \\
            + \infty & \text{otherwise} \ .
        \end{cases}
    \end{equation}
\end{definition}

\begin{remark}\label{rem:relativeEntropy_stdRelativeEntropy}
    The relation $\ssupp(\psi) \le \ssupp(\varphi)$ is equivalent to the condition $\Psi \in \supp(\varphi)$. Indeed, since $\supp(\psi) = \ran\bigl(\ssupp(\psi)\bigr) = [\comm{\MFM} \Psi]$ according to \cref{pro:operatorAlgebras_supportVectorFunctional}, it is clear that the former implies the latter. For the converse implication, let $\xi \in \supp(\psi)$ be arbitrary, \ie{}, $\xi = \dlim{\MH}{n \to \infty} \comm{A_n} \Psi$ for some $(\comm{A_n})_{n \in \N} \subset \comm{\MFM}$. As the assumption $\Psi \in \supp(\varphi)$ implies $\ssupp(\varphi) \Psi = \Psi$ (\cref{lem:operators_rangeProjection}), one also obtains
    \begin{equation*}
        \ssupp(\varphi) \xi = \lim_{n \to \infty} \ssupp(\varphi) \comm{A_n} \Psi = \lim_{n \to \infty} \comm{A_n} \ssupp(\varphi) \Psi = \xi \ ,
    \end{equation*}
    where in the second step $\ssupp(\varphi) \in \MFM$ was used, \cf{} \cref{para:operatorAlgebras_support}. This shows that $\xi \in \ran\bigl(\ssupp(\varphi)\bigr) = \supp(\varphi)$ by the aforementioned lemma, hence $\supp \psi \subset \supp \varphi$.
\end{remark}

The following representation of $\MS_\MFM^\mathrm{std}(\psi, \varphi)$ is given in \cite[Def. 3.1]{Araki77} as a definition.

\begin{lemma}\label{lem:relativeEntropy_spectralRepresentationRE}
    Let $\MFM$, $\psi, \varphi \in \NF{\MFM}$, and $\Psi, \Phi \in \MP$ be as above. Denote by $E_{\Psi, \Phi}$ the unique spectral measure associated with the positive self-adjoint relative modular operator $\Delta_{\Psi, \Phi}$. In case that $\ssupp(\psi) \le \ssupp(\varphi)$, the relative entropy $\MS_\MFM^\mathrm{std}(\psi, \varphi)$ can be written equivalently as
    \begin{equation}\label{eq:relativeEntropy_stdRESpectralIntegral}
        \MS_\MFM^\mathrm{std}(\psi, \varphi) = \int_{0}^{+\infty} \log(\lambda) \diff \braket{\Psi, E_{\Psi, \Phi}(\lambda) \Psi} \ .
    \end{equation}
\end{lemma}

\begin{Proof}
    Recall from \cref{pro:vonNeumann_propertiesRelativeModularOp} \ref{enu:vonNeumann_inverseRelModOp} the relation $\Delta_{\Psi, \Phi}^{-1} = J \Delta_{\Phi, \Psi} J$ which holds with restriction to $\supp(\Delta_{\Psi, \Phi})$, that is, $\Delta_{\Psi, \Phi}^{-1} \cdot \ssupp(\Delta_{\Psi, \Phi}) = J \Delta_{\Phi, \Psi} J \cdot \ssupp(\Delta_{\Psi, \Phi})$, where $\ssupp(\Delta_{\Psi, \Phi}) = \ssupp_\psi \comm{\ssupp}_\varphi$ is the support projection of the relative modular operator according to \cref{pro:vonNeumann_propertiesRelativeModularOp} \ref{enu:vonNeumann_supportRelModOp}. Taking the logarithm of $\Delta_{\Psi, \Phi}^{-1}$ and using properties of the functional calculus, one obtains (see \cref{lem:operators_funcCalcUnitaryConj} and \cite[Rem. 3.4]{Araki77})
    \begin{equation*}
        \log(\Delta_{\Psi, \Phi}) \cdot \ssupp(\Delta_{\Psi, \Phi}) = - J \log(\Delta_{\Phi, \Psi}) J \cdot \ssupp(\Delta_{\Psi, \Phi}) \ .
    \end{equation*}
    
    By the assumption $\ssupp(\psi) \le \ssupp(\varphi)$, it holds that $\Psi \in \supp(\varphi)$ according to \cref{rem:relativeEntropy_stdRelativeEntropy}. Therefore, $\ssupp(\varphi) \Psi = \Psi$ by \cref{lem:operators_rangeProjection}. Using also $\ssupp_\psi \comm{\ssupp}_\varphi = \ssupp(\psi) J \ssupp(\varphi) J$ (\cref{lem:vonNeumann_commMsupport}), it follows that $\ssupp(\Delta_{\Psi, \Phi}) \Psi = \Psi$ because $J \Psi = \Psi$ since $\Psi \in \NPC$, see \cref{pro:vonNeumann_propertiesNPC} \ref{enu:vonNeumann_JonNPC}. Therefore, by definition of the spectral integral in \cref{eq:relativeEntropy_stdRESpectralIntegral} and the identities found above, one computes
    \begin{align*}
        \int\nolimits_{0}^{+\infty} \log(\lambda) \diff \braket{\Psi, E_{\Psi, \Phi}(\lambda) \Psi} &= \vbraket[\big]{\Psi, \log(\Delta_{\Psi, \Phi}) \Psi} \\
        &= \vbraket[\big]{\Psi, \log(\Delta_{\Psi, \Phi}) \cdot \ssupp(\Delta_{\Psi, \Phi}) \Psi} \\[4pt]
        &= - \vbraket[\big]{\Psi, J \log(\Delta_{\Phi, \Psi}) J \cdot \ssupp(\Delta_{\Psi, \Phi}) \Psi} \\[4pt]
        &= - \vbraket[\big]{\Psi, \log(\Delta_{\Phi, \Psi}) \Psi} \\[4pt]
        &= \MS_\MFM^\mathrm{std}(\psi, \varphi) \ . \tag*{\qedhere}
    \end{align*}
\end{Proof}

\begin{remark}\label{rem:relativeEntropy_supportCondition}
    Using the spectral representation \eqref{eq:relativeEntropy_stdRESpectralIntegral}, it becomes apparent why the condition $\ssupp(\psi) \le \ssupp(\varphi)$ has to be imposed in order for $- \vbraket[\big]{\Psi, \log (\Delta_{\Phi, \Psi}) \Psi}$ to be well-defined: it was shown above that $\ssupp(\Delta_{\Psi, \Phi}) \Psi = \Psi$ is a consequence of the support condition, hence $\Psi \in \ker(\Delta_{\Psi, \Phi})^\perp$ (\cref{para:operators_supportPO} and \cref{lem:operators_rangeProjection}), and therefore $E_{\Psi, \Phi}(\set{0}) \, \Psi = 0$ because $E_{\Psi, \Phi}(\set{0})$ is the orthogonal projection onto $\ker(\Delta_{\Psi, \Phi})$ (\cref{pro:operators_eigenvalueSpectralMeasure}). Thus, the integral is well-defined at the lower end since the set $\{\lambda = 0\}$ has $E_{\Psi, \Phi}$-measure zero. If this condition is not satisfied, then $E_{\Psi, \Phi}(0) \Psi \neq 0$ in general, hence there is a divergent contribution to the integral implying $\MS_\MFM^\mathrm{std}(\psi, \varphi) = + \infty$.
\end{remark}

Before proving important properties of the Araki-Uhlmann relative entropy, it shall first be shown that $\MS_\MFM^\mathrm{std}$ is indeed a generalization of the Umegaki and classical relative entropies defined in \cref{sec:introduction_history}.

\begin{examples}\label{exa:relativeEntropy_specialCasesArakiUhlmann}
    \leavevmode
    \begin{enumerate}[env]
        \item Consider the commutative von Neumann algebra $\MFM = L^\infty(X, \mu)$ for a $\sigma$-finite measure space $(X, \Sigma, \mu)$. Recall from \cref{exa:vonNeumann_standardForm} \ref{enu:vonNeumann_exaStandardFormAbelian} that the standard form is given by $\bigl(\MFM, L^2(X, \mu), J, L^2(X, \mu)_+\bigr)$, and from \cref{exa:operatorAlgebras_normalStates} \ref{enu:operatorAlgebras_exaNormalStateAbelian} that $\NF{\MFM} \cong L^1(X, \mu)_+$ such that to every $\psi \in \NF{\MFM}$, one can associate the measure $\diff \nu_\psi = h_\psi \diff \mu$, $h_\psi \in L^1(X, \mu)_+$. Furthermore, \cref{exa:vonNeumann_relativeModularOperator} \ref{enu:vonNeumann_relativeModularOperatorAbelian} showed that for every pair of normal functionals $\psi, \varphi \in \NF{\MFM}$, the relative modular operator $\Delta_{\varphi, \psi}$ is given by the multiplication operator with the Radon-Nikodým derivative of $\nu_\varphi$ with respect to $\nu_\psi$, restricted to the support of $\psi$. Observe that the existence of $\Delta_{\varphi, \psi}$ implies that $\nu_\varphi \ll \nu_\psi$ \cite[p. 122]{Cohn13}, and that in case $\supp(\psi) \subset \supp(\varphi)$, one also has $\nu_\psi \ll \nu_\varphi$, hence it follows that \cite[p. 309]{Elstrodt18}
        \begin{equation*}
            \left(\od{\nu_\varphi}{\nu_\psi}\right)^{-1} = \od{\nu_\psi}{\nu_\varphi} \ .
        \end{equation*}
        
        Let $\Psi \ce \sqrt{h_\psi} \in L^2(X, \mu)_+$ be the vector representative of the functional $\psi$. With the above remarks, the relative entropy of $\psi$ with respect to $\varphi$ in case $\supp(\psi) \subset \supp(\varphi)$ takes the form
        \begin{equation*}
            \MS_\MFM^\mathrm{std}(\psi, \varphi) = - \vbraket[\big]{\Psi, \log(\Delta_{\varphi, \psi}) \Psi} = - \int_X \log\left(\od{\nu_\varphi}{\nu_\psi}\right) \, h_\psi \diff \mu = \int_X \log\left(\od{\nu_\psi}{\nu_\varphi}\right) \diff \nu_\psi = D(\nu_\psi, \nu_\varphi) \ .
        \end{equation*}
        That is, the Araki-Uhlmann relative entropy of two positive normal functionals on $\MFM = L^\infty(X, \mu)$ reduces to the Kullback-Leibler divergence from \cref{eq:introduction_KL}: $D(\mu, \nu) = \MS_{L^\infty}^\mathrm{std}(\omega_\mu, \omega_\nu)$.
        
        \item \label{enu:relativeEntropy_ArakiUhlmannTypeI} Next, consider the von Neumann algebra $\MFM = \BO(\MH)$ over a separable Hilbert space $\MH$. According to \cref{exa:vonNeumann_standardForm} \ref{enu:vonNeumann_exaStandardFormTypeI}, the standard form of $\MFM$ is given by $\bigl(\MFM, \HS(\MH), J, \HS(\MH)_+\bigr)$. For all positive normal functionals $\psi, \varphi \in \NF{\MFM}$, there exist positive trace-class operators $\rho_\psi, \rho_\varphi \in \NO(\MH)_+$ such that $\Psi = \sqrt{\rho_{\psi}}$ and $\Phi = \sqrt{\rho_\varphi}$ are their vector representatives in $\HS(\MH)_+$.
        
        In \cref{exa:vonNeumann_relativeModularOperator} \ref{enu:vonNeumann_relativeModularOperatorTypeI}, it was shown that if $\rho_\varphi = \sum_{i \in \N} \lambda_i \, P_i$ and $\rho_\psi = \sum_{j \in \N} \mu_j \, Q_j$ are the spectral decompositions of $\rho_\psi$ and $\rho_\varphi$, then the relative modular operator $\Delta_{\varphi, \psi}$ on $\HS(\MH)$ is given by the following expression:
        \begin{equation*}
            \Delta_{\varphi, \psi} = L_{\rho_\varphi} R_{\rho_\psi^{-1}} = \sum_{i,j=1}^{\infty} \lambda_i \inv{\mu_j} L_{P_i} R_{Q_j} \ .
        \end{equation*}
        Here, $L_T$ and $R_T$ denote the left and right multiplication operators with the operator $T$, respectively, and $\inv{\rho_\psi}$ is the generalized inverse. From the functional calculus, one obtains \cite[Exa. 5.1]{Schmüdgen12}
        \begin{align*}
            \log(\Delta_{\varphi, \psi}) &= \sum_{i,j=1}^{\infty} \log(\lambda_i \inv{\mu_j}) L_{P_i} R_{Q_j} = \sum_{i,j=1}^{\infty} \log(\lambda_i) L_{P_i} R_{Q_j} - \sum_{i,j=1}^{\infty} \log(\mu_j) L_{P_i} R_{Q_j} \\
            &= \sum_{i=1}^{\infty} \log(\lambda_i) L_{P_i} - \sum_{j=1}^{\infty} \log(\mu_j) R_{Q_j} = \log(\rho_\varphi) - \log(\rho_\psi) \ ,
        \end{align*}
        where it was used that $\sum_{j \in \N} R_{Q_j} = \id_{\HS(\MH)}$ and $\sum_{i \in \N} L_{P_i} = \id_{\HS(\MH)}$. With this result, the following expression for the relative entropy $\MS_\MFM^\mathrm{std}(\psi, \varphi)$ in case that $\supp(\rho_\psi) \subset \supp(\rho_\varphi)$ is found:
        \begin{align*}
            \MS_\MFM^\mathrm{std}(\psi, \varphi) &= - \vbraket[\big]{\Psi, \log(\Delta_{\Phi, \Psi}) \Psi} = - \tr\bigl(\rho_\psi^{1/2} \bigl(\log(\rho_\varphi) - \log(\rho_\psi)\bigr) \rho_\psi^{1/2}\bigr) \\
            &= \tr\bigl(\rho_\psi \log \rho_\psi - \rho_\psi \log \rho_\varphi\bigr) = S(\rho_\psi, \rho_\varphi) \ .
        \end{align*}
        That is, the Araki-Uhlmann relative entropy of two positive normal functionals on $\MFM = \BO(\MH)$ reduces to the Umegaki relative entropy from \cref{eq:introduction_Umegaki}: $S(\rho, \sigma) = \MS_{\BO(\MH)}^\mathrm{std}(\omega_\rho, \omega_\sigma)$.\qedhere
    \end{enumerate}
\end{examples}

\subsection{Fundamental Properties}

The subsequent properties of the relative entropy (\cref{pro:relativeEntropy_basicLowerBound,pro:relativeEntropy_independenceVecRep,pro:relativeEntropy_nonNegStdRelativeEntropy,pro:relativeEntropy_scaling}) are stated in \textsc{Araki}'s paper \cite{Araki77}, and the proofs given here are more detailed versions of the proofs given in that reference.

\begin{proposition}\label{pro:relativeEntropy_basicLowerBound}
    The relative entropy $\MS_\MFM^\mathrm{std}(\psi, \varphi)$ is well-defined. It takes a finite real value or $+ \infty$, and it satisfies the inequality
    \begin{equation}\label{eq:relativeEntropy_basicLowerBound}
        \MS_\MFM^\mathrm{std}(\psi, \varphi) \ge - \psi(\id_\MH) \log\left(\frac{\varphi\bigl(\ssupp(\psi)\bigr)}{\psi(\id_\MH)}\right) \ .
    \end{equation}
\end{proposition}

\begin{Proof}
    If $\ssupp(\psi) > \ssupp(\varphi)$, it holds that $\MS_\MFM^\mathrm{std}(\psi, \varphi) = + \infty$ by definition, and hence there is nothing to show. Thus, assume that $\ssupp(\psi) \le \ssupp(\varphi)$. The relative Tomita operator $S_{\Psi, \Phi}$ (\cref{para:vonNeumann_relativeTomita}) then satisfies $S_{\Psi, \Phi} \Phi = \ssupp(\varphi) \Psi = \Psi$ because $\Psi \in \supp(\varphi)$ by \cref{rem:relativeEntropy_stdRelativeEntropy}. Since also $J \Psi = \Psi$ by \cref{pro:vonNeumann_propertiesNPC} \ref{enu:vonNeumann_JonNPC}, it follows from \cref{pro:vonNeumann_propertiesRelativeModularOp} \ref{enu:vonNeumann_polarDecompRelTomita} that
    \begin{equation*}
        \Delta_{\Psi, \Phi}^{1/2} \, \Phi = J S_{\Psi, \Phi} \, \Phi = \Psi \ .
    \end{equation*}
    Next, recall that $\ssupp(\Delta_{\Psi, \Phi})$ is the orthogonal projection onto $\ker(\Delta_{\Psi, \Phi})^\perp$ (\cref{para:operators_supportPO}), and $E_{\Psi, \Phi}(\{0\})$ is the projection onto $\ker(\Delta_{\Psi, \Phi})$ (\cref{pro:operators_eigenvalueSpectralMeasure}). Therefore, $\ssupp(\Delta_{\Psi, \Phi}) = \id_\MH - E_{\Psi, \Phi}(\{0\}) = \id_\MH - \1_{\set{0}}(\Delta_{\Psi, \Phi})$ by \cref{pro:operators_functionalCalculus} \ref{enu:operators_spectralMeasureIndicatorFct}. Let $f : \R \longto \C \cup \set{+ \infty}$ be an $E_{\Psi, \Phi}$-almost everywhere finite Borel-measurable function. Since $\1_{\set{0}}$ is bounded, \cref{pro:operators_functionalCalculus} \ref{enu:operators_funcCalcCommutation} implies
    \begin{align*}
        f(\Delta_{\Psi, \Phi}) \, \ssupp(\Delta_{\Psi, \Phi}) &= f(\Delta_{\Psi, \Phi}) \bigl(\id_\MH - \1_{\set{0}}(\Delta_{\Psi, \Phi})\bigr) \\
        &= \bigl(\id_\MH - \1_{\set{0}}(\Delta_{\Psi, \Phi})\bigr) f\bigl(\Delta_{\Psi, \Phi}\bigr) = \ssupp(\Delta_{\Psi, \Phi}) f\bigl(\Delta_{\Psi, \Phi}\bigr) \ .
    \end{align*}
    
    With the above two observations, the fact that $\ssupp(\Delta_{\Psi, \Phi}) \Psi = \Psi$ (which was shown in the proof of \cref{lem:relativeEntropy_spectralRepresentationRE}), and \cref{pro:operators_functionalCalculus} \ref{enu:operators_funcCalcIP}, one can compute the following integral:
    \begin{align}\label{eq:relativeEntropy_auxiliaryIntegralBound}
        \begin{split}
            \int\nolimits_{0}^{+\infty} \inv{\lambda} \diff \braket{\Psi, E_{\Psi, \Phi}(\lambda) \Psi} &= \vbraket[\big]{\Psi, \Delta_{\Psi, \Phi}^{-1} \cdot \ssupp(\Delta_{\Psi, \Phi}) \Psi} = \braket{\Delta_{\Psi, \Phi}^{-1/2} \, \Psi, \ssupp(\Delta_{\Psi, \Phi}) \Delta_{\Psi, \Phi}^{-1/2} \, \Psi} \\
            &= \braket{\Phi, \ssupp(\Delta_{\Psi, \Phi}) \Phi} = \braket{\Phi, \ssupp_\psi \comm{\ssupp}_\varphi \Phi} = \braket{\Phi, \ssupp(\psi) \Phi} = \varphi\bigl(\ssupp(\psi)\bigr) \\[4pt]
            &\le \varphi\bigl(\ssupp(\varphi)\bigr) = \varphi(\id_\MH) \ .
        \end{split}
    \end{align}
    Using that $\log(\lambda) \le \lambda$ for all $\lambda > 0$, it follows that $- \log(\lambda) = \log(\inv{\lambda}) \le \inv{\lambda}$ and hence, using the spectral representation of $\MS_\MFM^\mathrm{std}(\psi, \varphi)$ from \cref{lem:relativeEntropy_spectralRepresentationRE}, one obtains the following bound:
    \begin{equation*}
        \MS_\MFM^\mathrm{std}(\psi, \varphi) = \int_{0}^{+\infty} \log(\lambda) \diff \braket{\Psi, E_{\Psi, \Phi}(\lambda) \Psi} \ge - \int_{0}^{+\infty} \inv{\lambda} \diff \braket{\Psi, E_{\Psi, \Phi}(\lambda) \Psi} \ge - \varphi(\id_\MH) \ .
    \end{equation*}
    This shows that the relative entropy $\MS_\MFM^\mathrm{std}(\psi, \varphi)$ is bounded from below, hence it is well-defined, and it takes real values or the value $+ \infty$.

    It remains to prove \cref{eq:relativeEntropy_basicLowerBound}. Assume still $\ssupp(\psi) \le \ssupp(\varphi)$ so that $\comm{\ssupp}_\varphi \Psi = \Psi$. Then it follows from the functional calculus as before that
    \begin{align*}
        \int_{0}^{+\infty} \diff \braket{\Psi, E_{\Psi, \Phi}(\lambda) \Psi} &= \braket{\Psi, \1_{(0, + \infty)}(\Delta_{\Psi, \Phi}) \Psi} \\
        &= \braket{\Psi, \ssupp(\Delta_{\Psi, \Phi}) \Psi} = \braket{\Psi, \ssupp_\psi \comm{\ssupp}_\varphi \Psi} = \psi\bigl(\ssupp(\psi)\bigr) = \psi(\id_\MH) \ .
    \end{align*}
    This implies that $\diff \braket{\Psi, E_{\Psi, \Phi} \Psi} / \psi(\id_\MH)$ is a probability measure on $(0, + \infty)$. Using that the logarithm is concave, one can apply Jensen's inequality \cite[Thm. VI.1.3]{Elstrodt18} and \cref{eq:relativeEntropy_auxiliaryIntegralBound} to find
    \begin{align*}
        \MS_\MFM^\mathrm{std}(\psi, \varphi) &= - \psi(\id_\MH) \int\nolimits_{0}^{+\infty} \log(\inv{\lambda}) \, \frac{1}{\psi(\id_\MH)} \diff \braket{\Psi, E_{\Psi, \Phi}(\lambda) \Psi} \\
        &\ge - \psi(\id_\MH) \log\biggl(\,\int\nolimits_{0}^{+\infty} \inv{\lambda} \, \frac{1}{\psi(\id_\MH)} \diff \braket{\Psi, E_{\Psi, \Phi}(\lambda) \Psi}\biggr) \\
        &\ge - \psi(\id_\MH) \log\left(\frac{\varphi\bigl(\ssupp(\psi)\bigr)}{\psi(\id_\MH)}\right) \ . \tag*{\qedhere}
    \end{align*}
\end{Proof}

\begin{proposition}\label{pro:relativeEntropy_independenceVecRep}
    The definition of the relative entropy $\MS_\MFM^\mathrm{std}(\psi, \varphi)$ does not depend on the choice of the natural positive cone $\NPC$ and the vector representatives $\Psi, \Phi \in \NPC$ of $\psi$ and $\varphi$.
\end{proposition}

\begin{Proof}
    Let $(\MFM, \MH, J^\prime, \NPC^\prime)$ be another standard form of the von Neumann algebra $\MFM$. By \cref{cor:vonNeumann_uniquenessStandardForm}, there exists a unique unitary operator $U \in \UO(\MH)$ such that (i) $U A U^\ast = A$ for all $A \in \MFM$ and (ii) $\NPC^\prime = U \NPC$. If $\Psi, \Phi \in \NPC$ are the vector representatives of $\psi, \varphi \in \NF{\MFM}$ in the cone $\NPC$, then $\Psi^\prime \ce U \Psi$ and $\Phi^\prime \ce U \Phi$ are the vector representatives in the cone $\NPC^\prime$ by virtue of property (ii). According to \cref{lem:vonNeumann_relativeModularOpUnitary}, it holds that $\Delta_{U \Phi, U \Psi} = U \Delta_{\Phi, \Psi} U^\ast$. With this relation, it follows that
    \begin{align*}
        - \vbraket[\big]{\Psi^\prime, \log (\Delta_{\Phi^\prime, \Psi^\prime}) \Psi^\prime} &= - \vbraket[\big]{U \Psi, \log (U \Delta_{\Phi, \Psi} U^\ast) U \Psi} \\
        &= - \vbraket[\big]{U \Psi, U \log (\Delta_{\Phi, \Psi}) U^\ast U \Psi} \\
        &= - \vbraket[\big]{\Psi, \log (\Delta_{\Phi, \Psi}) \Psi} \\
        &= \MS_\MFM^\mathrm{std}(\psi, \varphi) \ .
    \end{align*}
    In the second step, \cref{lem:operators_funcCalcUnitaryConj} was used. This shows that the relative entropy takes the same value for both the representatives $\Psi, \Phi \in \NPC$ and $\Psi^\prime, \Phi^\prime \in \NPC^\prime$ of the functionals $\psi, \varphi$.
\end{Proof}

\begin{proposition}\label{pro:relativeEntropy_nonNegStdRelativeEntropy}
    If $\psi(\id_\MH) = \varphi(\id_\MH) > 0$, then $\MS_\MFM^\mathrm{std}(\psi, \varphi) \ge 0$. Equality $\MS_\MFM^\mathrm{std}(\psi, \varphi) = 0$ holds true if and only if $\psi = \varphi$.
\end{proposition}

\begin{Proof}
    Using that $\varphi\bigl(\ssupp(\psi)\bigr) \le \varphi(\id_\MH)$ together with the inequality \eqref{eq:relativeEntropy_basicLowerBound} and the assumption $\psi(\id_\MH) = \varphi(\id_\MH)$, one obtains non-negativity of the relative entropy:
    \begin{equation*}
        \MS_\MFM^\mathrm{std}(\psi, \varphi) \ge - \psi(\id_\MH) \log\left(\frac{\varphi\bigl(\ssupp(\psi)\bigr)}{\psi(\id_\MH)}\right) \ge - \psi(\id_\MH) \log\left(\frac{\varphi(\id_\MH)}{\psi(\id_\MH)}\right) = 0 \ .
    \end{equation*}
    
    Regarding the case of equality, first assume that $\psi = \varphi$. Then $\Psi = \Phi$ and $\ssupp(\psi) = \ssupp(\varphi)$. The relative modular operator satisfies $\Delta_{\Psi, \Psi} \Psi = S_\Psi^\ast S_\Psi \Psi = \Psi$, that is, $\Psi \in \Eig(\Delta_\Psi, 1)$. From this and \cref{lem:operators_funcCalcEigenvalue}, it follows that $\MS_\MFM^\mathrm{std}(\psi, \varphi) = - \braket{\Psi, \log(\Delta_\Psi) \Psi} = - \braket{\Psi, \log(1) \Psi} = 0$. Conversely, if $\MS_\MFM^\mathrm{std}(\psi, \varphi) = 0$, then necessarily $\ssupp(\psi) \le \ssupp(\varphi)$ and $\varphi\bigl(\ssupp(\psi)\bigr) = \varphi(\id_\MH)$. (If $\varphi\bigl(\ssupp(\psi)\bigr) < \varphi(\id_\MH)$ were to hold true, then \cref{eq:relativeEntropy_basicLowerBound} would imply that $\MS_\MFM^\mathrm{std}(\psi, \varphi) > 0$, see above.) From the definition of the support projection, it follows that $\ssupp(\psi) = \ssupp(\varphi)$. In this case, $\MS_\MFM^\mathrm{std}(\psi, \varphi)$ coincides with the relative entropy $\MS_\MFN^\mathrm{std}(\psi, \varphi)$ on the algebra $\MFN = \ssupp_\psi \MFM \ssupp_\psi$ on which both $\psi$ and $\varphi$ are faithful \cite[Rem. 3.5]{Araki77}. Since the relative entropy between faithful states is strictly positive \cite[Eq. (1.3)]{Araki76}, it follows that $\psi = \varphi$.
\end{Proof}

\begin{proposition}\label{pro:relativeEntropy_scaling}
    For $\lambda, \mu > 0$, there holds
    \begin{equation*}
        \MS_\MFM^\mathrm{std}(\lambda \psi, \mu \varphi) = \lambda \MS_\MFM^\mathrm{std}(\psi, \varphi) - \lambda \psi(\id_\MH) \log\left(\frac{\mu}{\lambda}\right) \ .
    \end{equation*}
\end{proposition}

\begin{Proof}
    A scaling of the functional $\psi = \omega_\Psi$ by $\lambda$ induces a scaling of the vector representative $\Psi$ by $\lambda^{1/2}$. Therefore, using the identity $\Delta_{\mu \varphi, \lambda \psi} = \frac{\mu}{\lambda} \, \Delta_{\varphi, \psi}$ from \cref{lem:vonNeumann_relativeModularOpScaling}, it follows that
    \begin{align*}
        \MS_\MFM^\mathrm{std}(\lambda \psi, \mu \varphi) &= - \braket{\lambda^{1/2} \Psi, \log(\Delta_{\mu \varphi, \lambda \psi}) \lambda^{1/2} \Psi} \\
        &= - \lambda \, \vbraket[\bigg]{\Psi, \log\left(\frac{\mu}{\lambda} \, \Delta_{\Phi, \Psi}\right) \Psi} \\[3pt]
        &= - \lambda \, \braket{\Psi, \log(\Delta_{\Phi, \Psi}) \Psi} - \lambda \, \vbraket[\bigg]{\Psi, \log\left(\frac{\mu}{\lambda}\right) \Psi} \\
        &= \lambda \MS_\MFM^\mathrm{std}(\psi, \varphi) - \lambda \psi(\id_\MH) \log\left(\frac{\mu}{\lambda}\right) \ . \tag*{\qedhere}
    \end{align*}
\end{Proof}

The next inequality will be used in \cref{ch:perturbationTheory} in the perturbation theory of KMS-states. The proof presented here is an adaptation of \cite[Thm. 4.3 (4)]{DJP03}.

\begin{lemma}\label{lem:relativeEntropy_inequalityCenter}
    Let $A \in \SE{\MFM}$ be a self-adjoint element in the center $\MFM \cap \comm{\MFM}$, and assume that $\psi(\id_\MH) = 1$. Then the following inequality holds true:
    \begin{equation}\label{eq:relativeEntropy_inequalityCenter}
        \MS_\MFM^\mathrm{std}(\psi, \varphi) \ge \psi(A) - \log\bigl(\varphi(\ee^{A})\bigr) \ .
    \end{equation}
\end{lemma}

\begin{Proof}
    One may assume that $\ssupp_\psi \le \ssupp_\varphi$ since otherwise there is nothing to show. By the assumption on $A$ and \cref{lem:vonNeumann_relModOpCommutesCenter}, it follows that the operators $\ee^A$ and $\Delta_{\Phi, \Psi}$ commute with each other. Therefore, using the functional calculus for commuting operators \cite[Sect. 5.5]{Schmüdgen12}, one obtains
    \begin{equation*}
        \log(\Delta_{\Phi, \Psi}) + A - \log\bigl(\varphi(\ee^A \, \ssupp_\psi)\bigr) = \log(\Delta_{\Phi, \Psi}) + \log\bigl(\ee^A / \varphi(\ee^A \, \ssupp_\psi)\bigr) = \log\bigl(\Delta_{\Phi, \Psi} \, \ee^A / \varphi(\ee^A \, \ssupp_\psi)\bigr) \ .
    \end{equation*}
    From the well-known inequality $\log(x) \le x - 1$ for positive $x > 0$ (see \cref{ftn:perturabtionTheory_exponentialInequality} in \cref{ch:relativeEntropy} on \cpageref{ftn:perturabtionTheory_exponentialInequality} for a proof), it follows that
    \begin{equation*}
        \log\bigl(\Delta_{\Phi, \Psi} \, \ee^A / \varphi(\ee^A \, \ssupp_\psi)\bigr) \le \Delta_{\Phi, \Psi} \, \ee^A / \varphi(\ee^A \, \ssupp_\psi) - \id_\MH \ .
    \end{equation*}
    Thus, combining the definition \eqref{eq:relativeEntropy_stdRelativeEntropy} of $\MS_\MFM^\mathrm{std}$ and the previous two identities, one finds that
    \begin{align*}
        \MS_\MFM^\mathrm{std}(\psi, \varphi) - \psi(A) + \log\bigl(\varphi(\ee^{A} \, \ssupp_\psi)\bigr) &= - \vbraket[\big]{\Psi, \bigl(\log(\Delta_{\Phi, \Psi}) + A - \log \varphi(\ee^{A} \, \ssupp_\psi)\bigr) \Psi} \\
        &= - \vbraket[\big]{\Psi, \log\bigl(\Delta_{\Phi, \Psi} \, \ee^A / \varphi(\ee^A \, \ssupp_\psi)\bigr) \Psi} \\
        &\ge - \vbraket[\big]{\Psi, \bigl(\Delta_{\Phi, \Psi} \, \ee^A / \varphi(\ee^A \, \ssupp_\psi) - \id_\MH\bigr) \Psi} \\
        &= - \vnorm[\big]{\Delta_{\Phi, \Psi}^{1/2} \, \ee^{A/2} \, \Psi}^2 / \varphi(\ee^A \, \ssupp_\psi) + 1 \ .
    \end{align*}
    Using the relation $S_{\Phi, \Psi} = J \Delta_{\Phi, \Psi}^{1/2}$ (\cref{pro:vonNeumann_propertiesRelativeModularOp} \ref{enu:vonNeumann_polarDecompRelTomita}) and the definition \eqref{eq:vonNeumann_relativeTomitaOperator} of the relative Tomita operator $S_{\Phi, \Psi}$, one obtains
    \begin{equation*}
        \Delta_{\Phi, \Psi}^{1/2} \, \ee^{A/2} \, \Psi = J S_{\Phi, \Psi} \, \ee^{A/2} \, \Psi = J \, \ssupp_\psi \, \ee^{A/2} \, \Phi = J \, \ee^{A/2} \, \ssupp_\psi \Phi \ .
    \end{equation*}
    (Note that $\ee^{A / 2} \in \comm{\MFM}$ since $A \in \comm{\MFM}$ by assumption.) The norm of the this expression is given by $\norm{J \, \ee^{A/2} \, \ssupp_\psi \Phi}^2 = \norm{\ee^{A/2} \, \ssupp_\psi \Phi}^2 = \varphi(\ee^{A} \, \ssupp_\psi)$ because $J$ is anti-unitary.
    Inserting this result in the above inequality, the assertion follows:
    \begin{align*}
        \MS_\MFM^\mathrm{std}(\psi, \varphi) - \psi(A) + \log\bigl(\varphi(\ee^{A})\bigr) &= \MS_\MFM^\mathrm{std}(\psi, \varphi) - \psi(A) + \log\bigl(\varphi(\ee^{A} \, \ssupp_\varphi)\bigr) \\
        &\ge \MS_\MFM^\mathrm{std}(\psi, \varphi) - \psi(A) + \log\bigl(\varphi(\ee^{A} \, \ssupp_\psi)\bigr) \\
        &\ge - \norm{J \, \ee^{A/2} \, \ssupp_\psi \Phi}^2 / \varphi(\ee^A \, \ssupp_\psi) + 1 = 0 \ . \tag*{\qedhere}
    \end{align*}
\end{Proof}

\section{Spatial Form of the Relative Entropy}\label{sec:relativeEntropy_spatialForm}

Having introduced the relative entropy based on the relative modular operator, and having proved several properties of this functional, now, another definition for a relative entropy functional will be given which relies on the spatial derivative operator; it is taken from \cite[Eq. (5.1)]{OP04}.

\subsection{Definition and Relation to the Standard Form}

\begin{definition}[Relative entropy -- spatial form]\label{def:relativeEntropy_spaRelativeEntropy}
    Let $\MFM \subset \BO(\MH)$ be a von Neumann algebra acting on a Hilbert space $\MH$, let $\varphi \in \NS(\MFM)$ be an arbitrary normal state on $\MFM$, and let $\omega = \omega_\Omega \in \NS(\MFM)$ be a vector state on the algebra $\MFM$ induced by a cyclic vector $\Omega \in \MH$ for the algebra $\MFM$. The \bemph{Araki-Uhlmann relative entropy} $\MS_\MFM^\mathrm{spa}(\omega, \varphi)$ of the state $\omega$ with respect to the state $\varphi$ is defined to be
    \begin{equation}\label{eq:relativeEntropy_spaRelativeEntropy}
        \MS_\MFM^\mathrm{spa}(\omega, \varphi) \ce
        \begin{cases}
            - \vbraket[\big]{\Omega, \log \Delta(\varphi / \comm{\omega}) \Omega} & \text{if $\Omega \in \supp(\varphi)$} \ ,\\
            + \infty & \text{otherwise} \ .
        \end{cases}
    \end{equation}
    Here, $\comm{\omega} = \comm{\omega_\Omega}$ denotes the vector state induced by $\Omega$ on the commutant $\comm{\MFM}$, and $\Delta(\varphi / \comm{\omega})$ is the spatial derivative of $\varphi$ with respect to $\comm{\omega}$ which was introduced in \cref{para:vonNeumann_spatialDerivative}.
\end{definition}

\begin{remark}\label{rem:relativeEntropy_spaRelativeEntropy}
    Since $\Omega \in \MH$ is assumed to be cyclic for $\MFM$, it is separating for the commutant $\comm{\MFM}$ by \cref{pro:operatorAlgebras_cyclicSeparating} \ref{enu:operatorAlgebras_cyclicMseparatingCommM}. Therefore, it follows from \cref{cor:operatorAlgebras_faithfulSeparating} that the vector functional $\comm{\omega}$ is faithful on $\comm{\MFM}$. Hence, the spatial derivative $\Delta(\varphi / \comm{\omega})$ is well-defined. Similarly to \cref{rem:relativeEntropy_stdRelativeEntropy}, one shows that the condition $\Omega \in \supp(\varphi)$ is equivalent to $\ssupp(\omega) \le \ssupp(\varphi)$. Note that it is not clear yet whether the definition of $\MS_\MFM^\mathrm{spa}$ is independent of the vector representative $\Omega$ of the functionals $\omega$ and $\comm{\omega}$; this will be proved later with the help of \cref{thm:relativeEntropy_UhlmannMonotonicity}.
\end{remark}

The following proposition answers the obvious question regarding the relationship between the functionals $\MS_\MFM^\mathrm{spa}$ and $\MS_\MFM^\mathrm{std}$ on a von Neumann algebra $\MFM$.

\begin{proposition}\label{pro:relativeEntropy_equivalenceDefinitions}
    Let $\MFM$ be a von Neumann algebra in standard form $(\MFM, \MH, J, \NPC)$ and $\psi, \varphi \in \NF{\MFM}$ with vector representatives $\Psi, \Phi \in \NPC$. Assume that $\psi$ is faithful on $\MFM$, and define a normal functional $\comm{\psi} \in \NF{(\comm{\MFM})}$ by $\comm{\psi}(\comm{A}) \ce \psi(J (\comm{A})^\ast J)$ for all $\comm{A} \in \comm{\MFM}$. Then
    \begin{equation*}
        \MS_\MFM^\mathrm{std}(\psi, \varphi) = \MS_\MFM^\mathrm{spa}(\psi, \varphi) \ .
    \end{equation*}
\end{proposition}

\begin{Proof}
    First, note that the functional $\comm{\psi}$ is given by $\comm{\psi} = \comm{\omega}_{\Psi}$ on the commutant $\comm{\MFM}$:
    \begin{equation*}
        \comm{\psi}(\comm{A}) = \braket{\Psi, J (\comm{A})^\ast J \Psi} = \braket{(\comm{A})^\ast J \Psi, J \Psi} = \braket{\Psi, \comm{A} \Psi} = \comm{\omega_\Psi}(\comm{A}) \ .
    \end{equation*}
    Next, observe that since $\psi$ is assumed to be faithful, the vector representative $\Psi \in \NPC$ is cyclic and separating for $\MFM$ (\cref{rem:vonNeumann_cyclicSeparatingFaithful}). Therefore, $\Psi$ is also cyclic and separating for the commutant $\comm{\MFM}$ by \cref{pro:operatorAlgebras_cyclicSeparating}, hence $\comm{\psi} = \comm{\omega}_\Psi$ is faithful. Thus, the spatial derivative operator $\Delta(\varphi / \comm{\psi})$ is well-defined, and the claim immediately follows from \cref{pro:vonNeumann_spatialDerivativeRelativeModOp}.
\end{Proof}

\subsection{Uhlmann's Representation}

The following representation of the relative entropy, which was first given by \textsc{A. Uhlmann} in 1977 as a \emph{definition} of the relative entropy \cite[Eq. (45)]{Uhlmann77}, is very useful for proving certain properties of this functional; in particular, it will be used in the next section to prove \textsc{Uhlmann}'s monotonicity theorem. The following proof is a much more detailed version of the argument given in \cite[p. 80]{OP04}, and it corrects a mistake in this reference which claims that the convergence of \cref{eq:relativeEntropy_logLimit} is increasing for $\lambda > 1$.

\begin{ntheorem}[Uhlmann's representation of the relative entropy]\label{thm:relativeEntropy_UhlmannRepresentation}
    In the setting of \cref{def:relativeEntropy_spaRelativeEntropy}, the relative entropy $\MS_\MFM^\mathrm{spa}(\omega, \varphi)$ can be written as follows:
    \begin{equation}\label{eq:relativeEntropy_UhlmannRepresentation}
        \MS_\MFM^\mathrm{spa}(\omega, \varphi) = - \lim_{t \searrow 0} \frac{1}{t} \Bigl(\vnorm[\big]{\Delta(\varphi / \comm{\omega})^{t/2} \Omega}^2 - \norm{\Omega}^2\Bigr) \ .
    \end{equation}
\end{ntheorem}

\begin{Proof}
    \pp{1.} Assume first that $\Omega \in \supp(\varphi)$, and denote by $E_{\varphi, \comm{\omega}}$ the unique spectral measure associated with the operator $\Delta(\varphi / \comm{\omega})$. From \cref{pro:operators_functionalCalculus} \ref{enu:operators_funcCalcNorm}, it follows that
    \begin{align*}
        \vnorm[\big]{\Delta(\varphi / \comm{\omega})^{t/2} \Omega}^2 = \int_{0}^{+\infty} \lambda^t \diff \braket{\Omega, E_{\varphi, \comm{\omega}}(\lambda) \Omega} \tand \norm{\Omega}^2 = \int_{0}^{+\infty} \1_{[0, + \infty)} \diff \braket{\Omega, E_{\varphi, \comm{\omega}}(\lambda) \Omega} \ .
    \end{align*}

    \begin{nstatement}
        Let $\lambda \in (0, + \infty)$ be a positive number. Then the expression $\frac{\lambda^t - 1}{t}$, $t > 0$, converges monotonically decreasingly towards $\log(\lambda)$ as $t \searrow 0$:
        \begin{equation}\label{eq:relativeEntropy_logLimit}
            \log(\lambda) = \lim_{t \searrow 0} \frac{\lambda^t - 1}{t} \ .
        \end{equation}
    \end{nstatement}
    
    \begin{subproof}[of statement]
        It will be shown that (i) the limit of $\frac{\lambda^t - 1}{t}$ is actually given by $\log(\lambda)$, and (ii) that the convergence of this function to $\log(\lambda)$ is monotonically decreasing for all $\lambda \in (0, + \infty)$.
    
        (i) For every positive $\lambda > 0$, it holds that
        \begin{align*}
            \lim_{t \searrow 0} \frac{\lambda^t - 1}{t} = \lim_{t \searrow 0} \frac{\ee^{t \log(\lambda)} - \ee^{0 \cdot \log(\lambda)}}{t} = \odbound{}{t}{t=0} \, \ee^{t \log(\lambda)} = \log(\lambda) \ .
        \end{align*}
    
        (ii) For $\lambda \in (0, + \infty)$, define a function $f_\lambda : (0, + \infty) \longto \R$, $t \longmto f_\lambda(t)$, by
        \begin{equation*}
            f_{\lambda}(t) \ce
            \begin{cases}
                \displaystyle \int\nolimits_1^\lambda \xi^{t-1} \diff \xi & \text{if $\lambda \ge 1$} \ , \\[10pt]
                \displaystyle - \int\nolimits_\lambda^1 \xi^{t-1} \diff \xi & \text{if $\lambda < 1$} \ .
            \end{cases}
        \end{equation*}
        Observe that since $\int_1^\lambda \xi^{t-1} \diff \xi = \frac{1}{t} \, \xi^{t} \big|_1^\lambda = \frac{\lambda^t - 1}{t}$, it follows that $f_\lambda(t) = \frac{\lambda^t - 1}{t}$ for all $\lambda, t > 0$. Let $t_1, t_2 \in (0, + \infty)$ be arbitrary numbers such that $t_2 \ge t_1 > 0$. On the one hand, if $\lambda > 1$, then for all $\xi \in [1, \lambda]$ there holds $\log(\xi) \ge 0$, and hence
        \begin{equation*}
            \xi^{t_1 - 1} = \ee^{(t_1 - 1) \log(\xi)} = \ee^{t_1 \log(\xi)} \, \ee^{- \log(\xi)} \le \ee^{t_2 \log(\xi)} \, \ee^{- \log(\xi)} = \ee^{(t_2 - 1) \log(\xi)} = \xi^{t_2 - 1} \ .
        \end{equation*}
        On the other hand, if $0 < \lambda < 1$, then $\log(\xi) < 0$ for all $\xi \in [\lambda, 1]$; therefore
        \begin{equation*}
            - \xi^{t_1 - 1} = - \ee^{(t_1 - 1) \log(\xi)} = - \ee^{- t_1 \abs{\log(\xi)}} \, \ee^{- \log(\xi)} \le - \ee^{- t_2 \abs{\log(\xi)}} \, \ee^{- \log(\xi)} = - \ee^{(t_2 - 1) \log(\xi)} = - \xi^{t_2 - 1} \ .
        \end{equation*}
        These two inequalities show that for all $\lambda \in (0, + \infty)$, it holds that
        \begin{equation*}
            f_\lambda(t_1) \le f_\lambda(t_2) \ ,
        \end{equation*}
        hence $f_\lambda$ is a monotonically increasing function on $(0, + \infty)$. Therefore, since $ \lim_{t \searrow 0} f_\lambda(t) = \log(\lambda)$ by (i), it follows that $\frac{\lambda^t - 1}{t}$ converges monotonically decreasingly to $\log(\lambda)$ as $t \searrow 0$.
    \end{subproof}

    The derivation of \cref{eq:relativeEntropy_UhlmannRepresentation} may now be completed. Define a sequence $(t_n)_{n \in \N} \subset (0, 1]$ by setting $t_n \ce \frac{1}{n}$ for all $n \in \N$; this yields a monotonically decreasing sequence of positive numbers converging towards zero, that is, $t_{n + 1} \le t_n$ for all $n \in \N$ and $t_n \to 0$ as $n \to + \infty$; note also that $t_1 = 1$. Next, define for every $n \in \N$ a measurable function $f_n : (0, + \infty) \longto \R$ by
    \begin{equation*}
        f_n(\lambda) \ce f_{\lambda}(t_n) = \frac{\lambda^{t_n} - 1}{t_n} \ ,
    \end{equation*}
    where $f_\lambda$ was defined in part (ii) of the proof of the statement. From the properties of $f_\lambda$ shown there, it follows for all $\lambda \in (0, + \infty)$ that $f_{n+1}(\lambda) = f_\lambda(t_{n+1}) \le f_\lambda(t_n) = f_n(\lambda)$ and $f_n(\lambda) \to \log(\lambda)$ as $n \to + \infty$, that is, $(f_n)_{n \in \N}$ is a monotonically decreasing sequence converging pointwise to $\log(\lambda)$. Finally,
    \begin{equation*}
        \int_{0}^{+\infty} f_1(\lambda) \diff \braket{\Omega, E_{\varphi, \comm{\omega}}(\lambda) \Omega} = \int_{0}^{+\infty} (\lambda - 1) \diff \braket{\Omega, E_{\varphi, \comm{\omega}}(\lambda) \Omega} = \vnorm[\big]{\Delta(\varphi / \comm{\omega})^{1/2} \Omega}^2 - \norm{\Omega}^2 < + \infty
    \end{equation*}
    by \cref{pro:operators_functionalCalculus} \ref{enu:operators_funcCalcNorm}. The last expression is finite because $\Omega \in D(\MH, \comm{\omega}) \subset \dom\bigl(\Delta(\varphi / \comm{\omega})^{1 / 2}\bigr)$ according to the final remarks made in \cref{para:vonNeumann_spatialDerivative}. The previous observations show that the sequence $(f_n)_n$ satisfies all of the assumptions of the monotone convergence theorem for decreasing sequences \cite[p. 65]{Cohn13}; with it, one obtains the assertion of the theorem:
    \begin{align*}
        \MS_\MFM^\mathrm{spa}(\omega, \varphi) &= - \int_{0}^{+\infty} \log(\lambda) \diff \braket{\Omega, E_{\varphi, \comm{\omega}}(\lambda) \Omega} = - \int_{0}^{+\infty} \lim_{n \to \infty} f_n(\lambda) \diff \braket{\Omega, E_{\varphi, \comm{\omega}}(\lambda) \Omega} \\
        &= - \lim_{n \to \infty} \int_{0}^{+\infty} \frac{\lambda^{t_n} - 1}{t_n} \diff \braket{\Omega, E_{\varphi, \comm{\omega}}(\lambda) \Omega} = - \lim_{t \searrow 0} \frac{1}{t} \Bigl(\vnorm[\big]{\Delta(\varphi / \comm{\omega})^{t/2} \Omega}^2 - \norm{\Omega}^2\Bigr) \ .
    \end{align*}

    \pp{2.} (\cite[Thm. 4.3 (1)]{DJP03}) Next, consider the case that $\Omega \notin \supp(\varphi)$. Since according to \cref{para:vonNeumann_spatialDerivative} the support of the spatial derivative is given by $\supp\bigl(\Delta(\varphi, \comm{\omega})\bigr) = \supp(\varphi)$, it follows that $\Omega = \Omega_1 + \Omega_2$, where $\Omega_1 \perp \Omega_2$ and $\Omega_1 \in \ker\bigl(\Delta(\varphi / \comm{\omega})\bigr) \setminus \{0\}$. (Recall that the support of an operator is the projection onto the orthogonal complement of its kernel.) With this, one computes
    \begin{align*}
        \vnorm[\big]{\Delta(\varphi / \comm{\omega})^{t/2} \Omega}^2 - \norm{\Omega}^2 &= \vnorm[\big]{\Delta(\varphi / \comm{\omega})^{t/2} \Omega_2}^2 - \norm{\Omega_2}^2 - \norm{\Omega_1}^2
    \end{align*}
    using the spectral representation of the first term and the Pythagorean theorem \cite[Eq. (V.23)]{Werner18} for the second one. It then follows that
    \begin{equation*}
        - \lim_{t \searrow 0} \frac{1}{t} \Bigl(\vnorm[\big]{\Delta(\varphi / \comm{\omega})^{t/2} \Omega_2}^2 - \norm{\Omega_2}^2 - \norm{\Omega_1}^2\Bigr) = \MS^\mathrm{spa}(\omega_{\Omega_2}, \varphi) + \lim_{t \searrow 0} \frac{\norm{\Omega_1}^2}{t} = + \infty \ .
    \end{equation*}
    Thus, also in the case that $\Omega \notin \supp(\varphi)$, the value of the relative entropy $\MS^\mathrm{spa}(\omega, \varphi)$ agrees with the one of the expression \eqref{eq:relativeEntropy_UhlmannRepresentation}, hence the assertion is proved.
\end{Proof}

\section{Uhlmann's Monotonicity Theorem}\label{sec:relativeEntropy_monotonicity}

The goal of this section is to prove the very general monotonicity theorem for the relative entropy of \textsc{A. Uhlmann}. To this end, a technical interpolation lemma is crucially needed which, in turn, relies on certain auxiliary inequalities from the theory of operator monotone functions.

\subsection{An Interpolation Inequality}

The following inequality is well-known for matrices \cite[Thm. V.1.9]{Bhatia97} (originally proved by \textsc{K. Löwner} (1934) \cite{Löwner34} and \textsc{E. Heinz} (1951) \cite{Heinz51}), but it actually holds true for positive self-adjoint operators on Hilbert spaces as well. The proof can be found in \cite[Prop. 10.14]{Schmüdgen12}.

\begin{lemma}[Löwner-Heinz inequality]\label{lem:relativeEntropy_exponentialOperatorMonotone}
    Let $A$ and $B$ be two self-adjoint operators on a Hilbert space $\MH$. If $A \ge B \ge 0$ (in the sense of \cref{para:forms_orderRelation}), then also $A^t \ge B^t$ for all $t \in (0, 1)$.
\end{lemma}

Introduce the notation $f_t : [0, + \infty) \longto [0, + \infty)$, $\lambda \longmto \lambda^t$. The Löwner-Heinz inequality implies that $f_t$ is an \emph{operator monotone} function. (See \cite{Bhatia97, HiaiPetz14, Simon19} for a systematic treatment of these functions.) The next result was established by \textsc{F. Hansen} and \textsc{G. K. Pedersen} \cite{Hansen80, HansenPedersen82, HansenPedersen03} for \emph{any} operator monotone function. Proofs of this result can also be found in \cite[Thm. V.2.3]{Bhatia97} and \cite[Lem. 1.2]{OP04}.

\begin{lemma}[Bounded Hansen-Jensen-Pedersen inequality]\label{lem:relativeEntropy_contractionInequalityBO}
    Let $\MH$ be a Hilbert space, let $A \in \PE{\BO(\MH)}$ be a bounded positive self-adjoint operator, and let $K \in \BO(\MH)$ be a contraction, that is, an operator such that $\norm{K}_\mop \le 1$. Then
    \begin{equation}\label{eq:relativeEntropy_HJPinequality}
        \forall t \in (0, 1) \ : \ K^\ast f_t(A) K \le f_t(K^\ast A K) \ .
    \end{equation}
\end{lemma}

With the help of the previous two lemmata, one can prove a generalization of the Hansen-Jensen-Pedersen inequality which is valid for unbounded operators. The proof is a slightly modified version of an argument given by \textsc{D. Petz} in \cite[Thm. B]{Petz85}; he uses a similar result as \cref{lem:relativeEntropy_exponentialOperatorMonotone} (which he proves in \cite[Thm. A]{Petz85} but which first appeared in \cite[Lem. D]{ArakiMasuda82}) to establish the relevant inequality.

\begin{lemma}[Unbounded Hansen-Jensen-Pedersen inequality]\label{lem:relativeEntropy_contractionInequality}
    If $T$ is a positive self-adjoint operator on a Hilbert space $\MH$ and $K \in \BO(\MH)$ is a contraction, then
    \begin{equation*}
        \forall t \in (0, 1) \ : \ K^\ast f_t(T) K \le f_t(K^\ast T K) \ .
    \end{equation*}
\end{lemma}

\begin{Proof}
    \pp{1.} For every $n \in \N$, denote by $E_n \ce E_T([0, n]) = \1_{[0, n]}(T)$ the spectral projection onto $[0, n]$. Then the following relation (in the sense of \cref{para:forms_orderRelation}) holds true:
    \begin{equation*}
        K^\ast E_n T K \le K^\ast T K
    \end{equation*}
    
    To see this, first, define the closed operators $A \ce T^{1/2} E_n K$ and $B \ce T^{1/2} K$ (\cf{} \cite[Prop. 8.11]{Soltan18}), and observe that $K^\ast E_n T K = A^\ast A$ and $K^\ast T K = B^\ast B$, hence the claim is equivalent to $A^\ast A \le B^\ast B$. Secondly, recall that for every closed operator $C$ on $\MH$, it holds that $\dom(C) = \dom\bigl((C^\ast C)^{1/2}\bigr)$ and $\norm{C \xi} = \norm{(C^\ast C)^{1/2} \xi}$ for all $\xi \in \dom(C)$ \cite[Lem. 7.1]{Schmüdgen12}. Therefore, the relation $A^\ast A \le B^\ast B$ is equivalent to $\dom(B) \subset \dom(A)$ and $\norm{A \xi} \le \norm{B \xi}$ for all $\xi \in \dom(B)$.

    Let $\xi \in \dom(B) = \dom(T^{1/2} K) = \set{\eta \in \MH \, : \, K \eta \in \dom(T^{1/2})}$ be arbitrary. Then clearly $E_n K \xi \in \dom(T^{1/2})$, that is, $\xi \in \dom(T^{1/2} E_n K) = \dom(A)$. Furthermore, one computes
    \begin{align*}
        \norm{T^{1/2} E_n K \xi}^2 &\le \vnorm[\big]{E_n T^{1/2} K \xi}^2 + \vnorm[\big]{(\id_\MH - E_n) T^{1/2} K \xi}^2 \\
        &= \vnorm[\big]{(E_n + \id_\MH - E_n) T^{1/2} K \xi}^2 \\
        &= \vnorm[\big]{T^{1/2} K \xi}^2 \ ,
    \end{align*}
    where the fact that functions of $T$ commute with the spectral projections of $T$ (\cref{pro:operators_functionalCalculus} \ref{enu:operators_funcCalcCommutation}) and the Pythagorean theorem \cite[Eq. (V.23)]{Werner18} were used. This proves the claim.
    
    \pp{2.} Applying now the result of \cref{lem:relativeEntropy_exponentialOperatorMonotone}, one obtains $f_t(K^\ast E_n T K) \le f_t(K^\ast T K)$, where as before $f_t$ denotes the function $\lambda \longmto \lambda^t$, $t \in (0, 1)$. Therefore, using this inequality and \cref{lem:relativeEntropy_contractionInequalityBO} applied to the bounded operator $E_n T$, one obtains
    \begin{equation}\label{eq:relativeEntropy_almostHJPinequality}
        K^\ast f_t(E_n T) K \le f_t(K^\ast E_n T K) \le f_t(K^\ast T K) \ .
    \end{equation}
    
    \pp{3.} Finally, observe that the operator $f_t(E_n T)$ is given as follows, \cf{} \cref{pro:operators_functionalCalculus} \ref{enu:operators_funcCalcCommutation}:
    \begin{equation*}
        f_t(E_n T) = \bigl(f_t \circ (\1_{[0, n]} \id_{\sigma(T)})\bigr)(T) =
        \begin{cases}
            f_t(T) & \text{on $[0, n]$} \ , \\
            f_t(0) & \text{on $\sigma(T) \setminus [0, n]$} \ .
        \end{cases}
    \end{equation*}
    Thus, the left-hand side of \eqref{eq:relativeEntropy_almostHJPinequality} can be written as
    \begin{equation*}
        K^\ast f_t(E_n T) K = K^\ast f_t(T) E_n K + K^\ast f_t(0) (\id_\MH - E_n) K \ .
    \end{equation*}
    In the limit $n \to + \infty$, the projection $E_n$ converges strongly to $\id_\MH$ by properties of the functional calculus \cite[Thm. 2.20]{StratilaZsido19}. Thus, $K^\ast f_t(E_n T) K \to K^\ast f_t(T) K$ by the above identity. Therefore, taking this limit in \cref{eq:relativeEntropy_almostHJPinequality} proves the statement.
\end{Proof}

With these preliminary observations, the main result of this subsection can be proved. The structure of the following proof is the one outlined in \cite[Lem. 5.2]{OP04}.

\begin{lemma}[Interpolation inequality]\label{lem:relativeEntropy_interpolation}
    For $i \in \set{1, 2}$, let $A_i$ be a positive self-adjoint operator on a Hilbert space $\MH_i$. Assume that $T : \MH_1 \longto \MH_2$ is a bounded linear operator with the properties
    \begin{enumerate}[label=\normalfont(\greek*)]
        \item \label{enu:relativeEntropy_interpolationAssumption1} $T\bigl(\dom(A_1)\bigr) \subset \dom(A_2)$;
        \item \label{enu:relativeEntropy_interpolationAssumption2} $\forall \xi \in \dom(A_1) : \norm{A_2 T \xi}_{\MH_2} \le \norm{T}_\mathrm{op} \, \norm{A_1 \xi}_{\MH_1}$.
    \end{enumerate}
    Then for every $0 < t < 1$ and $\xi \in \dom(A_1^t)$, it holds that
    \begin{equation}\label{eq:relativeEntropy_interpolationInequality}
        \norm{A_2^t T \xi}_{\MH_2} \le \norm{T}_\mathrm{op} \, \norm{A_1^t \xi}_{\MH_1} \ .
    \end{equation}
\end{lemma}

\begin{Proof}
    Without loss of generality, one may assume that $\norm{T}_\mathrm{op} = 1$; otherwise just rescale the bounded operator $T$ by $\norm{T}_\mathrm{op}^{-1}$. The two assumptions \ref{enu:relativeEntropy_interpolationAssumption1} and \ref{enu:relativeEntropy_interpolationAssumption2} together, and the conclusion of the lemma are, respectively, equivalent to the two operator relations
    \begin{equation*}
        T^\ast A_2^2 T \le A_1^2 \tand T^\ast A_2^{2t} T \le A_1^{2t} \ ,
    \end{equation*}
    in the sense of \cref{para:forms_orderRelation}. Indeed, the relation $T^\ast A_2^2 T \le A_1^2$ is equivalent to $(A_2 T)^\ast (A_2 T) \le A_1^\ast A_1$, and, as argued at the beginning of the proof of \cref{lem:relativeEntropy_contractionInequality}, this relation is, itself, equivalent to $\dom(A_1) \subset \dom(A_2 T)$ and $\norm{A_2 T \xi}_{\MH_2} \le \norm{A_1 \xi}_{\MH_1}$ for all $\xi \in \dom(A_1)$. These two relations are evidently equivalent to the assumptions \ref{enu:relativeEntropy_interpolationAssumption1} and \ref{enu:relativeEntropy_interpolationAssumption1} in the case $\norm{T}_\mop = 1$. The argument showing that $T^\ast A_2^{2t} T \le A_1^{2t}$ is equivalent to the conclusion of the lemma proceeds analogously.   

    Using again the notation $f_t(\lambda) = \lambda^t$ for $t \in (0, 1)$, \cref{lem:relativeEntropy_exponentialOperatorMonotone}, and the hypothesis $T^\ast A_2^2 T \le A_1^2$, it follows that $f_t(T^\ast A_2^2 T) \le f_t(A_1^2)$ in the sense of the functional calculus, that is,
    \begin{equation*}
        (T^\ast A_2^2 T)^t \le A_1^{2t} \ .
    \end{equation*}
    Furthermore, since the operator $T \in \BO(\MH_1, \MH_2)$ is a contraction by assumption, \cref{lem:relativeEntropy_contractionInequality} implies that $T^\ast f_t(A_2^2) T \le f_t(T^\ast A_2^2 T)$, \ie{}, $T^\ast A_2^{2t} T \le (T^\ast A_2^2 T)^t$. This, together with the above inequality, shows 
    $T^\ast A_2^{2t} T \le A_1^{2t}$ which is the assertion of the lemma.
\end{Proof}

\begin{remark}
    There are other ways to obtain \cref{lem:relativeEntropy_interpolation}. On the one hand, one can use a theorem of \textsc{H. Triebel} \cite[Thm. 1.18.10]{Triebel78} and the Calderón-Lions interpolation theorem \cite[Thm. IX.20]{RS2} to prove \eqref{eq:relativeEntropy_interpolationInequality}; an argument along these lines is given in \textsc{Petz}' paper \cite[Prop. 1]{Petz86a}. Another different argument leading to \cref{eq:relativeEntropy_interpolationInequality}, based on analyticity properties and the so-called three-line theorem \cite[p. 33]{RS2}, is presented in \cite[Thm. A.2]{DJP03}.
\end{remark}

The next inequality is also of importance for proving the monotonicity theorem; essentially, it is a special case of a general inequality for 2-positive maps mentioned in \cref{subsec:operatorAlgebras_CP}. The proof given below follows \cite[Thm. E]{Petz85}, but once again more details are provided.

\begin{lemma}\label{lem:relativeEntropy_auxiliaryProjectionInequality}
    Let $\alpha : \MFA \longto \MFB$ be a unital 2-positive mapping between two unital $C^\ast$-algebras, and let $Q \in \MFA$ and $P \in \MFB$ be projections such that $P \le \alpha(Q)$. Then for all $A \in \MFA$, it holds that
    \begin{equation}\label{eq:relativeEntropy_auxiliaryProjectionInequality}
        \alpha(A) P \alpha(A)^\ast \le \alpha(A Q A^\ast) \ .
    \end{equation}
\end{lemma}

\begin{Proof}
    \pp{1.} For every $\lambda > 0$, define an operator $C_\lambda \ce Q + \lambda (\idm_\MFA - Q) \in \MFA$. Since $Q$ and $\idm_\MFA - Q$ are projections in $\MFA$, $C_\lambda$ is invertible with $C_\lambda^{-1} = Q + \frac{1}{\lambda} \, (\idm_\MFA - Q)$. Therefore, \cref{lem:operatorAlgebras_invertibleLowerBound} shows that $C_\lambda \ge \epsilon \, \idm_\MFA$ for some $\epsilon > 0$, and hence also $\alpha(C_\lambda) \ge \epsilon \, \idm_\MFB$; this implies that $\alpha(C_\lambda)$ is invertible, too. From \cref{pro:operatorAlgebras_2positivity}, it now follows that
    \begin{equation}\label{eq:relativeEntropy_proofProjectionInequality}
        \alpha(A) \alpha(C_\lambda)^{-1} \alpha(A)^\ast \le \alpha(A C_\lambda^{-1} A^\ast)
    \end{equation}
    for all $A \in \MFA$. In the limit $\lambda \to + \infty$, it holds that $C_\lambda^{-1} \to Q$ in the operator norm. Therefore, the right-hand side of the above inequality tends to $\alpha(A Q A^\ast)$ since positive mappings are uniformly continuous, see \cref{pro:operatorAlgebras_positiveMappingContinuous}.
    
    \pp{2.} To analyze the left-hand side of the inequality, first, note that $P$ and $\alpha(Q)$ commute with each other. To see this, assume that $\MFB \subset \BO(\MH)$ for some Hilbert space $\MH$ (\cref{thm:operatorAlgebras_structureCStar}), and note that from the assumption $P \le \alpha(Q) \le \idm_\MFB$, it follows that $\ran(P) \subset \ran(\alpha(Q))$, hence $\alpha(Q) \xi = \xi$ for all $\xi \in \ran(P)$ by \cref{lem:operators_rangeProjection}. From this, one obtains for any $\eta \in \MH$, which can be decomposed as $\eta = \eta_1 + \eta_2$ with $\eta_1 \in \ran(P)$, $\eta_2 \in \ran(P)^\perp$ (\cref{thm:operators_orthogonalProjection}), that $\alpha(Q) P \eta = \alpha(Q) \eta_1 = \eta_1$. This shows that $\alpha(Q) P = P$ on $\MH$, and hence, by taking adjoints, also $P \alpha(Q) = P$ since $\alpha$ is a $\ast$-preserving map according to \cref{pro:operatorAlgebras_positiveImpliesStarPreserving}. Thus, $\alpha(Q) P = P \alpha(Q)$. In particular, this property implies that $P$ also commutes with the element $\alpha(C_\lambda)$.
    
    \pp{3.} Next, observe that from the assumption $P \le \alpha(Q) \le \idm_\MFA$, the fact that $P (\idm_\MFB - P) = 0$, and the Schwarz inequality for 2-positive maps (\cref{pro:operatorAlgebras_2PmapsAreSchwarzMaps}), it follows that
    \begin{align*}
        P \alpha(C_\lambda) &= P \alpha(Q) + \lambda P \bigl(\idm_\MFB - \alpha(Q)\bigr) \le \alpha(Q)^2 + \lambda P (\idm_\MFB - P) \le \alpha(Q^2) \le \idm_\MFB \ .
    \end{align*}
    This shows that $\idm_\MFB - P \alpha(C_\lambda) = \idm_\MFB - \alpha(C_\lambda) P \ge 0$. Furthermore, note that  $\alpha(C_\lambda)^{-1} \ge 0$ as well since $\sigma(\alpha(C_\lambda)^{-1}) = \set[\big]{\inv{\lambda} \, : \, \lambda \in \sigma(\alpha(C_\lambda))}$ \cite[Prop. 2.2.3]{BR1} and $\alpha(C_\lambda) \ge 0$. \cref{lem:operatorAlgebras_productPositiveElements} thus implies
    \begin{align*}
        \alpha(C_\lambda)^{-1} - P = \alpha(C_\lambda)^{-1} \bigl(\idm_\MFB - \alpha(C_\lambda) P\bigr) \ge 0 \ ,
    \end{align*}
    showing that $P \le \alpha(C_\lambda)^{-1}$. Using this inequality in the result \eqref{eq:relativeEntropy_proofProjectionInequality} obtained above, one arrives, after taking the limit $\lambda \to + \infty$, at the desired inequality.
\end{Proof}

\subsection{Proof of Uhlmann's Theorem}

The necessary tools to prove the main theorem of this section are assembled. The proof is guided by \cite[Thm. 5.3]{OP04} and \cite[Thm. 2]{Petz86a} but contains more details. Originally, this theorem was proved by \textsc{A. Uhlmann} in the paper \cite[Prop. 18]{Uhlmann77} in 1977 using the language of abstract interpolation theory.

\begin{ntheorem}[Uhlmann's monotonicity theorem]\label{thm:relativeEntropy_UhlmannMonotonicity}
    For $i \in \set{1, 2}$, let $\MFM_i \subset \BO(\MH_i)$ be a von Neumann algebra acting on a Hilbert space $\MH_i$ and $\varphi_i, \omega_i \in \NF{(\MFM_i)}$. Assume that $\omega_i = \omega_{\Omega_i}$ is given by a cyclic vector $\Omega_i \in \MH_i$ for the algebra $\MFM_i$, and that $\omega_1(\id_{\MH_1}) = \omega_2(\id_{\MH_2})$. Let $\alpha : \MFM_1 \longto \MFM_2$ be a Schwarz mapping such that the following relations hold on $\MFM_1$:
    \begin{equation}\label{eq:relativeEntropy_UhlmannAssumptions}
        \omega_2 \circ \alpha \le \omega_1 \tand \varphi_2 \circ \alpha \le \varphi_1 \ ,
    \end{equation}
    that is, $\omega_2(\alpha(A)) \le \omega_1(A)$ and $\varphi_2(\alpha(A)) \le \varphi_1(A)$ for all $A \in \MFM_1$. Finally, assume that at least one of the following two conditions is satisfied in addition to the previous assumptions:
    \begin{enumerate}[label=\normalfont(U\arabic*)]
        \item \label{enu:relativeEntropy_assumptionU1} the vector $\Omega_i \in \MH_i$ is separating for the algebra $\MFM_i$;
        \item \label{enu:relativeEntropy_assumptionU2} the mapping $\alpha$ is 2-positive.
    \end{enumerate}
    Under these assumptions, the relative entropy satisfies the following inequality:
    \begin{equation}\label{eq:relativeEntropy_UhlmannMonotonicity}
        \MS_{\MFM_1}^\mathrm{spa}(\omega_1, \varphi_1) \le \MS_{\MFM_2}^\mathrm{spa}(\omega_2, \varphi_2) \ .
    \end{equation}
\end{ntheorem}

\begin{Proof}
    Let $\MD_1 \ce \MFM_1 \Omega_1 \subset \MH_1$ and observe that this defines a dense subset of the Hilbert space $\MH_1$ due to the assumption that $\Omega_1$ is cyclic for $\MFM_1$. Define an operator $T_0 : \MH_1 \supset \MD_1 \longto \MH_2$ by
    \begin{equation*}
        T_0 (A \Omega_1) \ce \alpha(A) \Omega_2 \quad (A \in \MFM_1) \ .
    \end{equation*}

    \begin{statement}\label{sta:relativeEntropy_auxiliaryContraction}
        The operator $T_0$ is well-defined and extends to a contraction $T : \MH_1 \longto \MH_2$.
    \end{statement}

    \begin{subproof}[of statement]
        Using the Schwarz inequality \eqref{eq:operatorAlgebras_SchwarzInequality} for $\alpha$ and the assumption \eqref{eq:relativeEntropy_UhlmannAssumptions} on the states $\omega_1$ and $\omega_2$, one finds the following inequality for all $A \in \MFM_1$:
        \begin{align*}
            \norm{\alpha(A) \Omega_2}_{\MH_2}^2 &= \omega_2 \bigl(\alpha(A)^\ast \alpha(A)\bigr) \le \omega_2 \bigl(\alpha(A^\ast A)\bigr) \le \omega_1(A^\ast A) = \norm{A \Omega_1}_{\MH_1}^2 \ .
        \end{align*}
        This shows that $T_0$ is a contraction, \ie{}, $\norm{T_0}_{\BO(\MD_1, \MH_2)} \le 1$. (Due to the linearity of Schwarz mappings, $T_0$ is clearly a linear operator on $\MD_1$.) Moreover, the above inequality shows that $T_0$ is well-defined: if $A, B \in \MFM_1$ such that $A \Omega_1 = B \Omega_1$, then $(A - B) \Omega_1 = 0$ and hence
        \begin{equation*}
            \norm{\alpha(A - B) \Omega_2}_{\MH_2}^2 \le \norm{(A - B) \Omega_1}_{\MH_1}^2 = 0 \ .
        \end{equation*}
        From linearity of $\alpha$, it now follows that $\bigl(\alpha(A) - \alpha(B)\bigr) \Omega_2 = 0$, thus $T_0(A \Omega_1) = T_0(B \Omega_1)$. This shows that $T_0$ is well-defined. By the bounded linear transformation theorem (\cref{thm:operators_BLT}), $T_0$ extends to a contraction $T : \MH_1 \longto \MH_2$ with $T|_{\MD_1} = T_0$ and $\norm{T}_\mop = \norm{T_0}_{\BO(\MD_1, \MH_2)} \le 1$.
    \end{subproof}

    Denote by $\Delta_i$ the spatial derivative $\Delta(\varphi_i / \comm{\omega_i})$, where $\comm{\omega_i}$ is the vector state on $\comm{\MFM_i}$ induced by $\Omega_i$. Since $\Omega_i \in \MH_i$ is cyclic for $\MFM_i$, $\Delta_i$ is well-defined, \cf{} \cref{rem:relativeEntropy_spaRelativeEntropy}. Furthermore, let $\MFq_i : \MH_i \supset \dom(\MFq_i) \longto \C$ be the quadratic form $\MFq_i \ce \MFq_{\varphi_i} = \varphi_i\bigl(\Theta^{\comm{\omega_i}}(\,\sbullet\,)\bigr)$ from which $\Delta_i$ was constructed, see \cref{lem:vonNeumann_quadraticFormSpatialDerivative}. Observe that the space $\MFM_i \Omega_i$ is a core for $\Delta_i^{1/2}$; this follows from the observations at the end of \cref{para:vonNeumann_spatialDerivative} and the fact that, according to \cref{lem:vonNeumann_linealVectorState}, $\dom(\MFq_i) = D(\MH, \comm{\omega_i}) = \MFM_i \Omega_i$.
    
    \begin{statement}\label{sta:relativeEntropy_UhlmannInterpolationAssumptions}
        It holds that $T \bigl(\dom(\Delta_1^{1/2})\bigr) \subset \dom(\Delta_2^{1/2})$ and $\vnorm[\big]{\Delta_2^{1/2} \, T \xi}_{\MH_2} \le \vnorm[\big]{\Delta_1^{1/2} \xi}_{\MH_1}$ for all $\xi \in \dom(\Delta_1^{1/2})$.
    \end{statement}
    
    \begin{subproof}[of statement]
        Let $\xi \in \dom(\Delta_1^{1/2})$ be arbitrary. By the previous observation (\cf{} \cref{def:operators_closed}), there exists a sequence $(\xi_n)_{n \in \N} \subset \MH_1$ given by $\xi_n = A_n \Omega_1$, $A_n \in \MFM_1$, such that $\xi_n \to \xi$ and $\Delta_1^{1/2} \xi_n \to \Delta_1^{1/2} \xi$ in $\MH_1$ as $n \to + \infty$. From \cref{eq:vonNeumann_spatialDerivativeForm} and the definition of the operator $T$ from above, it follows for all $n, m \in \N$ that
        \begin{align}\label{eq:relativeEntropy_normH2}
            \begin{split}
                \vnorm[\big]{\Delta_2^{1/2} \, T \xi_n - \Delta_2^{1/2} \, T \xi_m}_{\MH_2}^2 &= \vnorm[\big]{\Delta_2^{1/2} T\bigl((A_n - A_m) \Omega_1\bigr)}_{\MH_2}^2 = \vnorm[\big]{\Delta_2^{1/2} \alpha(A_n - A_m) \Omega_2}_{\MH_2}^2 \\
                &= \MFq_2\bigl(\alpha(A_n - A_m) \Omega_2\bigr) = \varphi_2\Bigl(\Theta^{\comm{\omega_2}}\bigl(\alpha(A_n - A_m) \Omega_2\bigr)\Bigr) \ .
            \end{split}
        \end{align}
        (Observe that $\alpha(A_n - A_m) \Omega_2 \in \MFM_2 \Omega_2 = D(\MH_2, \comm{\omega_2}) = \dom(\MFq_2)$, hence \cref{eq:vonNeumann_spatialDerivativeForm} could indeed be applied; similarly, it was used that $\Theta^{\comm{\omega_2}}$ is defined for elements in $D(\MH, \comm{\omega_2})$, see \cref{para:vonNeumann_intertwiningOperator}.) According to \cref{pro:vonNeumann_ThetaVectorState}, the operator $\Theta^{\comm{\omega_2}}$ takes the following form:
        \begin{equation}\label{eq:relativeEntropy_thetaForMonotonicity}
            \Theta^{\comm{\omega_2}}\bigl(\alpha(A_n - A_m) \Omega_2\bigr) = \alpha(A_n - A_m) \bigl[\comm{\MFM_2} \Omega_2\bigr] \alpha(A_n - A_m)^\ast \ .
        \end{equation}
        
        \pp{1.\ Case \ref{enu:relativeEntropy_assumptionU1}.} If $\Omega_2 \in \MH_2$ is separating for $\MFM_2$, \cref{pro:operatorAlgebras_cyclicSeparating} \ref{enu:operatorAlgebras_separatingMcyclicCommM} shows that it is cyclic for $\comm{\MFM_2}$. Therefore, $\bigl[\comm{\MFM_2} \Omega_2\bigr] = \id_{\MH_2}$ according to \cref{cor:operators_characterizationPerp} \ref{enu:operators_densePerp}. Using now the Schwarz inequality from \cref{lem:operatorAlgebras_basicPropertiesSchwarzMap} as well as the assumption \eqref{eq:relativeEntropy_UhlmannAssumptions} on $\varphi_1$ and $\varphi_2$, one finds that
        \begin{align*}
            \vnorm[\big]{\Delta_2^{1/2} \, T \xi_n - \Delta_2^{1/2} \, T \xi_m}_{\MH_2}^2 &= \varphi_2\bigl(\alpha(A_n - A_m) \alpha(A_n - A_m)^\ast\bigr) \\[2pt]
            &\le \varphi_2 \circ \alpha\bigl((A_n - A_m) (A_n - A_m)^\ast\bigr) \\[3pt]
            &\le \varphi_1\bigl((A_n - A_m) (A_n - A_m)^\ast\bigr) \ .
        \end{align*}
        Similarly, since $\Omega_1 \in \MH_1$ is separating for $\MFM_1$, one finds that $\Theta^{\comm{\omega_1}}(A \Omega_1) = A \bigl[\comm{\MFM_1} \Omega_1\bigr] A^\ast = A A^\ast$ for all $A \in \MFM_1$. Thus, the last expression obtained before can be recast into the form
        \begin{align*}
            \varphi_1\bigl((A_n - A_m) (A_n - A_m)^\ast\bigr) &= \varphi_1\Bigl(\Theta^{\comm{\omega_1}}\bigl((A_n - A_m) \Omega_1\bigr)\Bigr) = \MFq_1\bigl((A_n - A_m) \Omega_1\bigr) \\
            &= \vnorm[\big]{\Delta_1^{1/2} (A_n - A_m) \Omega_1}_{\MH_1}^2 = \vnorm[\big]{\Delta_1^{1/2} \xi_n - \Delta_1^{1/2} \xi_m}_{\MH_1}^2 \ .
        \end{align*}
        Combining the previous observations, one obtains the following inequality for all $n, m \in \N$:
        \begin{equation}\label{eq:relativeEntropy_UhlmannAuxiliaryIneq}
            \norm{\Delta_2^{1/2} \, T \xi_n - \Delta_2^{1/2} \, T \xi_m}_{\MH_2}^2 \le \norm{\Delta_1^{1/2} \xi_n - \Delta_1^{1/2} \xi_m}_{\MH_1}^2 \ .
        \end{equation}

        \pp{2.\ Case \ref{enu:relativeEntropy_assumptionU2}.} Define the projections $Q \ce [\comm{\MFM_1} \Omega_1] \equiv \ssupp(\omega_1) \in \MFM_1$ and $P \ce [\comm{\MFM_2} \Omega_2] \equiv \ssupp(\omega_2) \in \MFM_2$. From $\omega_1\bigl(\id_{\MH_1} - Q\bigr) = 0$, \cf{} \cref{eq:operatorAlgebras_supportProjection}, and the assumption $\omega_2 \circ \alpha \le \omega_1$, it follows that $\omega_2\bigl(\id_{\MH_2} - \alpha(Q)\bigr) = 0$. This shows that $P \le \alpha(Q)$ (by definition of the support projection $P$), so the assumptions of \cref{lem:relativeEntropy_auxiliaryProjectionInequality} are satisfied: with \cref{eq:relativeEntropy_auxiliaryProjectionInequality}, one obtains for the operator $\Theta^{\comm{\omega_2}}$ from \cref{eq:relativeEntropy_thetaForMonotonicity} that
        \begin{align*}
            \Theta^{\comm{\omega_2}}\bigl(\alpha(A_n - A_m) \Omega_2\bigr) &= \alpha(A_n - A_m) P \alpha(A_n - A_m)^\ast \\[3pt]
            &\le \alpha\bigl((A_n - A_m) \, Q \, (A_n - A_m)^\ast\bigr) \\[2pt]
            &= \alpha\bigl(\Theta^{\comm{\omega_1}}\bigl((A_n - A_m) \Omega_1\bigr)\bigr) \ .
        \end{align*}
        Using this inequality, one can now estimate \cref{eq:relativeEntropy_normH2} for all $n, m \in \N$ as follows:
        \begin{align*}
            \vnorm[\big]{\Delta_2^{1/2} \, T \xi_n - \Delta_2^{1/2} \, T \xi_m}_{\MH_2}^2 &= \varphi_2\Bigl(\Theta^{\comm{\omega_2}}\bigl(\alpha(A_n - A_m) \Omega_2\bigr)\Bigr) \\
            &\le \varphi_2 \circ \alpha\Bigl(\Theta^{\comm{\omega_1}}\bigl((A_n - A_m) \Omega_1\bigr)\Bigr) \\
            &\le \varphi_1 \Bigl(\Theta^{\comm{\omega_1}}\bigl((A_n - A_m) \Omega_1\bigr)\Bigr) = \vnorm[\big]{\Delta_1^{1/2} \xi_n - \Delta_1^{1/2} \xi_m}_{\MH_1}^2 \ ,
        \end{align*}
        where in the third step the assumption $\varphi_2 \circ \alpha \le \varphi_1$ from \cref{eq:relativeEntropy_UhlmannAssumptions} was used. This is precisely the same inequality as obtained before in \cref{eq:relativeEntropy_UhlmannAuxiliaryIneq}.
        
        \pp{3. Conclusion.} It was shown that assuming either \ref{enu:relativeEntropy_assumptionU1} or \ref{enu:relativeEntropy_assumptionU2}, the inequality \eqref{eq:relativeEntropy_UhlmannAuxiliaryIneq} holds true; with its help, the proof of the statement can be completed. Since the sequence $(\Delta_1^{1/2} \xi_n)_{n \in \N} \subset \MH_1$ converges to $\Delta_1^{1/2} \xi$ by construction, it is a Cauchy sequence in $\MH_1$. Therefore, \eqref{eq:relativeEntropy_UhlmannAuxiliaryIneq} shows that $(\Delta_2^{1/2} T \xi_n)_{n \in \N} \subset \MH_2$ has the Cauchy property as well, hence it converges to some element $\eta \in \MH_2$. Furthermore,
        \begin{equation*}
            (T \xi_n)_{n \in \N} = \bigl(\alpha(A_n) \Omega_2\bigr)_{n \in \N} \subset \MFM_2 \Omega_2 = \dom(\MFq_2) \subset \dom(\Delta_2^{1/2})
        \end{equation*}
        as discussed above. It holds that $T \xi_n \to T \xi$ in $\MH_2$ as $n \to + \infty$ because $\xi_n \to \xi$ in $\MH_1$ by assumption and boundedness of $T$, \cf{} \cref{sta:relativeEntropy_auxiliaryContraction}. Since the operator $\Delta_2^{1/2}$ is, in particular, closed, the convergence of the sequences $(T\xi_n)_n$ and $(\Delta_2^{1/2} T \xi_n)_n$ implies that $T \xi \in \dom(\Delta_2^{1/2})$ and $\Delta_2^{1/2} T \xi = \eta$ \cite[Prop. 1.4]{Schmüdgen12}. As $\xi \in \dom(\Delta_1^{1/2})$ was arbitrary, the former shows that
        \begin{equation*}
            T \bigl(\dom(\Delta_1^{1/2})\bigr) \subset \dom(\Delta_2^{1/2})
        \end{equation*}
        which is the first assertion of the statement. On the other hand, one also finds
        \begin{equation*}
            \norm{\Delta_2^{1/2} \, T \xi}_{\MH_2} \le \norm{\Delta_1^{1/2} \xi}_{\MH_1}
        \end{equation*}
        by noting that the argument leading to \cref{eq:relativeEntropy_UhlmannAuxiliaryIneq} can be repeated analogously with $\xi_m \equiv 0$, giving $\norm{\Delta_2^{1/2} \, T \xi_n}_{\MH_2}^2 \le \norm{\Delta_1^{1/2} \xi_n}_{\MH_1}^2$ for all $n \in \N$, and then taking the limit $n \to + \infty$ and using the convergences discussed before. This is the second assertion of the statement, hence its proof is complete.
    \end{subproof}

    \cref{sta:relativeEntropy_UhlmannInterpolationAssumptions} shows that the assumptions of \cref{lem:relativeEntropy_interpolation} are satisfied for the operators $A_i = \Delta_i^{1/2}$ on $\MH_i$ and $T \in \BO(\MH_1, \MH_2)$ as constructed in \cref{sta:relativeEntropy_auxiliaryContraction}. Therefore, it follows that
    \begin{equation*}
        \norm{\Delta_2^{t/2} \, T \xi}_{\MH_2} \le \norm{\Delta_1^{t/2} \xi}_{\MH_1}
    \end{equation*}
    for all $t \in (0, 1)$ and $\xi \in \dom(\Delta_1^{1/2})$. Choosing $\xi = \id_{\MH_1} \Omega_1$, which is possible since $\MFM_1 \Omega_1 \subset \dom(\Delta_1^{1/2})$, yields $T \xi = \Omega_2$ because $\alpha$ is unital. Furthermore, by assumption it holds that $\norm{\Omega_1}_{\MH_1} = \norm{\Omega_2}_{\MH_2}$. Therefore, the above inequality implies
    \begin{equation}\label{eq:relativeEntropy_UhlmannAuxiliaryIneq2}
        \frac{1}{t} \Bigl(\norm{\Delta_2^{t/2} \Omega_2}_{\MH_2}^2 - \norm{\Omega_2}_{\MH_2}^2\Bigr) \le \frac{1}{t} \Bigl(\norm{\Delta_1^{t/2} \Omega_1}_{\MH_1}^2 - \norm{\Omega_1}_{\MH_1}^2\Bigr)
    \end{equation}
    for all $0 < t < 1$. Multiplying both sides by minus one, taking the limit $t \searrow 0$, and using Uhlmann's representation \eqref{eq:relativeEntropy_UhlmannRepresentation} of the relative entropy from \cref{thm:relativeEntropy_UhlmannRepresentation} finally yields the desired \cref{eq:relativeEntropy_UhlmannMonotonicity}.
\end{Proof}

\begin{remarks}\label{rem:relativeEntropy_UhlmannMonotonicity}
    \leavevmode
    \begin{enumerate}[env]
        \item \label{enu:relativeEntropy_monotonicitySTD} The assertion of the theorem also holds true for the relative entropy $\MS_\MFM^\mathrm{std}$ in standard form which can be shown by an almost analogous proof, see \cite[Thm. 4.4]{DJP03} for details.

        \item The two different assumptions \ref{enu:relativeEntropy_assumptionU1} and \ref{enu:relativeEntropy_assumptionU2} might be interesting in different situations. On the one hand, for the local algebras of quantum field theory, the Reeh-Schlieder theorem \cite[Sect. II.5.3]{Haag96}, \cite[p. 26]{HollandsSanders18}, \cite{ReehSchlieder61} provides cyclic and separating vectors. On the other hand, many of the Schwarz mappings that appear in applications are, in fact, completely positive (\cf{} \cref{subsec:operatorAlgebras_CP}).
        
        \item The proof of \cref{thm:relativeEntropy_UhlmannMonotonicity} given in \cite[Thm. 5.3]{OP04} seems to contain some inaccuracies: first, it is not very precise regarding the necessary assumptions on $\Omega_1, \Omega_2$ and $\alpha : \MFM_1 \longto \MFM_2$. While it is clear that $\Omega_i$ has to be cyclic for the algebra $\MFM_i$, $i \in \set{1,2}$, in order for the spatial derivative to be well-defined, in \cite{OP04} it is not explicitly assumed that $\Omega_i$ is separating for $\MFM_i$ as well. (This assumption is mentioned explicitly in \cite[Thm. 2]{Petz86a}, corresponding to \ref{enu:relativeEntropy_assumptionU1}, whereas \cite[Thm. 4]{Petz85} assumes $\alpha$ to be 2-positive, corresponding to \ref{enu:relativeEntropy_assumptionU2}.) Secondly, it appears that \cite{OP04} implicitly assumes the von Neumann algebras to be in standard form; this somewhat contradicts using the spatial form of the relative entropy in the first place. The proof of \cref{thm:relativeEntropy_UhlmannMonotonicity} provided here is truly valid in any representation of the von Neumann algebras, with the cost of having to pose additional assumptions on the vector representatives of $\omega_1, \omega_2$ or on the Schwarz mapping $\alpha$.
    \end{enumerate}
\end{remarks}

\subsection{Some Consequences of Monotonicity}

In the following, a couple of consequences of \cref{thm:relativeEntropy_UhlmannMonotonicity} shall be derived; they are all mentioned in \cite[p. 82]{OP04}. By \cref{rem:relativeEntropy_UhlmannMonotonicity} \ref{enu:relativeEntropy_monotonicitySTD}, these corollaries are also valid for $\MS_\MFM^\mathrm{std}$.

\begin{corollary}\label{cor:relativeEntropy_monotonicitySubalgebra}
    Let $\MFM_1$ and $\MFM_2$ be two von Neumann algebras acting on a Hilbert space $\MH$ such that $\MFM_1 \subset \MFM_2$, and let $\varphi, \omega \in \NF{(\MFM_2)}$ be two normal functionals on $\MFM_2$ such that $\omega$ is induced by a cyclic vector $\Omega \in \MH$ for the algebra $\MFM_1$. Then
    \begin{equation}\label{eq:relativeEntropy_monotonicitySubalgebra}
        \MS_{\MFM_1}^\mathrm{spa}\bigl(\omega\big|_{\MFM_1}, \varphi\big|_{\MFM_1}\bigr) \le \MS_{\MFM_2}^\mathrm{spa}(\omega, \varphi) \ .
    \end{equation}
\end{corollary}

\begin{Proof}
    Let $\omega_1 \ce \omega\big|_{\MFM_1}$ and $\varphi_1 \ce \varphi\big|_{\MFM_1}$ denote the restrictions of the the states $\omega, \varphi$ to the subalgebra $\MFM_1 \subset \MFM_2$. Define $\alpha : \MFM_1 \longto \MFM_2$ to be the inclusion $A \longmto A$ which is a completely positive mapping satisfying $\omega\bigl(\alpha(A)\bigr) = \omega_1(A)$ and $\varphi\bigl(\alpha(A)\bigr) = \varphi_1(A)$ for all $A \in \MFM_1$. This shows that the assumptions of \cref{thm:relativeEntropy_UhlmannMonotonicity} are satisfied, hence $\MS_{\MFM_1}^\mathrm{spa}(\omega_1, \varphi_1) \le \MS_{\MFM_2}^\mathrm{spa}(\omega, \varphi)$.
\end{Proof}

\begin{corollary}\label{cor:relativeEntropy_basicLowerBoundRevisited}
    Let $\MFM \subset \BO(\MH)$ be a von Neumann algebra and $\varphi, \omega \in \NF{\MFM}$ be two normal functionals with $\omega$ being induced by a cyclic vector $\Omega \in \MH$ such that $\Omega \in \supp(\varphi)$. Then
    \begin{equation*}
        \MS_\MFM^\mathrm{spa}(\omega, \varphi) \ge - \omega(\id_\MH) \log\left(\frac{\varphi(\id_\MH)}{\omega(\id_\MH)}\right) \ .
    \end{equation*}
\end{corollary}

\begin{Proof}
    Define the subalgebra $\MFM_0 \ce \set{\id_\MH} \subset \MFM$ and note that $\varphi|_{\MFM_0} \equiv \varphi(\id_\MH)$ and $\omega|_{\MFM_0} \equiv \omega(\id_\MH)$. If $\1 : \MFM \longto \C$, $A \longmto 1$, denotes the constant functional, one can write $\varphi_0 \ce \varphi|_{\MFM_0} = \varphi(\id_\MH) \, \1$ and $\omega_0 \ce \omega|_{\MFM_0} = \omega(\id_\MH) \, \1$. Using \cref{eq:relativeEntropy_monotonicitySubalgebra} and the scaling property of the relative entropy from \cref{pro:relativeEntropy_scaling} (which is proved analogously for $\MS_\MFM^\mathrm{spa}$), it follows that
    \begin{align*}
        \MS_{\MFM}^{\mathrm{spa}}(\omega, \varphi) &\ge \MS_{\MFM_0}^{\mathrm{spa}}(\omega_0, \varphi_0) = \omega(\id_\MH) \, \mathord{\underbrace{\MS_{\MFM_0}^{\mathrm{spa}}(\1, \1)}_{=0}} - \omega(\id_\MH) \log\left(\frac{\varphi(\id_\MH)}{\omega(\id_\MH)}\right) \\
        &= - \omega(\id_\MH) \log\left(\frac{\varphi(\id_\MH)}{\omega(\id_\MH)}\right) \ . \tag*{\qedhere}
    \end{align*}
\end{Proof}

\begin{corollary}\label{cor:relativeEntropy_nonnegSpaRelativeEntropy}
    Let $\MFM \subset \BO(\MH)$ be a von Neumann algebra and $\varphi, \omega \in \NF{\MFM}$ with $\omega$ being induced by a cyclic vector $\Omega \in \MH$. Then $\MS_\MFM^\mathrm{spa}(\omega, \varphi) \ge 0$.
\end{corollary}

\begin{corollary}\label{cor:relativeEntropy_independenceSpaRelEntVectorRepr}
    Let $\MFM \subset \BO(\MH)$ be a von Neumann algebra, let $\varphi \in \NF{\MFM}$ be arbitrary, and let $\omega = \omega_\Omega \in \NF{\MFM}$ be induced by a cyclic vector $\Omega \in \MH$. Then the the relative entropy $\MS_\MFM^\mathrm{spa}(\omega, \varphi)$ does not depend on the choice of $\Omega$.
\end{corollary}

\begin{Proof}
    Suppose that $\Omega^\prime \in \MH$ is another cyclic vector for $\MFM$ such that $\omega = \omega_{\Omega^\prime}$. Then, since $\omega_{\Omega} = \omega_{\Omega^\prime}$, it follows from \cref{thm:relativeEntropy_UhlmannMonotonicity} by choosing $\alpha : \MFM \longto \MFM$, $A \longmto A$, that $\MS_\MFM^{\mathrm{spa}}(\omega_{\Omega}, \varphi) \le \MS_\MFM^{\mathrm{spa}}(\omega_{\Omega^\prime}, \varphi)$ and also $\MS_\MFM^{\mathrm{spa}}(\omega_{\Omega^\prime}, \varphi) \le \MS_\MFM^{\mathrm{spa}}(\omega_{\Omega}, \varphi)$, hence $\MS_\MFM^{\mathrm{spa}}(\omega_{\Omega}, \varphi) = \MS_\MFM^{\mathrm{spa}}(\omega_{\Omega^\prime}, \varphi)$.
\end{Proof}

\section{Monotonicity on the Hilbert-Space Level}\label{sec:relativeEntropy_vectorMonotonicity}

\cref{thm:relativeEntropy_UhlmannMonotonicity} is the most general formulation of the monotonicity property of the relative entropy; this is highlighted, for example, by the various corollaries mentioned above. In this section, further consequences of Uhlmann's theorem shall be studied, namely monotonicity properties of $\MS_\MFM^\mathrm{spa}(\omega, \varphi)$ for the case in which $\varphi = \omega_\Phi$ is a functional induced by some vector $\Phi \in \MH$, too. (Recall from \cref{def:relativeEntropy_spaRelativeEntropy} that this is not a necessary assumption in the definition of the \emph{spatial form} of the relative entropy.)

\begin{para}[Motivation]
    In applications, this setting is interesting for the following reasons. \emph{First}, one usually does not consider completely arbitrary normal states on a von Neumann algebra, but rather (faithful) normal states induced by (cyclic and separating) vectors; therefore, it is interesting to examine how the monotonicity of the relative entropy is reflected in terms of the inducing vectors, \ie{}, investigate what kind of mapping on the Hilbert-space level increases or decreases the relative entropy when applied to the vector representatives. \emph{Second}, if one is given a vector on a Hilbert space which induces a state on a certain von Neumann algebra (\eg{}, the vacuum vector or coherent vectors, to which much attention has been paid recently regarding the computation of relative entropies \cite{CGP19, Galanda23, Hollands20, Longo19}), it might be handy to have a more specialized, concrete monotonicity inequality at hand, rather than having to rely on the very general but abstract \cref{thm:relativeEntropy_UhlmannMonotonicity}. \emph{Third}, in the applications mentioned before, it is often the case that one starts with a given vector $\Omega \in \MH$ inducing a state, and then transforms it via an operator $V \in \BO(\MH)$ to another vector $\Omega^\prime \ce V \Omega$ which is then interpreted to induce a new, transformed state. Therefore, it might be interesting to see the monotonicity property of the relative entropy explicitly in terms of this transformation $V$.

    The strategy for obtaining monotonicity results along these lines will be the following: given von Neumann algebras $\MFM_1 \subset \MFM_2 \subset \BO(\MH)$ and an element $V \in \MFM_2$ with certain properties, a Schwarz mapping $\alpha : \MFM_1 \longto \MFM_2$ will be constructed satisfying the additional assumptions of \cref{thm:relativeEntropy_UhlmannMonotonicity}; thus, the proof of the proposed results will be reduced to applying Uhlmann's much more general theorem.
\end{para}

\begin{notation}
    Let $\MFM \subset \BO(\MH)$ be a von Neumann algebra and $\Psi_1, \Psi_2 \in \MH$ be some vectors such that $\Psi_1$ is cyclic for $\MFM$. Define a functional $\MR_\MFM : \MH \times \MH \longto \R \cup \{+ \infty\}$ by
    \begin{equation*}
        \MR_\MFM(\Psi_1, \Psi_2) \ce \MS_\MFM^\mathrm{spa}(\omega_{\Psi_1}, \omega_{\Psi_2}) \ .
    \end{equation*}
\end{notation}

\subsection{A First Result}

The following proposition contains the first result in the direction described above.

\begin{proposition}\label{pro:relativeEntropy_vectorMonotonicityIsometry}
    Let $\MFM_1$ and $\MFM_2$ be von Neumann algebras such that $\MFM_1 \subset \MFM_2$, and assume that both algebras act on a Hilbert space $\MH$. Furthermore, let $V \in \MFM_2$ be an isometry, and let $\Omega, \Phi \in \MH$ be two vectors such that $V \Omega$ and $\Omega$ are cyclic for the algebras $\MFM_1$ and $\MFM_2$, respectively. Then the following inequality is satisfied:
    \begin{equation}\label{eq:relativeEntropy_vectorMonotonicityIsometry}
        \MR_{\MFM_1}(V \Omega, V \Phi) \le \MR_{\MFM_2}(\Omega, \Phi) \ .
    \end{equation}
\end{proposition}

\begin{Proof}
    Denote by $\omega_1 = \omega_{V \Omega}$ and $\varphi_1 = \omega_{V \Phi}$ the vector functionals on the algebra $\MFM_1$ induced by $V \Omega$ and $V \Phi$, and let $\omega_2 = \omega_{\Omega}$ and $\varphi_2 = \omega_{\Phi}$ be the vector functionals on $\MFM_2$ corresponding to $\Omega$ and $\Phi$. Note that since $V \Omega$ and $\Omega$ are assumed to be cyclic for $\MFM_1$ and $\MFM_2$, respectively, the spatial derivatives $\Delta(\varphi_1 / \comm{\omega_1})$ and $\Delta(\varphi_2 / \comm{\omega_2})$ are well-defined, hence are the relative entropies $\MS_{\MFM_1}^\mathrm{spa}(\omega_1, \varphi_1)$ and $\MS_{\MFM_2}^\mathrm{spa}(\omega_2, \varphi_2)$.

    Let $\iota : \MFM_1 \longto \MFM_2$, $A \longmto A$, denote the inclusion mapping which is a well-defined continuous injective $\ast$-homomorphism by the assumption $\MFM_1 \subset \MFM_2$. Define $\alpha : \MFM_1 \longto \MFM_2$ by
    \begin{equation}\label{eq:relativeEntropy_SchwarzMapFromLinOp}
        \alpha(A) \ce V^\ast \iota(A) V \quad (A \in \MFM_1) \ .
    \end{equation}
    Observe that this defines a Schwarz mapping between the von Neumann algebras $\MFM_1$ and $\MFM_2$: first, $\alpha$ is a completely positive mapping by \cref{exa:operatorAlgebras_CP} \ref{enu:operatorAlgebras_exaConjugationCP}, hence it satisfies the Schwarz inequality according to \cref{pro:operatorAlgebras_2PmapsAreSchwarzMaps}, and also assumption \ref{enu:relativeEntropy_assumptionU2} of \cref{thm:relativeEntropy_UhlmannMonotonicity}. Second, $\alpha$ is unital because $\alpha(\id_\MH) = V^\ast V = \id_\MH$ by the isometry property of $V$.

    Next, observe that the following identity holds true for the functional $\omega_2$ composed with the Schwarz mapping $\alpha$: for all $A \in \MFM_1$,
    \begin{align*}
        \omega_2 \circ \alpha(A) = \braket{\Omega, \alpha(A) \Omega} = \braket{\Omega, V^\ast \iota(A) V \Omega} = \braket{V \Omega, A \, V \Omega} = \omega_1(A) \ .
    \end{align*}
    Analogously, one also obtains $\varphi_2 \circ \alpha(A) = \varphi_1(A)$. Furthermore, it holds that $\norm{V \Omega} = \norm{\Omega}$ due to the isometry property, hence $\omega_1(\id_\MH) = \omega_2(\id_\MH)$. Therefore, all assumptions of \cref{thm:relativeEntropy_UhlmannMonotonicity} are satisfied, and it follows that $\MS_{\MFM_1}^\mathrm{spa}(\omega_1, \varphi_1) \le \MS_{\MFM_2}^\mathrm{spa}(\omega_2, \varphi_2)$. In the notation introduced above, this is exactly the inequality \eqref{eq:relativeEntropy_vectorMonotonicityIsometry}, hence the proposition is proved.
\end{Proof}

In the following, the assumptions on the vectors $\Omega, \Phi \in \MH$ and on the operator $V \in \BO(\MH)$ shall be discussed; the form of the mapping $\alpha$ from \cref{eq:relativeEntropy_SchwarzMapFromLinOp} will be commented on below.

\begin{remarks}
    \leavevmode
    \begin{enumerate}[env]
        \item The assumption that $V \Omega \in \MH$ should be cyclic for $\MFM_1$ appears to be rather strong since it restricts the admissible isometries $V$ in the above proposition. One may modify the assumptions on $\Omega$ and $V$ slightly as to still guarantee cyclicity of $V \Omega$ for $\MFM_1$; see \cref{eq:relativeEntropy_vectorMonotonicityIsometry2} below.
        
        \item The isometry property of the operator $V$ was only used to establish that the Schwarz mapping defined in \eqref{eq:relativeEntropy_SchwarzMapFromLinOp} is unital. Going back to \cref{thm:relativeEntropy_UhlmannMonotonicity}, one recognizes that this property of $\alpha$ is used in two places. First, in the case of assumption \ref{enu:relativeEntropy_assumptionU2}, to apply \cref{lem:relativeEntropy_auxiliaryProjectionInequality} which requires the mapping $\alpha$ to be unital. Second, to derive \cref{eq:relativeEntropy_UhlmannAuxiliaryIneq2} which uses that the operator $T : \MH_1 \longto \MH_2$, $T(A \Omega_1) = \alpha(A) \Omega_2$, defined in the proof of the theorem satisfies $T(\id_{\MH_1} \Omega_1) = \alpha(\id_{\MH_1}) \Omega_2 = \Omega_2$. Looking at this last expression, one may relax the assumption on $\alpha$, at least under the assumption \ref{enu:relativeEntropy_assumptionU1} of \cref{thm:relativeEntropy_UhlmannMonotonicity} which does not rely on \cref{eq:relativeEntropy_auxiliaryProjectionInequality}, as follows:
        \begin{itemize}[labelindent=\parindent, leftmargin=*]
            \item {\itshape The assertion of Uhlmann's theorem, assuming \ref{enu:relativeEntropy_assumptionU1}, still holds true if $\alpha$ satisfies the Schwarz inequality and the identity $\alpha(\id_{\MH_1}) \Omega_2 = \Omega_2$.}
        \end{itemize}
        One class of mappings $V$, more general than isometries as considered in the previous proposition but still interesting enough for applications, for which this relaxed condition on $\alpha$ might be feasible to verify is considered below in \cref{pro:relativeEntropy_vectorMonotonicityPartialIsometry}.

        \item The properties of the vector $\Phi$ also have a strong influence on the usefulness of the inequality \eqref{eq:relativeEntropy_vectorMonotonicityIsometry}: recall from \cref{def:relativeEntropy_spaRelativeEntropy} that $\Omega \in \supp(\varphi)$ is a necessary condition for $\MS_\MFM^\mathrm{spa}(\omega_\Omega, \varphi)$ to be finite. Therefore, in the context of \cref{pro:relativeEntropy_vectorMonotonicityIsometry}, the cases $V \Omega \in \supp(\omega_{V \Phi})$ and $\Omega \in \supp(\omega_{\Phi})$ are relevant. By \cref{pro:operatorAlgebras_supportVectorFunctional}, this is to say that
        \begin{equation*}
            V \Omega \in \clos_{\ndot}(\comm{\MFM_1} V \Phi) \tand \Omega \in \clos_{\ndot}(\comm{\MFM_2} \Phi) \ .
        \end{equation*}
        It is noteworthy to observe that if $V \Phi$ and $\Phi$ are separating for the algebras $\MFM_1$ and $\MFM_2$, respectively, then the above conditions are automatically satisfied by \cref{pro:operatorAlgebras_cyclicSeparating} \ref{enu:operatorAlgebras_separatingMcyclicCommM}.
    \end{enumerate}
\end{remarks}

Before turning to some modifications of the assumptions of \cref{pro:relativeEntropy_vectorMonotonicityIsometry}, some immediate consequences of \cref{eq:relativeEntropy_vectorMonotonicityIsometry} shall be recorded. The first one is nothing but the monotonicity statement of \cref{cor:relativeEntropy_monotonicitySubalgebra} for the special case of restricting vector functionals to a subalgebra.

\begin{corollary}\label{cor:relativeEntropy_vectorMonotonicitySubalgebra}
    Let $\MFM_1 \subset \MFM_2$ be von Neumann algebras acting on a Hilbert space $\MH$. If $\Omega, \Phi \in \MH$ are arbitrary vectors such that $\Omega$ is cyclic for the algebra $\MFM_1$, then
    \begin{equation*}
        \MR_{\MFM_1}(\Omega, \Phi) \le \MR_{\MFM_2}(\Omega, \Phi) \ .
    \end{equation*}
\end{corollary}

\begin{corollary}\label{cor:relativeEntropy_vectorMonotonicityOneAlgebra}
    Let $\MFM \subset \BO(\MH)$ be a von Neumann algebra, let $V \in \MFM$ be an isometry, and let $\Omega, \Phi \in \MH$ be arbitrary vectors such that $V \Omega$ and $\Omega$ are cyclic for $\MFM$. Then
    \begin{equation*}
        \MR_\MFM(V \Omega, V \Phi) \le \MR_\MFM(\Omega, \Phi) \ .
    \end{equation*}
\end{corollary}

\subsection{Modification of the Assumptions}

The next two propositions concern modifications of the setting of \cref{pro:relativeEntropy_vectorMonotonicityIsometry}.

\begin{proposition}\label{eq:relativeEntropy_vectorMonotonicityIsometry2}
    Let $\MFM_1 \subset \MFM_2$ be von Neumann algebras on a Hilbert space $\MH$, let $V \in \MFM_1$ be an isometry, and let $\Omega, \Phi \in \MH$ be two vectors such that $\Omega$ is cyclic for the algebra $\MFM_1$. Then
    \begin{equation*}
        \MR_{\MFM_1}(V \Omega, V \Phi) \le \MR_{\MFM_2}(\Omega, \Phi) \ .
    \end{equation*}
\end{proposition}

\begin{Proof}
    The argument establishing the inequality is identical to that of \cref{pro:relativeEntropy_vectorMonotonicityIsometry}, the only difference being that it has to be verified whether $V \Omega$ and $\Omega$ are cyclic for the algebras $\MFM_1$ and $\MFM_2$, respectively. The latter property is clear since $\Omega$ being cyclic for the subalgebra $\MFM_1$ implies that it is also cyclic for the larger algebra $\MFM_2$. For the former, let $\xi \in \MH$ be arbitrary. By cyclicity of $\Omega$ for $\MFM_1$, there exists a sequence $(B_n)_{n \in \N} \subset \MFM_1$ such that $\xi = \dlim{\MH}{n \to \infty} B_n \Omega$. For every $n \in \N$, define an operator $A_n \ce B_n V^\ast \in \MFM_1$. From the isometry property of $V$, it follows that $\xi = \dlim{\MH}{n \to \infty} B_n V^\ast V \Omega = \dlim{\MH}{n \to \infty} A_n V \Omega$, showing cyclicity of $V \Omega$ since $\xi \in \MH$ was arbitrary.
\end{Proof}

\begin{proposition}\label{pro:relativeEntropy_vectorMonotonicityPartialIsometry}
    Let $\MFM_1 \subset \MFM_2$ be von Neumann algebras on a Hilbert space $\MH$. Furthermore, let $V \in \MFM_2$ be a partial isometry, and let $\Omega, \Phi \in \MH$ such that $\Omega$ is an element of the initial subspace of $V$, and $V \Omega$ and $\Omega$ are cyclic and separating for the algebras $\MFM_1$ and $\MFM_2$, respectively. Then
    \begin{equation*}
        \MR_{\MFM_1}(V \Omega, V \Phi) \le \MR_{\MFM_2}(\Omega, \Phi) \ .
    \end{equation*}
\end{proposition}

\begin{Proof}
    The argument is again analogous to that of \cref{pro:relativeEntropy_vectorMonotonicityIsometry} with the following changes: since the vectors $V \Omega, \Omega \in \MH$ are now assumed to be cyclic and separating for the respective algebra, assumption \ref{enu:relativeEntropy_assumptionU1} of \cref{thm:relativeEntropy_UhlmannMonotonicity} is satisfied. Moreover, the mapping $\alpha(A) \ce V^\ast \iota(A) V$ is still completely positive, hence it satisfies the Schwarz inequality. Finally, one has that
    \begin{equation*}
        \alpha(\id_\MH) \Omega = V^\ast V \Omega = \Omega
    \end{equation*}
    because by assumption, $\Omega$ is in the initial subspace of $V$ onto which $V^\ast V$ is the unique orthogonal projection (\cf{} \cref{def:operators_partialIsometry}), hence \cref{lem:operators_rangeProjection} applies. The validity of the relations $\omega_2 \circ \alpha = \omega_1$ and $\varphi_2 \circ \alpha = \varphi_1$ is not influenced by the specific properties of $V$, hence all the assumptions of Uhlmann's \cref{thm:relativeEntropy_UhlmannMonotonicity} are satisfied, and the asserted inequality follows.
\end{Proof}

If the element $V \in \MFM_2$ is not only an isometry but a unitary operator, it turns out that there is a \enquote{reversed} form of \cref{eq:relativeEntropy_vectorMonotonicityIsometry}; this is the content of the next result.

\begin{proposition}\label{pro:relativeEntropy_vectorMonotonicityUnitary}
    Let $\MFM_1 \subset \MFM_2$ be von Neumann algebras on a Hilbert space $\MH$, let $U \in \MFM_2$ be a unitary operator, and let $\Omega, \Phi \in \MH$ be two vectors such that $\Omega$ is cyclic for $\MFM_1$. Then
    \begin{equation*}
        \MR_{\MFM_1}(\Omega, \Phi) \le \MR_{\MFM_2}(U \Omega, U \Phi) \ .
    \end{equation*}
\end{proposition}

\begin{Proof}
    Let $\omega_1$ and $\varphi_1$ be the vector functionals on $\MFM_1$ induced by $\Omega$ and $\Phi$, and let $\omega_2$ and $\varphi_2$ be the vector functionals on $\MFM_2$ coming from $U \Omega$ and $U \Phi$, respectively. From the assumption on $\Omega$ and $U^\ast U = \id_\MH$, it follows that $U \Omega$ is cyclic for $\MFM_2$ (see the proof of \cref{pro:relativeEntropy_vectorMonotonicityPartialIsometry}). This shows that both of the relative entropies are well-defined.
    
    Consider now the mapping $\alpha : \MFM_1 \longto \MFM_2$ given by $\alpha(A) \ce U \iota(A) U^\ast$ for all $A \in \MFM_1$, where $\iota : \MFM_1 \longhookrightarrow \MFM_2$ is the inclusion. Since $U$ is assumed to be unitary, it follows that $\alpha$ is unital. Moreover, invoking \cref{exa:operatorAlgebras_CP} \ref{enu:operatorAlgebras_exaConjugationCP} once more, $\alpha$ is completely positive, hence satisfies the Schwarz inequality; in fact, one can directly compute that
    \begin{equation*}
        \alpha(A^\ast A) = U A^\ast A U^\ast = U A^\ast U^\ast U A U^\ast = \alpha(A)^\ast \alpha(A)
    \end{equation*}
    for all $A \in \MFM_1$. Finally, there holds
    \begin{equation*}
        \omega_2 \circ \alpha(A) = \braket{U \Omega, \alpha(A) \, U \Omega} = \braket{U^\ast U \Omega, A \, U^\ast U \Omega} = \braket{\Omega, A \Omega} = \omega_1(A) \ ,
    \end{equation*}
    and similarly $\varphi_2 \circ \alpha = \varphi_1$. The assumptions of Uhlmann's \cref{thm:relativeEntropy_UhlmannMonotonicity} are satisfied; therefore, $\MS_{\MFM_1}^\mathrm{spa}(\omega_1, \varphi_1) \le \MS_{\MFM_2}^\mathrm{spa}(\omega_2, \varphi_2)$ which is the asserted inequality.
\end{Proof}

Combining \cref{cor:relativeEntropy_vectorMonotonicityOneAlgebra} and \cref{pro:relativeEntropy_vectorMonotonicityUnitary} yields the following result which is actually a special case of invariance of the relative entropy under $\ast$-automorphisms \cite[Eq. (6.6)]{Araki76}.

\begin{corollary}\label{cor:relativeEntropy_invarianceAutomorphisms}
    Let $\MFM$ be a von Neumann algebra acting on a Hilbert space $\MH$, let $U \in \MFM$ be a unitary operator, and let $\Omega, \Phi \in \MH$ be vectors such that $\Omega$ is cyclic for $\MFM$. Then
    \begin{equation*}
        \MR_\MFM(U \Omega, U \Phi) = \MR_\MFM(\Omega, \Phi) \ .
    \end{equation*}
\end{corollary}

\begin{Proof}
    First, observe that both relative entropies are well-defined due to the assumption on $\Omega$. (The argument is analogous to that given in the proof of \cref{pro:relativeEntropy_vectorMonotonicityPartialIsometry}.) Now, on the one hand, choosing $\MFM_1 = \MFM_2 = \MFM$ in \cref{pro:relativeEntropy_vectorMonotonicityUnitary} shows that $\MR_{\MFM}(\Omega, \Phi) \le \MR_{\MFM}(U \Omega, U \Phi)$. On the other hand, \cref{cor:relativeEntropy_vectorMonotonicityOneAlgebra} implies $\MR_\MFM(U \Omega, U \Phi) \le \MR_\MFM(\Omega, \Phi)$. Thus, $\MR_\MFM(U \Omega, U \Phi) \le \MR_\MFM(\Omega, \Phi) \le \MR_{\MFM}(U \Omega, U \Phi)$ which proves equality.
\end{Proof}

Having discussed modifications on the assumptions on the vectors $\Omega$ and $\Phi$, and on the operator $V \in \MFM_2$, now, the structure of the mapping $\alpha : \MFM_1 \longto \MFM_2$ defined in \cref{eq:relativeEntropy_SchwarzMapFromLinOp}, and the assumptions on the algebras $\MFM_1$ and $\MFM_2$ shall be commented on.

\begin{remarks}\label{rem:relativeEntropy_vectorMonotonicityAssumptions2}
    \leavevmode
    \begin{enumerate}[env]
        \item \label{enu:relativeEntropy_vectorMonotonicityAQFT} The assumption $\MFM_1 \subset \MFM_2$, which might seem artificial or very restrictive at first glance, appears to be rather natural in the context of algebraic quantum field theory \cite{Araki09, BDFY15, Haag96}. There, one considers nets $\MFA : \MO \longmto \MFA(\MO)$ of unital $C^\ast$-algebras $\MFA(\MO)$ for every open and bounded region $\MO \subset M$ of the spacetime manifold $M$. Each algebra $\MFA(\MO)$ is interpreted as the algebra of observables of the region $\MO$, \ie{}, it should contain the observables which can be measured within $\MO$.
        
        The local nets have to satisfy certain axioms, one of them being the axiom of \bemph{isotony} \cite[p. 2]{BDFY15}: if $\MO_1 \subset \MO_2$ for two open and bounded regions $\MO_1, \MO_2 \subset M$, then there exists an inclusion $\iota : \MFA(\MO_1) \longhookrightarrow \MFA(\MO_2)$ which implies that $\MFA(\MO_1) \subset \MFA(\MO_2)$. Choosing a representation $\pi : \MFA(\MO) \longto \BO(\MH)$ of the local algebras on a Hilbert space $\MH$, it follows that $\pi\bigl(\MFA(\MO_1)\bigr) \subset \pi\bigl(\MFA(\MO_2)\bigr)$, hence
        \begin{equation*}
            \bicomm{\pi\bigl(\MFA(\MO_1)\bigr)} \subset \bicomm{\pi\bigl(\MFA(\MO_2)\bigr)}
        \end{equation*}
        since the operation of forming the commutant is order-reversing \cite[Rem. 6.2]{Moretti19}. With this observation, the setting of \cref{pro:relativeEntropy_vectorMonotonicityIsometry} is recovered: suppose one is interested in calculating the relative entropy $\MR_{\MFM_1}(V \Omega, V \Phi)$ on a local algebra $\MFM_1 \ce \bicomm{\pi\bigl(\MFA(\MO_1)\bigr)}$, where $\Omega, \Phi \in \MH$ are elements in the representation space and $V \in \bicomm{\pi\bigl(\MFA(\MO_2)\bigr)} \ec \MFM_2$ is an isometry. If this turns out to be difficult or not analytically possible at all, one can consider the relative entropy $\MR_{\MFM_2}(\Omega, \Phi)$ on the larger algebra $\MFM_2$ which one might be able to compute. Then, by virtue of \cref{eq:relativeEntropy_vectorMonotonicityIsometry}, one obtains at least an upper bound for $\MR_{\MFM_1}(V \Omega, V \Phi)$.

        \item \label{enu:relativeEntropy_dilationSchwarz} A crucial step in the proofs of the above propositions consisted in invoking \cref{exa:operatorAlgebras_CP} \ref{enu:operatorAlgebras_exaConjugationCP} which showed that the conjugation of a $\ast$-homomorphism is a completely positive mapping. In fact, Stinespring's dilation theorem \cite[Thm. 1]{Stinespring55}, \cite[Thm. 1.2.7]{Stormer13} implies that all completely positive mappings on $C^\ast$-algebras are of this form. It turns out that there is a generalization of this theorem to von Neumann algebras \cite[Thm. 2.10]{Longo18}, \cite[Thm. 2.3]{BMSS12}, \cite[Ch. V, App. B, Cor. 9]{Connes94} which shows that \cref{pro:relativeEntropy_vectorMonotonicityIsometry} actually examines the implications of Uhlmann's \cref{thm:relativeEntropy_UhlmannMonotonicity} for the special case of a completely positive map $\alpha : \MFM_1 \longto \MFM_2$. Therefore, one might further investigate the manifestation of monotonicity of the relative entropy on the Hilbert-space level in the case of a general Schwarz mapping. To this end, it would be desirable to obtain a dilation theorem for Schwarz maps.
    \end{enumerate}
\end{remarks}

A final generalization of \cref{pro:relativeEntropy_vectorMonotonicityIsometry} that shall be considered here is relaxing the assumption $\MFM_1 \subset \MFM_2$. This is the content of the next proposition.

\begin{proposition}\label{pro:vectorMonotonicity_contraction}
    Let $\MFM_1$ and $\MFM_2$ be von Neumann algebras on a Hilbert space $\MH$. Assume that there exists a $\ast$-homomorphism $\rho : \MFM_1 \longto \MFM_2$ such that $A - \rho(A) \ge 0$ in $\BO(\MH)$ for all $A \in \MFM_1$. Let $V \in \MFM_2$ be an isometry and $\Phi, \Omega \in \MH$ such that $V \Omega$ and $\Omega$ are cyclic for $\MFM_1$ and $\MFM_2$, respectively. Then it holds that
    \begin{equation*}
        \MR_{\MFM_1}(V \Omega, V \Phi) \le \MR_{\MFM_2}(\Omega, \Phi) \ .
    \end{equation*}
\end{proposition}

\begin{Proof}
    Define $\omega_1 = \omega_{V \Omega}$, $\varphi_1 = \omega_{V \Phi}$ on $\MFM_1$ and $\omega_2 = \omega_{\Omega}$, $\varphi_2 = \omega_{\Phi}$ on $\MFM_2$, and consider the mapping $\alpha : \MFM_1 \longto \MFM_2$ given by $\alpha(A) \ce V^\ast \rho(A) V$ for all $A \in \MFM_1$. It still satisfies the Schwarz inequality since it is completely positive. Moreover, for all $A \in \MFM_1$ it holds that
    \begin{equation*}
        \omega_2 \circ \alpha(A) = \braket{V \Omega, \rho(A) V \Omega} \le \braket{V \Omega, A \, V \Omega} = \omega_1(A) \ ,
    \end{equation*}
    and analogously $\varphi_2 \circ \alpha \le \varphi_1$. The claim now follows again from Uhlmann's \cref{thm:relativeEntropy_UhlmannMonotonicity}.
\end{Proof}

\chapter{Perturbation Theory in Operator Algebras}\label{ch:perturbationTheory}

In this chapter, another aspect in the theory of operator algebras shall be illuminated: the perturbation of automorphism groups and KMS-states. This is a cornerstone in the theory of operator algebras and important for applications in quantum field theory and quantum statistical mechanics. In \cref{sec:perturbationTheory_dynamicalSystems}, the necessary background material is discussed, namely the \emph{Liouvillian} of a \emph{$W^\ast$-dynamical system} and the important notion of \emph{KMS-states}. \cref{sec:perturbationTheory_perturbation} gives a brief introduction to general perturbation theory and provides some results for unbounded perturbations. After that, \cref{sec:perturbationTheory_perturbationKMS} focuses on \emph{perturbations of KMS-states} using the framework developed in \cite{DJP03}. Finally, these methods are used in \cref{sec:perturbationTheory_twoSidedBogoliubov} to extend the so-called \emph{two-sided Bogoliubov inequality} to arbitrary von Neumann algebras.

\section{Dynamical Systems, Liouvillians, and KMS-States}\label{sec:perturbationTheory_dynamicalSystems}

Before the study of perturbation theory can commence, some technical background will be introduced in this section: continuous one-parameter groups of $\ast$-automorphisms, their associated Liouvillians, and the notion of KMS-states on von Neumann algebras.

\subsection{\texorpdfstring{The Liouvillian of a $W^\ast$-Dynamical System}{The Liouvillian of a W*-Dynamical System}}

\begin{definition}[$W^\ast$-dynamical system]\label{def:perturbationTheory_WStarDS}
    Let $\tau : \R \longto \Aut(\MFM)$, $t \longmto \tau_t$, be a one-parameter group of $\ast$-automorphisms of a von Neumann algebra $\MFM$, that is, $\tau_0 = \id_\MFM$ and $\tau_{s+t} = \tau_s \circ \tau_t$ for all $s, t \in \R$. The group $(\tau_t)_{t \in \R}$ is said to be \bemph{pointwise $\sigma$-weakly continuous} iff for every $A \in \MFM$, $t \longmto \tau_t(A)$ is continuous in the $\sigma$-weak topology on $\MFM$ \cite[p. 453]{DJP03}. In this case, the pair $(\MFM, \tau)$ is called a \bemph{$W^\ast$-dynamical system}, and the mapping $\tau$ is called \bemph{$W^\ast$-dynamics} on $\MFM$.
\end{definition}

\begin{examples}\label{exa:perturbationTheory_oneParameterGroups}
    \leavevmode
    \begin{enumerate}[env]
        \item \label{enu:perturbationTheory_exaFQS} Consider a \bemph{finite quantum system}, that is, the von Neumann algebra $\MFM = \BO(\MH)$ over a complex separable Hilbert space $\MH$, and let $H : \MH \supset \dom(H) \longto \MH$ be the self-adjoint \bemph{Hamiltonian} of $\MFM$. Define a one-parameter group of $\ast$-automorphisms $\tau = (\tau_t)_{t \in \R}$ of $\MFM$ by
        \begin{equation*}
            \tau_t(A) \ce \ee^{\ii t H} A \, \ee^{- \ii t H} \quad (t \in \R,\, A \in \MFM) \ .
        \end{equation*}
        One can show that $\tau$ is pointwise $\sigma$-weakly continuous \cite[Exa. 4.16]{AJP06}, hence $(\MFM, \tau)$ is a $W^\ast$-dynamical system. Furthermore, it holds that $\tau$ is strongly continuous on $\MFM$ if and only if the Hamiltonian $H$ is \emph{bounded} \cite[p. 133]{AJP06}. In this case, the infinitesimal generator of $\tau$, that is, the operator $\delta : \MFM \supset \dom(\delta) \longto \MFM$ satisfying $\odbound{}{t}{t=0} \, \tau_t(A) = \delta(A)$ \cite[Lem. II.1.3]{EngelNagel00}, is given by
        \begin{equation*}
            \delta(A) \ce \ii \, [H, A] \com A \in \dom(\delta) = \MFM \ .
        \end{equation*}

        \item Let $\MFA$ be a $C^\ast$-algebra. A $\ast$-automorphism $\alpha \in \Aut(\MFA)$ of $\MFA$ is called \bemph{inner} iff there exists a unitary element $U \in \MFA$ such that $\alpha(A) = U A \, U^\ast$ for all $A \in \MFA$. Suppose that $\MFA \subset \BO(\MH)$. If one is given a one-parameter group $(U_t)_{t \in \R}$ of unitary operators $U_t \in \UO(\MH)$, then one obtains a one-parameter group $\tau = (\tau_t)_{t \in \R}$ of $\ast$-automorphisms on $\MFA$ by setting \cite[p. 133]{AJP06}
        \begin{equation*}
            \tau_t(A) \ce U_t \, A \, U_t^\ast \quad (t \in \R,\, A \in \MFA) \ .
        \end{equation*}
        Groups of $\ast$-automorphisms of the above type are called \bemph{spatial}; one also says that the group $\tau : \R \longto \Aut(\MFA)$ is \bemph{implemented} by the group $U : \R \longto \UO(\MH)$.\qedhere
    \end{enumerate}
\end{examples}

\begin{para}[Standard Liouvillian]\label{para:perturbationTheory_standardLiouvillian}
    (\cite[pp. 456 f.]{DJP03}, \cite[pp. 149 f.]{AJP06})
    Let $(\MFM, \tau)$ be a $W^\ast$-dynamical system. \cref{thm:vonNeumann_existenceStandardForm} implies that there is a standard form $(\MFM, \MH, J, \NPC)$, and \cref{cor:vonNeumann_unitaryReprAut} yields for every $t \in \R$ a unique unitary operator $U_t \in \UO(\MH)$ such that $\tau_t(A) = U_t \, A \, U_t^\ast$ and $U_t \NPC \subset \NPC$. One can show that the group $(U_t)_{t \in \R}$ is strongly continuous \cite[p. 149]{AJP06}, hence Stone's theorem \cite[Thm. 6.2]{Schmüdgen12} implies that there exists a unique self-adjoint operator $L : \MH \supset \dom(L) \longto \MH$, called the \bemph{standard Liouvillian} of $\tau$, such that $U_t = \ee^{\ii t L}$ for every $t \in \R$. Moreover, it follows from the properties of the standard implementation of $\ast$-automorphisms given in \cref{cor:vonNeumann_unitaryReprAut} that
    \begin{enumerate}
        \item \label{enu:perturbationTheory_standardLiouvillianJ} $J L + L J = 0$;
        
        \item \label{enu:perturbationTheory_standardLiouvillianVectorRepr} $L \xi_\psi = 0$ for all $\tau$-invariant $\psi \in \NF{\MFM}$.
    \end{enumerate}
\end{para}

\begin{proposition}[\protect{\cite[Thm. 2.9]{DJP03}, \cite[Prop. 4.46]{AJP06}}]\label{pro:perturbationTheory_standardLiouvillian}
    The standard Liouvillian is the unique self-adjoint operator $L$ on $\MH$ such that for all $t \in \R$ and $A \in \MFM$, there holds
    \begin{enumerate}[num]
        \item $\ee^{\ii t L} \NPC \subset \NPC$;
        \item $\tau_t(A) = \ee^{\ii t L} A \, \ee^{- \ii t L}$.
    \end{enumerate}
\end{proposition}

\begin{examples}\label{exa:perturbationTheory_omegaLiouvillian}
    \leavevmode
    \begin{enumerate}[env]
        \item \label{enu:perturbationTheory_exaLiouvillianFQS} (\cite[Exa. 4.51]{AJP06}) Consider the finite quantum system from \cref{exa:perturbationTheory_oneParameterGroups} \ref{enu:perturbationTheory_exaFQS}. According to \cref{exa:vonNeumann_standardForm} \ref{enu:vonNeumann_exaStandardFormTypeI}, the standard form of $\MFM$ is given by $\bigl(\MFM, \HS(\MH), J, \HS(\MH)_+\bigr)$. The standard Liouvillian $L : \HS(\MH) \supset \dom(L) \longto \HS(\MH)$ is determined by the identity
        \begin{equation}\label{eq:perturbationTheory_exaLiouvillianFQS}
            \ee^{\ii t L} X = \ee^{\ii t H} X \, \ee^{- \ii t H} \com X \in \HS(\MH) \ .
        \end{equation}
        This can be seen as follows: on the one hand, the above relation for $L$ results in a unitary implementation of the dynamics $\tau_t$ since
        \begin{equation*}
            \ee^{\ii t L} A \, \ee^{- \ii t L} X = \ee^{\ii t H} \bigl(A \, \ee^{- \ii t L} X\bigr) \, \ee^{- \ii t H} = \ee^{\ii t H} \bigl(A \, (\ee^{- \ii t H} X \, \ee^{\ii t H})\bigr) \, \ee^{- \ii t H} = \ee^{\ii t H} A \, \ee^{- \ii t H} X \ .
        \end{equation*}
        On the other hand, the operator $\ee^{\ii t L}$ acting as in \eqref{eq:perturbationTheory_exaLiouvillianFQS} preserves the natural positive cone $\HS(\MH)_+$. Therefore, by \cref{pro:perturbationTheory_standardLiouvillian}, this determines $L$ uniquely.

        \item \label{enu:perturbationTheory_exaModularGroup} (\cite[p. 455]{DJP03}) Let $\MFM \subset \BO(\MH)$ be a $\sigma$-finite von Neumann algebra with cyclic and separating vector $\Omega \in \MH$, and standard form representation $(\MFM, \MH, J, \NPC)$ induced by $\Omega$ (\cf{} \cref{para:vonNeumann_NPC}). Let $\Delta_\Omega$ be the modular operator of $(\MFM, \Omega)$. Then the modular group (see \cref{def:vonNeumann_modularGroup}) $\sigma_t : \MFM \longto \MFM$,
        \begin{equation*}
            \sigma_t(A) = \Delta_\Omega^{\ii t} A \, \Delta_\Omega^{- \ii t} \quad (t \in \R,\, A \in \MFM) \ ,
        \end{equation*}
        is a $W^\ast$-dynamics on $\MFM$ with standard Liouvillian $L_\Omega = \log(\Delta_\Omega)$. To see this, note that \cref{pro:vonNeumann_propertiesNPC} \ref{enu:vonNeumann_DeltaNPC} implies that $\ee^{\ii t L_\Omega} \NPC = \Delta_\Omega^{\ii t} \NPC = \NPC$. Moreover, for every $t \in \R$ and $A \in \MFM$, one finds that
        \begin{align*}
            \ee^{\ii t L_\Omega} A \, \ee^{- \ii t L_\Omega} = \ee^{\ii t \log \Delta_\Omega} A \, \ee^{- \ii t \log \Delta_\Omega} = \Delta_\Omega^{\ii t} A \, \Delta_\Omega^{- \ii t} = \tau_t(A) \ .
        \end{align*}
        Thus, according to \cref{pro:perturbationTheory_standardLiouvillian}, $L_\Omega = \log(\Delta_\Omega)$ must be the unique standard Liouvillian of the $W^\ast$-dynamical system $(\MFM, \sigma)$.\qedhere
    \end{enumerate}
\end{examples}

\subsection{KMS-States on von Neumann Algebras}

In the following, thermal equilibrium states for systems with infinitely many degrees of freedom will be introduced. Their mathematical description reveals a deep connection with modular theory.

\begin{para}[Motivation]\label{para:perturbationTheory_motivationKMS}
    The \emph{Gibbs variational principle} \cite[Prop. 1.10]{OP04} shows that on a finite quantum system $\MFM = \BO(\MH)$ at inverse temperature\footnote{In this text, Boltzmann's constant $k_B$ is set equal to one; thus, the entropy is dimensionless, and temperature is measured in units of Joule.} $\beta = 1 / T$ with Hamiltonian $H \in \SE{\MFM}$ such that $\ee^{- \beta H} \in \NO(\MH)$, there exists a unique normal state $\varphi_{\beta} \in \NS(\MFM)$, called \emph{canonical Gibbs state} and given by
    \begin{align*}
        \rho_\beta = \frac{1}{Z_\beta} \, \ee^{- \beta H} \ , \twhere Z_\beta = \tr(\ee^{- \beta H}) \ ,
    \end{align*}
    which minimizes the \emph{free energy} functional $F_\beta(\varphi) = \varphi(H) + \inv{\beta} \tr(\rho_\varphi \log \rho_\varphi)$ and is hence a thermal equilibrium state \cite[p. 40]{Dongen17}. (The specific form of $\rho_\beta$ is implied by the \emph{principle of maximum entropy} and the conservation of internal energy \cite[p. 242]{GustafsonSigal20}.)

    For systems with infinitely many degrees of freedom, canonical Gibbs states do not exist because the operator $\ee^{- \beta H}$ will not be trace-class since $H$ has continuous spectrum \cite[p. 4]{BuchholzFredenhagen23}, \cite[pp. 206 f.]{Haag96}, \cite[p. 597]{Moretti18}. There is, however, an alternative description of equilibrium states based on their strong relationship with the time evolution \cite[p. 169]{AJP06}. To make this precise, consider the $W^\ast$-dynamics $\tau_t(A) = \ee^{\ii t H} A \, \ee^{- \ii t H}$ on $\MFM$, \cf{} \cref{exa:perturbationTheory_oneParameterGroups} \ref{enu:perturbationTheory_exaFQS}, and define for arbitrary $A, B \in \MFM$ the function \cite[p. 21]{Hiai21}, \cite[p. 19]{HollandsSanders18}
    \begin{equation}\label{eq:perturabtionTheory_thermalGreenFunction}
        F_{\beta, A, B}(z) \ce \frac{1}{Z_\beta} \, \tr\bigl(\ee^{- \beta H} A \, \ee^{\ii z H} B \, \ee^{- \ii z H}\bigr) \quad (z \in \C)
    \end{equation}
    which is analytic in the strip $S_\beta \ce \set{z \in \C \ : \ 0 < \Im(z) < \beta}$. For $z = t \in \R$, that is, $\Im(z) = 0$, there holds $F_{\beta, A, B}(t) = \varphi_\beta\bigl(A \tau_t(B)\bigr)$, and for $z = t + \ii \beta$, that is, $\Im(z) = \beta$, one finds
    \begin{align*}
        F_{\beta, A, B}(t + \ii \beta) &= \frac{1}{Z_\beta} \, \tr(\ee^{- \beta H} A \, \ee^{- \beta H} \ee^{\ii t H} B \, \ee^{- \ii t H} \ee^{\beta H}) \\
        &= \frac{1}{Z_\beta} \tr\bigl(\ee^{\ii t H} B \, \ee^{- \ii t H} A \, \ee^{- \beta H}\bigr) \\
        &= \varphi_\beta\bigl(\tau_t(B) A\bigr) \ .
    \end{align*}
    Thus, in a certain sense, the function $F_{\beta, A, B}$ encodes the \emph{non-commutativity} of the product $A \tau_t(B)$ in the thermal state $\varphi_\beta$ \cite[p. 169]{AJP06}.
    
    The boundary behavior of \eqref{eq:perturabtionTheory_thermalGreenFunction} and its analyticity properties were first pointed out by \textsc{R. Kubo}, \textsc{P. C. Martin}, and \textsc{J. Schwinger} \cite{Kubo57, MartinSchwinger59}. The following generalization of this observation to arbitrary normal states on von Neumann algebras is due to \textsc{R. Haag}, \textsc{N. M. Hugenholtz}, and \textsc{M. Winnink} \cite{HHW67} who postulated these conditions to be the defining property of equilibrium states \cite[pp. 200 ff.]{Haag96}.
\end{para}

\begin{definition}[KMS-state]\label{def:perturbationTheory_KMS}
    Let $(\MFM, \tau)$ be a $W^\ast$-dynamical system and $\beta > 0$. A normal state $\omega \in \NS(\MFM)$ is called a \bemph{$(\tau, \beta)$-KMS-state} iff for all $A, B \in \MFM$ there exists a function $F_{\beta, A, B} : \C \longto \C$ which is analytic in the strip $S_\beta = \set{z \in \C \, : \, 0 < \Im(z) < \beta}$ and continuous on its closure, and which satisfies the following \emph{KMS-boundary conditions} for all $t \in \R$:
    \begin{equation*}
        F_{\beta, A, B}(t) = \omega\bigl(A \tau_t(B)\bigr) \tand F_{\beta, A, B}(t + \ii \beta) = \omega\bigl(\tau_t(B) A\bigr) \ .
    \end{equation*}
\end{definition}

A thorough discussion of KMS-states, their properties, and arguments in favor of using them as equilibrium states can be found in \cite[Ch. 5]{BR2} or \cite[Sect. V.1 \& V.3]{Haag96}. Here, only a few results which are needed below in the development of perturbation theory shall be mentioned. To begin with, the following definition, taken from \cite[p. 460]{DJP03}, introduces a class of elements $A \in \MFM$ for which $t \longmto \tau_t(A)$ extends to an entire analytic function on $\C$. A proof of the subsequent important \cref{thm:perturbationTheory_analyticApproximation} can be found in \cite[Prop. 2.5.22]{BR1}, \cite[Lem. 2.13]{Hiai21}, or \cite[p. 170]{AJP06}.

\begin{definition}[$\tau$-entire element]\label{def:perturbationTheory_tauAnalytic}
    Let $(\MFM, \tau)$ be a $W^\ast$-dynamical system. An element $A \in \MFM$ is called \bemph{$\tau$-entire} iff there exists a function $f : \C \longto \MFM$ such that
    \begin{enumerate}[num]
        \item $f(t) = \tau_t(A)$ for all $t \in \R$;
        \item $\C \owns z \longmto \varphi\bigl(f(z)\bigr) \in \C$ is analytic for all $\varphi \in \pdual{\MFM}$.
    \end{enumerate}
    If these conditions are satisfied, one also writes $f(z) \equiv \tau_z(A)$ for $z \in \C$. The set of all $\tau$-entire elements of $\MFM$ will be denoted by $\MFM^\tau$.
\end{definition}

\begin{ntheorem}[Approximation by $\tau$-entire elements]\label{thm:perturbationTheory_analyticApproximation}
    Let $(\MFM, \tau)$ be a $W^\ast$-dynamical system. Then the set of all $\tau$-entire elements $\MFM^\tau$ forms a $\tau$-invariant $\sigma$-weakly and strongly dense $\ast$-subalgebra of $\MFM$.
\end{ntheorem}

The next proposition gives a characterization of the KMS-condition which is often more convenient to use than \cref{def:perturbationTheory_KMS}; it indicates that the dynamics $\tau$ measures the non-tracial character of a KMS-state \cite[p. 78]{BR2}. The proof can be found in \cite[Prop. 5.3.7]{BR2} or \cite[Thm. 5.4]{AJP06}.

\begin{proposition}\label{pro:perturbationTheory_differentCharacterizationKMS}
    Let $(\MFM, \tau)$ be a $W^\ast$-dynamical system and $\beta > 0$. A normal state $\omega \in \NS(\MFM)$ is a $(\tau, \beta)$-KMS-state if and only if there exists a $\sigma$-weakly dense $\tau$-invariant $\ast$-subalgebra $\mathfrak{D} \subset \MFM^\tau$ of $\tau$-entire elements of $(\MFM, \tau)$ such that for all $A, B \in \mathfrak{D}$:
    \begin{equation*}
        \omega\bigl(A \tau_{\ii \beta}(B)\bigr) = \omega(BA) \ .
    \end{equation*}
\end{proposition}

The following profound theorem \cite[Thm. 5.3.10]{BR2}, \cite[Thm. 2.14]{Hiai21}, originally due to \textsc{M. Takesaki} \cite[Thm. 13.1 \& 13.2]{Takesaki70}, relates the KMS-condition with the modular group.

\begin{ntheorem}[Takesaki]\label{thm:perturbationTheory_Takesaki}
    Let $\MFM \subset \BO(\MH)$ be a $\sigma$-finite von Neumann algebra with cyclic and separating vector $\Omega \in \MH$ and corresponding faithful normal functional $\omega = \omega_\Omega \in \NF{\MFM}$, and let $\Delta$, $J$ be the modular data and $(\sigma_t^\omega)_{t \in \R}$ be the modular automorphism group associated with $(\MFM, \Omega)$. Then the following assertions hold true:
    \begin{enumerate}
        \item \label{enu:perturbationTheory_modularKMS} $\omega$ satisfies the $(\sigma, - 1)$-KMS-condition (with respect to the strip $\set{z \in \C \, : \, \beta < \Im(z) < 0}$).
        
        \item \label{enu:perturbationTheory_uniquenessModularGroup} The modular group is uniquely determined as the $W^\ast$-dynamics on $\MFM$ for which $\omega$ satisfies the KMS-condition at $\beta = - 1$.
    \end{enumerate}
\end{ntheorem}

\begin{example}\label{exa:perturbationTheory_KMSfiniteQuantumSystem}
    (\cite[Prop. 4.7]{Attal06}, \cite[Prop. 9.25]{Landsman17}, \cite[Exa. 5.5]{AJP06})
    Consider the finite quantum system $(\MFM, \tau)$ from \cref{exa:perturbationTheory_oneParameterGroups} \ref{enu:perturbationTheory_exaFQS}. Let $\beta > 0$ and $\omega$ be an arbitrary $(\tau, \beta)$-KMS-state on $\MFM$ given by a density matrix $\rho_\omega \in \DM(\MH)$, and assume that the Hamiltonian $H$ is chosen such that $\ee^{- \beta H} \in \NO(\MH)$.
    
    For any two vectors $\xi, \eta \in \MH$, consider the operator $A \ce \braket{\xi, \,\sbullet\ } \, \eta \in \MFM$, and let $B \in \MFM^\tau$ be an arbitrary $\tau$-entire element. Observe that
    \begin{equation*}
        \omega\bigl(A \tau_{\ii \beta}(B)\bigr) = \tr\bigl(A \tau_{\ii \beta}(B) \rho_\omega\bigr) = \braket{\xi, \tau_{\ii \beta}(B) \rho_\omega \eta} \tand \omega(BA) = \tr\bigl(\rho_\omega B A\bigr) = \braket{\xi, \rho_\omega B \eta} \ .
    \end{equation*}
    The characterization of KMS-states from \cref{pro:perturbationTheory_differentCharacterizationKMS} now implies that
    \begin{equation*}
        \braket{\xi, \tau_{\ii \beta}(B) \rho_\omega \eta} = \braket{\xi, \rho_\omega B \eta} \ .
    \end{equation*}
    Since $\xi, \eta \in \MH$ were arbitrary, this shows that the operator relation $\tau_{\ii \beta}(B) \rho_\omega = \rho_\omega B$ holds true for all $\tau$-entire elements $B$ in $\MFM$; by definition of $\tau$, it can be written as
    \begin{equation*}
        B (\ee^{\beta H} \rho_\omega) = (\ee^{\beta H} \rho_\omega) B \ .
    \end{equation*}
    Since $B \in \MFM^\tau$ was arbitrary and $\MFM^\tau$ lies strongly dense in $\MFM$ by \cref{thm:perturbationTheory_analyticApproximation}, it follows that the above relation is valid for all $B \in \MFM$. Therefore, one may conclude that $\ee^{\beta H} \rho_\omega \in \comm{\MFM} = \C \cdot \id_\MH$, hence there exists $c \in \C$ such that $\rho_\omega = c \, \ee^{- \beta H}$. Since $\rho_\omega$ must be a density matrix, it follows that
    \begin{equation*}
        \rho_\omega = \frac{\ee^{- \beta H}}{\tr(\ee^{- \beta H})} \ .
    \end{equation*}
    This argument, together with the computations in \cref{para:perturbationTheory_motivationKMS}, shows that the $W^\ast$-dynamical system $(\MFM, \tau)$ admits a $(\tau, \beta)$-KMS-state if and only if $\ee^{- \beta H}$ is a trace-class operator, and that in this case the KMS-state is given by the canonical Gibbs state.
\end{example}

The next result, taken from \cite[Thm. 2.13]{DJP03}, gives a characterization of the KMS-condition expressed on the Hilbert-space level using the Liouvillian of the dynamics.

\begin{proposition}\label{pro:perturbationTheory_LiouvillianKMS}
    Let $(\MFM, \MH, J, \MP)$ be a von Neumann algebra in standard form and $\Omega \in \MP$ be a unit vector. Furthermore, let $\tau$ be a $W^\ast$-dynamics on $\MFM$ with standard Liouvillian $L$. Then $\omega_\Omega$ is a $(\tau, \beta)$-KMS state if and only if $\MFM \Omega \subset \dom(\ee^{- \beta L / 2})$ and
    \begin{equation}\label{eq:perturbationTheory_LiouvillianKMS}
        \forall A \in \MFM \ : \ \ee^{- \beta L / 2} A \Omega = J A^\ast \Omega \ .
    \end{equation}
    In this case, if $\Omega$ is cyclic and $\Delta_\Omega$ denotes the corresponding modular operator, then
    \begin{equation}\label{eq:perturbationTheory_modularOperatorLiouvillian}
        \Delta_\Omega = \ee^{- \beta L} \ .
    \end{equation}
\end{proposition}

This section shall be concluded by providing certain results regarding the convergence of $W^\ast$-dynamics and their associated Liouvillians, invariant states, and KMS-states. The proof of the following proposition can be found in \cite[Thm. 2.14]{DJP03}.

\begin{proposition}\label{pro:perturbationTheory_convergence}
    Consider a von Neumann algebra $\MFM$ in standard form $(\MFM, \MH, J, \NPC)$.
    \begin{enumerate}
        \item Let $(\tau_n)_{n \in \N}$ be a sequence of $W^\ast$-dynamics $\tau_n = (\tau_{n, t})_{t \in \R}$ on $\MFM$, let $L_n$ denote the standard Liouvillian of $\tau_n$, and assume that $L$ is a self-adjoint operator on $\MH$ such that $L_n \to L$ in the strong resolvent sense (\cref{def:operators_SRConvergence}). Then
        \begin{equation*}
            \tau_t(A) \ce \ee^{\ii t L} A \, \ee^{- \ii t L} \quad (t \in \R,\, A \in \MFM)
        \end{equation*}
        defines a $W^\ast$-dynamics on $\MFM$ with standard Liouvillian given by $L$.

        \item Assume, in addition to the previous statement, that for every $n \in \N$, $\omega_n \in \NF{\MFM}$ is a $\tau_n$-invariant normal functional with standard vector representative $\Omega_n \in \NPC$, and that $\dlim{w}{n \to \infty} \Omega_n = \Omega$. Then $\Omega \in \NPC$ and the vector functional $\omega_\Omega$ is $\tau$-invariant.
        
        \item Let $\beta > 0$. In the situation of the previous statement, assume additionally that for every $n \in \N$, $\omega_n$ is a $(\tau_n, \beta)$-KMS-state, and that $\Omega \neq 0$. Then $\omega_{\Omega / \norm{\Omega}}$ is a $(\tau, \beta)$-KMS-state.
    \end{enumerate}
\end{proposition}

\section{Foundations of Perturbation Theory}\label{sec:perturbationTheory_perturbation}

This section begins with the presentation of perturbation theory in operator algebras. The focus lies on unbounded perturbations since they are especially interesting in the context of the two-sided Bogoliubov inequality, and bounded perturbations are more commonly treated in the existing literature.

\subsection{Bounded Perturbations}

\begin{para}[Motivation]
    Consider a finite quantum system $\MFM = \BO(\MH)$ with bounded self-adjoint Hamiltonian $H_0 \in \SE{\MFM}$ generating the dynamics $\tau_t(A) = \ee^{\ii t H_0} A \, \ee^{- \ii t H_0}$. It was mentioned in \cref{exa:perturbationTheory_oneParameterGroups} \ref{enu:perturbationTheory_exaFQS} that the infinitesimal generator of $\tau$ is given by $\delta_0(A) = \ii \, [H_0, A]$ for all $A \in \MFM$.
    
    A natural question to ask in this context is how the dynamics change if $H_0$ is \emph{perturbed} by some bounded self-adjoint operator $V \in \SE{\MFM}$; physically speaking, this corresponds to adding an interaction to the Hamiltonian $H_0$ describing a free quantum system. In this case, it follows that the generator $\delta$ for the dynamics coming from the Hamiltonian $H = H_0 + V$ is given by
    \begin{equation}\label{eq:perturbationTheory_perturbedGenerator}
        \delta(A) = \delta_0(A) + \ii \, [V, A] \quad (A \in \MFM) \ .
    \end{equation}
    For more general $W^\ast$-dynamical systems, where there is no fixed Hamiltonian which could be perturbed, one may take \eqref{eq:perturbationTheory_perturbedGenerator} as a \emph{definition} for the perturbed dynamics since this identity does not require any particular form for the unperturbed generator $\delta_0$, and it is a purely algebraic identity, hence perfectly suited for the operator-algebraic setting. This approach will be outlined below.

    Bounded perturbation theory in operator algebras was initiated by \textsc{H. Araki} in 1973 \cite{Araki73a, Araki73b, Araki73c}; his work was inspired by the perturbative expansions used in quantum electrodynamics. A different approach to the theory was developed by \textsc{M. J. Donald} in \cite{Donald90} which deals with semi-bounded perturbations; a review can be found in \cite[Ch. 12]{OP04}.
\end{para}

\begin{para}[Perturbed dynamics]\label{para:perturbationTheory_localPerturbations}
    (\cite[pp. 147 f.]{BR2}, \cite[Sect. 4.8]{AJP06})
    Let $(\MFM, \tau)$ be a $W^\ast$-dynamical system and $\delta_0 : \MFM \supset \dom(\delta_0) \longto \MFM$ be the infinitesimal generator of $\tau$. A \bemph{bounded perturbation} of the dynamics $\tau$ is obtained by perturbing $\delta_0$ with $V \in \SE{\MFM}$ as follows:
    \begin{equation*}
        \delta_V \ce \delta_0 + \ii \, [V, \,\sbullet\,] \twith \dom(\delta_V) = \dom(\delta_0) \ .
    \end{equation*}
    Indeed, using the Hille-Yosida theorem \cite[Thm. 3.2.50]{BR1}, \cite[Prop. 4.15]{AJP06}, it follows that $\delta_V$ is the generator of a pointwise $\sigma$-weakly continuous one-parameter group of $\ast$-automorphisms $\tau^V : \R \longto \Aut(\MFM)$ on $\MFM$ \cite[p. 165]{AJP06}, called the \bemph{perturbed dynamics}.

    One can show that for every $t \in \R$ and $A \in \MFM$, the perturbed dynamics $\tau_t^V(A)$ can be expressed in terms of the following series expansion, where the integrals exist in the $\sigma$-weak operator topology and define a norm-convergent series of bounded operators \cite[Prop. 5.4.1]{BR2}, \cite[p. 165]{AJP06}:
    \begin{align}\label{eq:perturbationTheory_DysonRobinson}
        \begin{split}
            \tau_t^V(A) &= \tau_t(A) \\
            &\quad + \sum_{n=1}^{\infty} \ii^n \int_{0}^{t} \diff t_1 \int_{0}^{t_1} \diff t_2 \, \dotsb \int_{0}^{t_{n-1}} \diff t_n \, \Bigl[\tau_{t_n}(V), \bigl[\dotsb, \bigl[\tau_{t_2}(V), [\tau_{t_1}(V), \tau_t(A)]\bigr] \dotsb\bigr]\Bigr] \ .
        \end{split}
    \end{align}

    An important tool in perturbation theory is the following one-parameter family $\expan_\tau^V : \R \longto \MFM$, $t \longmto \expan_\tau^V(t)$, of elements in $\MFM$, called \bemph{Araki-Dyson expansional}, which is defined by
    \begin{equation}\label{eq:perturbationTheory_ArakiDyson}
        \expan_\tau^V(t) \ce \idm_\MFM + \sum_{n=1}^{\infty} \ii^n \int_{0}^{t} \diff t_1 \int_{0}^{t_1} \diff t_2 \, \dotsb \int_{0}^{t_{n-1}} \diff t_n \, \tau_{t_n}(V) \dotsm \tau_{t_2}(V) \tau_{t_1}(V) \ .
    \end{equation}
    One can show that $\bigl(\expan_\tau^V(t)\bigr)_{t \in \R}$ forms a one-parameter family of unitary elements in $\MFM$, and that the following relationship with the perturbed dynamics holds true \cite[Prop. 5.4.1]{BR2}:
    \begin{equation*}
        \tau_t^V(A) = \expan_\tau^V(t) \, \tau_t(A) \, \expan_\tau^V(t)^\ast \ .
    \end{equation*}
    
    If $\MFM \subset \BO(\MH)$ and the dynamics $\tau$ is unitarily implemented on $\MH$, the same follows for the perturbed dynamics $\tau^V$, and also the cocycle $\expan_\tau^V$ is given in terms of a concrete expression.
\end{para}

\begin{proposition}[\protect{\cite[Cor. 5.4.2]{BR2}}]\label{pro:perturbationTheory_spatialPerturbation}
    Assume that $\tau$ is given by $\tau_t(A) = U_t A U_t^\ast$, where $t \longmto U_t = \ee^{\ii t L}$ is a strongly continuous unitary group on $\MH$ with $L : \MH \supset \dom(L) \longto \MH$ a self-adjoint operator. Then the perturbed dynamics $\tau^V$ and the expansional $\expan_\tau^V$ are given for all $t \in \R$ and $A \in \MFM$ by
    \begin{equation}\label{eq:perturbationTheory_dynamicsSpatialPerturbation}
        \tau_t^V(A) = \ee^{\ii t (L + V)} A \, \ee^{- \ii t (L + V)} \tand \expan_\tau^V(t) = \ee^{\ii t (L + V)} \, \ee^{- \ii t L} \ .
    \end{equation}
\end{proposition}

\subsection{Analytic and Unbounded Perturbations}\label{subsec:perturbationTheory_analyticUnboundedPerturbations}

\begin{para}[Analytic perturbations]\label{para:perturbationTheory_analyticPerturbations}
    (\cite[Sect. 3.2]{DJP03})
    Let $(\MFM, \tau)$ be a $W^\ast$-dynamical system and $V$ be a self-adjoint, $\tau$-entire element of $\MFM$. In analogy to the expansion \eqref{eq:perturbationTheory_DysonRobinson} for the perturbed dynamics $\tau^V$, one can define for all $z \in \C$ and $A \in \MFM^\tau$:
    \begin{equation*}
        \tau_z^V(A) \ce \tau_z(A) + \sum_{n=1}^{\infty} (\ii z)^n \int_{0}^{1} \diff s_1 \int_{0}^{s_1} \diff s_2 \, \dotsb \int_{0}^{s_{n-1}} \diff s_n \, \Bigl[\tau_{s_n z}(V), \bigl[\dotsb, [\tau_{s_1 z}(V), \tau_z(A)] \dotsb\bigr]\Bigr] \ .
    \end{equation*}
    Note that by assumption, the elements $\tau_{s_n z}(V)$ and $\tau_z(A)$ are well-defined. Therefore, it immediately follows that the perturbed dynamics $\tau^V$ possesses the same analytic elements as $\tau$, that is, $\MFM^{\tau^V} = \MFM^\tau$. Furthermore, one can define the Araki-Dyson expansionals for all $V \in \MFM^\tau \cap \SE{\MFM}$ and $z \in \C$ by
    \begin{equation*}
        \expan_\tau^V(z) \ce \sum_{n=0}^{\infty} (\ii z)^n \int_{0}^{1} \diff s_1 \int_{0}^{s_1} \diff s_2 \, \dotsb \int_{0}^{s_{n-1}}  \diff s_n \, \tau_{s_n z}(V) \dotsm \, \tau_{s_2 z}(V) \tau_{s_1 z}(V) \ .
    \end{equation*}
    Both of the above series converge uniformly in the norm topology for all $z$ in compact sets, and they define analytic functions with values in $\MFM$ \cite[p. 461]{DJP03}.
\end{para}

\begin{proposition}[\protect{\cite[Thm. 3.2]{DJP03}}]\label{pro:perturbationTheory_analyticPerturbations}
    Let $\MFM \subset \BO(\MH)$ and assume that there is a self-adjoint operator $L$ on $\MH$ such that $\tau_t(B) = \ee^{\ii t L} B \, \ee^{- \ii t L}$ for all $B \in \MFM$. Furthermore, let $V \in \MFM^\tau \cap \SE{\MFM}$, $A \in \MFM^\tau$, and $z, z_1, z_2 \in \C$ be arbitrary. Then the following assertions hold true:
    \begin{enumerate}
        \item $E_\tau^V(z) \in \MFM^\tau$.
        
        \item $\tau_z^V(A) = \ee^{\ii z (L + Q)} A \, \ee^{- \ii z (L + Q)}$.
        
        \item \label{enu:perturbationTheory_complexSpatialExpansional} $\expan_\tau^V(z) = \ee^{\ii z (L + Q)} \, \ee^{- \ii z L}$.

        \item \label{enu:perturbationTheory_complexExpansionalInverse} $\expan_\tau^V(z)^{-1} = \expan_\tau^V(\ol{z})^\ast = \tau_z\bigl(\expan_\tau^V(-z)\bigr)$.
        
        \item \label{enu:perturbationTheory_complexCocycle} $\expan_\tau^V(z_1 + z_2) = \expan_\tau^V(z_1) \, \tau_{z_1}\bigl(\expan_\tau^V(z_2)\bigr)$.
    \end{enumerate}
\end{proposition}

\begin{para}[Unbounded perturbations]\label{para:perturbationTheory_unboundedPerturbations}
    (\cite[Sect. 3.3]{DJP03})
    Let $\MFM \subset \BO(\MH)$ be a von Neumann algebra and $\tau$ be a $W^\ast$-dynamics on $\MFM$ implemented by a self-adjoint operator $L$ on $\MH$ (\ie{}, $\tau_t(A) = \ee^{\ii t L} A \, \ee^{- \ii t L}$ for all $t \in \R$, $A \in \MFM$). Consider a self-adjoint operator $V \in \AE{\MFM}$ affiliated with the algebra $\MFM$ (\cf{} \cref{def:operatorAlgebras_affiliatedElement}). The following property shall be assumed in the sequel:
    \begin{enumerate}[label=\normalfont(A\arabic*)]
        \item \label{enu:perturbationTheory_assumptionEssSA} $L + V$ is essentially self-adjoint on $\dom(L) \cap \dom(V)$.
    \end{enumerate}
    Based on this assumption, one may construct the perturbed dynamics also for unbounded perturbations, as the next result illustrates. Its proof is an expanded version of \cite[Thm. 3.3]{DJP03}. In the following, denote the self-adjoint closure $\ol{L + V}$ (\cf{} \cref{para:operators_closure}) also by $L + V$ for simplicity.
\end{para}

\begin{proposition}\label{pro:perturbationTheory_unboundedPerturbations}
    In the setting described above, define for all $t \in \R$ and $A \in \MFM$:
    \begin{equation}\label{eq:perturbationTheory_unboundedPerturbedDynamics}
        \tau_t^V(A) \ce \ee^{\ii t (L + V)} A \, \ee^{- \ii t (L + V)} \ .
    \end{equation}
    Then $t \longmto \tau_t^V$ is a $W^\ast$-dynamics on $\MFM$ which, in the case that the operator $V \in \SE{\MFM}$ is bounded, coincides with the expression \eqref{eq:perturbationTheory_DysonRobinson}.
\end{proposition}

\begin{Proof}
    Let $A \in \MFM$ be arbitrary. By assumption \ref{enu:perturbationTheory_assumptionEssSA}, one may apply the Trotter product formula (\cref{thm:operators_TrotterProductFormula}) to obtain
    \begin{align*}
        \tau_t^V(A) &= \ee^{\ii t (L + V)} A \, \ee^{- \ii t (L + V)} = \dlim{so}{n \to \infty} \bigl(\ee^{\ii t L / n} \, \ee^{\ii t V / n}\bigr)^n A \, \bigl(\ee^{- \ii t V / n} \, \ee^{- \ii t L / n}\bigr)^n \ .
    \end{align*}
    Since $\ee^{\ii t L}$ is unitary, one can insert $n$ identity operators $\ee^{- \ii k t L / n} \, \ee^{\ii k t L / n}$, $k \in \set{1, \dotsc, n}$, in between the $n$ factors of $\ee^{\ii t L / n} \, \ee^{\ii t V / n}$ in the first term of the above limit:
    \begin{align*}
        \bigl(\ee^{\ii t L / n} \, \ee^{\ii t V / n}\bigr)^n &= \bigl(\ee^{\ii t L / n} \, \ee^{\ii t V / n}\bigr) \bigl(\ee^{-\ii t L / n} \, \ee^{\ii t L / n}\bigr) \bigl(\ee^{\ii t L / n} \, \ee^{\ii t V / n}\bigr) \bigl(\ee^{- 2 \ii t L / n} \, \ee^{2 \ii t L / n}\bigr) \\
        &\qquad \dotsm \bigl(\ee^{\ii t L / n} \, \ee^{\ii t V / n}\bigr) \bigl(\ee^{- (n-1) \ii t L / n} \, \ee^{(n - 1) \ii t L / n}\bigr) \bigl(\ee^{\ii t L / n} \, \ee^{\ii t V / n}\bigr) \bigl(\ee^{- \ii t L} \, \ee^{\ii t L}\bigr) \\[3pt]
        &= \tau_{t/n}(\ee^{\ii t V / n}) \, \tau_{2t/n}(\ee^{\ii t V / n}) \dotsm \tau_{(n-1)t/n}(\ee^{\ii t V / n}) \, \tau_{t}(\ee^{\ii t V / n}) \, \ee^{\ii t L} \\[3pt]
        &= \biggl[\,\prod_{k=1}^{n} \tau_{k t / n}(\ee^{\ii t V / n})\biggr] \, \ee^{\ii t L} \ .
    \end{align*}
    For the third term in the expression for $\tau_t^V(A)$, one may proceed similarly (but this time inserting the identity expression in a \enquote{descending} order) and write
    \begin{align*}
        \bigl(\ee^{- \ii t V / n} \, \ee^{- \ii t L / n}\bigr)^n &= \bigl(\ee^{- \ii t L} \, \ee^{\ii t L}\bigr) \bigl(\ee^{- \ii t V / n} \, \ee^{- \ii t L / n}\bigr) \bigl(\ee^{- (n-1) \ii t L / n} \, \ee^{(n - 1) \ii t L / n}\bigr) \bigl(\ee^{- \ii t V / n} \, \ee^{- \ii t L / n}\bigr) \\
        &\qquad \dotsm \bigl(\ee^{- \ii 2 t L / n} \, \ee^{\ii 2 t L / n}\bigr) \bigl(\ee^{- \ii t V / n} \, \ee^{- \ii t L / n}\bigr) \bigl(\ee^{- \ii t L / n} \, \ee^{\ii t L / n}\bigr) \bigl(\ee^{- \ii t V / n} \, \ee^{- \ii t L / n}\bigr) \\[3pt]
        &= \ee^{- \ii t L} \, \tau_{t}(\ee^{- \ii t V / n}) \, \tau_{(n-1)t/n}(\ee^{- \ii t V / n}) \dotsm \tau_{t/n}(\ee^{- \ii t V / n}) \\[3pt]
        &= \ee^{- \ii t L} \, \biggl[\,\prod_{k=0}^{n-1} \tau_{(n - k) t / n}(\ee^{- \ii t V / n})\biggr] \ .
    \end{align*}
    Therefore, combining these expression with the limit representation, it follows that
    \begin{equation*}
        \tau_t^V(A) = \dlim{so}{n \to \infty} \biggl[\,\prod_{k=1}^{n} \tau_{k t / n}(\ee^{\ii t V / n})\biggr] \, \ee^{\ii t L} A \, \ee^{- \ii t L} \, \biggl[\,\prod_{k=0}^{n-1} \tau_{(n - k) t / n}(\ee^{- \ii t V / n})\biggr] \ .
    \end{equation*}
    It holds that $\ee^{\pm \ii t V / n} \in \MFM$ because $V$ is affiliated with $\MFM$, hence \cref{pro:operatorAlgebras_characterizationAffiliation} applies. Thus, since von Neumann algebras are strongly closed, it follows that $\tau_t^V(A) \in \MFM$. This, together with \cref{exa:perturbationTheory_oneParameterGroups} \ref{enu:perturbationTheory_exaFQS}, shows that $t \longmto \tau_t^V$ is a $W^\ast$-dynamics on $\MFM$, establishing the first assertion. The second statement now follows immediately from \cref{pro:perturbationTheory_spatialPerturbation}.
\end{Proof}

\begin{para}[Expansionals for unbounded perturbations]\label{para:perturbationTheory_unboundedExpansional}
    (\cite[p. 462]{DJP03})
    Consider the situation described in \cref{para:perturbationTheory_unboundedPerturbations}, in particular, assume that \ref{enu:perturbationTheory_assumptionEssSA} holds true. Based on \cref{eq:perturbationTheory_dynamicsSpatialPerturbation}, define the Araki-Dyson expansional $t \longmto \expan_\tau^V(t)$ for an unbounded self-adjoint $V \in \AE{\MFM}$ by
    \begin{equation}\label{eq:perturbationTheory_unboundedExpansional}
        \expan_\tau^V(t) \ce \ee^{\ii t (L + V)} \, \ee^{- \ii t L} \ .
    \end{equation}
    Using again the Trotter product formula, \cref{thm:operators_TrotterProductFormula}, and performing a similar computation as in the proof of \cref{pro:perturbationTheory_unboundedPerturbations}, it follows that $\expan_\tau^V(t)$ takes the form
    \begin{align*}
        \expan_\tau^V(t) &= \dlim{so}{n \to \infty} \bigl(\ee^{\ii t L / n} \, \ee^{\ii t V / n}\bigr)^n \, \ee^{- \ii t L} = \dlim{so}{n \to \infty} \biggl[\,\prod_{k=1}^{n} \tau_{k t / n}(\ee^{\ii t V / n})\biggr] \ ,
    \end{align*}
    hence $\expan_\tau^V(t) \in \MFM$ for all $t \in \R$. Using the explicit definition \eqref{eq:perturbationTheory_unboundedExpansional} of the expansionals, one can verify the statements of \cref{pro:perturbationTheory_analyticPerturbations} also in the unbounded case \cite[Thm. 3.4]{DJP03}.
\end{para}

\subsection{Perturbation of Liouvillians}

\begin{para}[Assumptions]\label{para:perturbationTheory_perturbationLiouvillian}
    (\cite[Sect. 3.4]{DJP03})
    Let $(\MFM, \tau)$ be a $W^\ast$-dynamical system and $V \in \AE{\MFM}$ be a self-adjoint operator affiliated with $\MFM$. Assume that the algebra $\MFM$ is represented in standard form $(\MFM, \MH, J, \NPC)$, and that $L$ is the standard Liouvillian of $\tau$ (\cf{} \cref{para:perturbationTheory_standardLiouvillian}). Define
    \begin{equation*}
        L_V \ce L + V - J V J \ ,
    \end{equation*}
    and suppose that the following assumption is satisfied in addition to \ref{enu:perturbationTheory_assumptionEssSA} from \cref{para:perturbationTheory_unboundedPerturbations}:
    \begin{enumerate}[label=\normalfont(A\arabic*)]
        \setcounter{enumi}{1}
        \item \label{enu:perturbationTheory_assumptionPLiouvillian} $L_V$ is essentially self-adjoint on $\dom(L) \cap \dom(V) \cap \dom(J V J)$.
    \end{enumerate}
    The following proposition regarding the operator $L_V$ is essential for the next section. The proof is a slightly expanded version of \cite[Thm. 3.5]{DJP03}. As before, the closure $\ol{L_V}$ will be denoted simply by $L_V$.
\end{para}

\begin{proposition}\label{pro:perturbationTheory_perturbationLiouvillian}
    Under the assumptions \ref{enu:perturbationTheory_assumptionEssSA} and \ref{enu:perturbationTheory_assumptionPLiouvillian}, it follows that $L_V$ is the standard Liouvillian of the perturbed dynamics $\tau^V$ from \eqref{eq:perturbationTheory_unboundedPerturbedDynamics}.
\end{proposition}

\begin{Proof}
    It will be verified that the operator $L_V$ satisfies the two properties of \cref{pro:perturbationTheory_standardLiouvillian}. To this end, first note that from the functional calculus for the self-adjoint operator $V$ (see \cref{lem:operators_funcCalcUnitaryConj}), the fact that $\ee^{- \ii t V} \in \MFM$, and \cref{def:vonNeumann_standardForm} \ref{enu:vonNeumann_standardFormTomita}, it follows that
    \begin{equation}\label{eq:perturbationTheory_perturbationCommutant}
        \ee^{\ii t J V J} = J \, \ee^{- \ii t V} J \in \comm{\MFM} \ .
    \end{equation}
    Moreover, by definition it holds that $\dom(L + V) = \dom(L) \cap \dom(V)$, hence $\dom(L + V) \cap \dom(J V J) = \dom(L) \cap \dom(V) \cap \dom(J V J)$. Assumption \ref{enu:perturbationTheory_assumptionPLiouvillian} implies that $L_V$ is essentially self-adjoint on $\dom(L + V) \cap \dom(J V J)$. Therefore, one may apply the Trotter product formula from \cref{thm:operators_TrotterProductFormula} to write
    \begin{equation*}
        \ee^{\ii t L_V} = \dlim{so}{n \to \infty} \bigl(\ee^{\ii t (L + V) / n} \, \ee^{- \ii t J V J / n}\bigr)^n \ .
    \end{equation*}
    From this limit representation, the above observation \eqref{eq:perturbationTheory_perturbationCommutant}, and the fact that $\tau^V$, as defined in \cref{eq:perturbationTheory_unboundedPerturbedDynamics}, is indeed a $W^\ast$-dynamics (\cref{pro:perturbationTheory_unboundedPerturbations}), it follows for all $A \in \MFM$ that
    \begin{align*}
        \ee^{\ii t L_V} A \, \ee^{- \ii t L_V} &= \dlim{so}{n \to \infty} \bigl(\ee^{\ii t (L + V) / n} \, \ee^{- \ii t J V J / n}\bigr)^n A \, \bigl(\ee^{\ii t J V J / n} \, \ee^{- \ii t (L + V) / n}\bigr)^n \\
        &= \dlim{so}{n \to \infty} \biggl[\,\prod_{k=1}^{n-1} \ee^{\ii t (L + V) / n} \, \ee^{- \ii t J V J / n}\biggr] \, \ee^{\ii t (L + V) / n} \underbrace{\bigl(\ee^{- \ii t J V J / n} A \, \ee^{\ii t J V J / n}\bigr)}_{=A} \, \ee^{- \ii t (L + V) / n} \\
        &\qquad \times \biggl[\,\prod_{k=1}^{n-1} \ee^{\ii t J V J / n} \, \ee^{- \ii t (L + V) / n}\biggr] \\
        &= \dlim{so}{n \to \infty} \biggl[\,\prod_{k=1}^{n-2} \ee^{\ii t (L + V) / n} \, \ee^{- \ii t J V J / n}\biggr] \, \ee^{2 \ii t (L + V) / n} A \, \ee^{- 2 \ii t (L + V) / n} \\
        &\qquad \times \biggl[\,\prod_{k=1}^{n-2} \ee^{\ii t J V J / n} \, \ee^{- \ii t (L + V) / n}\biggr] \\
        &\vdotswithin{=} \\
        &= \dlim{so}{n \to \infty} \ee^{n \ii t (L + V) / n} A \, \ee^{- n \ii t (L + V) / n} = \ee^{\ii t (L + V)} A \, \ee^{- \ii t (L + V)} = \tau_t^V(A) \ .
    \end{align*}
    This establishes the second property of \cref{pro:perturbationTheory_standardLiouvillian}. To see also the first one, observe that since the operators $\ee^{\ii t V}$ and $\ee^{\ii t J V J}$ commute by \eqref{eq:perturbationTheory_perturbationCommutant} and the assumption on $V$, the same follows for $V$ and $J V J$ \cite[Thm. VIII.13]{RS1}, hence one can write
    \begin{equation*}
        \ee^{\ii t (V - J V J)} = \ee^{\ii t V} \ee^{- \ii t J V J} = \ee^{\ii t V} J \, \ee^{\ii t V} J = \ee^{\ii t V} j(\ee^{\ii t V}) \ .
    \end{equation*}
    Therefore, $\ee^{\ii t (V - J V J)} \, \NPC \subset \NPC$ by \cref{def:vonNeumann_standardForm} \ref{enu:vonNeumann_standardFormJConjugationNPC}. Furthermore, from \cref{pro:perturbationTheory_standardLiouvillian} it follows that $\ee^{\ii t L} \, \NPC \subset \NPC$ as $L$ is the standard Liouvillian of $\tau$. Since $\dom(L) \cap \dom(V - J V J) = \dom(L) \cap \dom(V) \cap \dom(J V J)$, the operator $L_V$ is essentially self-adjoint on $\dom(L) \cap \dom(V - J V J)$, and hence one can apply Trotter's formula \eqref{eq:operators_TrotterProductFormula} once more to write
    \begin{equation}\label{eq:perturbationTheory_perturbedStdLiouvillianTrotter}
        \ee^{\ii t L_V} = \dlim{so}{n \to \infty} \bigl(\ee^{\ii t L / n} \, \ee^{\ii t (V - J V J) / n}\bigr)^n \ .
    \end{equation}
    This, together with the previous observations regarding $\NPC$ and the fact that the latter is a closed subset of $\MH$, implies that $\ee^{\ii t L_V} \NPC \subset \NPC$. Therefore, $L_V$ also satisfies the first property in \cref{pro:perturbationTheory_standardLiouvillian}, hence it must be the unique standard Liouvillian of the dynamics $\tau^V$.
\end{Proof}

The following auxiliary lemma, whose proof is inspired by \cite[Thm. 4.70]{AJP06}, will be used below in the perturbation theory of KMS-states.

\begin{lemma}\label{lem:perturbationTheory_LJVJessSA}
    Assume that \ref{enu:perturbationTheory_assumptionEssSA} is satisfied. Then the operator $L - J V J$ is essentially self-adjoint on $\dom(L) \cap \dom(J V J)$.
\end{lemma}

\begin{Proof}
    For brevity, introduce the notation $W \ce J V J$. Then it holds that
    \begin{equation*}
        \dom(W) = J \dom(V) \ .
    \end{equation*}
    To see this, first let $\eta \in \dom(W) = \set{\xi \in \MH : J \xi \in \dom(V)}$. Since $J^2 = \id_\MH$ (\cref{pro:vonNeumann_propertiesModularData} \ref{enu:vonNeumann_propertiesJ}), it follows that $\eta = J (J \eta)$ with $J \eta \in \dom(V)$, hence $\eta \in J \dom(V)$. Conversely, if $\eta \in J \dom(V)$, then $\eta = J \xi$ for some $\xi \in \dom(V)$. Therefore, $J \eta = J^2 \xi = \xi \in \dom(V)$, and so $\eta \in \dom(W)$. This shows equality of the domains. Similarly, there holds the relation
    \begin{equation*}
        \dom(L) = J \dom(L) \ .
    \end{equation*}
    For its proof, recall the relation $J L + L J = 0$ from \cref{para:perturbationTheory_standardLiouvillian} \ref{enu:perturbationTheory_standardLiouvillianJ}. Then, if $\xi \in \dom(L)$, it follows that $J L \xi = - L J \xi$, hence $J \xi \in \dom(L)$ and $\xi = J (J \xi) \in J \dom(L)$. Conversely, $\eta \in J \dom(L)$ implies $\eta = J \xi$ for $\xi \in \dom(L)$, hence $L \eta = L J \xi = - J L \xi$. This shows that $\eta \in \dom(L)$, and thus equality of the domains. Furthermore, from the relation $J L + L J = 0$ one also obtains that
    \begin{equation*}
        L - W = - J L J - J V J = - J (L + V) J \ ,
    \end{equation*}
    where the domain of this operator is given by
    \begin{equation*}
        \dom(L - W) = \dom(L) \cap \dom(W) = J \bigl(\dom(L) \cap \dom(V)\bigr) \ ,
    \end{equation*}
    and if $\eta = J \xi \in \dom(L - W)$ for $\xi \in \dom(L + V)$, then the action of the operator can be written as $(L - W) \eta = - J (L + V) J \eta = - J (L + V) \xi$. By assumption \ref{enu:perturbationTheory_assumptionEssSA} from \cref{para:perturbationTheory_unboundedPerturbations}, the operator $L + V$ is essentially self-adjoint on $\dom(L) \cap \dom(V)$. Since $J$ is a bounded, self-adjoint, bijective operator, this shows that $L - W$ is essentially self-adjoint on $\dom(L) \cap \dom(W)$ \cite[Prop. 8.11]{Soltan18}.
\end{Proof}

The following properties of the complex expansional will also be needed in the sequel.

\begin{proposition}[\protect{\cite[Thm. 3.6]{DJP03}}]\label{pro:perturbationTheory_expansionalPerturbedLiouvillian}
    Let $V \in \MFM^\tau \cap \SE{\MFM}$. Then for all $z \in \C$,
    \begin{equation}\label{eq:perturbationTheory_expansionalPerturbedLiouvillian}
        \expan_\tau^V(z) = \ee^{\ii z L_V} \, \ee^{- \ii z (L - J V J)} \tand \ee^{\ii z L_V} = J \, \expan_\tau^V(\ol{z}) J \, \ee^{\ii z L} \, \expan_\tau^V(- z)^{-1} \ .
    \end{equation}
\end{proposition}

\section{Perturbation of KMS-States}\label{sec:perturbationTheory_perturbationKMS}

In this section, given a \emph{faithful} $(\tau, \beta)$-KMS-state $\omega \in \NS(\MFM)$ and a perturbation $V$, a state $\omega_V \in \NS(\MFM)$ shall be constructed which is a KMS-state for the \emph{perturbed} dynamics $\tau^V$. In the first step, this will be done for analytic perturbations $V \in \SE{\MFM} \cap \MFM^\tau$, and then the framework will be extended successively to bounded self-adjoint $V \in \SE{\MFM}$ and unbounded self-adjoint $V \in \AE{\MFM}$ affiliated with $\MFM$.

\subsection{Bounded Perturbations}\label{subsec:perturbationTheory_boundedPerturbationKMS}

Below, the main results for bounded perturbations are collected; they are taken from \cite[Thm. 5.1]{DJP03}, and their proofs will be given in \cref{subsec:perturbationTheory_analyticProofs,subsec:perturbationTheory_boundedProofs}. Originally, these results were obtained by \textsc{H. Araki} in 1973 \cite{Araki73b, Araki73c}. However, the method of proof employed in \cref{subsec:perturbationTheory_boundedProofs} was developed by \textsc{J. Derezi\'{n}ski}, \textsc{V. Jak\v{s}i\'{c}}, and \textsc{C.-A. Pillet} in \cite[Sect. 5]{DJP03}.

\begin{proposition}[Existence of bounded perturbations]\label{pro:perturbationTheory_boundedPerturbedVector}
    Let $\beta > 0$, $(\MFM, \MH, J, \NPC)$ be a standard form, $\tau$ be a $W^\ast$-dynamics on $\MFM$ with standard Liouvillian $L$, and $\omega$ be a faithful $(\tau, \beta)$-KMS-state with vector representative $\Omega \in \NPC$. Then for every $V \in \SE{\MFM}$, it holds that
    \begin{equation*}
        \Omega \in \dom(\ee^{- \beta (L + V) / 2}) \ .
    \end{equation*}
\end{proposition}

\begin{definition}[Perturbed KMS-state]\label{def:perturbationTheory_boundedPerturbedVector}
    In the notation of the previous proposition, define the following vector $\Omega_V \in \MH$ and state $\omega_V \in \SS(\MFM)$:
    \begin{equation*}
        \Omega_V \ce \ee^{- \beta (L + V) / 2} \, \Omega \tand \omega_V \ce \frac{1}{\norm{\Omega_V}^2} \, \braket{\Omega_V, \,\sbullet\ \Omega_V} \ .
    \end{equation*}
\end{definition}

\begin{ntheorem}[Properties of $\Omega_V$ and $\omega_V$ for bounded $V$]\label{thm:perturbationTheory_boundedPerturbation}
    Consider the situation described in \cref{pro:perturbationTheory_boundedPerturbedVector}. Then for every $V \in \SE{\MFM}$, the following properties of the perturbed vector $\Omega_V$ and the perturbed state $\omega_V$ hold true:
    \begin{enumerate}
        \item \label{enu:perturbationTheory_boundedPerturbationA} $\Omega_V \in \MP$.

        \item \label{enu:perturbationTheory_boundedPerturbationB} $\Omega_V$ is cyclic and separating for $\MFM$.
        
        \item \label{enu:perturbationTheory_boundedPerturbationC} The state $\omega_V$ is a $(\tau^V, \beta)$-KMS-state on $\MFM$.
        
        \item \label{enu:perturbationTheory_boundedPerturbationD} The modular operator of $(\MFM, \Omega_V)$ satisfies $\log(\Delta_{\Omega_V}) = - \beta L_V$.
        
        \item \label{enu:perturbationTheory_boundedPerturbationE} $\log(\Delta_{\Omega_V, \Omega}) = \log(\Delta_{\Omega}) - \beta V$.
        
        \item \label{enu:perturbationTheory_boundedPerturbationF} $\log(\Delta_{\Omega, \Omega_V}) = \log(\Delta_{\Omega_V}) + \beta V$.
        
        \item \label{enu:perturbationTheory_REidentity1} $\MS_\MFM^\mathrm{std}(\omega, \omega_V) = \log\bigl(\norm{\Omega_V}^2\bigr) + \beta \omega(V)$.
        
        \item \label{enu:perturbationTheory_REidentity2} $\MS_\MFM^\mathrm{std}(\omega_V, \omega) = - \log\bigl(\norm{\Omega_V}^2\bigr) - \beta \omega_V(V)$.
    \end{enumerate}
\end{ntheorem}

Another very important identity involving the perturbed vector $\Omega_V$ is the famous \emph{Golden-Thompson inequality} which is a crucial element for the proof of the above results.

\begin{proposition}[Golden-Thompson inequality]\label{pro:perturbationTheory_boundedGT}
    Consider the situation of \cref{pro:perturbationTheory_boundedPerturbedVector}. For every $V \in \SE{\MFM}$, the following inequality holds true:
    \begin{equation}\label{eq:perturbationTheory_boundedGT}
        \norm{\Omega_V} \le \norm{\ee^{- \beta V / 2} \Omega} \ .
    \end{equation}
\end{proposition}

\begin{remark}\label{rem:perturbationTheory_GT}
    The Golden-Thompson inequality, which was first discussed in \cite{Golden65, Symanzik65, Thompson65} and extended to the form \eqref{eq:perturbationTheory_boundedGT} by \textsc{H. Araki} \cite[Thm. 2]{Araki73b}, is an important tool in statistical mechanics \cite{Bebiano04}, \cite[Sect. 8.1]{Simon05}. It is usually formulated as follows: for two positive self-adjoint operators $A$ and $B$ on a Hilbert space $\MH$ such that $B$ is relatively $A$-bounded with $A$-bound strictly less than one \cite[Def. 8.1]{Schmüdgen12}, and such that $\ee^{-A}, \ee^{-B} \in \NO(\MH)$, there holds \cite[Thm. 1]{BreiteneckerGrümm72}
    \begin{equation*}
        \tr\bigl(\ee^{- (A + B)}\bigr) \le \tr(\ee^{-A} \, \ee^{-B}) \ .
    \end{equation*}
\end{remark}

The plan for the proof of the above results, which will be followed in the next two subsections, can be summarized as follows:
\begin{equation*}
    \begin{tikzcd}[row sep=large]
        \text{\ref{pro:perturbationTheory_boundedPerturbedVector} for analytic $V$} \arrow[Rightarrow]{r} & \text{\ref{thm:perturbationTheory_boundedPerturbation} for analytic $V$} \arrow[Rightarrow]{r} & \text{\ref{pro:perturbationTheory_boundedGT} for analytic $V$} \arrow[Rightarrow]{d} \\
        \text{\ref{pro:perturbationTheory_boundedGT} for bounded $V$} & \arrow[Rightarrow]{l} \text{\ref{thm:perturbationTheory_boundedPerturbation} for bounded $V$} & \arrow[Rightarrow]{l} \text{\ref{pro:perturbationTheory_boundedPerturbedVector} for bounded $V$}
    \end{tikzcd}
\end{equation*}

\subsection{Proof for Analytic Perturbations}\label{subsec:perturbationTheory_analyticProofs}

Let $V \in \SE{\MFM} \cap \MFM^\tau$, and let $\tau^V$ and $\expan_\tau^V$ be the perturbed dynamics and Araki-Dyson expansional from \cref{para:perturbationTheory_analyticPerturbations}, respectively. Recall that according to \cref{pro:perturbationTheory_perturbationLiouvillian}, $L_V = L + V - J V J$ is the standard Liouvillian of $\tau^V$. The following proofs are taken from \cite[Sect. 5.2]{DJP03}, except the one of \cref{thm:perturbationTheory_boundedPerturbation} \ref{enu:perturbationTheory_boundedPerturbationF} which is independent.

\begin{Proof}[of \cref{pro:perturbationTheory_boundedPerturbedVector}]
    The standard Liouvillian $L$ of $\tau$ satisfies $L \Omega = 0$, \cf{} \cref{para:perturbationTheory_standardLiouvillian} \ref{enu:perturbationTheory_standardLiouvillianVectorRepr}, hence $\ee^{- \ii t L} \Omega = \Omega$ for all $t \in \R$ by \cref{lem:operators_funcCalcEigenvalue}. Thus, using the spatial representation $\expan_\tau^V(t) = \ee^{\ii t (L + V)} \, \ee^{- \ii t L}$ of the Araki-Dyson expansional from \cref{eq:perturbationTheory_dynamicsSpatialPerturbation}, one obtains
    \begin{equation*}
        \expan_\tau^V(t) \Omega = \ee^{\ii t (L + V)} \, \ee^{- \ii t L} \Omega = \ee^{\ii t (L + V)} \Omega \ .
    \end{equation*}
    According to \cref{para:perturbationTheory_analyticPerturbations}, the expansional $\expan_\tau^V(t)$ possesses an analytic continuation to an  entire function $\C \owns z \longmto \expan_\tau^V(z) \in \MFM$. Hence, it follows that $\Omega \in \dom(\ee^{\ii z (L + V)})$ and $\expan_\tau^V(z) \Omega = \ee^{\ii z (L + V)} \Omega$ for all $z \in \C$. In particular, choosing $z = \ii \beta / 2$, this shows that
    \begin{equation}\label{eq:perturbationTheory_perturbedVectorFromExpansional}
        \Omega_V \ce \expan_\tau^V(\ii \beta / 2) \, \Omega = \ee^{- \beta (L + V) / 2} \, \Omega
    \end{equation}
    is a well-defined vector in $\MH$, so it must be an element of the domain of the operator $\ee^{- \beta (L + V) / 2}$.
\end{Proof}

\begin{Proof}[of \cref{thm:perturbationTheory_boundedPerturbation}]
    \tAd{} \ref{enu:perturbationTheory_boundedPerturbationA}. Consider the element $\expan_\tau^V(\ii \beta / 2) \in \MFM^\tau$. From \cref{pro:perturbationTheory_analyticPerturbations} \ref{enu:perturbationTheory_complexExpansionalInverse} and \ref{enu:perturbationTheory_complexCocycle}, it follows by choosing $z_1 = z_2 = \ii \beta / 4$ that
    \begin{equation*}
        \expan_\tau^V(\ii \beta / 2) = \expan_\tau^V(\ii \beta / 4 + \ii \beta / 4) = \expan_\tau^V(\ii \beta / 4) \, \tau_{\ii \beta / 4}\bigl(\expan_\tau^V(\ii \beta / 4)\bigr) = \expan_\tau^V(\ii \beta / 4) \, \tau_{\ii \beta / 2}\bigl(\expan_\tau^V(\ii \beta / 4)^\ast\bigr) \ .
    \end{equation*}
    Combining this identity with \cref{eq:perturbationTheory_perturbedVectorFromExpansional}, one obtains
    \begin{align*}
        \Omega_V &= \expan_\tau^V(\ii \beta / 4) \, \tau_{\ii \beta / 2}\bigl(\expan_\tau^V(\ii \beta / 4)^\ast\bigr) \Omega \\
        &= \expan_\tau^V(\ii \beta / 4) \, \ee^{- \beta L / 2} \, \expan_\tau^V(\ii \beta / 4)^\ast \, \Omega \\
        &= \expan_\tau^V(\ii \beta / 4) \, J \, \expan_\tau^V(\ii \beta / 4) J \Omega \ .
    \end{align*}
    In the second step, it was used that $\ee^{\beta L / 2} \Omega = \Omega$, and in the third, \cref{pro:perturbationTheory_LiouvillianKMS} was employed. This shows that $\Omega_V = \expan_\tau^V(\ii \beta / 4) \, j\bigl(\expan_\tau^V(\ii \beta / 4)\bigr) \Omega$, hence $\Omega_V \in \NPC$ according to \cref{pro:vonNeumann_propertiesNPC} \ref{enu:vonNeumann_jNPC}.

    \tAd{} \ref{enu:perturbationTheory_boundedPerturbationB}. $\expan_\tau^V(\ii \beta / 2)$ is an invertible element in $\MFM$ by \cref{pro:perturbationTheory_analyticPerturbations} \ref{enu:perturbationTheory_complexExpansionalInverse}. Therefore, since $\Omega$ is cyclic and separating for $\MFM$ by assumption (recall that $\omega$ is faithful, hence \cref{rem:vonNeumann_cyclicSeparatingFaithful} applies), it follows from \cref{eq:perturbationTheory_perturbedVectorFromExpansional} that $\Omega_V$ is cyclic and separating, too.

    \tAd{} \ref{enu:perturbationTheory_boundedPerturbationC} \& \ref{enu:perturbationTheory_boundedPerturbationD}. From the second identity in \eqref{eq:perturbationTheory_expansionalPerturbedLiouvillian}, one obtains by choosing $z = \ii \beta / 2$ that
    \begin{equation*}
        \ee^{- \beta L_V / 2} = J \, \expan_\tau^V(- \ii \beta / 2) \, J \, \ee^{- \beta L / 2} \, \expan_\tau^V(- \ii \beta / 2)^{-1} \ .
    \end{equation*}
    \cref{pro:perturbationTheory_LiouvillianKMS} implies that $\MFM \Omega \subset \dom(\ee^{- \beta L / 2})$ since $\omega = \omega_\Omega$ is a $(\tau, \beta)$-KMS-state, hence $\MFM \Omega \subset \dom(\ee^{- \beta L_V / 2})$ as well by the above identity. Since $\MFM \Omega_V = \MFM \Omega$ according to \eqref{eq:perturbationTheory_perturbedVectorFromExpansional}, one can conclude that $\MFM \Omega_V \subset \dom(\ee^{- \beta L_V / 2})$. Furthermore, for all $A \in \MFM$ one computes
    \begin{spreadlines}{4pt}
        \begin{align*}
            \ee^{- \beta L_V / 2} A \Omega_V &= J \, \expan_\tau^V(- \ii \beta / 2) \, J \, \ee^{- \beta L / 2} \, \bigl(\expan_\tau^V(- \ii \beta / 2)^{-1} A \, \expan_\tau^V(\ii \beta / 2)\bigr) \Omega \\
            &= J \, \expan_\tau^V(- \ii \beta / 2) \, J^2 \bigl(\expan_\tau^V(- \ii \beta / 2)^{-1} A \, \expan_\tau^V(\ii \beta / 2)\bigr)^\ast \Omega \\
            &= J \, \expan_\tau^V(- \ii \beta / 2) \, \expan_\tau^V(- \ii \beta / 2)^{-1} A^\ast \, \expan_\tau^V(\ii \beta / 2) \, \Omega \\
            &= J A^\ast \Omega_V \ ,
        \end{align*}
    \end{spreadlines}
    where in the second line \cref{eq:perturbationTheory_LiouvillianKMS}, and in the third line \cref{pro:perturbationTheory_analyticPerturbations} \ref{enu:perturbationTheory_complexExpansionalInverse} was used. This shows that both conditions of \cref{pro:perturbationTheory_LiouvillianKMS} are satisfied, hence $\omega_V$ is a $(\tau^V, \beta)$-KMS-state. In light of \ref{enu:perturbationTheory_boundedPerturbationB} proved above, the aforementioned proposition also shows that $\Delta_{\Omega_V} = \ee^{- \beta L_V}$.

    \tAd{} \ref{enu:perturbationTheory_boundedPerturbationE}. Using the definition \eqref{eq:vonNeumann_relativeTomitaOperator} of the relative Tomita operator $S_{\Omega_V, \Omega}$ and the definition \eqref{eq:vonNeumann_preclosedTomitaOperator} of the Tomita operator $S_\Omega$ (noting that both are well-defined because $\Omega, \Omega_V \in \NPC$ are cyclic and separating), as well as the definition \eqref{eq:perturbationTheory_perturbedVectorFromExpansional} of the perturbed vector, one obtains for all $A \in \MFM$:
    \begin{equation*}
        S_{\Omega_V, \Omega} A \Omega = A^\ast \Omega_V = A^\ast \, \expan_\tau^V(\ii \beta / 2) \, \Omega = \bigl(\expan_\tau^V(\ii \beta / 2)^\ast A\bigr)^\ast \Omega = S_\Omega \, \expan_\tau^V(\ii \beta / 2)^\ast A \Omega \ .
    \end{equation*}
    This implies that $S_{\Omega_V, \Omega} = S_\Omega \, \expan_\tau^V(\ii \beta / 2)^\ast$. Therefore, using the definitions \eqref{eq:vonNeumann_relativeModularOperator} and \eqref{eq:vonNeumann_modularOperator} of the relative modular operator $\Delta_{\Omega_V, \Omega}$ and modular operator $\Delta_\Omega$, it follows that
    \begin{align*}
        \Delta_{\Omega_V, \Omega} &= \bigl(S_\Omega \, \expan_\tau^V(\ii \beta / 2)^\ast\bigr)^\ast S_\Omega \, \expan_\tau^V(\ii \beta / 2)^\ast \\
        &= \expan_\tau^V(\ii \beta / 2) \, S_\Omega^\ast S_\Omega \, \expan_\tau^V(\ii \beta / 2)^\ast \\
        &= \expan_\tau^V(\ii \beta / 2) \, \Delta_\Omega \, \expan_\tau^V(\ii \beta / 2)^\ast \ .
    \end{align*}
    Introducing the identity $\Delta_\Omega = \ee^{- \beta L}$ from \cref{pro:perturbationTheory_LiouvillianKMS} and writing, as before, $\expan_\tau^V(\ii \beta / 2) = \ee^{- \beta (L + V) / 2} \, \ee^{\beta L / 2}$ according to \cref{pro:perturbationTheory_analyticPerturbations} \ref{enu:perturbationTheory_complexSpatialExpansional}, one furthermore obtains
    \begin{align*}
        \Delta_{\Omega_V, \Omega} &= \bigl(\ee^{- \beta (L + V) / 2} \, \ee^{\beta L / 2} \, \ee^{- \beta L / 2}\bigr) \bigl(\ee^{- \beta L / 2} \, \ee^{\beta L / 2} \, \ee^{- \beta (L + V) / 2}\bigr) = \ee^{- \beta (L + V)} \ .
    \end{align*}
    Thus, taking the logarithm of this expression and inserting again \cref{eq:perturbationTheory_modularOperatorLiouvillian} yields the assertion:
    \begin{align*}
        \log(\Delta_{\Omega_V, \Omega}) &= - \beta (L + V) = \log(\Delta_\Omega) - \beta V \ .
    \end{align*}

    \tAd{} \ref{enu:perturbationTheory_boundedPerturbationF}. Proceeding similarly as in the previous part of the proof, it follows from \cref{eq:perturbationTheory_perturbedVectorFromExpansional} and \cref{pro:perturbationTheory_analyticPerturbations} \ref{enu:perturbationTheory_complexExpansionalInverse} that
    \begin{align*}
        S_{\Omega, \Omega_V} A \Omega_V &= A^\ast \Omega = A^\ast \, \expan_\tau^V(\ii \beta / 2)^{-1} \Omega_V = A^\ast \, \expan_\tau^V(- \ii \beta / 2)^\ast \Omega_V = S_{\Omega_V} \, \expan_\tau^V(- \ii \beta / 2) A \Omega_V
    \end{align*}
    for all $A \in \MFM$, hence $S_{\Omega, \Omega_V} = S_{\Omega_V} \, \expan_\tau^V(- \ii \beta / 2)$. Therefore, using the first identity in \cref{eq:perturbationTheory_expansionalPerturbedLiouvillian} and the relation $\Delta_{\Omega_V} = \ee^{- \beta L_V}$ proved above in \ref{enu:perturbationTheory_boundedPerturbationD}, one computes
    \begin{align*}
        \Delta_{\Omega, \Omega_V} &= \expan_\tau^V(- \ii \beta / 2)^\ast \, S_{\Omega_V}^\ast S_{\Omega_V} \, \expan_\tau^V(- \ii \beta / 2) \\
        &= \bigl(\ee^{\beta L_V / 2} \, \ee^{- \beta (L - J V J) / 2}\bigr)^\ast \Delta_{\Omega_V} \bigl(\ee^{\beta L_V / 2} \, \ee^{- \beta (L - J V J) / 2}\bigr) \\
        &= \ee^{- \beta (L - J V J) / 2} \, \ee^{\beta L_V / 2} \, \ee^{- \beta L_V} \, \ee^{\beta L_V / 2} \, \ee^{- \beta (L - J V J) / 2} \\
        &= \ee^{- \beta (L - J V J)} \ .
    \end{align*}
    Since the perturbed Liouvillian $L_V$ is given by $L_V = L + V - J V J$, one has $L - J V J = L_V - V$, and this operator is essentially self-adjoint on its natural domain by \cref{lem:perturbationTheory_LJVJessSA}, thus $\Delta_{\Omega, \Omega_V} = \ee^{- \beta (L_V - V)}$; taking the logarithm yields the desired expression:
    \begin{equation*}
        \log(\Delta_{\Omega, \Omega_V}) = - \beta L_V + \beta V = \log(\Delta_{\Omega_V}) + \beta V \ .
    \end{equation*}

    \tAd{} \ref{enu:perturbationTheory_REidentity1} \& \ref{enu:perturbationTheory_REidentity2}. Define the operator $\wt{V} \ce V + \inv{\beta} \log\bigl(\norm{\Omega_V}^2\bigr) \, \id_\MH$ which is still self-adjoint and $\tau$-entire. According to \cref{eq:perturbationTheory_perturbedVectorFromExpansional}, the $\wt{V}$-perturbation of the vector $\Omega$ is given by
    \begin{equation*}
        \Omega_{\wt{V}} = \ee^{- \beta (L + \wt{V}) / 2} \, \Omega = \ee^{- \log \norm{\Omega_V}} \, \ee^{- \beta (L + V)} \Omega = \frac{\Omega_V}{\norm{\Omega_V}} \ ,
    \end{equation*}
    and hence the corresponding state $\omega_{\wt{V}}$ (see \cref{def:perturbationTheory_boundedPerturbedVector}) takes the form
    \begin{equation*}
        \omega_{\wt{V}}(A) = \frac{1}{\vnorm[\big]{\Omega_{\wt{V}}}^2} \, \vbraket[\big]{\Omega_{\wt{V}}, A \, \Omega_{\wt{V}}} = \braket{\Omega_{\wt{V}}, A \, \Omega_{\wt{V}}} = \frac{1}{\norm{\Omega_V}^2} \, \braket{\Omega_V, A \, \Omega_V} = \omega_V(A)
    \end{equation*}
    for all $A \in \MFM$. Employing the identity $\log(\Delta_{\Omega_V, \Omega}) = \log(\Delta_{\Omega}) - \beta V$ derived above in \ref{enu:perturbationTheory_boundedPerturbationE}, one obtains a corresponding relation for $\log(\Delta_{\Omega_{\tilde{V}}, \Omega})$:
    \begin{equation*}
        \log(\Delta_{\Omega_{\tilde{V}}, \Omega}) = \log(\Delta_{\Omega}) - \beta \wt{V} = \log(\Delta_{\Omega}) - \beta V - \log\bigl(\norm{\Omega_V}^2\bigr) \ .
    \end{equation*}
    This and the fact that $\omega_V = \omega_{\wt{V}}$ imply the following identity for the relative entropy $\MS_\MFM^\mathrm{std}(\omega, \omega_V)$:
    \begin{align*}
        \MS_\MFM^\mathrm{std}(\omega, \omega_V) &= \MS_\MFM^\mathrm{std}(\omega, \omega_{\wt{V}}) = - \braket{\Omega, \log (\Delta_{\Omega_{\tilde{V}}, \Omega}) \Omega} \\[3pt]
        &= - \braket{\Omega, \log(\Delta_{\Omega}) \Omega} + \beta \braket{\Omega, V \Omega} + \log\bigl(\norm{\Omega_V}^2\bigr) \braket{\Omega, \Omega} \\[3pt]
        &= \beta \omega(V) + \log\bigl(\norm{\Omega_V}^2\bigr) \ .
    \end{align*}
    To get from the second to the third line, it was used that $\Delta_\Omega \Omega = \Omega$ (\cref{pro:vonNeumann_propertiesModularData} \ref{enu:vonNeumann_actionModDatOmega}), hence $\log(\Delta_\Omega) \Omega = 0$ by \cref{lem:operators_funcCalcEigenvalue}, and that the vector $\Omega$ is normalized by assumption. This proves the first identity. For the second one, note that with $\wt{V}$ as above, it follows from \ref{enu:perturbationTheory_boundedPerturbationF} that
    \begin{equation*}
        \log (\Delta_{\Omega, \Omega_{\tilde{V}}}) = \log(\Delta_{\Omega_{\tilde{V}}}) + \beta \wt{V} = \log(\Delta_{\Omega_{\tilde{V}}}) + \beta V + \log\bigl(\norm{\Omega_V}^2\bigr) \ .
    \end{equation*}
    Therefore, the relative entropy takes the form
    \begin{align*}
        \MS_\MFM^\mathrm{std}(\omega_V, \omega) &= \MS_\MFM^\mathrm{std}(\omega_{\wt{V}}, \omega) = - \vbraket[\big]{\Omega_{\wt{V}}, \log (\Delta_{\Omega, \Omega_{\tilde{V}}}) \Omega_{\wt{V}}} \\[3pt]
        &= - \vbraket[\big]{\Omega_{\wt{V}}, \log(\Delta_{\Omega_{\tilde{V}}}) \Omega_{\wt{V}}} - \beta \vbraket[\big]{\Omega_{\wt{V}}, V \Omega_{\wt{V}}} - \log\bigl(\norm{\Omega_V}^2\bigr) \vbraket[\big]{\Omega_{\wt{V}}, \Omega_{\wt{V}}} \\[3pt]
        &= - \beta \omega_{\wt{V}}(V) - \log\bigl(\norm{\Omega_V}^2\bigr) \ . \tag*{\qedhere}
    \end{align*}
\end{Proof}

\begin{Proof}[of \cref{pro:perturbationTheory_boundedGT}]
    Denote by $\MFN$ the commutative von Neumann subalgebra of $\MFM$ generated by the perturbation $V$. According to monotonicity of the relative entropy under restriction of the functionals to a subalgebra (\cref{cor:relativeEntropy_monotonicitySubalgebra}), it holds that
    \begin{equation*}
        \MS_{\MFN}^\mathrm{std}\bigl(\omega_V\big|_{\MFN}, \omega\big|_{\MFN}\bigr) \le \MS_\MFM^\mathrm{std}(\omega_V, \omega) \ .
    \end{equation*}
    Moreover, because $- \beta V \in \MFN \cap \comm{\MFN}$ is a self-adjoint element in the center of $\MFN$ (since $\MFN$ is Abelian, $\MFN \subset \comm{\MFN}$), one can employ the inequality \eqref{eq:relativeEntropy_inequalityCenter} from \cref{lem:relativeEntropy_inequalityCenter} to obtain
    \begin{equation*}
        \MS_{\MFN}^\mathrm{std}\bigl(\omega_V\big|_{\MFN}, \omega\big|_{\MFN}\bigr) + \beta \omega_V(V) \ge - \log \omega(\ee^{- \beta V}) \ .
    \end{equation*}
    Combining the above two inequalities with \cref{thm:perturbationTheory_boundedPerturbation} \ref{enu:perturbationTheory_REidentity2} for analytic perturbations, it follows that
    \begin{spreadlines}{4pt}
        \begin{align*}
            \log\bigl(\norm{\Omega_V}^2\bigr) &= - \MS_\MFM^\mathrm{std}(\omega_V, \omega) - \beta \omega_V(V) \\
            &\le - \MS_{\MFN}^\mathrm{std}\bigl(\omega_V\big|_{\MFN}, \omega\big|_{\MFN}\bigr) - \beta \omega_V(V) \\
            &\le \log \omega(\ee^{- \beta V}) \\
            &= \log\bigl(\norm{\ee^{- \beta V / 2} \Omega}^2\bigr) \ .
        \end{align*}
    \end{spreadlines}
    Exponentiating this inequality and taking the square root, \cref{eq:perturbationTheory_boundedGT} follows.
\end{Proof}

\subsection{Proof for Bounded Perturbations}\label{subsec:perturbationTheory_boundedProofs}

Now, following \cite[Sect. 5.3]{DJP03}, the assertions of \cref{subsec:perturbationTheory_boundedPerturbationKMS} shall be proved for $V \in \SE{\MFM}$. The general strategy of the proofs consists in taking a sequence $(V_n)_{n \in \N} \subset \MFM^\tau \cap \SE{\MFM}$ converging strongly to $V$, which exists by \cref{thm:perturbationTheory_analyticApproximation}, and then using certain limit arguments to reduce the problem to the situation of \cref{subsec:perturbationTheory_analyticProofs}. In this spirit, first an auxiliary result is needed which is stated in \cite[Lem. 5.1]{DJP03}.

\begin{lemma}\label{lem:perturbationTheory_boundedConvergenceLiouvillian}
    Let $(V_n)_{n \in \N} \subset \MFM^\tau \cap \SE{\MFM}$ be a sequence of self-adjoint, $\tau$-entire elements converging strongly to $V \in \SE{\MFM}$. Then the following two properties are satisfied:
    \begin{enumerate}
        \item \label{enu:perturbationTheory_convergenceUnperturbedL} $L + V_n \to L + V$ in the strong resolvent sense.
        
        \item \label{enu:perturbationTheory_convergencePerturbedL} $L_{V_n} \to L_V$ in the strong resolvent sense.
    \end{enumerate}
\end{lemma}

\begin{Proof}
    Only \ref{enu:perturbationTheory_convergenceUnperturbedL} shall be proved in detail; the argument to establish \ref{enu:perturbationTheory_convergencePerturbedL} is similar (see \cref{lem:perturbationTheory_unboundedConvergenceLiouvillian} below). Using the second resolvent identity \cite[p. 34]{Oliveira09}, \cite[Thm. 5.13]{Weidmann80}, it follows for all $\xi \in \MH$ that
    \begin{align*}
        \vnorm[\big]{(L + V_n - \ii)^{-1} \xi - (L + V - \ii)^{-1} \xi} &= \vnorm[\big]{(L + V_n - \ii)^{-1} (V - V_n) (L + V - \ii)^{-1} \xi} \\[3pt]
        &\le \vnorm[\big]{(L + V_n - \ii)^{-1}}_\mop \, \vnorm[\big]{(V - V_n) (L + V - \ii)^{-1} \xi} \ .
    \end{align*}
    Since $L + V_n$ is self-adjoint, it holds that $\norm{(L + V_n - \ii)^{-1}}_\mop \le 1$ \cite[Thm. 2.2.17]{Oliveira09}. Furthermore, as $V_n \to V$ in the strong operator topology by assumption, the right-hand side in the above inequality tends to zero as $n \to + \infty$, showing that $(L + V_n - \ii)^{-1} \to (L + V - \ii)^{-1}$ strongly. By \cref{def:operators_SRConvergence}, this implies $L + V_n \to L + V$ in the strong resolvent sense.
\end{Proof}

\begin{Proof}[of \cref{pro:perturbationTheory_boundedPerturbedVector}]
    From \cref{thm:perturbationTheory_analyticApproximation}, it follows that there exists a sequence $(V_n)_{n \in \N} \subset \SE{\MFM} \cap \MFM^\tau$ of self-adjoint, $\tau$-entire elements such that $V_n \to V$ strongly. This implies $L + V_n \to L + V$ in the strong resolvent sense by \cref{lem:perturbationTheory_boundedConvergenceLiouvillian} \ref{enu:perturbationTheory_convergenceUnperturbedL}. Moreover, it holds that
    \begin{equation}\label{eq:perturbationTheory_GTconvergence}
        \lim_{n \to \infty} \ee^{- \beta V_n / 2} \Omega = \ee^{- \beta V / 2} \Omega
    \end{equation}
    since $V_n \to V$ also in the strong resolvent sense (\cref{pro:operators_SOimpliesSR}), hence \cref{pro:operators_characterizationSRConvergence} applies. Therefore, there exists a constant $C > 0$ such that for all $n \in \N$,
    \begin{equation*}
        \vnorm[\big]{\ee^{- \beta V_n / 2} \Omega} \le C \ .
    \end{equation*}
    Furthermore, using the Golden-Thompson inequality \eqref{eq:perturbationTheory_boundedGT} for analytic perturbations,
    \begin{equation*}
        \norm{\Omega_{V_n}} \le \vnorm[\big]{\ee^{-\beta V_n / 2} \Omega} \ .
    \end{equation*}
    Combining these two inequalities yields $\norm{\Omega_{V_n}} \le C$. From \cref{pro:operators_SRimpliesW} applied to $T_n \equiv \ee^{- \beta (L + V_n) / 2}$, $T \equiv \ee^{- \beta (L + V) / 2}$, $\Omega_n \equiv \Omega$, noting that $T_n \to T$ in the strong resolvent sense by \cref{lem:perturbationTheory_boundedConvergenceLiouvillian} \ref{enu:perturbationTheory_convergenceUnperturbedL} and \cref{lem:operators_SRimpliesFSR} and $\norm{\ee^{- \beta (L + V_n) / 2} \Omega} \le C$ according to \eqref{eq:perturbationTheory_perturbedVectorFromExpansional}, one obtains $\Omega \in \dom(\ee^{- \beta (L + V) / 2})$ and
    \begin{equation*}
        \dlim{w}{n \to \infty} \ee^{- \beta (L + V_n) / 2} \Omega = \ee^{- \beta (L + V) / 2} \Omega \ . \tag*{\qedhere}
    \end{equation*}
\end{Proof}

\begin{Proof}[of \cref{thm:perturbationTheory_boundedPerturbation}]
    As in the proof of \cref{pro:perturbationTheory_boundedPerturbedVector}, one may assume, by virtue of \cref{thm:perturbationTheory_analyticApproximation}, that there exists a sequence $(V_n)_{n \in \N} \subset \SE{\MFM} \cap \MFM^\tau$ of self-adjoint, $\tau$-entire elements such that $V_n \to V$ in the strong operator topology. Note that, as shown above,
    \begin{equation}\label{eq:perturbationTheory_weakConvergenceOmegaV}
        \dlim{w}{n \to \infty} \Omega_{V_n} = \Omega_V \ .
    \end{equation}
    
    \tAd{} \ref{enu:perturbationTheory_boundedPerturbationA}. According to \cref{subsec:perturbationTheory_analyticProofs}, it holds that $\Omega_{V_n} \in \NPC$ for all $n \in \N$. Since 
    $\NPC$ is weakly closed by \cref{pro:vonNeumann_propertiesNPC} \ref{enu:vonNeumann_NPCclosedConvex} and Mazur's theorem \cite[Cor. 2.11]{Voigt20}, \eqref{eq:perturbationTheory_weakConvergenceOmegaV} implies the assertion.

    \tAd{} \ref{enu:perturbationTheory_boundedPerturbationB}. Define 
    $P \ce \id_\MH - \ssupp(\omega_V) \in \PO(\MFM)$. It holds that $P \Omega_V = 0$ and $\tau_t^V(P) = P$. To see the latter identity, note that $\tau^V$ is implemented by $\ee^{\ii t L_V}$ (\cref{pro:perturbationTheory_perturbationLiouvillian}), and that $\ee^{- \ii t L_V} \Omega_V = \dlim{w}{n \to \infty} \ee^{- \ii t L_{V_n}} \Omega_{V_n}$ according to \cref{lem:perturbationTheory_boundedConvergenceLiouvillian} \ref{enu:perturbationTheory_convergencePerturbedL}, \cref{pro:operators_characterizationSRConvergence}, and \cref{eq:perturbationTheory_weakConvergenceOmegaV}. Furthermore, since $L_{V_n} \Omega_{V_n} = 0$ by \cref{thm:perturbationTheory_boundedPerturbation} \ref{enu:perturbationTheory_boundedPerturbationD} for analytic perturbations and \cref{pro:vonNeumann_propertiesModularData} \ref{enu:vonNeumann_actionModDatOmega}, it follows from \cref{lem:operators_funcCalcEigenvalue} that $\ee^{- \ii t L_V} \Omega_V = \Omega_V$. Therefore, $\omega_V$ is $\tau^V$-invariant, and hence $\tau_t^V(\ssupp(\omega_V)) = \ssupp(\omega_V)$ \cite[p. 454]{DJP03}. For $z \in \C$, define a vector
    \begin{equation*}
        \Omega(z) \ce \ee^{- z (L + V)} \Omega \in \MH \ .
    \end{equation*}
    One can show that the function $z \longmto \Omega(z)$ is analytic in the strip $0 < \Re(z) < \beta / 2$ and norm continuous on its closure \cite[Prop. A.1]{DJP03}, \cite[Lem. 3]{Araki73b}. Furthermore, $\Omega(\beta / 2) = \Omega_V$ by \cref{def:perturbationTheory_boundedPerturbedVector} and
    \begin{equation*}
        \ee^{\ii t (L + V)} P \Omega(\ii t + \beta / 2) = \ee^{\ii t (L + V)} P \, \ee^{- \ii t (L + V)} \Omega(\beta / 2) = \tau_t^V(P) \Omega_V = P \Omega_V = 0 \ .
    \end{equation*}
    This shows that $P \Omega(\ii t + \beta / 2) = 0$ for all $t \in \R$. The three-line theorem \cite[Lem. 1.1.2]{BerghLöfström76}, \cite[Thm. 3.7]{Conway78} implies $P \Omega(z) = 0$ for every $z$ in the strip $0 \le \Re(z) \le \beta / 2$. In particular, $P \Omega(0) = P \Omega = 0$. Since $\Omega$ is separating for $\MFM$ by the assumption on $\omega$ (\cf{} \cref{rem:vonNeumann_cyclicSeparatingFaithful}), it follows that $P = 0$. Therefore, $\ssupp(\omega_V) = \id_\MH$, and hence $\Omega_V$ is separating for $\MFM$ as well, see \cref{pro:operatorAlgebras_normalFaithful} and \cref{cor:operatorAlgebras_faithfulSeparating}. Finally, since $\Omega_V \in \NPC$ by assertion \ref{enu:perturbationTheory_boundedPerturbationA}, using \cref{lem:vonNeumann_vectorsNPC} one concludes that $\Omega_V$ is cyclic for $\MFM$.

    \tAd{} \ref{enu:perturbationTheory_boundedPerturbationC}. It was shown in \cref{subsec:perturbationTheory_analyticProofs} that for every $n \in \N$, the state $\omega_{V_n}$ is a $(\tau^{V_n}, \beta)$-KMS-state. In particular, this implies that $\omega_{V_n}$ is $\tau^{V_n}$-invariant (see \cite[Prop. 5.3.3]{BR2}, \cite[Lem. 2.12]{Hiai21} or \cite[Thm. 5.3]{AJP06}). Moreover, $L_{V_n}$ is the standard Liouvillian of $\tau^{V_n}$ (\cref{pro:perturbationTheory_perturbationLiouvillian}), and it holds that $L_{V_n} \to L_V$ in the strong resolvent sense according to \cref{lem:perturbationTheory_boundedConvergenceLiouvillian} \ref{enu:perturbationTheory_convergencePerturbedL}. Therefore, all assumptions of \cref{pro:perturbationTheory_convergence} are satisfied; it implies that $\omega_{\Omega_V / \norm{\Omega_V}} = \omega_V$ is a $(\tau^V, \beta)$-KMS-state.

    \tAd{} \ref{enu:perturbationTheory_boundedPerturbationD}. As in the analytic case, the assertion follows from the statements \ref{enu:perturbationTheory_boundedPerturbationB} and \ref{enu:perturbationTheory_boundedPerturbationC}, as well as from \cref{pro:perturbationTheory_LiouvillianKMS}.

    \tAd{} \ref{enu:perturbationTheory_boundedPerturbationE} \& \ref{enu:perturbationTheory_boundedPerturbationF}. For all $n \in \N$, it holds that $\log(\Delta_{\Omega_{V_n}, \Omega}) = \log(\Delta_{\Omega}) - \beta V_n = - \beta (L + V_n)$ as shown in \cref{subsec:perturbationTheory_analyticProofs}. From \cref{lem:perturbationTheory_boundedConvergenceLiouvillian} \ref{enu:perturbationTheory_convergenceUnperturbedL}, it follows that $L + V_n \to L + V$ in the strong resolvent sense, hence $\Delta_{\Omega_{V_n}, \Omega} = \ee^{- \beta (L + V_n)} \to \ee^{- \beta (L + V)}$ in the strong resolvent sense by \cref{lem:operators_SRimpliesFSR}. With this observation, \cref{lem:vonNeumann_convergenceRelModOp} implies that $\log(\Delta_{\Omega_V, \Omega}) = - \beta (L + V) = \log(\Delta_\Omega) - \beta V$. The proof of the other identity proceeds analogously.

    \tAd{} \ref{enu:perturbationTheory_REidentity1} \& \ref{enu:perturbationTheory_REidentity2}. These assertions follow from the previously shown identities \ref{enu:perturbationTheory_boundedPerturbationE} and \ref{enu:perturbationTheory_boundedPerturbationF} in exactly the same way as in the analytic case.
\end{Proof}

\begin{Proof}[of \cref{pro:perturbationTheory_boundedGT}]
    As before, let $(V_n)_{n \in \N} \subset \SE{\MFM} \cap \MFM^\tau$ be a sequence of self-adjoint, $\tau$-entire elements such that $V_n \to V$ strongly. It was argued above in the proof of \cref{pro:perturbationTheory_boundedPerturbedVector} that in this case, $\lim_{n \to \infty} \ee^{- \beta V_n / 2} \Omega = \ee^{- \beta V / 2} \Omega$, \cf{} \cref{eq:perturbationTheory_GTconvergence}. This implies
    \begin{equation*}
        \lim_{n \to \infty} \norm{\ee^{- \beta V_n / 2} \Omega} = \norm{\ee^{- \beta V / 2} \Omega} \ .
    \end{equation*}
    Furthermore, since $\dlim{w}{n \to \infty} \Omega_{V_n} = \Omega_V$ by \cref{eq:perturbationTheory_weakConvergenceOmegaV}, it follows that
    \begin{equation*}
        \norm{\Omega_V} \le \liminf_{n \to \infty} \norm{\Omega_{V_n}}
    \end{equation*}
    since the norm of $\MH$ is weakly lower semi-continuous (\cf{} \cref{def:forms_semiContinuity}) by the Hahn-Banach theorem.\footnotemark
    \footnotetext{Let $(X, \ndot)$ be a Banach space, let $(x_n)_{n \in \N} \subset X$ be a weakly convergent sequence with weak limit $x \in X \setminus \{0\}$, and let $f \in \cdual{X}$ such that $\norm{f}_\mop = 1$ and $\norm{x} = f(x)$ which exists by the Hahn-Banach theorem \cite[Cor. III.1.6]{Werner18}. Then $\norm{x} = f(x) = \lim_{n \to \infty} f(x_n) \le \liminf_{n \to \infty} \norm{f}_\mop \norm{x_n} = \liminf_{n \to \infty} \norm{x_n}$.}
    Applying the Golden-Thompson inequality \eqref{eq:perturbationTheory_boundedGT} for analytic perturbations to the perturbed vector $\Omega_{V_n}$, one obtains
    \begin{equation*}
        \norm{\Omega_{V_n}} \le \norm{\ee^{- \beta V_n / 2} \Omega} \ .
    \end{equation*}
    Combining these three identities, the asserted inequality follows:
    \begin{equation*}
        \norm{\Omega_V} \le \liminf_{n \to \infty} \norm{\Omega_{V_n}} \le \liminf_{n \to \infty} \norm{\ee^{- \beta V_n / 2} \Omega} = \lim_{n \to \infty} \norm{\ee^{- \beta V_n / 2} \Omega} = \norm{\ee^{- \beta V / 2} \Omega} \ . \tag*{\qedhere}
    \end{equation*}
\end{Proof}

\subsection{Unbounded Perturbations}\label{subsec:perturbationTheory_unboundedPerturbationKMS}

The results of \cref{subsec:perturbationTheory_boundedPerturbationKMS} shall now be generalized to unbounded self-adjoint perturbations $V \in \AE{\MFM}$ affiliated with $\MFM$. The strategy is similar as before: $V$ will be approximated by a suitable sequence $(V_n)_{n \in \N} \subset \SE{\MFM}$ such that one can reduce the problem to an application of the results for bounded perturbations. This framework was developed by \textsc{Derezi\'{n}ski}, \textsc{Jak\v{s}i\'{c}}, and \textsc{Pillet} in \cite[Sect. 5.5]{DJP03}.

\begin{para}[Assumptions]\label{para:perturbationTheory_assumptionsUnboundedPKMS}
    Let $(\MFM, \MH, J, \MP)$ be a von Neumann algebra represented in standard form, let $\tau$ be a $W^\ast$-dynamics on $\MFM$ with standard Liouvillian $L$, and let $\beta > 0$ and $\omega \in \NS(\MFM)$ be a $(\tau, \beta)$-KMS-state with vector representative $\Omega \in \NPC$. Consider a self-adjoint operator $V \in \AE{\MFM}$ affiliated with $\MFM$, and assume that the following two properties are satisfied:
    \begin{enumerate}[label=\normalfont(A\arabic*)]
        \item \label{enu:perturbationTheory_A1} $L + V$ is essentially self-adjoint on $\dom(L) \cap \dom(V)$.
        
        \item \label{enu:perturbationTheory_A2} $L_V \ce L + V - J V J$ is essentially self-adjoint on $\dom(L) \cap \dom(V) \cap \dom(J V J)$.
    \end{enumerate}
    It was shown in \cref{pro:perturbationTheory_unboundedPerturbations} that under the first assumption, $\tau_t^V(A) = \ee^{\ii t (L + V)} A \, \ee^{- \ii t (L + V)}$ defines a $W^\ast$-dynamics on $\MFM$. Assuming additionally \ref{enu:perturbationTheory_assumptionPLiouvillian}, \cref{pro:perturbationTheory_perturbationLiouvillian} showed that $L_V = L + V - J V J$ is the standard Liouvillian of $\tau^V$. To construct KMS-states also for unbounded perturbations, an additional assumption is required \cite[p. 479]{DJP03}:
    \begin{enumerate}[resume, label=\normalfont(A\arabic*)]
        \item \label{enu:perturbationTheory_A3} $\norm{\ee^{- \beta V / 2} \Omega} < + \infty$.
    \end{enumerate}
    Introduce the following class of \enquote{DJP-perturbations} to simplify notation:
    \begin{equation*}
        \MFS_1(\MFM, L, \Omega) \ce \left\{V \in \AE{\MFM} \ : \ \parbox{0.6\textwidth}{\centering $V$ self-adjoint and satisfies assumptions \ref{enu:perturbationTheory_A1}, \ref{enu:perturbationTheory_A2}, and \ref{enu:perturbationTheory_A3} with respect to $L$ and $\Omega$}\right\} \ .
    \end{equation*}

    To obtain an analogue of \cref{thm:perturbationTheory_boundedPerturbation}, first an approximating sequence $(V_n)_{n \in \N} \subset \SE{\MFM}$ for $V$ will be constructed, and then a corresponding version of \cref{lem:perturbationTheory_boundedConvergenceLiouvillian} will be proved.
\end{para}

\begin{para}[Approximation of unbounded operators]\label{para:perturbationTheory_unboundedApproximation}
    (\cite[p. 479]{DJP03})
    Let $V : \MH \supset \dom(V) \longto \MH$ be an unbounded self-adjoint operator affiliated with the von Neumann algebra $\MFM$. Denote by $E_V$ the unique spectral measure associated with $V$. For all $n \in \N$, define an operator
    \begin{equation*}
        V_n \ce \1_{[-n, n]}(V) V \ .
    \end{equation*}
    From \cref{pro:operators_functionalCalculus}, it follows that $V_n$ can be written as
    \begin{equation}\label{eq:perturbationTheory_boundedApproxSequence}
        V_n = \1_{[-n, n]}(V) \, \id_{\sigma(V)}(V) = \bigl(\1_{[-n, n]} \cdot \id_{\sigma(V)}\bigr)(V) = \int_{\sigma(V) \cap [-n, n]} \lambda \diff E_V(\lambda) \ .
    \end{equation}
    Define a function $f_n : \R \longto \R$ by setting $f_n \ce \1_{[-n, n]} \cdot \id_{\sigma(V)}$. It holds that $\norm{f_n}_\infty = n$, hence $V_n = f_n(V)$ is a bounded operator. In particular, since $V$ is affiliated with $\MFM$, it follows that $V_n \in \SE{\MFM}$ for all $n \in \N$ by \cref{pro:operatorAlgebras_characterizationAffiliation}.
\end{para}

\begin{lemma}\label{lem:perturbationTheory_convergenceUnboundedApproximation}
    Let $\xi \in \dom(V)$ be arbitrary. Then $\lim_{n \to \infty} V_n \xi = V \xi$.
\end{lemma}

\begin{Proof}
    For $n \in \N$, define the function $g_n \ce \1_{\R \setminus [-n, n]} \cdot \id_{\sigma(V)} = \id_{\sigma(V)} - f_n$ on $\R$. Let $\mu_\xi^V$ denote the positive measure on the Borel $\sigma$-algebra of $\sigma(V)$ given by $A \longmto \braket{\xi, E_V(A) \xi}$. From properties of the functional calculus for the operator $V$ (\cref{pro:operators_functionalCalculus}), it follows that
    \begin{align*}
        \norm{V \xi - V_n \xi}^2 &= \vnorm[\big]{(\id_{\sigma(V)} - f_n)(V) \, \xi}^2 = \int_{\sigma(V)} \abs{g_n(\lambda)}^2 \, \diff \mu_\xi^V(\lambda) = \int_{\sigma(V)} \1_{\R \setminus [-n, n]}(\lambda) \, \abs{\lambda}^2 \, \diff \mu_\xi^V(\lambda) \ .
    \end{align*}
    Since $\xi \in \dom(V)$, it holds that $\int_{\sigma(V)} \abs{\lambda}^2 \diff \mu_\xi^V(\lambda) < + \infty$ (\cref{def:operators_spectralIntegral}). Moreover, as $\1_{\R \setminus [-n, n]}$ converges pointwise to the zero function as $n \to + \infty$ and $\abs{g_n}^2 \le \abs{\id_{\sigma(V)}}^2$ for all $n \in \N$ on $\sigma(V)$, it follows from Lebesgue's dominated convergence theorem \cite[Thm. 2.4.5]{Cohn13} that
    \begin{equation*}
        \norm{V \xi - V_n \xi}^2 = \int_{\sigma(V)} \1_{\R \setminus [-n, n]} \, \abs{\lambda}^2 \, \diff \mu_\xi^V(\lambda) \longto 0 \ . \tag*{\qedhere}
    \end{equation*}
\end{Proof}

The first two assertions of the next auxiliary result are from \cite[Thm. 5.6]{DJP03}.

\begin{lemma}\label{lem:perturbationTheory_unboundedConvergenceLiouvillian}
    Let $V \in \MFS_1(\MFM, L, \Omega)$ and $(V_n)_{n \in \N}$ be an approximating sequence as in \cref{para:perturbationTheory_unboundedApproximation}.
    \begin{enumerate}
        \item \label{enu:perturbationTheory_unbdConvergenceUnperturbedL} $L + V_n \to L + V$ in the strong resolvent sense.
        
        \item \label{enu:perturbationTheory_unbdConvergencePerturbedL} $L_{V_n} \to L_V$ in the strong resolvent sense.
        
        \item \label{enu:perturbationTheory_unbdConvergenceSum} $- L_{V_n} + V_n \to - L_V + V$ in the strong resolvent sense.
    \end{enumerate}
\end{lemma}

\begin{Proof}
    The proof of assertion \ref{enu:perturbationTheory_unbdConvergenceUnperturbedL} is similar to the one for \cref{lem:perturbationTheory_boundedConvergenceLiouvillian} \ref{enu:perturbationTheory_convergenceUnperturbedL}.

    \tAd{} \ref{enu:perturbationTheory_unbdConvergencePerturbedL} Define the set $\MD_0 \ce \dom(L) \cap \dom(V) \cap \dom(J V J)$. According to assumption \ref{enu:perturbationTheory_A2}, the operator $L_V$ is essentially self-adjoint on $\MD_0$. Furthermore, it holds that $L_{V_n} \xi \to L_V \xi$ for all $\xi \in \MD_0$ by \cref{lem:perturbationTheory_convergenceUnboundedApproximation}. Therefore, the assertion follows from an application of \cref{pro:operators_coreSOimpliesSR}.

    \tAd{} \ref{enu:perturbationTheory_unbdConvergenceSum}. By definition of the operators $L_{V_n}$ and $L_V$, it holds that $- L_{V_n} + V_n = - L + J V_n J$ and $- L_V + V = - L + J V J$. According to \cref{lem:perturbationTheory_LJVJessSA}, the second operator is essentially self-adjoint on $\MD_1 \ce \dom(L) \cap \dom(J V J) = J \bigl(\dom(L) \cap \dom(V)\bigr)$. Furthermore, \cref{lem:perturbationTheory_convergenceUnboundedApproximation} implies that for all $\xi = J \eta \in \MD_1$, there holds $(- L_{V_n} + V_n) \xi = - L \xi + J V_n \eta \to - L \xi + J V \eta = (- L_V + V) \xi$. Hence, the claim follows again from \cref{pro:operators_coreSOimpliesSR}.
\end{Proof}

The following theorem is taken from \cite[Thm. 5.5]{DJP03}, except assertion \ref{enu:perturbationTheory_unboundedRelModOp2} which is independent.

\begin{ntheorem}[Properties of $\Omega_V$ and $\omega_V$ for unbounded $V$]\label{thm:perturbationTheory_unboundedPerturbation}
    Let $\beta > 0$ be arbitrary, let $\MFM$ be a von Neumann algebra represented in standard form $(\MFM, \MH, J, \NPC)$, let $\tau$ be a $W^\ast$-dynamics on $\MFM$ with standard Liouvillian $L$, and let $\omega$ be a faithful $(\tau, \beta)$-KMS-state with vector representative $\Omega \in \NPC$. Then for every $V \in \MFS_1(\MFM, L, \Omega)$, the following properties are satisfied.
    \begin{enumerate}
        \item \label{enu:perturbationTheory_unboundedPerturbationA} $\Omega \in \dom(\ee^{- \beta (L + V) / 2})$. Define $\Omega_V$ and $\omega_V$ as in \cref{def:perturbationTheory_boundedPerturbedVector}.
        \item $\Omega_V \in \MP$.
        \item \label{enu:perturbationTheory_unboundedPerturbationC} $\Omega_V$ is cyclic and separating for $\MFM$.
        \item \label{enu:perturbationTheory_unboundedPerturbationD} $\omega_V$ is a $(\tau^V, \beta)$-KMS-state.
        \item \label{enu:perturbationTheory_unboundedPerturbationE} $\log(\Delta_{\Omega_V}) = - \beta L_V$.
        \item \label{enu:perturbationTheory_unboundedRelModOp1} $\log(\Delta_{\Omega_V, \Omega}) = \log(\Delta_\Omega) - \beta V$.
        \item \label{enu:perturbationTheory_unboundedRelModOp2} $\log(\Delta_{\Omega, \Omega_V}) = \log(\Delta_{\Omega_V}) + \beta V$.
    \end{enumerate}
\end{ntheorem}

\begin{Proof}
    \tAd{} \ref{enu:perturbationTheory_unboundedPerturbationA}. For $n \in \N$, let the bounded operator $V_n = f_n(V) \in \SE{\MFM}$ be defined as in \cref{para:perturbationTheory_unboundedApproximation}. According to assumption \ref{enu:perturbationTheory_A3}, it holds that $\Omega \in \dom(\ee^{- \beta V / 2})$. Therefore, employing essentially the same argument as in the proof of \cref{lem:perturbationTheory_convergenceUnboundedApproximation}, it follows that
    \begin{equation*}
        \lim_{n \to \infty} \ee^{- \beta V_n / 2} \Omega = \ee^{- \beta V / 2} \Omega \ .
    \end{equation*}
    This shows that there exists $C > 0$ such that for all $n \in \N$, $\norm{\ee^{- \beta V_n / 2} \Omega} \le C$. From the Golden-Thompson inequality \eqref{eq:perturbationTheory_boundedGT} for bounded perturbations, one obtains
    \begin{equation*}
        \norm{\Omega_{V_n}} \le \vnorm[\big]{\ee^{-\beta V_n / 2} \Omega} \ .
    \end{equation*}
    Therefore, $\norm{\Omega_{V_n}} \le C$, and applying \cref{pro:operators_SRimpliesW} gives $\Omega \in \dom(\ee^{- \beta (L + V) / 2})$ as well as
    \begin{equation*}
        \dlim{w}{n \to \infty} \ee^{- \beta (L + V_n) / 2} \Omega = \ee^{- \beta (L + V) / 2} \Omega \ .
    \end{equation*}

    By virtue of \cref{lem:perturbationTheory_unboundedConvergenceLiouvillian}, all of the remaining assertions follow from their bounded versions in precisely the same way as the corresponding statements of \cref{thm:perturbationTheory_boundedPerturbation} followed from their analytic versions. This shall be illustrated explicitly only for \ref{enu:perturbationTheory_unboundedRelModOp2}.
    
    \tAd{} \ref{enu:perturbationTheory_unboundedRelModOp2} For every $n \in \N$, it holds that $\log(\Delta_{\Omega, \Omega_{V_n}}) = \log(\Delta_{\Omega_{V_n}}) + \beta V_n = - \beta L_{V_n} + \beta V_n$ according to \cref{thm:perturbationTheory_boundedPerturbation} \ref{enu:perturbationTheory_boundedPerturbationD} and \ref{enu:perturbationTheory_boundedPerturbationF}. Furthermore, \cref{lem:perturbationTheory_unboundedConvergenceLiouvillian} \ref{enu:perturbationTheory_unbdConvergenceSum} shows that $- L_{V_n} + V_n \to - L_V + V$ in the strong resolvent sense, hence $\Delta_{\Omega, \Omega_{V_n}} = \ee^{- \beta L_{V_n} + \beta V_n} \to \ee^{- \beta L_V + \beta V}$ in the strong resolvent sense by \cref{lem:operators_SRimpliesFSR}. Since $\Omega_{V_n} \to \Omega_V$ weakly in $\MH$ by \ref{enu:perturbationTheory_unboundedPerturbationA}, \cref{lem:vonNeumann_convergenceRelModOp} and assertion \ref{enu:perturbationTheory_unboundedPerturbationE} imply that $\log(\Delta_{\Omega, \Omega_V}) = - \beta L_V + \beta V = \log(\Delta_{\Omega_V}) + \beta V$. (Note that $\ssupp(\omega_{V_n}) = \id_\MH$ and $\ssupp(\omega_V) = \id_\MH$ since the vectors $\Omega_{V_n}, \Omega_V$ are separating, hence assumption \ref{enu:vonNeumann_convergenceRelModOpA2} of \cref{lem:vonNeumann_convergenceRelModOp} is also satisfied.)
\end{Proof}

\begin{para}[Additional assumptions]\label{para:perturbationTheory_additionalAssumptions}
    The generalization of the properties stated in \cref{thm:perturbationTheory_boundedPerturbation} \ref{enu:perturbationTheory_REidentity1} and \ref{enu:perturbationTheory_REidentity2} is a bit more subtle. In fact, another assumption is necessary to obtain the latter identity, while the former can be obtained using only \ref{enu:perturbationTheory_A1} -- \ref{enu:perturbationTheory_A3}. Before introducing this additional assumption, an important consequence of \ref{enu:perturbationTheory_A3} shall be discussed.
    
    Let $V \in \AE{\MFM}$ be a self-adjoint operator affiliated with $\MFM$. Then $V$ can be decomposed into the sum $V = V_+ + V_-$, where $V_+$ is the \emph{positive part} of $V$ and $V_-$ is the \emph{negative part} of $V$ which are given by the expression
    \begin{align*}
        V_+ &\ce \1_{(0, + \infty)}(V) V = \bigl(\1_{(0, + \infty)} \cdot \id_{\sigma(V)}\bigr)(V) \ , \\
        V_- &\ce \1_{(- \infty, 0]}(V) V = \bigl(\1_{(- \infty, 0]} \cdot \id_{\sigma(V)}\bigr)(V) \ .
    \end{align*}
    Consider $\beta > 0$ and a $(\tau, \beta)$-KMS-state $\omega$ on $\MFM$ as given in \cref{thm:perturbationTheory_unboundedPerturbation}. Using \cref{pro:operators_functionalCalculus} \ref{enu:operators_funcCalcIP}, it follows that
    \begin{align*}
        \beta \omega(V_-) = \vbraket[\big]{\Omega, \bigl(\beta \1_{(- \infty, 0]} \cdot \id_{\sigma(V)}\bigr)(V) \Omega} = \int_{\sigma(V) \cap (- \infty, 0]} \beta \lambda \diff \mu_\Omega^V(\lambda) \ .
    \end{align*}
    One can estimate the integral by using that $x > - \ee^{- x}$ for all $x \in \R$ which is, in essence, a consequence of Bernoulli's inequality.\footnotemark
    \footnotetext{Indeed, from the limit representation $\ee^x = \lim_{n \to \infty} (1 + x / n)^n$ and the Bernoulli inequality $(1 + x)^n \ge 1 + nx$, which is valid for all $x \ge -1$ and $n \in \N$ \cite[§ 3, Thm. 2]{Forster16}, it follows that $\ee^{x} \ge 1 + x$ for all $x \in \R$. (The validity of this inequality for $x < - 1$ is clear.) Hence, substituting $x \to - x$, one obtains $\ee^{-x} \ge 1 - x$, and so $- \ee^{-x} \le - 1 + x < x$.\label{ftn:perturabtionTheory_exponentialInequality}}
    With this, it follows that
    \begin{align*}
        \beta \omega(V_-) &> - \int_{\sigma(V) \cap (- \infty, 0]} \ee^{- \beta \lambda} \diff \mu_\Omega^V(\lambda) \ge - \int_{\sigma(T)} \ee^{- \beta \lambda} \diff \mu_\Omega^V(\lambda) \ ,
    \end{align*}
    where it was used that $\int_{\sigma(T) \cap (0, + \infty)} \ee^{- \beta \lambda} \diff \mu_\Omega^V(\lambda) \ge 0$ since $\ee^{- \beta \lambda} \ge 0$ and $\mu_\Omega^V$ is a positive measure. Finally, rewriting the last spectral integral from above, one obtains
    \begin{align*}
        \beta \omega(V_-) &> - \int_{\sigma(T)} \ee^{- \beta \lambda} \diff \mu_\Omega^V(\lambda) = - \braket{\Omega, \ee^{- \beta V} \Omega} = - \norm{\ee^{- \beta V / 2} \Omega}^2 > - \infty \ .
    \end{align*}
    Note that the last inequality is true because of the assumption \ref{enu:perturbationTheory_A3} from \cref{para:perturbationTheory_assumptionsUnboundedPKMS}. The condition $\beta \omega(V_-) > - \infty$ will be essential to establish the formula for the relative entropy of the unperturbed state $\omega$ with respect to the perturbed state $\omega_V$. Similarly, to obtain the relative entropy for the states in reversed order, the following additional assumption is required:
    \begin{enumerate}[label=\normalfont(A\arabic*)]
        \setcounter{enumi}{3}
        \item \label{enu:perturbationTheory_A4} $\beta \omega_V(V_+) < + \infty$.
    \end{enumerate}
    To simplify the notation again, introduce the following class of \enquote{Bogoliubov perturbations}:
    \begin{equation*}
        \MFS_2(\MFM, L, \Omega) \ce \left\{V \in \AE{\MFM} \ : \ \parbox{0.6\textwidth}{\centering $V$ self-adjoint and satisfies assumptions \ref{enu:perturbationTheory_A1}, \ref{enu:perturbationTheory_A2}, \ref{enu:perturbationTheory_A3}, and \ref{enu:perturbationTheory_A4} with respect to $L$ and $\Omega$}\right\} \ .
    \end{equation*}

    Note that the additional assumption \ref{enu:perturbationTheory_A4} does not appear in \cite{DJP03} and, in fact, the latter does not state the second identity for the relative entropy $\MS_\MFM^\mathrm{std}(\omega_V, \omega)$ contained in the next proposition.
\end{para}

\begin{proposition}\label{pro:perturbationTheory_relEntUnboundedPerturbedState}
    Let $\beta > 0$ be arbitrary, $(\MFM, \MH, J, \NPC)$ be a von Neumann algebra in standard form, $\tau$ be a $W^\ast$-dynamics on $\MFM$ with Liouvillian $L$, and $\omega$ be a faithful $(\tau, \beta)$-KMS-state with vector representative $\Omega \in \NPC$. The following identities hold true for all $V \in \MFS_2(\MFM, L, \Omega)$:
    \begin{equation*}
        \MS_\MFM^\mathrm{std}(\omega, \omega_V) = \beta \omega(V) + \log\bigl(\norm{\Omega_V}^2\bigr) \tand \MS_\MFM^\mathrm{std}(\omega_V, \omega) = - \beta \omega_V(V) - \log\bigl(\norm{\Omega_V}^2\bigr) \ .
    \end{equation*}
\end{proposition}

\begin{Proof}
    By the assumptions \ref{enu:perturbationTheory_A3} and \ref{enu:perturbationTheory_A4}, it follows that $\beta \omega(V) = \braket{\Omega, \beta V \Omega}$ is a finite number or $+ \infty$, and similarly $\beta \omega_V(V) = \braket{\Omega_V, \beta V \Omega_V} / \norm{\Omega_V}^2$ is a finite number or $- \infty$. This shows that the right-hand sides in the above equations cannot be equal to $- \infty$ which is necessary to guarantee non-negativity of the relative entropy (\cref{cor:relativeEntropy_nonnegSpaRelativeEntropy}). With this in mind, the identities follow from \cref{thm:perturbationTheory_unboundedPerturbation} \ref{enu:perturbationTheory_unboundedRelModOp1} and \ref{enu:perturbationTheory_unboundedRelModOp2} in the same way as \cref{thm:perturbationTheory_boundedPerturbation} \ref{enu:perturbationTheory_REidentity1} and \ref{enu:perturbationTheory_REidentity2} were proved for analytic perturbations.
\end{Proof}

As the identity for $\MS_\MFM^\mathrm{std}(\omega, \omega_V)$ in the previous proposition, the Golden-Thompson inequality for unbounded perturbations is also part of \cite[Thm. 5.5]{DJP03}.

\begin{proposition}[Golden-Thompson inequality]\label{pro:perturbationTheory_unboundedGT}
    For $V \in \MFS_1(\MFM, L, \Omega)$, there holds
    \begin{equation*}
        \norm{\Omega_V} \le \norm{\ee^{- \beta V / 2} \Omega} \ .
    \end{equation*}
\end{proposition}

\begin{Proof}
    Since $\lim_{n \to \infty} \ee^{- \beta V_n / 2} \Omega = \ee^{- \beta V / 2} \Omega$ and $\dlim{w}{n \to \infty} \Omega_{V_n} = \Omega_V$ as observed in the proof of \cref{thm:perturbationTheory_unboundedPerturbation} \ref{enu:perturbationTheory_unboundedPerturbationA}, where $(V_n)_{n \in \N}$ is given by \cref{eq:perturbationTheory_boundedApproxSequence}, one can employ the same proof strategy as used for bounded perturbations in \cref{pro:perturbationTheory_boundedGT} which, in particular, relies on using the Golden-Thompson inequality \eqref{eq:perturbationTheory_boundedGT} for bounded perturbations.
\end{Proof}

\section{The Two-sided Bogoliubov Inequality}\label{sec:perturbationTheory_twoSidedBogoliubov}

Finally, the two-sided Bogoliubov inequality, which was recently proved for quantum mechanical systems \cite{Reible22} and applied to the problem of determining finite-size effects in molecular simulations \cite{Reible23, DS24}, can be generalized to arbitrary von Neumann algebras, relying on the perturbation theory of KMS-states developed above. First, the quantum-mechanical case will be reviewed quickly following \cite{Reible22}.

\subsection{The Quantum-Mechanical Case}\label{subsec:perturbationTheory_twoSidedBogoliubovQM}

\begin{para}[Mathematical setup]\label{para:perturbationTheory_setupBogoliubovQM}
    Consider a finite quantum system $\MFM = \BO(\MH)$ at inverse temperature $\beta > 0$ confined to a region $\Sigma \subset \R^3$ with Hamiltonian $H : \MH \supset \dom(H) \longto \MH$ such that $\ee^{- \beta H} \in \NO(\MH)$. One may characterize the system in terms of the canonical Gibbs state (\cf{} \cref{para:perturbationTheory_motivationKMS})
    \begin{equation}\label{eq:perturbationTheory_fullGibbs}
        \rho = \frac{1}{Z} \, \ee^{- \beta H} \com Z = \tr(\ee^{- \beta H}) \ .
    \end{equation}
    Suppose that the region $\Sigma$ is divided into $d \in \N$ subregions $\Sigma_k \subset \R^3$, $1 \le k \le d$, such that $\Sigma = \bigcup_{k=1}^{d} \Sigma_k$. Then each $\Sigma_k$ represents a smaller subsystem of the total system $\Sigma$ which shall be described independently from the other ones. To implement this idea, one decomposes the the Hamiltonian $H$ into the sum
    \begin{equation*}
        H = H_0 + U \ ,
    \end{equation*}
    where $H_0$ and $U$ are self-adjoint operators such that $H$ is again self-adjoint on its domain. Physically, the operator $H_0$ represents the Hamiltonian of the non-interacting subsystems, and it takes the form
    \begin{equation*}
        H_0 = \sum_{k=1}^{d} H_0^{(k)} \ ,
    \end{equation*}
    with $H_0^{(k)}$ being the self-adjoint Hamiltonian describing only the $k$-th subsystem located in $\Sigma_k$. The operator $U$ mediates the interaction between the $d$ different subsystems. Assume that also $\ee^{- \beta H_0}$ is trace-class and define the Gibbs state
    \begin{equation}\label{eq:perturbationTheory_freeGibbs}
        \rho_0 = \frac{1}{Z_0} \, \ee^{-\beta H_0} \com Z_0 = \tr(\ee^{- \beta H_0})
    \end{equation}
    representing the thermal state of the uncoupled subsystems.
    
    To quantify the thermodynamic difference between the full system and the collection of uncoupled systems, one defines the \bemph{relative free energy} (or: \emph{interface energy})
    \begin{equation}\label{eq:perturbationTheory_interfaceEnergy}
        \Delta F \ce - \inv{\beta} \log\left(\frac{Z}{Z_0}\right) \ .
    \end{equation}
    This is nothing but the difference between the free energy $F = - \inv{\beta} \log(Z)$ of the full system and the free energy $F_0 = - \inv{\beta} \log(Z_0)$ of the uncoupled subsystems. Physically, $\Delta F$ represents the free energy difference associated with the partitioning of the large system into smaller, independent subsystems.
\end{para}

\begin{proposition}[\protect{Two-sided Bogoliubov inequality \cite[Thm. 4.1]{Reible22}}]\label{pro:perturbationTheory_BogoliubovQM}
    In the situation described above, assume that $\EV_{\rho_0}[U] \ce \tr(\rho_0 U) < + \infty$ and $\EV_\rho[U] \ce \tr(\rho U) < + \infty$. Then
    \begin{equation}\label{eq:perturbationTheory_BogoliubovQM}
        \EV_\rho[U] \le \Delta F \le \EV_{\rho_0}[U] \ .
    \end{equation}
\end{proposition}

\begin{Proof}
    Consider the Umegaki relative entropy $S$, defined in \cref{eq:introduction_Umegaki}, between the Gibbs states of the uncoupled and coupled system. One computes \cite[pp. 8 f.]{Reible22}
    \begin{align*}
        S(\rho_0, \rho) = \beta \, \EV_{\rho_0}[U] - \beta \Delta F \tand S(\rho, \rho_0) = - \beta \, \EV_{\rho}[U] + \beta \Delta F \ .
    \end{align*}
    From non-negativity of $S$ (see \cref{cor:relativeEntropy_nonnegSpaRelativeEntropy}), it follows that $\Delta F \le \EV_{\rho_0}[U]$, which is the upper bound, and $\EV_{\rho}[U] \le \Delta F$, which is the lower bound.
\end{Proof}

\begin{remarks}\label{rem:perturbationTheory_BogoliubovQM}
    \leavevmode
    \begin{enumerate}[env]
        \item The upper bound is the famous \emph{Peierls-Bogoliubov inequality} \cite{Peierls38, Bogoliubov58, Symanzik65} which is an important tool in statistical mechanics \cite{Bebiano04}, \cite[Sect. 2.5]{Ruelle99}, \cite[Sect. 8.3]{Simon05}, and which, for two self-adjoint operators $A$ and $B$ such that $\ee^A, \ee^{A + B} \in \NO(\MH)$, is usually written \cite[Eq. (2.14)]{Carlen10}, \cite[Cor. 3.14]{OP04}
        \begin{equation*}
            \log\left(\frac{\tr(\ee^{A + B})}{\tr(\ee^A)}\right) \ge \frac{\tr(B \, \ee^A)}{\tr(\ee^A)} \ .
        \end{equation*}
        Indeed, choosing $A = - \beta H_0$ and $B = - \beta U$, the above inequality becomes $\log(Z / Z_0) \ge - \beta \, \EV_{\rho_0}[U]$ which is equivalent to the upper bound in \cref{eq:perturbationTheory_BogoliubovQM}.
        
        \item The two-sided Bogoliubov inequality can be applied in numerical simulations to define and estimate the error due to finite-size effects \cite[p. 9]{Reible22}; a general computational protocol highlighting the utility of \eqref{eq:perturbationTheory_BogoliubovQM} is discussed in \cite[Sect. 5]{Reible22}. In particular, the inequality can be used to justify the simulation of a small subsystem instead of the computationally unfeasible total system \cite{Reible23}; see also \cite{DS17, DS24}.
    \end{enumerate}
\end{remarks}

The inequalities in \eqref{eq:perturbationTheory_BogoliubovQM} can be sharpened by using the upper bound in this equation and the Golden-Thompson trace inequality (\cf{} \cref{rem:perturbationTheory_GT}). The following result is proved in \cite[Sect. 4.1]{Reible22}.

\begin{proposition}[\protect{\cite[Cor. 4.5]{Reible22}}]\label{pro:perturbationTheory_variationalBoundsQM}
    Under the assumptions of \cref{pro:perturbationTheory_BogoliubovQM}, it holds that
    \begin{equation}\label{eq:perturbationTheory_variationalBoundsQM}
        \sup_{V \ge 0} \Bigl\{\EV_\rho[U - V] - \inv{\beta} \log \tr\bigl(\ee^{\log \rho_0 - \beta V}\bigr)\Bigr\} = \Delta F = \inf_{\gamma \in \DM(\MH)} \Bigl\{\EV_\gamma[U] + \inv{\beta} S(\gamma, \rho_0)\Bigr\} \ .
    \end{equation}
    The supremum is taken over all densely defined, positive self-adjoint operators $V$ on $\MH$, and the infimum ranges over all density matrices $\gamma$ on $\MH$. In particular, the bounds of \cref{pro:perturbationTheory_BogoliubovQM} are a special case of \eqref{eq:perturbationTheory_variationalBoundsQM} obtained by taking $V = 0$ and $\gamma = \rho_0$.
\end{proposition}

\subsection{Generalization to Arbitrary von Neumann Algebras}\label{subsec:perturbationTheory_twoSidedBogoliubovVN}

A natural question to ask is whether the two-sided Bogoliubov inequality \eqref{eq:perturbationTheory_BogoliubovQM} can be extended to quantum systems of infinitely many degrees of freedom, that is, to arbitrary von Neumann algebras. A first step towards this goal consists of defining an analogue of the quantity $\Delta F$ from \cref{eq:perturbationTheory_interfaceEnergy}.

\begin{definition}[Relative free energy]\label{def:perturabtionTheory_relativeFreeEnergy}
    Let $(\MFM, \MH, J, \MP)$ be a von Neumann algebra in standard form, let $\tau$ be a $W^\ast$-dynamics on $\MFM$ with standard Liouvillian $H_0$, and let $\beta > 0$ and $\omega_0$ be a faithful $(\tau, \beta)$-KMS-state on $\MFM$ with vector representative $\Omega_0 \in \NPC$. Given a Bogoliubov perturbation $U \in \MFS_2(\MFM, H_0, \Omega_0)$, define the \bemph{relative free energy} of the perturbed state $\omega_U$ with respect to the state $\omega_0$ to be
    \begin{equation}\label{eq:perturbationTheory_relativeFreeEnergy}
        \MF(\omega_U, \omega_0) \ce \omega_U(U) + \inv{\beta} \MS_\MFM^\mathrm{std}(\omega_U, \omega_0) \ .
    \end{equation}
\end{definition}

\begin{remark}\label{rem:perturabtionTheory_relativeFreeEnergy}
    Consider a finite quantum system $\MFM = \BO(\MH)$ with Hamiltonian $H_0$ at inverse temperature $\beta > 0$, let $\omega_0 = \tr(\rho_0 \ \sbullet\,)$ be the canonical Gibbs state for $H_0$, and let $\omega_U = \tr(\rho \ \sbullet\,)$ be the corresponding one for the Hamiltonian $H = H_0 + U$, see \cref{eq:perturbationTheory_fullGibbs,eq:perturbationTheory_freeGibbs}. In the proof of \cref{pro:perturbationTheory_BogoliubovQM}, it was mentioned that
    \begin{equation*}
        \Delta F = \tr(\rho U) + \inv{\beta} \tr\bigl(\rho (\log \rho - \log \rho_0)\bigr) = \omega_U(U) + \inv{\beta} S(\rho, \rho_0) \ ,
    \end{equation*}
    where $S(\rho, \rho_0)$ is the Umegaki relative entropy. Since the right-hand side of this expression can be defined in arbitrary von Neumann algebras, this justifies calling \eqref{eq:perturbationTheory_relativeFreeEnergy} \enquote{relative free energy}.
    
    A definition of a relative free energy in operator algebras similar to \eqref{eq:perturbationTheory_relativeFreeEnergy}, justified by physical principles and the analogy to finite quantum systems, was already suggested by \textsc{M. J. Donald} in 1987 \cite[Eq. (1.2)]{Donald87a}; he considers a relative free energy of an \emph{arbitrary} normal state with respect to an equilibrium state, both defined with respect to the \emph{same} Hamiltonian of the quantum system.
\end{remark}

The next proposition is the main result of this section.

\begin{proposition}[Two-sided Bogoliubov inequality in von Neumann algebras]\label{pro:perturbationTheory_BogoliubovVN}
    Let $(\MFM, \MH, J, \MP)$ be a von Neumann algebra in standard form, let $\tau$ be a $W^\ast$-dynamics on $\MFM$ with standard Liouvillian $H_0$, and let $\beta > 0$ and $\omega_0$ be a faithful $(\tau, \beta)$-KMS-state on $\MFM$ with vector representative $\Omega_0 \in \NPC$. For all $U \in \MFS_2(\MFM, H_0, \Omega_0)$, the following two-sided inequality holds true:
    \begin{equation}\label{eq:perturbationTheory_BogoliubovVN}
        \omega_U(U) \le \MF(\omega_U, \omega_0) \le \omega_0(U) \ .
    \end{equation}
\end{proposition}

\begin{Proof}
    From the second identity in \cref{pro:perturbationTheory_relEntUnboundedPerturbedState}, $\MS_\MFM^\mathrm{std}(\omega_U, \omega_0) = - \beta \omega_U(U) - \log \norm{\Omega_U}^2$, and \cref{def:perturabtionTheory_relativeFreeEnergy}, it follows that the relative free energy can also be written as
    \begin{equation}\label{eq:perturbationTheory_perturbedVectorRFE}
        \MF(\omega_U, \omega_0) = - \inv{\beta} \log\bigl(\norm{\Omega_U}^2\bigr) \ .
    \end{equation}
    Using that the Araki-Uhlmann relative entropy $\MS_\MFM^\mathrm{std}(\omega_U, \omega_0)$ is non-negative, \cf{} \cref{pro:relativeEntropy_nonNegStdRelativeEntropy}, it follows that $\omega_U(U) \le \MF(\omega_U, \omega_0)$ which is the lower bound. Similarly, combining the first identity from \cref{pro:perturbationTheory_relEntUnboundedPerturbedState}, $\MS_\MFM^\mathrm{std}(\omega_0, \omega_U) = \beta \omega_0(U) + \log \norm{\Omega_U}^2$, with non-negativity of $\MS_\MFM^\mathrm{std}(\omega_0, \omega_U)$, one obtains $\MF(\omega_U, \omega_0) \le \omega_0(U)$ which is the upper bound.
\end{Proof}

\begin{remarks}\label{rem:perturbationTheory_boundedBogoliubov}
    \leavevmode
    \begin{enumerate}[env]
        \item The inequality \eqref{eq:perturbationTheory_BogoliubovVN} still holds true if $U \in \MFS_2(\MFM, H_0, \Omega_0)$ is replaced by a bounded perturbation $U \in \SE{\MFM}$. This can also be proved directly by employing \cref{thm:perturbationTheory_boundedPerturbation}.
        
        \item As in the quantum-mechanical case (\cref{pro:perturbationTheory_BogoliubovQM}), the upper bound $- \inv{\beta} \log \norm{\Omega_U}^2 \allowbreak\le \omega_0(U)$ of \cref{eq:perturbationTheory_BogoliubovVN} is a version of the \emph{Peierls-Bogoliubov inequality}. For operator algebras, it was established by \textsc{H. Araki} in \cite[Thm. 1]{Araki73b} for bounded perturbations. In its typical form, this inequality is written as
        \begin{align*}
            \ee^{- \beta \omega_0(U)} \le \norm{\Omega_U}^2 \ .
        \end{align*}
        It is proved for unbounded perturbations $U \in \MFS_1(\MFM, H_0, \Omega_0)$ also in \cite[Thm. 5.5 (8)]{DJP03}.
        
        The lower bound $\omega_U(U) \le - \inv{\beta} \log \norm{\Omega_U}^2$, although obtained in a straightforward fashion, does not seem to have been discussed in the literature so far. Usually, rather the \emph{Golden-Thompson inequality}, see \cref{pro:perturbationTheory_unboundedGT} and \cref{rem:perturbationTheory_GT}, which can be written in the more familiar form
        \begin{align*}
            \norm{\Omega_U}^2 \le \omega_0(\ee^{- \beta U}) \ ,
        \end{align*}
        is investigated. Note that together, the Peierls-Bogoliubov and Golden-Thompson inequality yield the following two-sided bound for the relative free energy $\MF(\omega_U, \omega_0)$:
        \begin{equation}\label{eq:perturbationTheory_twoSidedGTPB}
            - \inv{\beta} \log \omega_0(\ee^{- \beta U}) \le \MF(\omega_U, \omega_0) \le \omega_0(U) \ .
        \end{equation}
        Thus, it needs to be stressed that the lower bound in the two-sided Bogoliubov inequality \eqref{eq:perturbationTheory_BogoliubovVN} does not reproduce the known Golden-Thompson inequality, but is actually a new estimate. Of course, the virtue of the bounds in \cref{eq:perturbationTheory_twoSidedGTPB} is that they require only knowledge of the given state $\omega_0$, and not of the perturbed state $\omega_U$. However, for many applications in physics, \cref{eq:perturbationTheory_BogoliubovVN} seems to be more interesting, see \cite{DS17, DS24, Reible22, Reible23}.

        \item One can rewrite the expression \eqref{eq:perturbationTheory_perturbedVectorRFE} as follows:
        \begin{align*}
            - \inv{\beta} \log\bigl(\norm{\Omega_U}\bigr)^2 &= - \inv{\beta} \log \vbraket[\big]{\ee^{- \beta (H_0 + U) / 2} \Omega_0, \ee^{- \beta (H_0 + U) / 2} \Omega_0} \\
            &= - \inv{\beta} \log \vbraket[\big]{\Omega_0, \ee^{- \beta (H_0 + U)} \Omega_0} \ .
        \end{align*}
        If one defines $H \ce H_0 + U$ to be the \enquote{physical Liouvillian} of the system, then one obtains
        \begin{equation*}
            \MF(\omega_U, \omega_0) = - \inv{\beta} \log\bigl(\norm{\Omega_U}\bigr)^2 = - \inv{\beta} \log\bigl(\braket{\Omega_0, \ee^{- \beta H} \Omega_0}\bigr) \ .
        \end{equation*}
        Comparing this expression with some of the works of \textsc{R. Longo} \cite[p. 471]{Longo97}, \cite[Eq. (30)]{Longo18}, \cite[p. 116]{LongoXu18}, where the right-hand side of the above expression appears naturally and is termed \emph{relative} or \emph{incremental free energy}, reinforces the interpretation of the quantity $\MF(\omega_U, \omega_0)$ as the relative free energy between the states $\omega_U$ and $\omega_0$.
    \end{enumerate}
\end{remarks}

In the following, it will be shown that \cref{pro:perturbationTheory_BogoliubovVN} can be considered a proper generalization of the two-sided Bogoliubov inequality for quantum systems, \cref{pro:perturbationTheory_BogoliubovQM}, to unbounded perturbations of the \enquote{free Liouvillian} $H_0$ on general von Neumann algebras $\MFM$.

\begin{corollary}\label{cor:perturbationTheory_boundedBogoliubovImpliesQMBogoliubov}
    Let $\MFM = \BO(\MH)$ and $\beta > 0$, let $\tau$ be a $W^\ast$-dynamics on $\MFM$ with standard Liouvillian $H_0$ such that $\ee^{- \beta H_0} \in \NO(\MH)$, let $\omega_0$ be a faithful $(\tau, \beta)$-KMS-state, and let $U \in \MFS_2(\MFM, H_0, \Omega_0)$ be an unbounded perturbation such that $\ee^{- \beta (H_0 + U)} \in \NO(\MH)$. Then there exist $\rho_0, \rho \in \DM(\MH)$ satisfying the inequality
    \begin{equation*}
        \EV_\rho[U] \le \Delta F \le \EV_{\rho_0}[U] \ .
    \end{equation*}
\end{corollary}

\begin{Proof}
    Since $H_0$ is the standard Liouvillian of $\tau$, it holds that $\tau_t(A) = \ee^{\ii t H_0} A \, \ee^{- \ii t H_0}$ (\cref{pro:perturbationTheory_standardLiouvillian}). Therefore, as $\omega_0$ is a $(\tau, \beta)$-KMS-state on $\MFM$, \cref{exa:perturbationTheory_KMSfiniteQuantumSystem} implies that $\omega_0$ is given by the density matrix $\rho_0$ from \cref{eq:perturbationTheory_freeGibbs}. Furthermore, recall that the perturbed dynamics takes the form $\tau_t^U(A) = \ee^{\ii t (H_0 + U)} A \, \ee^{- \ii t (H_0 + U)}$ (\cref{pro:perturbationTheory_unboundedPerturbations}). Thus, as the perturbed state is a $(\tau^U, \beta)$-KMS-state according to \cref{thm:perturbationTheory_unboundedPerturbation} \ref{enu:perturbationTheory_unboundedPerturbationD}, it follows again from \cref{exa:perturbationTheory_KMSfiniteQuantumSystem} that $\omega_U$ is given by the density matrix $\rho$ from \cref{eq:perturbationTheory_fullGibbs} with $H = H_0 + U$.
    
    As was already discussed in \cref{rem:perturabtionTheory_relativeFreeEnergy} (see also \cref{exa:relativeEntropy_specialCasesArakiUhlmann} \ref{enu:relativeEntropy_ArakiUhlmannTypeI} for the proof that $\MS_\MFM^\mathrm{std} = S$ if $\MFM = \BO(\MH)$), the relative free energy $\MF(\omega_U, \omega_0)$ in the present situation is given by $\Delta F$ from \cref{eq:perturbationTheory_interfaceEnergy}. Therefore, \cref{pro:perturbationTheory_BogoliubovVN} implies the asserted inequality.
\end{Proof}

\subsection{Towards Variational Bounds}\label{subsec:perturbationTheory_variationalBoundsVN}

In the following, it will be investigated whether the variational bounds for the relative free energy stated in \cref{pro:perturbationTheory_variationalBoundsQM} can be extended to the operator-algebraic setting. It turns out that for bounded perturbations, corresponding results are well-known. To extend these to more general unbounded perturbations, first two auxiliary results will be proved. They are both extensions of known result due to \textsc{H. Araki} (\cf{} \cite[Eq. (3.6)]{Araki73b} and \cite[Eq. (4.28)]{Araki77}, respectively).

\begin{lemma}\label{lem:perturbationTheory_generalUnboundedPerturbedVector}
    Let $(\MFM, \MH, J, \NPC)$ be a von Neumann algebra in standard form, let $\Phi \in \NPC$ be separating, and let $L \ce \log(\Delta_\Phi)$. If $V \in \MFS_1(\MFM, L, \Phi)$, then $\Phi \in \dom(\ee^{- \beta (L + V) / 2})$, hence one may define the following vector: 
    \begin{align*}
        \Phi_V \ce \ee^{- \beta (L + V) / 2} \, \Phi \ .
    \end{align*}
\end{lemma}

\begin{Proof}
    One can prove this result in an analogous way as the existence of the perturbed $(\tau, \beta)$-KMS-vector for unbounded $V$ in \cref{thm:perturbationTheory_unboundedPerturbation} \ref{enu:perturbationTheory_unboundedPerturbationA}.

    Let $\varphi = \omega_\Phi$ be the faithful normal functional induced by $\Phi$, and let $\tau = \sigma^\varphi$ be the modular automorphism group. According to \cref{exa:perturbationTheory_omegaLiouvillian} \ref{enu:perturbationTheory_exaModularGroup}, $L$ is the standard Liouvillian of $(\MFM, \tau)$. With this observation, one establishes the assertion of the lemma first for analytic $V \in \MFM^\tau$ as in \cref{subsec:perturbationTheory_analyticProofs} by noting that the proof of the latter did not use any KMS-specific properties.
    
    Next, one shows that the claim is true for bounded $V \in \SE{\MFM}$ as was done in \cref{subsec:perturbationTheory_boundedProofs}. It is important to note that the argument requires the Golden-Thompson inequality for analytic perturbations whose proof relied on an expression for the relative entropy, \cref{thm:perturbationTheory_boundedPerturbation} \ref{enu:perturbationTheory_REidentity2}, specific for KMS-states. However, the Golden-Thompson inequality can be proved for arbitrary normal functionals and bounded perturbations without relying on this technique \cite[Thm. 2]{Araki73b}, hence the argument proceeds as mentioned.
    
    Finally, one establishes the assertion for self-adjoint $V \in \AE{\MFM}$ satisfying \ref{enu:perturbationTheory_A1} -- \ref{enu:perturbationTheory_A3} by an approximation through bounded operators as in \cref{thm:perturbationTheory_unboundedPerturbation} \ref{enu:perturbationTheory_unboundedPerturbationA} which makes use of the Golden-Thompson inequality for bounded perturbations which, as mentioned before, is available for non-KMS-states.
\end{Proof}

\begin{lemma}\label{lem:perturbationTheory_generalPerturbedRelModOp}
    Let $\varphi, \psi \in \NS(\MFM)$ be two faithful normal states on $\MFM$, denote by $\Phi, \Psi \in \NPC$ their vector representatives, and let $L = \log(\Delta_\Phi)$ and $V \in \MFS_1(\MFM, L, \Phi)$. Assume that the operator $\log(\Delta_{\Phi, \Psi}) - \beta V$ is essentially self-adjoint on $\MD \ce \dom(\log \Delta_{\Phi, \Psi}) \cap \dom(V)$. Then
    \begin{equation}\label{eq:perturbationTheory_generalPerturbedRelModOp}
        \log(\Delta_{\Phi_V, \Psi}) = \log(\Delta_{\Phi, \Psi}) - \beta V \ .
    \end{equation}
\end{lemma}

\begin{Proof}
    Let $(V_n)_{n \in \N} \subset \SE{\MFM}$ be a sequence of bounded self-adjoint elements of $\MFM$ approximating the operator $V$ as constructed in \cref{para:perturbationTheory_unboundedApproximation}.
    
    \emph{First}, note that for all $\xi \in \MD$, it holds that $\lim_{n \to \infty} V_n \xi = V \xi$ according to \cref{lem:perturbationTheory_convergenceUnboundedApproximation}. Therefore, the assumption on $\log(\Delta_{\Phi, \Psi}) - \beta V$ together with \cref{pro:operators_coreSOimpliesSR} implies that
    \begin{equation*}
        \log(\Delta_{\Phi, \Psi}) - \beta V_n \longto \log(\Delta_{\Phi, \Psi}) - \beta V
    \end{equation*}
    in the strong resolvent sense as $n \to + \infty$. \emph{Second}, observe that by the last identity in the proof of \cref{thm:perturbationTheory_unboundedPerturbation} \ref{enu:perturbationTheory_unboundedPerturbationA} (which, as mentioned in the proof of \cref{lem:perturbationTheory_generalUnboundedPerturbedVector}, also holds for non-KMS-states), one has $\dlim{w}{n \to \infty} \Phi_{V_n} = \Phi_V$. Furthermore, since $V_n \in \SE{\MFM}$ is bounded, it is well-known that $\Phi_{V_n}$, $n \in \N$, is a cyclic and separating element of the natural positive cone $\NPC$. (See \cite[Prop. 4.1 \& Cor. 4.4]{Araki73c} together with the identity \cite[Eq. (3.6)]{Araki73b}, or \cite[p. 221 \& Cor. 12.7]{OP04} together with the remarks regarding the normalization convention of the perturbed functional in \cite[Exa. 3.3]{Donald90}.) Therefore, since $\NPC$ is weakly closed (\cf{} \cpageref{eq:perturbationTheory_weakConvergenceOmegaV}), it follows that $\Phi_V \in \NPC$. \emph{Third}, for bounded perturbations, the following identity is known to hold true \cite[Eq. (4.28)]{Araki77}, \cite[Eq. (3.6) \& Thm. A.7]{Donald90}, \cite[Thm. 12.6]{OP04}:
    \begin{equation*}
        \log(\Delta_{\Phi_{V_n}, \Psi}) = \log(\Delta_{\Phi, \Psi}) - \beta V_n \ .
    \end{equation*}
    Since, as argued above, the right-hand side converges to $\log(\Delta_{\Phi, \Psi}) - \beta V$ in the strong resolvent sense, \cref{lem:vonNeumann_convergenceRelModOp} implies \cref{eq:perturbationTheory_generalPerturbedRelModOp}.
\end{Proof}

The next proposition contains the Gibbs variational principle for unbounded perturbations. The proof presented here is an adaptation to the unbounded setting of an argument due to \textsc{D. Petz} who established the following result for bounded perturbations \cite[Prop. 1 \& Cor. 2]{Petz88}. In light of the previous \cref{lem:perturbationTheory_generalPerturbedRelModOp}, introduce the following final class of \enquote{Gibbs perturbations}:
\begin{equation*}
    \MFS_3(\MFM, L, \Phi, \Psi) \ce \left\{V \in \MFS_1(\MFM, L, \Phi) \ : \ \parbox{0.4\textwidth}{\centering $\log(\Delta_{\Phi, \Psi}) - \beta V$ is essentially self-adjoint on $\dom(\log \Delta_{\Phi, \Psi}) \cap \dom(V)$}\right\} \ .
\end{equation*}

\begin{proposition}[Gibbs variational principle]\label{pro:perturbationTheory_unboundedGibbsVP}
    Let $\varphi, \psi \in \NS(\MFM)$ be two faithful normal states on $\MFM$, let $\Phi, \Psi \in \NPC$ be their vector representatives, and let $L = \log(\Delta_\Phi)$ and $V \in \MFS_3(\MFM, L, \Phi, \Psi)$. Then
    \begin{equation}\label{eq:perturbationTheory_generalizedPBinequality}
        - \inv{\beta} \log\bigl(\norm{\Phi_V}^2\bigr) \le \psi(V) + \inv{\beta} \MS_\MFM^\mathrm{std}(\psi, \varphi) \ ,
    \end{equation}
    and equality holds true if and only if $\psi = \omega_{\Phi_V} / \norm{\Phi_V}^2$. Moreover,
    \begin{equation}\label{eq:perturbationTheory_generalizedGibbsVP}
        - \inv{\beta} \log\bigl(\norm{\Phi_V}^2\bigr) = \inf_{\gamma \in \Xi(\MFM, V, \Phi)} \Bigl(\gamma(V) + \inv{\beta} \MS_\MFM^\mathrm{std}(\gamma, \varphi)\Bigr) \ ,
    \end{equation}
    where the infimum is taken over the set
    \begin{equation*}
    \Xi(\MFM, V, \Phi) \ce \set[\Big]{\omega \in \NS(\MFM) \ : \ \text{$\omega = \omega_\Omega$ faithful and $\log(\Delta_{\Phi, \Omega}) - \beta V$ ess. self-adjoint}} \ .
    \end{equation*}
\end{proposition}

\begin{Proof}
    \pp{1.} Write $\varphi^V \ce \omega_{\Phi_V} = \braket{\Phi_V, \,\sbullet\ \Phi_V}$,\footnotemark
    \footnotetext{The non-normalized perturbed functional $\varphi^V$ should not be confused with the \emph{state} $\omega_V$ considered in \cref{def:perturbationTheory_boundedPerturbedVector}.}
    and note that $\Phi_V$ is cyclic and separating for $\MFM$. (This follows from the fact that, as mentioned above, $\Phi_V$ is cyclic and separating for bounded $V$, hence it can be proved as in \cref{thm:perturbationTheory_unboundedPerturbation} \ref{enu:perturbationTheory_unboundedPerturbationC}.) Therefore, $\MS_\MFM^\mathrm{std}(\psi, \varphi^V)$ is given by the usual expression. From the assumption on the operator $\log(\Delta_{\Phi, \Psi}) - \beta V$ and \cref{lem:perturbationTheory_generalPerturbedRelModOp}, it follows that
    \begin{align*}
        \MS_\MFM^\mathrm{std}(\psi, \varphi^V) &= - \vbraket[\big]{\Psi, \log(\Delta_{\Phi_V, \Psi}) \Psi} = - \vbraket[\big]{\Psi, \bigl(\log(\Delta_{\Phi, \Psi}) - \beta V\bigr) \Psi} = \MS_\MFM^\mathrm{std}(\psi, \varphi) + \beta \psi(V) \ .
    \end{align*}
    (This identity is well-known for bounded perturbations \cite[Thm. 3.10]{Araki77}, \cite[Cor. 5.9]{Donald90}, \cite[Cor. 12.8]{OP04}.) Using the fundamental inequality for $\MS_\MFM^\mathrm{std}$ from \cref{pro:relativeEntropy_basicLowerBound}, one obtains
    \begin{equation}\label{eq:perturbationTheory_basicLowerBound}
        \MS_\MFM^\mathrm{std}(\psi, \varphi^V) \ge - \psi(\id_\MH) \log\left(\frac{\varphi^V\bigl(\ssupp(\psi)\bigr)}{\psi(\id_\MH)}\right) = - \log\bigl(\varphi^V(\id_\MH)\bigr) \ .
    \end{equation}
    Writing $\varphi^V(\id_\MH) = \norm{\Phi_V}^2$, the relation for the relative entropy $\MS_\MFM^\mathrm{std}(\psi, \varphi_V)$ derived above together with the previous inequality imply that
    \begin{equation*}
        - \log\bigl(\norm{\Phi_V}^2\bigr) \le \beta \psi(V) + \MS_\MFM^\mathrm{std}(\psi, \varphi) \ .
    \end{equation*}

    \pp{2.} For the statement regarding equality in \eqref{eq:perturbationTheory_generalizedPBinequality}, recall from \cref{cor:relativeEntropy_basicLowerBoundRevisited} that the right-hand side of \cref{eq:perturbationTheory_basicLowerBound} can be interpreted as the relative entropy $\MS_{\MFM_0}^\mathrm{std}(\psi|_{\MFM_0}, \varphi^V|_{\MFM_0})$ on the subalgebra $\MFM_0 = \set{\id_\MH}$, and hence the inequality \eqref{eq:perturbationTheory_basicLowerBound} is just a special instance of \cref{cor:relativeEntropy_monotonicitySubalgebra}. In light of this fact, it follows from a theorem of \textsc{D. Petz} \cite[Thm. 4]{Petz86b} that equality in \eqref{eq:perturbationTheory_basicLowerBound} is obtained if and only if
    \begin{equation*}
        \Delta_{\varphi^V, \psi}^{\ii t} \, \Delta_{\psi}^{- \ii t} = \Delta_{\varphi^V|_{\MFM_0}, \psi|_{\MFM_0}}^{\ii t} \, \Delta_{\psi|_{\MFM_0}}^{- \ii t} = \left(\frac{\varphi^V(\id_\MH)}{\psi(\id_\MH)}\right)^{\ii t} \quad (t \in \R) \ ,
    \end{equation*}
    where it was used that $\Delta_{\psi|_{\MFM_0}} = \id_\MH$ because $\psi|_{\MFM_0}$ is a tracial state (\cf{} \cref{exa:vonNeumann_modularData} \ref{enu:vonNeumann_exaModOpLInf} and the proof of \cref{cor:relativeEntropy_basicLowerBoundRevisited}), and that $\Delta_{\varphi^V|_{\MFM_0}, \psi|_{\MFM_0}} = \varphi^V(\id_\MH) / \psi(\id_\MH)$ according to \cref{lem:vonNeumann_relativeModularOpScaling}. (In \cite[Thm. 4]{Petz86b}, the first equality in the above equation is stated as an equality of \emph{Connes' cocycle derivatives}, \cf{} \cite[Sect. 6.2 \& 10.2]{Hiai21} for definitions, which was avoided here for simplicity.) It is noted in \cite[p. 346]{Petz88} that the above identity implies $\varphi^V = \lambda \psi$ for some $\lambda > 0$, namely $\lambda = \varphi^V(\id_\MH) / \psi(\id_\MH) = \varphi^V(\id_\MH)$. This shows that equality in \eqref{eq:perturbationTheory_generalizedPBinequality} is attained only for the state $\psi = \varphi^V / \norm{\Phi_V}^2$.

    \pp{3.} Finally, observe that on the one hand, the inequality \eqref{eq:perturbationTheory_generalizedPBinequality} implies that
    \begin{equation*}
        - \inv{\beta} \log\bigl(\norm{\Phi_V}^2\bigr) \le \inf_{\gamma \in \Xi(\MFM, V, \Phi)} \Bigl(\gamma(V) + \inv{\beta} \MS_\MFM^\mathrm{std}(\gamma, \varphi)\Bigr)
    \end{equation*}
    since by assumption on $\psi$ and $V$, it holds that $\psi \in \Xi(\MFM, V, \Phi)$. On the other hand, the fact that equality in \eqref{eq:perturbationTheory_generalizedPBinequality} is obtained if and only if $\psi$ takes the form $\wt{\psi} \ce \varphi^V / \norm{\Phi_V}^2$ shows that
    \begin{equation*}
        - \inv{\beta} \log\bigl(\norm{\Phi_V}^2\bigr) = \wt{\psi}(V) + \inv{\beta} \MS_\MFM^\mathrm{std}\bigl(\wt{\psi}, \varphi\bigr) \ge \inf_{\gamma \in \Xi(\MFM, V, \Phi)} \Bigl(\gamma(V) + \inv{\beta} \MS_\MFM^\mathrm{std}(\gamma, \varphi)\Bigr) \ .
    \end{equation*}
    Together, the previous two inequalities prove the variational principle in \cref{eq:perturbationTheory_generalizedGibbsVP}.
\end{Proof}

\begin{remark}\label{rem:perturbationTheory_unboundedGibbsVP}
    There is a different approach to perturbation theory in operator algebras, going back to \textsc{S. Sakai} \cite{Sakai87} and \textsc{M. J. Donald} \cite{Donald90}, which deals with semi-bounded perturbations which are modeled as so-called extended-valued self-adjoint operators, see \cite[Ch. 12]{OP04} and \cite[Sect. 7.1]{Hiai21}. In this framework, given a normal state $\varphi \in \NS(\MFM)$ and a generalized operator $V$, one \emph{defines} the perturbed state $\varphi^V$ as the unique minimizer of the Gibbs variational principle \eqref{eq:perturbationTheory_generalizedGibbsVP}; this definition is essentially based on the fact that $\rho_0 = \ee^{- \beta H_0} / \tr(\ee^{- \beta H_0})$ is the unique minimizer of the free energy according to the Gibbs variational principle \cite[Prop. 1.10]{OP04}.
    
    Thus, in this alternative approach to perturbation theory, the Gibbs variational principle is already satisfied by construction. In light of this, it is appropriate to remark that the argument presented above has the advantage of relying on the more constructive dynamical approach to perturbation theory presented in \cref{sec:perturbationTheory_perturbation,sec:perturbationTheory_perturbationKMS}, and that it does not employ the abstract concept of extended-valued operators.
\end{remark}

Next, a variational expression for the relative entropy in terms of bounded perturbations $V \in \SE{\MFM}$, also due to \textsc{D. Petz} \cite[Thm. 9]{Petz88}, shall be proved. This variational principle is the dual form of the Gibbs variational principle \eqref{eq:perturbationTheory_generalizedGibbsVP}, and it is known in the literature as the \emph{Donsker-Varadhan principle}. The proof of this result presented here uses methods from convex analysis and follows \cite[p. 228]{OP04}.

\begin{proposition}[Petz]\label{pro:perturbationTheory_PetzVP}
    Let $\varphi, \psi \in \NS(\MFM)$ be two faithful normal states on $\MFM$ with $\varphi$ given by $\Phi \in \NPC$. Then the Araki-Uhlmann relative entropy of $\psi$ and $\varphi$ can be written as
    \begin{align}\label{eq:perturbationTheory_PetzVP}
        \MS_\MFM^\mathrm{std}(\psi, \varphi) = \sup_{V \in \SE{\MFM}} \Bigl(- \beta \psi(V) - \log\bigl(\norm{\Phi_V}^2\bigr)\Bigr) \ .
    \end{align}
\end{proposition}

\begin{Proof}
    Let $E \ce \NS(\MFM)$ and $F \ce \SE{\MFM}$, and define the duality pairing $E \times F \owns (\omega, V) \longmto \braket{\omega, V} \ce - \beta \omega(V)$. Equip $E$ and $F$ with the weak topologies $\sigma(E, F)$ and $\sigma(F, E)$, respectively (\cf{} \cref{para:topologicalVectorSpaces_weakTopology}), and consider the function $g : E \longto \R \cup \{+ \infty\}$ given by $g(\omega) \ce \MS_\MFM^\mathrm{std}(\omega, \varphi)$. It holds that the \emph{Legendre-Fenchel conjugate} \cite[p. 255 f.]{NiculescuPersson18} $g^\ast : F \longto \R \cup \{+ \infty\}$ of $g$ is given for all $V \in F$ by
    \begin{align*}
        g^\ast(V) &\ce \sup_{\omega \in E} \Bigl(\braket{\omega, V} - g(\omega)\Bigr) = \sup_{\omega \in E} \Bigl(- \beta \omega(V) - \MS_\MFM^\mathrm{std}(\omega, \varphi)\Bigr) \\
        &= - \inf_{\omega \in E} \Bigl(\beta \omega(V) + \MS_\MFM^\mathrm{std}(\omega, \varphi)\Bigr) = \log\bigl(\norm{\Phi_V}^2\bigr) \ ,
    \end{align*}
    where in the last step the Gibbs variational principle \eqref{eq:perturbationTheory_generalizedGibbsVP} (see also \cite[Cor. 2]{Petz88}) was used. Now, since the function $g$ is lower semi-continuous and convex due to the respective properties of the relative entropy \cite[Thm. 3.7 \& 3.8]{Araki77}, \cite[Cor. 5.12]{OP04}, it follows from the Fenchel-Hörmander-Moreau theorem \cite[Thm. 6.1.2]{NiculescuPersson18} that $g$ is equal to its \emph{biconjugate} $\dadj{g}$, \ie{}, for all $\psi \in E$ there holds
    \begin{equation*}
        \MS_\MFM^\mathrm{std}(\psi, \varphi) = g(\psi) = \dadj{g}(\psi) \ce \sup_{V \in F} \Bigl(\braket{\psi, V} - g^\ast(V)\Bigr) = \sup_{V \in \SE{\MFM}} \Bigl(- \beta \psi(V) - \log\bigl(\norm{\Phi_V}^2\bigr)\Bigr) \ . \tag*{\qedhere}
    \end{equation*}
\end{Proof}

With the help of the previous proposition and the Gibbs variational principle, one can now establish the Donsker-Varadhan principle for unbounded perturbations.

\begin{proposition}[Donsker-Varadhan principle]\label{pro:perturabtionTheory_unboundedDVVP}
    Let $\varphi, \psi \in \NS(\MFM)$ be two faithful normal states on $\MFM$ given by $\Phi, \Psi \in \NPC$, and let $L = \log(\Delta_\Phi)$. Then
    \begin{equation}\label{eq:perturabtionTheory_unboundedDVVP}
        \inv{\beta} \MS_\MFM^\mathrm{std}(\psi, \varphi) = \sup_{V \in \MFS_3(\MFM, L, \Phi, \Psi)} \Bigl(- \psi(V) - \inv{\beta} \log\bigl(\norm{\Phi_V}^2\bigr)\Bigr) \ .
    \end{equation}
\end{proposition}

\begin{Proof}
    From the inequality \eqref{eq:perturbationTheory_generalizedPBinequality} in \cref{pro:perturbationTheory_unboundedGibbsVP}, it immediately follows that
    \begin{equation*}
        \inv{\beta} \MS_\MFM^\mathrm{std}(\psi, \varphi) \ge \sup_{V \in \MFS_3(\MFM, L, \Phi, \Psi)} \Bigl(- \psi(V) - \inv{\beta} \log\bigl(\norm{\Phi_V}^2\bigr)\Bigr) \ .
    \end{equation*}
    On the other hand, since $\SE{\MFM} \subset \MFS_3(\MFM, L, \Phi, \Psi)$, the variational expression \eqref{eq:perturbationTheory_PetzVP} from \cref{pro:perturbationTheory_PetzVP} for bounded perturbations shows that
    \begin{equation*}
        \inv{\beta} \MS_\MFM^\mathrm{std}(\psi, \varphi) = \sup_{V \in \SE{\MFM}} \Bigl(- \psi(V) - \inv{\beta} \log\bigl(\norm{\Phi_V}^2\bigr)\Bigr) \le \sup_{V \in \MFS_3(\MFM, L, \Phi, \Psi)} \Bigl(- \psi(V) - \inv{\beta} \log\bigl(\norm{\Phi_V}^2\bigr)\Bigr) \ .
    \end{equation*}
    Together, these two inequalities prove the assertion.
\end{Proof}

Combining the Gibbs variational principle from \cref{pro:perturbationTheory_unboundedGibbsVP} and the Donsker-Varadhan principle from the previous \cref{pro:perturabtionTheory_unboundedDVVP}, one obtains variational bounds for the relative free energy.

\begin{corollary}[Variational bounds for the relative free energy]\label{cor:perturbationTheory_unboundedVariationalBounds}
    Let $(\MFM, \MH, J, \MP)$ be a von Neumann algebra in standard form, let $\tau$ be a $W^\ast$-dynamics on $\MFM$ with standard Liouvillian $H_0$, and let $\omega_0$ be a faithful $(\tau, \beta)$-KMS-state given by $\Omega_0 \in \NPC$. For all $U \in \MFS_2(\MFM, H_0, \Omega_0)$, it holds that
    \begin{align}\label{eq:perturbationTheory_unboundedVariationalBounds}
        \begin{split}
            &\sup_{V \in \MFS_3(\MFM, H_0, \Omega_0, \Omega_U)} \Bigl(\omega_U(U - V) - \inv{\beta} \log\bigl(\norm{\Omega_V}^2\bigr)\Bigr) \\[4pt]
            &\qquad = \MF(\omega_U, \omega_0) = \inf_{\psi \in \Xi(\MFM, U, \Omega_0)} \Bigl(\psi(U) + \inv{\beta} \MS_\MFM^\mathrm{std}(\psi, \omega_0)\Bigr) \ .
        \end{split}
    \end{align}
    In particular, there is the following family of two-sided bounds for the relative free energy $\MF(\omega, \omega_0)$ which holds for all $V \in \MFS_3(\MFM, L, \Omega_0, \Omega_U)$ and $\psi \in \Xi(\MFM, U, \Omega_0)$:
    \begin{equation}\label{eq:perturbationTheory_unboundedFamilyBounds}
        \omega_U(U - V) - \inv{\beta} \log\bigl(\norm{\Omega_V}^2\bigr) \le \MF(\omega_U, \omega_0) \le \psi(U) + \inv{\beta} \MS_\MFM^\mathrm{std}(\psi, \omega_0) \ .
    \end{equation}
\end{corollary}

\begin{Proof}
    First, choosing $\varphi = \omega_0$ in \eqref{eq:perturbationTheory_generalizedGibbsVP} and noting that $\MF(\omega_U, \omega_0) = - \inv{\beta} \log \norm{\Omega_U}^2$ immediately yields the second line of \cref{eq:perturbationTheory_unboundedVariationalBounds}. Next, combining the definition of the relative free energy from \cref{eq:perturbationTheory_relativeFreeEnergy} with \eqref{eq:perturabtionTheory_unboundedDVVP} gives
    \begin{align*}
        \MF(\omega_U, \omega_0) &= \omega_U(U) + \sup_{V \in \MFS_3(\MFM, H_0, \Omega_0, \Omega_U)} \Bigl(- \omega_U(V) - \inv{\beta} \log\bigl(\norm{\Omega_V}^2\bigr)\Bigr) \\
        &= \sup_{V \in \MFS_3(\MFM, H_0, \Omega_0, \Omega_U)} \Bigl(\omega_U(U - V) - \inv{\beta} \log\bigl(\norm{\Omega_V}^2\bigr)\Bigr) \ . \tag*{\qedhere}
    \end{align*}
\end{Proof}

\begin{remark}
    The variational principles in \cref{eq:perturbationTheory_unboundedVariationalBounds} are the generalization of \cref{pro:perturbationTheory_variationalBoundsQM} to unbounded perturbations on a general von Neumann algebra. Note also that by choosing $V = 0$ and $\psi = \omega_0$ in \cref{eq:perturbationTheory_unboundedFamilyBounds}, one recovers the bounds of \cref{pro:perturbationTheory_BogoliubovVN}. Thus, the above corollary generalizes the two-sided Bogoliubov inequality \eqref{eq:perturbationTheory_BogoliubovVN}.
\end{remark}

\chapter{Conclusion}\label{ch:conclusion}

\section{Summary}

Utilizing the mathematically involved theory of operator algebras and, in particular, Tomita-Takesaki modular theory and its relatives, the abstract Araki-Uhlmann relative entropy of two normal functionals was defined for arbitrary von Neumann algebras. A selection of important properties of this functional were discussed, most prominently Uhlmann's monotonicity theorem which states that the relative entropy is non-increasing under Schwarz mappings between two von Neumann algebras. A complete and detailed proof of this assertion under different assumptions on the involved normal functionals and the Schwarz mapping was given, and certain insights from this proof were used to obtain further monotonicity inequalities for the relative entropy. Namely, the situation of two von Neumann algebras $\MFM_1 \subset \MFM_2 \subset \BO(\MH)$ and two vector functionals $\omega_\Omega, \omega_\Phi \in \NF{(\MFM_2)}$ was studied, and it was investigated under what Hilbert-space transformations $V \in \BO(\MH)$, applied to the vector representatives $\Omega, \Phi \in \MH$, the relative entropy is monotonic. A number of corresponding inequalities were obtained, all being derived from Uhlmann's ultimate monotonicity inequality, under varying assumptions on the operator $V$ and the vectors $\Omega$ and $V \Omega$, see \cref{pro:relativeEntropy_vectorMonotonicityIsometry,eq:relativeEntropy_vectorMonotonicityIsometry2,pro:relativeEntropy_vectorMonotonicityPartialIsometry,pro:relativeEntropy_vectorMonotonicityUnitary,pro:vectorMonotonicity_contraction}. In particular, special instances of certain well-known results for the relative entropy, \emph{viz.}, monotonicity under restriction of the functionals to subalgebras and invariance under $\ast$-automorphisms of the underlying algebra, were re-obtained in this framework, see \cref{cor:relativeEntropy_vectorMonotonicitySubalgebra,cor:relativeEntropy_invarianceAutomorphisms}.

Furthermore, perturbation theory in operator algebras was developed to some extend, focusing on unbounded perturbations of KMS-states. Results from this area obtained fairly recently in the literature \cite{DJP03} were presented and re-proved in great detail, and even some slight extensions were suggested (\cref{lem:perturbationTheory_unboundedConvergenceLiouvillian} \ref{enu:perturbationTheory_unbdConvergenceSum}, \cref{thm:perturbationTheory_unboundedPerturbation} \ref{enu:perturbationTheory_unboundedRelModOp2}, and \cref{pro:perturbationTheory_relEntUnboundedPerturbedState}). The culmination of the known methods from the literature as well as of the extensions are certain identities for the Araki-Uhlmann relative entropy between perturbed and unperturbed KMS-states; they were already used in the past to prove, \eg{}, the Golden-Thompson inequality in von Neumann algebras for bounded perturbations. Most importantly for the present text, this framework, and the relative entropy identities in particular, were employed to extend the two-sided Bogoliubov inequality, which was previously established for quantum systems (that is, \enquote{type I von Neumann algebras}), to unbounded perturbations of equilibrium states on general von Neumann algebras (\cref{pro:perturbationTheory_BogoliubovVN}). It was shown that these generalized bounds reproduce the known inequalities from \cite[Thm. 4.1]{Reible22} in case that the von Neumann algebra is $\BO(\MH)$ and the perturbation is some unbounded self-adjoint operator on $\MH$ (\cref{cor:perturbationTheory_boundedBogoliubovImpliesQMBogoliubov}). Moreover, under further technical assumptions and by extending certain well-known identities from bounded perturbation theory to unbounded perturbations (\cref{lem:perturbationTheory_generalUnboundedPerturbedVector,lem:perturbationTheory_generalPerturbedRelModOp}), it was possible to prove analoga of the Gibbs and Donsker-Varadhan variational principles for unbounded perturbations (\cref{pro:perturbationTheory_unboundedGibbsVP,pro:perturabtionTheory_unboundedDVVP}). This facilitated a generalization of the known variational bounds for the relative free energy \cite[Sect. 4.1]{Reible22} to unbounded perturbations on von Neumann algebras (\cref{cor:perturbationTheory_unboundedVariationalBounds}).

\section{Outlook}

Since the investigations in this text were conducted in a purely abstract mathematical setting, the main question that should be settled in future work is the applicability of the derived monotonicity inequalities to concrete situations, for example, to estimate relative entropies in certain models of algebraic quantum field theory, \cf{} \cref{rem:relativeEntropy_vectorMonotonicityAssumptions2} \ref{enu:relativeEntropy_vectorMonotonicityAQFT}. Theoretically, it is plausible that there are situations where one is given von Neumann algebras $\MFM_1 \subset \MFM_2 \subset \BO(\MH)$ and states $\omega_{V \Omega}, \omega_{V \Phi} \in \NF{(\MFM_1)}$, where $\Omega, \Phi \in \MH$ and $V \in \BO(\MH)$ with suitable properties, for which the relative entropy shall be computed. Under certain technical conditions, which can be realized in applications, the bounds of \cref{sec:relativeEntropy_vectorMonotonicity} then allow to estimate this relative entropy by the one between $\omega_{\Omega}$ and $\omega_{\Phi}$ on the larger algebra. Similarly, it would be interesting to examine situations in which the two-sided Bogoliubov inequality, \cf{} \cref{eq:perturbationTheory_BogoliubovVN,eq:perturbationTheory_unboundedFamilyBounds}, can be applied in quantum field theory. It has already been established in \cite{Reible23, DS24} that it is a useful inequality for certain applications in quantum statistical mechanics, hence one may investigate whether conceptually similar situations can be approached in the setting of quantum field theory.

Aside from this question, one may also further analyze the theoretical setup leading to the results of this text. Regarding the monotonicity inequalities, it would be interesting to examine whether one can prove inequalities of the form \eqref{eq:relativeEntropy_vectorMonotonicityIsometry} also for projection operators, or if one can relax the assumption in \cref{pro:relativeEntropy_vectorMonotonicityPartialIsometry} about the vector $\Omega$ being contained in the initial subspace of the partial isometry $V$. This would be a step towards monotonicity results for non-unital Schwarz mappings which are of interest in the context of operational quantum physics \cite{Davies1970,Kraus1971} and measurement theory in quantum field theory \cite{FewsterVerch20, FewsterVerch23}. There might be a connection to the problem raised in \cref{rem:relativeEntropy_vectorMonotonicityAssumptions2} \ref{enu:relativeEntropy_dilationSchwarz}, namely deriving a dilation theorem for Schwarz mappings between von Neumann algebras, in the spirit of the theorem of Stinespring for completely positive mappings between $C^\ast$-algebras. Furthermore, it might be interesting to investigate whether similar inequalities are valid if the operator $V$ is substituted by an unbounded operator, no longer an element of the von Neumann algebra, but affiliated with it. This would greatly increase the situations in which such estimates are applicable. Finally, one might try to eliminate using Uhlmann's theorem in the proofs of the inequalities of \cref{sec:relativeEntropy_vectorMonotonicity} and instead find independent proofs. (Possibly along similar lines as the discussion in \cite{Witten18}.) This would allow a broader audience, which does not necessarily have all the technical background from von Neumann algebra theory, to be introduced to monotonicity results for the Araki-Uhlmann relative entropy.

With regard to the two-sided Bogoliubov inequality and the related variational principles, one should investigate whether the variational bounds \eqref{eq:perturbationTheory_unboundedVariationalBounds} for the relative free energy can be proved under relaxed assumptions on the perturbation. In the present situation, in order to arrive at the final \cref{cor:perturbationTheory_unboundedVariationalBounds}, one has to take the supremum and infimum over somewhat unwieldy sets of operators and normal functionals; this was necessary in order for the inequality \eqref{eq:perturbationTheory_generalizedPBinequality} to be satisfied which was central to proving the variational principles. Ideally, one would hope to prove \cref{eq:perturbationTheory_unboundedVariationalBounds} for the case in which the supremum ranges over the set $\MFS_2(\MFM, H_0, \Omega_0)$ and the infimum over all normal states $\NS(\MFM)$ since these seem to be the most natural classes of operators and functionals in light of the variational bounds known from the quantum-mechanical case. In particular, it would be ideal to eliminate the technical assumption on the operator $\log(\Delta_{\Phi, \Psi}) - \beta V$ which first appeared in \cref{lem:perturbationTheory_generalPerturbedRelModOp}, and which then propagated to the variational principles in \cref{pro:perturbationTheory_unboundedGibbsVP,pro:perturabtionTheory_unboundedDVVP}, resulting in the use of the cumbersome classes of perturbations $\MFS_3(\MFM, L, \Phi, \Psi)$ and of normal states $\Xi(\MFM, V, \Phi)$. Similarly, one might try to dispense with the assumption of faithfulness of the involved normal states which would considerably increase generality of the results.

\appendix
\addtocontents{toc}{\protect\setcounter{tocdepth}{0}}

\chapter{Topological Vector Spaces}\label{app:topologicalVectorSpaces}

This first appendix collects some selected notions from general topology and the theory of locally convex spaces, in particular, the concept of weak topologies, which are used in different places throughout the main part of the text.

\section{Preliminaries from General Topology}

\begin{para}[Closure]\label{para:topologicalVectorSpaces_topology}
    (\cite[pp. 1 f.]{Voigt20})
    Let $X \neq \emptyset$ be a non-empty set and $\MFT \subset \MFP(X)$ be a topology on $X$. Denote the system of closed sets in $(X, \MFT)$ by $\MFC_\MFT \ce \set{C \subset X \ : \ X \setminus C \in \MFT}$. For an arbitrary subset $A \subset X$, one defines the \bemph{closure} of $A$ to be the set
    \begin{equation*}
        \clos_\MFT(A) \ce \bigcap \set{C \in \MFC_\MFT \ : \ A \subset C} \ .
    \end{equation*}
    The closure of $A$ is also commonly denoted by $\ol{A}^\MFT$ or $\ol{A}$. It is clear that the closure of $A$ is the smallest closed set containing $A$.
\end{para}

\begin{para}[Comparison of topologies]\label{para:topologicalVectorSpaces_comparisonTopologies}
    (\cite[Def. 2.2.5]{Waldmann14})
    Let $\MFT_1$ and $\MFT_2$ be two topologies on a set $X$. $\MFT_2$ is said to be \bemph{finer} than $\MFT_1$, and $\MFT_1$ is said to be \bemph{coarser} than $\MFT_2$, iff $\MFT_1 \subset \MFT_2$. For every subset $A \subset X$ and every topological space $(Y, \MFS)$, it holds that \cite[p. 432]{Werner18}
    \begin{equation*}
        \clos_{\MFT_2}(A) \subset \clos_{\MFT_1}(A) \tand C^0\bigl((X, \MFT_1), Y\bigr) \subset C^0\bigl((X, \MFT_2), Y\bigr) \ .
    \end{equation*}
    
    Let $\MST(X) \subset \MFP\bigl(\MFP(X)\bigr)$ be the collection of topologies on a set $X$. Define a binary relation $\le$ on $\MST(X)$ by setting $\MFT_1 \le \MFT_2$ iff $\MFT_1 \subset \MFT_2$. It holds that $(\MST(X), \le)$ is a partially ordered set.
\end{para}

\begin{definition}[\protect{Initial topology \cite[Def. 3.12]{Querenburg01}}]\label{def:topologicalVectorSpaces_initialTopology}
    Let $X$ be a non-empty set, let $I$ be an index set, and for each $i \in I$ let $(X_i, \MFT_i)$ be a topological space and $f_i : X \longto X_i$ be a map. A topology $\MFT$ on $X$ is called \bemph{initial topology} with respect to $(f_i)_{i \in I}$ iff it satisfies the following universal property:
    \begin{enumerate}[label=\normalfont(IT)]
        \item \label{enu:topologicalVectorSpaces_universalPropertyIT} For an arbitrary topological space $(Y, \MFS)$, a mapping $g : Y \longto X$ is continuous if and only if $f_i \circ g$ is continuous for every $i \in I$.
    \end{enumerate}
\end{definition}

The initial topology is characterized as follows.

\begin{ntheorem}[\protect{Characterization of the initial topology \cite[Thm. 3.13]{Querenburg01}}]\label{thm:topologicalVectorSpaces_initialTopology}
    Let $X \neq \emptyset$ be a set, let $I$ be an index set, and for each $i \in I$ let $(X_i, \MFT_i)$ be a topological space and $f_i : X \longto X_i$ be a map. There exists a unique initial topology $\MFT$ with respect to $(f_i)_{i \in I}$ on $X$, and it is the coarsest topology for which the family of maps $(f_i)_{i \in I}$ is continuous. In particular, a net $(x_k)_{k \in N}$ converges in $(X, \MFT)$ to $x \in X$ if and only if $f_i(x_k) \to f_i(x)$ in $(X_i, \MFT_i)$ for every $i \in I$.
\end{ntheorem}

\section{Linear and Weak Topologies}

\begin{definition}[Topological vector space]\label{def:topologicalVectorSpaces_topVS}
    (\cite[p. 4]{Voigt20})
    Let $E$ be a vector space over a field $\K$ and $\MFT$ be a topology on the underlying set $E$. $\MFT$ is called a \bemph{vector space topology} (or: \emph{linear topology}), and the pair $(E, \MFT)$ is called a \bemph{topological vector space}, iff the vector space operations $(x, y) \longmto x + y$ and $(\lambda, x) \longmto \lambda x$ are continuous with respect to $\MFT$. The (\bemph{topological}) \bemph{dual space} $\cdual{(E, \MFT)}$ of $E$ (short: $\cdual{E}$) is defined to be the space of all $\MFT$-continuous linear functionals on $E$.
\end{definition}

\begin{proposition}[\protect{\cite[Thm. 1.5]{Voigt20}}]\label{pro:topologicalVectorSpaces_ITTVS}
    Let $E$ be a vector space, let $\set{(E_i, \MFT_i)}_{i \in I}$ be a family of topological vector spaces, and for every $i \in I$ let $f_i : E \longto E_i$ be a linear map. Then the initial topology $\MFT$ on $E$ with respect to the family $(f_i)_{i \in I}$ is a linear topology.
\end{proposition}

\begin{proposition}[\protect{\cite[p. 7]{Voigt20}}]\label{pro:topologicalVectorSpaces_LCS}
    Let $E$ be a $\K$-vector space and $P$ be a set of semi-norms on $E$. Then the initial topology $\MFT_P$ on $E$ with respect to the mappings $(\id_p)_{p \in P}$, where $\id_p : E \longmto (E, p)$, $x \longmto x$, is a vector space topology on $E$ called the \emph{topology generated} by $P$. In particular, it holds that a net $(x_i)_{i \in I}$ in $E$ converges to $x \in E$ if and only if
    \begin{equation*}
        \forall p \in P \ : \ \lim_{i \in I} p(x_i - x) = 0 \ .
    \end{equation*}
\end{proposition}

\begin{para}[Weak topologies]\label{para:topologicalVectorSpaces_weakTopology}
    (\cite[p. 6]{Voigt20})
    A \bemph{dual pair} $\braket{E, F}$ consists of two $\K$-vector spaces $E$ and $F$ and a bilinear mapping $b \equiv \braket{\cdot, \cdot} : E \times F \longto \K$, called the \emph{dual pairing}. The map $b$ induces the two mappings
    \begin{align*}
        b_1 : E \longto \Hom(F, \K), \ x \longmto b_1(x) \ce \braket{x, \,\sbullet\,} \ , \\
        b_2 : F \longto \Hom(E, \K), \ y \longmto b_2(y) \ce \braket{\,\sbullet\,, y} \ .
    \end{align*}

    The \bemph{weak topology} $\sigma(E, F)$ on the space $E$ with respect to the dual pair $\braket{E, F}$ is defined to be the initial topology with respect to the family $(b_2(y))_{y \in F}$; the weak topology $\sigma(F, E)$ on $F$ is defined analogously as the initial topology with respect to $(b_1(x))_{x \in E}$. According to \cref{pro:topologicalVectorSpaces_ITTVS}, they both are vector space topologies.
\end{para}

\begin{proposition}[\protect{\cite[Thm. 1.8]{Voigt20}}]\label{pro:topologicalVectorSpaces_dualSpaceWeakTopology}
    Let $\braket{E, F}$ be a dual pair. Then $\cdual{\bigl(E, \sigma(E, F)\bigr)} = b_2(F)$.
\end{proposition}

\begin{example}\label{exa:topologicalVectorSpaces_theWeakTop}
    (\cite[Sect. V.1]{Conway85})
    Let $E$ be a vector space and $\Lambda \subset \Hom(E, \K)$ be a linear subspace of the algebraic dual space of $E$. Define a bilinear mapping $b : E \times \Lambda \longto \K$ by setting $b(x, \varphi) \ce \varphi(x)$. Then the weak topology $\sigma(E, \Lambda)$ on $E$ with respect to $\braket{E, \Lambda}$ is given as the initial topology of $(\varphi)_{\varphi \in \Lambda}$.
    
    Note that for every $\varphi \in \Lambda$, the functional $\abs{\varphi} : E \longto [0, + \infty)$, $x \longmto \abs{\varphi(x)}$, is a semi-norm on $E$. It holds that the topology $\MFT_{P_\Lambda}$ generated by the family $P_\Lambda \ce \{\abs{\varphi} \, : \, \varphi \in \Lambda\}$ and the weak topology $\sigma(E, \Lambda)$ coincide. Furthermore, by \cref{pro:topologicalVectorSpaces_dualSpaceWeakTopology},
    \begin{equation*}
        \cdual{\bigl(E, \sigma(E, \Lambda)\bigr)} = \Lambda \ .
    \end{equation*}
    
    Let $(E, \MFT)$ be a topological vector space. Then the topology $\sigma(E, \cdual{E})$ is called \bemph{the weak topology} of $E$. It is the coarsest linear topology on $E$ that yields the same dual space as $\MFT$.
\end{example}

\begin{para}[Weak-$\ast$ topology]\label{para:topologicalVectorSpaces_weakStar}
    (\cite[Sect. V.1]{Conway85})
    Let $(E, \MFT)$ be a topological vector space, let $x \in E$ be a point and let $\varphi \in \cdual{E}$ be a continuous linear functional. Consider the map $\iota(x) : \cdual{E} \longto \C$ given by $\iota(x)(\varphi) \ce \varphi(x)$. The weak topology on $\cdual{E}$ with respect to $\Lambda \ce \set{\iota(x) \, : \, x \in E} \subset \Hom(\cdual{E}, \K)$ and $\braket{\varphi, \iota(x)} \ce \iota(x)(\varphi)$ is called the \bemph{weak-$\ast$ topology} on $\cdual{E}$ and denoted by $\sigma(\cdual{E}, E)$. Note that according to \cref{exa:topologicalVectorSpaces_theWeakTop}, $\sigma(\cdual{E}, E)$ is generated by the family of semi-norms $\set{\abs{\iota(x)} \, : \, x \in E}$ on $\cdual{E}$.
\end{para}

\chapter{Hilbert Space Operators}\label{app:operators}

In this appendix, several definitions and propositions from (bounded and unbounded) operator theory on Hilbert space shall be collected as they are heavily used in the main part of this text. Throughout, let $\MH \equiv (\MH, \bdot)$ denote a complex Hilbert space, let $\ndot = \bdot^{1/2}$ be the norm on $\MH$, and let $\bigl(\BO(\MH), \ndot_\mop\bigr)$ be the Banach space of bounded linear operators on $\MH$ with the operator norm.

\section{Bounded Operators}

\begin{ntheorem}[\protect{Bounded linear extension \cite[Thm. II.1.5]{Werner18}}]\label{thm:operators_BLT}
    Let $D \subset E$ be a dense linear subspace of a normed space $E$, let $F$ be a Banach space, and let $T \in \BO(D, F)$ be a bounded linear operator. Then there exists a unique extension $\wh{T} \in \BO(E, F)$ of $T$, that is, a bounded linear operator such that $\wh{T}|_D = T$, which satisfies $\norm{\wh{T}}_\mop = \norm{T}_\mop$.
\end{ntheorem}

\begin{definition}[\protect{Trace-class and Hilbert-Schmidt operators}]\label{def:operators_NOHS}
    (\cite[Def. 26.1]{BlanchardBrüning15})
    Let $T \in \BO(\MH)$ be a bounded linear operator and $\abs{T} \ce (T^\ast T)^{1/2}$, defined via the functional calculus. $T$ is called a \bemph{trace-class operator} iff there exists an orthonormal basis $(e_i)_{i \in I} \subset \MH$ such that
    \begin{equation*}
        \norm{T}_\mathrm{tr} \ce \sum_{i \in I} \vbraket[\big]{e_i, \abs{T} e_i} = \sum_{i \in I} \vnorm[\big]{\abs{T}^{1/2} e_i}^2 < + \infty \ .
    \end{equation*}
    The set of all trace-class operators is denoted by $\NO(\MH)$. Similarly, $T$ is called a \bemph{Hilbert-Schmidt operator}, and the set of all Hilbert-Schmidt operators is denoted by $\HS(\MH)$, iff there exists an orthonormal basis $(e_i)_{i \in I} \subset \MH$ such that
    \begin{equation*}
        \norm{T}_\mathrm{HS} \ce \norm{T^\ast T}_\mathrm{tr}^{1/2} = \sumbra{\bigg}{\sum_{i \in I} \norm{Te_i}^2}^{1/2} = \sumbra{\bigg}{\sum_{i \in I} \braket{e_i, T^\ast T e_i}}^{1/2} < + \infty \ .
    \end{equation*}
    It is evident that $T$ is a Hilbert-Schmidt operator if and only if $T^\ast T$ is a trace-class operator. Furthermore, $\ndot_\mathrm{tr}$ and $\ndot_\mathrm{HS}$ are independent of the choice of the orthonormal basis \cite[Lem. 26.1]{BlanchardBrüning15}.
\end{definition}

\begin{para}[Trace and Hilbert-Schmidt inner product]\label{para:operators_traceHS}
    (\cite[Thm. 26.1 \& Cor. 26.3]{BlanchardBrüning15})
    On the space of trace-class operators, the \bemph{trace} $\tr : \NO(\MH) \longto \K$ is a well-defined bounded linear functional \cite[Cor. 26.3]{BlanchardBrüning15} which is given for \emph{any} orthonormal basis $(e_i)_{i \in I} \subset \MH$ by
    \begin{equation*}
        \tr(T) \ce \sum_{i \in I} \braket{e_i, T e_i} \com T \in \NO(\MH) \ .
    \end{equation*}
    Moreover, one can define an inner product on the space of Hilbert-Schmidt operators:
    \begin{equation}\label{eq:operators_HSinnerProduct}
        \braket{S, T}_\mathrm{HS} \ce \tr(S^\ast T) = \sum_{i \in I} \braket{S e_i, T e_i} \com S, T \in \HS(\MH) \ .
    \end{equation}
\end{para}

\begin{proposition}[\protect{\cite[Thm. 26.1 \& 26.2]{BlanchardBrüning15}}]\label{pro:operators_NOHS}
    $\bigl(\NO(\MH), \ndot_\mathrm{tr}\bigr)$ and $\bigl(\HS(\MH), \ndot_\mathrm{HS}\bigr)$ are normed vector spaces with $\norm{T^\ast}_\mathrm{tr} = \norm{T}_\mathrm{tr}$ and $\norm{S^\ast}_\mathrm{HS} = \norm{S}_\mathrm{HS}$ for all $T \in \NO(\MH)$ and $S \in \HS(\MH)$, and they are two-sided $\ast$-ideals in $\BO(\MH)$. Furthermore, $\bigl(\NO(\MH), \ndot_\mathrm{tr}\bigr)$ is a Banach space and $\bigl(\HS(\MH), \bdot_\mathrm{HS}\bigr)$ is a Hilbert space.
\end{proposition}

\begin{definition}[Density matrix]\label{def:operators_densityMatrix}
    (\cite[Def. 26.2]{BlanchardBrüning15})
    A trace-class operator $\rho \in \NO(\MH)$ is called a \bemph{density matrix} iff $\rho$ is self-adjoint ($\rho^\ast = \rho$), positive ($\braket{\xi, \rho \xi} \ge 0$ for all $\xi \in \MH)$ and normalized ($\norm{\rho}_\mathrm{tr} = \tr(\rho) = 1$). The space of all density matrices on $\MH$ is denoted by $\DM(\MH)$.
\end{definition}

\begin{definition}[Orthogonal projection]\label{def:operators_orthogonalProjections}
    An operator $P \in \BO(\MH)$ is called an \bemph{orthogonal projection} (or: \emph{projection} for short) iff $P^2 = P$ and $\ker(P) \perp \ran(P)$, equivalently, if $P^2 = P$ and $P^\ast = P$ \cite[Thm. V.5.9]{Werner18}. The set of all orthogonal projections on $\MH$ is denoted by $\PO(\MH)$.
\end{definition}

\begin{lemma}[Range of a projection]\label{lem:operators_rangeProjection}
    Let $P \in \PO(\MH)$. Then $\eta \in \ran(P)$ if and only if $P \eta = \eta$.
\end{lemma}

\begin{Proof}
    If $\eta \in \ran(P)$, then there exists $\xi \in \MH$ such that $\eta = P \xi$. Since $P$ is idempotent, it follows that $P \eta = P^2 \xi = P \xi = \eta$. The converse implication is clear.
\end{Proof}

\begin{ntheorem}[\protect{of the orthogonal projection \cite[Thm. V.3.4]{Werner18}}]\label{thm:operators_orthogonalProjection}
    Let $U \subset \MH$ be a closed non-empty subspace. Then there exists a unique non-zero orthogonal projection $P = P_U \in \PO(\MH) \setminus \{0\}$ such that $\ran(P) = U$ and $\ker(P) = U^\perp$. Moreover, the operator $Q \ce \id_\MH - P$ is also an orthogonal projection, and its range is given by $\ran(Q) = \ran(P)^\perp$. Finally, the Hilbert space decomposes into an orthogonal direct sum $\MH = U \oplus U^\perp$.
\end{ntheorem}

\begin{remark}\label{rem:operators_orthogonalProjectionClosedSubspace}
    It follows from \cref{thm:operators_orthogonalProjection} that the set $\PO(\MH)$ and the set of closed subspaces of $\MH$ are in one-to-one correspondence by identifying every orthogonal projection $P$ with the closed subspace $U \ce \ran(P) \subset \MH$, and \emph{vice versa} every closed subspace with the unique orthogonal projection $P_U$ onto it \cite[Prop. 23.1 (d)]{BlanchardBrüning15}. This justifies using the following notation for any subspace $U \subset \MH$:
    \begin{equation*}
        [U] \ce P_{\clos_{\ndot}(U)} \in \PO(\MH) \ .
    \end{equation*}
\end{remark}

\begin{corollary}\label{cor:operators_characterizationPerp}
    Let $U \subset \MH$ be a subspace.
    \begin{enumerate}
        \item \label{enu:operators_closurePerp} $\clos_{\ndot}(U) = (U^\perp)^\perp$.
        
        \item \label{enu:operators_densePerp} The following assertion are equivalent:
        \begin{enumerate}[equiv]
            \item $U$ lies dense in $\MH$;
            \item $U^\perp = \{0\}$;
            \item $[U] = \id_\MH$.
        \end{enumerate}
    \end{enumerate}
\end{corollary}

\begin{Proof}
    Assertion \ref{enu:operators_closurePerp} is a well-known consequence of the theorem of the orthogonal projection \cite[Cor. V.3.5]{Werner18}. The first equivalence (i) $\iff$ (ii) in \ref{enu:operators_densePerp} follows from  \ref{enu:operators_closurePerp} by noting that $\MH = \{0\}^\perp$, and the equivalence (i) $\iff$ (iii) can be seen by using that $\eta = [U] \eta$ iff $\eta \in \ran([U])$ according to \cref{lem:operators_rangeProjection}.
\end{Proof}

\begin{lemma}\label{lem:operators_projectionRestrictedID}
    Let $U \subset \MH$ be a closed subspace and $P \in \PO(\MH)$ be an orthogonal projection such that $P|_U = \id_\MH$. If $U^\perp \subset \ker(P)$, then $P = [U]$.
\end{lemma}

\begin{Proof}
    The identity $P|_U = \id_\MH$ implies that $U \subset \ran(P)$ by \cref{lem:operators_rangeProjection}. Since $\ker(P) = \ran(P)^\perp$ (\cref{thm:operators_orthogonalProjection}), it follows from the assumption that $U^\perp \subset \ran(P)^\perp$.
    
    Let $\eta \in \ran(P)$ be arbitrary. Then the previous observation implies that $\eta \perp \xi$ for all $\xi \in U^\perp$, hence $\eta \in (U^\perp)^\perp = \clos_{\ndot}(U) = U$ by \cref{cor:operators_characterizationPerp} \ref{enu:operators_closurePerp} and closedness of $U$. Thus, $\ran(P) = U$ which is equivalent to $P = [U]$.
\end{Proof}

\begin{para}[Order relation on $\PO(\MH)$]\label{para:operators_orderPO}
    (\cite[p. 57]{Moretti19})
    For orthogonal projections $P, Q \in \PO(\MH)$, define the relation $P \le Q$ to hold true if and only if $\ran(P) \subset \ran(Q)$. This gives a partial order on $\MSP(\MH)$ due to the respective properties of the set-theoretic inclusion, \cf{} \cref{rem:operators_orthogonalProjectionClosedSubspace}.
\end{para}

\begin{para}[Support projections]\label{para:operators_supportPO}
    (\cite[Sect. 2.13]{StratilaZsido19})
    Let $T \in \BO(\MH)$. The following projections related to $T$ are of interest:
    \begin{equation*}
        \nsupp(T) \ce \bigl[\ker(T)\bigr] \com \lsupp(T) \ce \bigl[\clos_{\ndot} \ran(T)\bigr] \tand \rsupp(T) \ce \id_\MH - \nsupp(T) = \bigl[\ker(T)^\perp\bigr] \ .
    \end{equation*}
    Since $\ker(T) = \ran(T^\ast)^\perp$ \cite[Thm. V.5.2 (g)]{Werner18}, it follows that $\rsupp(T) = \lsupp(T^\ast)$, and in case that $T$ is self-adjoint, the \bemph{support} (\bemph{projection}) of $T$ is defined to be
    \begin{equation*}
        \ssupp(T) \ce \lsupp(T) = \rsupp(T) \ .
    \end{equation*}
\end{para}

\begin{definition}[Partial isometry]\label{def:operators_partialIsometry}
    (\cite[p. 137]{Schmüdgen12}, \cite[p. 27]{StratilaZsido19})
    Let $\MH_1$ and $\MH_2$ be two Hilbert spaces, and let $\MG_1 \subset \MH_1$ and $\MG_2 \subset \MH_2$ be closed subspaces. A linear operator $V : \MH_1 \longto \MH_2$ that maps the space $\MG_1$ isometrically onto the space $\MG_2$ and annihilates $\MG_1^\perp$ is called a \bemph{partial isometry}. The space $\MG_1$ is called the \bemph{initial space} of $V$, and $\MG_2$ is called the \bemph{final space} of $V$. It holds that $\rsupp(V) = V^\ast V$ is the projection onto $\MG_1$ and $\lsupp(V) = V V^\ast$ is the projection onto $\MG_2$.
\end{definition}

\section{Closed and Closable Operators}

\begin{definition}\label{def:operators_closed}
    (\cite[Def. 1.3 \& 1.5]{Schmüdgen12})
	Let $T : \MH \supset \dom(T) \longto \MH$ be a linear operator. $T$ is called \bemph{closed} iff its graph $\gr(T) \ce \set{(\xi, T \xi) \, : \, \xi \in \dom(T)}$ is closed in $\MH \oplus \MH$, and it is called \bemph{closable} iff there exists a closed linear operator $S : \MH \supset \dom(S) \longto \MH$ such that $S \supset T$, \ie{}, $\dom(S) \supset \dom(T)$ and $S \xi = T \xi$ for all $\xi \in \dom(T)$. A linear subspace $\MD \subset \dom(T)$ is called a \bemph{core} for the operator $T$ iff for every $\xi \in \dom(T)$, there exists a sequence $(\xi_n)_{n \in \N} \subset \MD$ such that $\xi_n \to \xi$ and $T \xi_n \to T \xi$ in $\MH$.
\end{definition}

\begin{para}[Closure of an operator]\label{para:operators_closure}
    (\cite[pp. 6 f.]{Schmüdgen12})
    Let $T$ be a closable operator and $\overline{T} : \MH \supset \dom(\ol{T}) \longto \MH$ be the closed operator defined by the relation $\gr(\ol{T}) = \ol{\gr(T)}$. This operator, called the \bemph{closure} of $T$, is the smallest closed extension of $T$ with respect to the operator inclusion $\subset$. By construction, $\dom(\ol{T})$ consists of all those $\xi \in \MH$ for which there exists a sequence $(\xi_n)_{n \in \N} \subset \dom(T)$ such that $\xi = \lim_{n \to \infty} \xi_n$ and $\lim_{n \to \infty} T \xi_n$ exists in $\MH$; one defines $\ol{T} \xi \ce \lim_{n \to \infty} T \xi_n$.
\end{para}

\begin{para}[Adjoint operators]
    (\cite[p. 8]{Schmüdgen12})
    Let $T : \MH \supset \dom(T) \longto \MH$ be a densely defined linear operator, \ie{}, $\dom(T) \subset \MH$ lies dense. The \bemph{adjoint operator} $T^\ast : \MH \supset \dom(T^\ast) \longto \MH$ of $T$ is defined by
    \begin{align*}
        \dom(T^\ast) &\ce \set[\Big]{\xi \in \MH \ : \ \exists \zeta \in \MH \ \forall \eta \in \dom(T) : \braket{\xi, T \eta} = \braket{\zeta, \eta}} \ ,\\[4pt]
        T^\ast \xi &\ce \zeta \tfor \xi \in \dom(T^\ast) \ .
    \end{align*}
\end{para}

\begin{para}[Anti-linear operators]\label{para:operators_antiLinear}
    (\cite[Sect. 9.35]{StratilaZsido19})
    An \bemph{anti-linear operator} is a mapping $T : \MH \supset \dom(T) \longto \MH$ which satisfies for all $\xi, \eta \in \dom(T)$ and $\alpha, \beta \in \C$
    \begin{equation*}
        T(\alpha \xi + \beta \eta) = \ol{\alpha} \, T \xi + \ol{\beta} \, T \eta \ .
    \end{equation*}
    The \bemph{adjoint operator} of $T$ is the anti-linear operator $T^\ast : \MH \supset \dom(T^\ast) \longto \MH$ defined by
    \begin{align*}
        \dom(T^\ast) &\ce \set[\Big]{\xi \in \MH \ : \ \exists \zeta \in \MH \ \forall \eta \in \dom(T) : \braket{x, T \eta} = \braket{\eta, \zeta}} \ ,\\[4pt]
        T^\ast \xi &\ce \zeta \tfor \xi \in \dom(T^\ast) \ .
    \end{align*}
\end{para}

\begin{proposition}[Properties of $T^\ast$ and $\ol{T}$ \protect{\cite[Prop. 1.6 \& Thm. 1.8]{Schmüdgen12}}]\label{pro:operators_propertiesAdjoint}
    Let $T : \MH \supset \dom(T) \longto \MH$ be a densely defined linear operator. The following properties obtain:
    \begin{enumerate}
        \item \label{enu:operators_adjointClosed} $T^\ast : \MH \supset \dom(T^\ast) \longto \MH$ is a closed linear operator.
        
        \item \label{enu:operators_ranTperp} $\ran(T)^\perp = \ker(T^\ast)$.
        
        \item \label{enu:operators_closability} $T$ is closable if and only if $T^\ast$ is densely defined.
        
        \item \label{enu:operators_adjointClosure} If $T$ is closable, then $(\,\ol{T}\,)^\ast = T^\ast$ and $\ol{T} = \dadj{T}$

        \item \label{enu:operators_inverseAdjoint} Suppose that $\ker(T) = \set{0}$ and that $\ran(T)$ lies dense in $\MH_2$. Then $T^\ast$ is invertible with $(T^\ast)^{-1} = (T^{-1})^\ast$.
        
        \item \label{enu:operators_inverseClosure} Assume that $T$ is closable and injective. Then $\inv{T}$ is closable if and only if $\ker(\ol{T}) = \{0\}$. In this case, there holds $(\,\ol{T}\,)^{-1} = \ol{(\inv{T})}$.
    \end{enumerate}
\end{proposition}

\begin{ntheorem}[\protect{Polar decomposition of closed operators \cite[Thm. 7.2]{Schmüdgen12}}]\label{thm:operators_polarDecompositionCO}
    Let $T : \MH \supset \dom(T) \longto \MH$ be a densely defined, closed linear operator. Then there exists a partial isometry $U_T : \MH \longto \MH$ with initial space $\ker(T)^\perp = \clos_{\ndot} \ran(T^\ast) = \clos_{\ndot} \ran(\abs{T})$, where $\abs{T} \ce (T^\ast T)^{1/2}$, and final space $\ker(T^\ast)^\perp = \clos_{\ndot} \ran(T)$ such that
    \begin{equation*}
        T = U_T \abs{T} \ .
    \end{equation*}
    This decomposition is unique. Moreover, if $T$ is an anti-linear operator, then the partial isometry $U_T$ is anti-linear as well \cite[p. 263]{StratilaZsido19}.
\end{ntheorem}

\section{Self-adjoint Operators}

\begin{ntheorem}[\protect{Spectral thereom for self-adjoint operators \cite[Thm. 5.7]{Schmüdgen12}}]
    Let $T : \MH \supset \dom(T) \longto \MH$ be a self-adjoint operator. Then there exists a unique spectral measure $E_T : \BO(\R) \longto \PO(\MH)$ on the Borel $\sigma$-algebra $\BO(\R)$ such that
    \begin{equation*}
        T = \int\nolimits_\R \lambda \diff E_T(\lambda) \ .
    \end{equation*}
\end{ntheorem}

\begin{definition}[Spectral integral]\label{def:operators_spectralIntegral}
    (\cite[Sect. 4.3 \& 5.3]{Schmüdgen12})
    Let $T$ be a self-adjoint operator on $\MH$ with spectral measure $E_T$, and let $f : \R \longto \C \cup \set{+ \infty}$ be an $E_T$-almost everywhere finite Borel-measurable function. The \bemph{spectral integral} of $f$ with respect to $E_T$ is denoted by
    \begin{equation*}
        \I_T(f) \equiv f(T) \ce \int\nolimits_\R f(\lambda) \diff E_T(\lambda) \ .
    \end{equation*}
    For $\xi, \eta \in \MH$, let $\mu_{\xi, \eta}^T : \MB(\R) \longto [0, + \infty]$ be the complex measure $\mu_{\xi, \eta}^T(A) \ce \braket{\xi, E_T(A) \eta}$ on $\MB(\R)$, and set $\mu_\xi^T(A) \ce \mu_{\xi, \xi}^T(A)$. Then the domain of the spectral integral $\I_T(f)$ is given by
    \begin{equation*}
        \dom\bigl(\I_T(f)\bigr) \ce \set[\bigg]{\xi \in \MH \ : \ \int\nolimits_\R \abs{f(\lambda)}^2 \diff \mu_{\xi}^T(\lambda) < + \infty} \ .
    \end{equation*}
\end{definition}

\begin{proposition}[\protect{Properties of the functional calculus \cite[Thm. 5.9]{Schmüdgen12}}]\label{pro:operators_functionalCalculus}
    Let $T$ be a self-adjoint operator on $\MH$, let $f, g : \R \longto \C \cup \set{+ \infty}$ be two $E_T$-almost everywhere finite Borel-measurable functions, and let $\xi, \eta \in \dom\bigl(\I_T(f)\bigr)$. Then the following properties hold true:
    \begin{enumerate}
        \item \label{enu:operators_funcCalcIP} $\displaystyle \braket{\xi, f(T) \eta} = \int\nolimits_\R f(\lambda) \diff \mu_{\xi, \eta}^T(\lambda)$.
        
        \item \label{enu:operators_funcCalcNorm} $\displaystyle \norm{f(T) \xi}^2 = \int\nolimits_\R \abs{f(\lambda)}^2 \diff \mu_\xi^T(\lambda)$.
        
        \item \label{enu:operators_funcCalcCommutation} $(fg)(T) = \ol{f(T) g(T)}$. If at least one of the functions is bounded, then $(fg)(T) = f(T) g(T)$ \cite[p. 231]{StratilaZsido19}, and in this case, the operators $f(T)$ and $g(T)$ commute with each other.
        
        \item \label{enu:operators_spectralMeasureIndicatorFct} $\1_A(T) = E_T(A)$ for all $A \in \MB(\R)$.
    \end{enumerate}
\end{proposition}

\begin{proposition}[\protect{\cite[Prop. 5.10]{Schmüdgen12}}]\label{pro:operators_eigenvalueSpectralMeasure}
    Let $T$ be a self-adjoint operator with spectral measure $E_T$. A number $\lambda \in \R$ is an eigenvalue of $T$ if and only if $E_T(\set{\lambda}) \neq 0$, and in this case $E_T(\set{\lambda}) \in \PO(\MH)$ is the orthogonal projection onto the eigenspace $\Eig(T, \lambda)$.
\end{proposition}

\begin{lemma}\label{lem:operators_funcCalcUnitaryConj}
    Let $T$ be a self-adjoint operator on $\MH$ and $U \in \UO(\MH)$ be a unitary operator. For a Borel measurable function $\varphi : \sigma(T) \longto \C$, it holds that
    \begin{equation*}
        \varphi(U T U^\ast) = U \varphi(T) U^\ast \ .
    \end{equation*}
\end{lemma}

\begin{Proof}
    By the multiplication operator form of the spectral theorem \cite[Thm. 16.5]{Kaballo14}, \cite[Thm. VIII.4]{RS1}, there exists a measure space $(X, \Sigma, \nu)$, a measurable function $g : X \longto \C$, and a unitary operator $V : \MH \longto L^2(X, \nu)$ such that $T = V^\ast M_g V$. Define the unitary $W \ce V U^\ast : \MH \longto L^2(X, \nu)$. One can then write $U T U^\ast = W^\ast M_g W$. From the functional calculus \cite[p. 426]{Kaballo14}, it now follows that
    \begin{equation*}
        \varphi(U T U^\ast) = \varphi(W^\ast M_g W) = W^\ast M_{\varphi \circ g} W = U V^\ast M_{\varphi \circ g} V U^\ast = U \varphi(T) U^\ast \ . \tag*{\qedhere}
    \end{equation*}
\end{Proof}

\begin{lemma}\label{lem:operators_funcCalcEigenvalue}
    Let $T$ be a self-adjoint operator on $\MH$ and $\Omega \in \Eig(T, \lambda)$ be an eigenvector of $T$ corresponding to the eigenvalue $\lambda \in \R$, \ie{}, $T \Omega = \lambda \Omega$. Then for every $E_T$-almost everywhere finite Borel-measurable function $f : \R \longto \C \cup \set{+ \infty}$, it holds that $f(T) \Omega = f(\lambda) \Omega$.
\end{lemma}

\begin{Proof}
    First, note that from \cref{lem:operators_rangeProjection} and \cref{pro:operators_eigenvalueSpectralMeasure}, it follows that $E_T(\{\lambda\}) \Omega = \Omega$. Next, recall that if $A \cap B = \emptyset$ for $A, B \in \MB(\R)$, then $E_T(A \cup B) = E_T(A) + E_T(B)$ and $E_T(A) E_T(B) = E_T(B) E_T(A) = 0$ \cite[Def. 4.2 \& Lem. 4.3]{Schmüdgen12}. Hence, if $\lambda \in A$ for some $A \in \MB(\R)$, then one obtains
    \begin{align*}
        E_T(A) \Omega = E_T\bigl(A \setminus \{\lambda\}\bigr) \Omega + E_T(\{\lambda\}) \Omega = E_T(A \setminus \{\lambda\}) E_T(\{\lambda\}) \Omega + \Omega = \Omega \ .
    \end{align*}
    Similarly, if $\lambda \notin A$, it holds that $E_T(A) \Omega = E_T(A) E_T(\{\lambda\}) \Omega = 0$.

    Let $\xi \in \MH$ be arbitrary and consider the complex measure $A \longmto \braket{\xi, E_T(A) \Omega} = \mu_{\xi, \Omega}^T(A)$ on $\MB(\R)$. By the previous two observations, it holds that $\mu_{\xi, \Omega}^T(A) = \braket{\xi, \Omega} \, \delta_\lambda(A)$, where $A \longmto \delta_\lambda(A) = \1_{A}(\lambda)$ is the Dirac measure at $\lambda$. Using \cref{pro:operators_functionalCalculus} \ref{enu:operators_funcCalcIP}, one obtains
    \begin{equation*}
        \braket{\xi, f(T) \Omega} = \int_\R f(t) \diff \braket{\xi, E_T(t) \Omega} = \int_\R f(t) \, \braket{\xi, \Omega} \diff \delta_\lambda(t) = \braket{\xi, f(\lambda) \Omega} \ .
    \end{equation*}
    Since $\xi \in \MH$ was arbitrary, this shows that $f(T) \Omega = f(\lambda) \Omega$.
\end{Proof}

\begin{ntheorem}[\protect{Trotter product formula \cite[Thm. VIII.31]{RS1}, \cite[Thm. A.1]{DJP03}}]\label{thm:operators_TrotterProductFormula}
    Let $A$ and $B$ be two self-adjoint operators on $\MH$, and assume that the operator sum $A + B$ is essentially self-adjoint on $\dom(A) \cap \dom(B)$. Then the following limit exists in the strong operator topology on $\BO(\MH)$:
    \begin{equation}\label{eq:operators_TrotterProductFormula}
        \ee^{\ii t (A + B)} = \dlim{so}{n \to \infty} \bigl(\ee^{\ii t A / n} \, \ee^{\ii t B / n}\bigr)^n \ .
    \end{equation}
\end{ntheorem}

\begin{definition}[Strong resolvent convergence]\label{def:operators_SRConvergence}
    (\cite[Def. 10.1.1]{Oliveira09})
    Let $T_n$ ($n \in \N$) and $T$ be self-adjoint operators on $\MH$. One says that the sequence $(T_n)_{n \in \N}$ converges to $T$ in the \bemph{strong resolvent sense} iff $(T_n - \ii)^{-1} \to (T - \ii)^{-1}$ in the strong operator topology, that is,
    \begin{equation*}
        \forall \xi \in \MH \ : \ \lim_{n \to \infty} \vnorm[\big]{(T_n - \ii)^{-1} \xi - (T - \ii)^{-1} \xi} = 0 \ .
    \end{equation*}
\end{definition}

\begin{proposition}[\protect{\cite[Prop. 10.1.8 \& 10.1.9]{Oliveira09}, \cite[Thm. VIII.20 \& VIII.21]{RS1}}]\label{pro:operators_characterizationSRConvergence}
    The following assertions are equivalent:
    \begin{enumerate}[equiv]
        \item $T_n \to T$ in the strong resolvent sense;
        \item $\ee^{- \ii t T_n} \to \ee^{- \ii t T}$ strongly for all $t \in \R$;
        \item $f(T_n) \to f(T)$ strongly for all bounded continuous $f : \R \longto \C$.
    \end{enumerate}
\end{proposition}

\begin{proposition}[\protect{\cite[Prop. 10.1.13]{Oliveira09}}]\label{pro:operators_SOimpliesSR}
    Let $T_n, T \in \BO(\MH)$ be self-adjoint. If $T_n \to T$ in the strong operator topology, then $T_n \to T$ in the strong resolvent sense.
\end{proposition}

\begin{proposition}[\protect{\cite[Prop. 10.1.18]{Oliveira09}, \cite[Thm. VIII.25]{RS1}, \cite[Prop. A.3]{DJP03}}]\label{pro:operators_coreSOimpliesSR}
    Let $T$ and $(T_n)_{n \in \N}$ be self-adjoint operators on a Hilbert space $\MH$. Suppose that $\MD \subset \MH$ is a subspace contained in $\dom(T)$ and $\dom(T_n)$ for all $n \in \N$, and that $T$ is essentially self-adjoint on $\MD$. Furthermore, assume that $\lim_{n \to \infty} T_n \xi = T \xi$ for all $\xi \in \MD$. Then it follows that $T_n \to T$ in the strong resolvent sense.
\end{proposition}

\begin{proposition}[\protect{\cite[Prop. A.4]{DJP03}}]\label{pro:operators_SRimpliesW}
    Let $(T_n)_{n \in \N}$ be a sequence of self-adjoint operators on a Hilbert space $\MH$, and assume that there exists another self-adjoint operator $T$ such that $T_n \to T$ in the strong resolvent sense. Furthermore, suppose that $(\Omega_n)_{n \in \N} \subset \MH$ and $\Omega \in \MH$ are vectors such that $\Omega_n \to \Omega$ weakly, and that there exists a constant $C \ge 0$ such that for all $n \in \N$, $\norm{T_n \Omega_n} \le C$. Then $\Omega \in \dom(T)$, $\dlim{w}{n \to \infty} T_n \Omega_n$ exists, and $T \Omega = \dlim{w}{n \to \infty} T_n \Omega_n$.
\end{proposition}

\begin{remark}\label{rem:operators_SRimpliesW}
    (\cite[p. 482]{DJP03})
    Suppose that $\dlim{w}{n \to \infty} T_n \Omega_n$ exists. Then it follows from the uniform boundedness principle that $(T_n \Omega_n)_{n \in \N} \subset \MH$ is bounded \cite[Cor. IV.2.3]{Werner18}, \ie{}, there is a constant $C \ge 0$ such that $\norm{T_n \Omega_n} \le C$ for all $n \in \N$. Therefore, one may replace the assumption $\norm{T_n \Omega_n} \le C$ in the above proposition by the existence of $\dlim{w}{n \to \infty} T_n \Omega_n$.
\end{remark}

\begin{lemma}\label{lem:operators_SRimpliesFSR}
    Let $(T_n)_{n \in \N}$ be a sequence of self-adjoint operators on $\MH$ converging in the strong resolvent sense to a self-adjoint operator $T$. Then for any continuous real-valued function $g : \R \longto \R$, it holds that $g(T_n) \to g(T)$ as $n \to + \infty$ in the strong resolvent sense.\footnote{The author wishes to thank Professor Dr.~Jan Derezi\'{n}ski for clarifying this result and for discussing its proof.}
\end{lemma}

\begin{Proof}
    Let $z \in \C \setminus \R$ be arbitrary, and observe that the function $h_z : \R \longto \C$, $\lambda \longmto \bigl(g(\lambda) - z\bigr)^{-1}$, is continuous and bounded. Indeed, since clearly $\abs{g(\lambda) - z} \ge \abs{\Im(z)}$, one obtains
    \begin{equation*}
        \abs{h_z(\lambda)} \le \frac{1}{\abs{\Im(z)}} \ .
    \end{equation*}
    Therefore, since $T_n \to T$ in the strong resolvent sense, it follows from \cref{pro:operators_characterizationSRConvergence} that $h_z(T_n) \to h_z(T)$ as $n \to + \infty$ in the strong operator topology. But $h_\ii(T_n) = (g(T_n) - \ii)^{-1}$ and $h_\ii(T) = (g(T) - \ii)^{-1}$, hence the previous statement is equivalent to convergence $g(T_n) \to g(T)$ in the strong resolvent sense.
\end{Proof}

\chapter{Quadratic Forms and Self-adjoint Operators}\label{app:forms}

In this final appendix, some definitions and results from the theory of positive quadratic forms will be outlined. Again, $\MH \equiv (\MH, \braket{\cdot, \cdot})$ shall denote a complex Hilbert space, and $\ndot = \bdot^{1/2}$ the norm induced by the inner product.

\section{Closed and Closable Quadratic Forms}

\begin{definition}[Quadratic form]\label{def:forms_quadraticForm}
    (\cite[Def. A.9]{Hiai21})
    A \bemph{positive quadratic form} on $\MH$ is a mapping $\MFq : \MH \supset \dom(\MFq) \longto [0, + \infty)$ which satisfies for all $\xi, \eta \in \dom(\MFq)$ and $\lambda \in \C$ the following identities:
    \begin{equation*}
        \MFq(\lambda \xi) = \abs{\lambda}^2 \MFq(\xi) \tand \MFq(\xi + \eta) + \MFq(\xi - \eta) = 2 \MFq(\xi) + 2 \MFq(\eta) \ .
    \end{equation*}
    The subspace $\dom(\MFq) \subset \MH$ is called the \bemph{domain} of $\MFq$; if it lies dense in $\MH$, then the quadratic form $\MFq$ is said to be \bemph{densely defined}.
\end{definition}

\begin{definition}[Closed and closable form]\label{def:forms_closedClosableForm}
    (\cite[Def. A.9]{Hiai21}, \cite[Def. 10.2]{Schmüdgen12})
    A positive quadratic form $\MFq$ on $\MH$ is called \bemph{closed} iff $(\xi_n)_{n \in \N} \subset \dom(\MFq)$, $\xi \in \dom(\MFq)$, $\lim_{n \to \infty} \xi_n = \xi$ and $\lim_{n, k \to \infty} \MFq(\xi_n - \xi_k) = 0$ implies that $\xi \in \dom(\MFq)$ and $\lim_{n \to \infty} \MFq(\xi_n - \xi) = 0$. The form $\MFq$ is said to be \bemph{closable} iff there exists a closed positive quadratic form $\MFt$ on $\MH$ which is an extension of $\MFq$, that is, $\dom(\MFq) \subset \dom(\MFt)$ and $\MFt(\xi) = \MFq(\xi)$ for all $\xi \in \dom(\MFq)$.
\end{definition}

\begin{definition}[Lower semi-continuity]\label{def:forms_semiContinuity}
    (\cite[p. 223]{Schmüdgen12})
    Let $(X, d)$ be a metric space. An extended real-valued function $f : X \longto \R \cup \{+ \infty\}$ is called \bemph{lower semi-continuous} iff for every sequence $(x_n)_{n \in \N}$ in $X$ which converges to some element $x \in X$, there holds
    \begin{equation*}
        f(x) \le \liminf_{n \to \infty} f(x_n) \ .
    \end{equation*}
\end{definition}

\begin{proposition}[\protect{Characterization of closable forms \cite[Thm. A.12]{Hiai21}, \cite[Prop. 10.3]{Schmüdgen12}}]\label{pro:forms_closableForms}
    Let $\MFq$ be a positive quadratic form on $\MH$. The following assertions are equivalent:
    \begin{enumerate}[equiv]
        \item $\MFq$ is closable.
        \item $\MFq$ is lower semi-continuous on $\dom(\MFq)$.
        \item For every sequence $(\xi_n)_{n \in \N} \subset \dom(\MFq)$ satisfying $\dlim{\MH}{n \to \infty} \xi_n = 0$ and $\lim_{n, m \to \infty} \MFq(\xi_n - \xi_m) = 0$, there holds $\lim_{n \to \infty} \MFq(\xi_n) = 0$.
    \end{enumerate}
\end{proposition}

\begin{para}[Closure of a form]\label{para:forms_closure}
    (\cite[p. 224]{Hiai21}, \cite[p. 224]{Schmüdgen12})
    Let $\MFq$ be a closable positive quadratic form. Define $\dom(\ol{\MFq}) \ce \set[\big]{\xi \in \MH \, : \, \exists (\xi_n)_{n \in \N} \subset \dom(\MFq) \ \text{s.t.} \ \norm{\xi_n - \xi} \to 0 \ \text{and} \ \MFq(\xi_n - \xi_m) \to 0} \subset \MH$ and the mapping $\ol{\MFq} : \MH \supset \dom(\ol{\MFq}) \longto [0, + \infty)$ by $\ol{\MFq}(\xi) \ce \lim_{n \to \infty} \MFq(\xi_n)$ for $\xi \in \dom(\ol{\MFq})$. This is a well-defined positive quadratic form. Furthermore, one can show that $\ol{\MFq}$ is closed, and that it is the smallest closed extension of the form $\MFq$. Therefore, one calls $\ol{\MFq}$ the \bemph{closure} of $\MFq$.
\end{para}

\section{The Form Representation Theorem}

\begin{para}[Order relation for self-adjoint operators]\label{para:forms_orderRelation}
    (\cite[Lem. A.1]{Hiai21}, \cite[p. 230]{Schmüdgen12})
    Let $T$ and $S$ be two positive self-adjoint linear operators on $\MH$. Define the relation $S \le T$ to be satisfied iff
    \begin{equation}\label{eq:forms_orderRelation}
        \dom\bigl(T^{1/2}\bigr) \subset \dom\bigl(S^{1/2}\bigr) \tand \vnorm[\big]{S^{1/2} \xi} \le \vnorm[\big]{T^{1/2} \xi} \quad \text{for all $\xi \in \dom\bigl(T^{1/2}\bigr)$} \ .
    \end{equation}
\end{para}

\begin{ntheorem}[\protect{Form representation theorem \cite[Thm. A.11]{Hiai21}, \cite[Thm. 10.7]{Schmüdgen12}}]\label{thm:forms_formReprThm}
    Let $\MFq : \MH \supset \dom(\MFq) \longto [0, + \infty)$ be a densely defined positive quadratic form. Then $\MFq$ is closed if and only if there exists a positive self-adjoint operator $Q_\MFq : \MH \supset \dom(Q_\MFq) \longto \MH$ such that
    \begin{equation}\label{eq:forms_formReprThm}
        \dom\bigl(Q_\MFq^{1/2}\bigr) = \dom(\MFq) \tand \forall \xi \in \dom\bigl(Q_\MFq^{1/2}\bigr) \, : \, \MFq(\xi) = \vnorm[\big]{Q_\MFq^{1/2} \xi}^2 \ .
    \end{equation}
\end{ntheorem}

\begin{para}[Operator associated with a positive form]\label{para:forms_fromFormToOp}
    (\cite[Rem. A.13]{Hiai21})
    Let $\MFq : \MH \supset \dom(\MFq) \longto [0, + \infty)$ be a densely defined positive quadratic form. If $\MFq$ is lower semi-continuous, then \cref{pro:forms_closableForms} implies that the closure $\ol{\MFq}$ (\cf{} \cref{para:forms_closure}) is well-defined. Applying \cref{thm:forms_formReprThm} to this closed form, one obtains a positive self-adjoint operator $T \ce Q_{\ol{\MFq}} : \MH \supset \dom(T) \longto \MH$ satisfying \eqref{eq:forms_formReprThm}. Since the closure $\ol{\MFq}$ is, in particular, an extension of $\MFq$, it follows that
    \begin{equation}\label{eq:forms_fromFormToOp}
        \dom\bigl(T^{1/2}\bigr) \supset \dom(\MFq) \tand \MFq(\xi) = \vnorm[\big]{T^{1/2} \xi} \com \xi \in \dom(\MFq) \ .
    \end{equation}
    Moreover, it holds that $T$ is the \emph{largest} positive self-adjoint linear operator, with respect to the order relation \eqref{eq:forms_orderRelation}, which satisfies the properties \eqref{eq:forms_fromFormToOp}, and, furthermore, $\dom(\MFq)$ is a core for the self-adjoint operator $T^{1/2}$ (\cref{def:operators_closed}). To see this, assume that $S$ were another positive self-adjoint operator on $\MH$ satisfying the relations \eqref{eq:forms_fromFormToOp}. For every $\xi \in \dom(T^{1/2}) = \dom(\ol{\MFq})$, there exists a sequence $(\xi_n)_{n \in \N} \subset \dom(\MFq)$ such that $\norm{\xi_n - \xi} \to 0$ and $\MFq(\xi_n) \to \ol{\MFq}(\xi)$ as $n \to + \infty$ (\cref{para:forms_closure}). Therefore, from the lower semi-continuity of $\MFq$, it follows that
    \begin{align*}
        \vnorm[\big]{S^{1/2} \xi}^2 = \MFq(\xi) \le \liminf_{n \to \infty} \MFq(\xi_n) = \ol{\MFq}(\xi) = \vnorm[\big]{T^{1/2} \xi} \ ,
    \end{align*}
    hence $S \le T$ according to \cref{para:forms_orderRelation}. For the second statement, note that by definition of the closure $\ol{\MFq}$ (\cref{para:forms_closure}) also $\MFq(\xi_n - \xi_m) \to 0$ as $n, m \to + \infty$. Therefore,
    \begin{align*}
        \lim_{n \to \infty} \norm{\xi - \xi_n} + \lim_{n \to \infty} \norm{T^{1/2} (\xi - \xi_n)} = \lim_{n \to \infty} \ol{\MFq}(\xi - \xi_n) = \lim_{n, m \to \infty} \MFq(\xi_m - \xi_n) = 0 \ .
    \end{align*}
\end{para}

\addtocontents{toc}{\protect\setcounter{tocdepth}{3}}

\backmatter
\printbibliography[heading=bibintoc]

\end{document}